\documentclass[10pt,reqno]{amsart}
\usepackage{graphicx}
\graphicspath{ {./images/} }
\usepackage{xargs}
\usepackage{geometry}
\usepackage{fancyhdr} 
\usepackage{lastpage} 
\usepackage{extramarks} 
\usepackage[most]{tcolorbox} 
\usepackage{xcolor} 
\usepackage[hidelinks]{hyperref} 
\usepackage[permil]{overpic}
\usepackage{pict2e} 
\usepackage{listings}
\usepackage{amssymb}
\usepackage{amsthm}
\theoremstyle{plain}
\usepackage[
backend=biber,
style=numeric,
sorting=ynt
]{biblatex}
\usepackage{xargs}
\usepackage{float}
\usepackage{lipsum}
\usepackage{geometry}
\usepackage{fancyhdr} 
\usepackage{lastpage} 
\usepackage{extramarks} 
\usepackage[most]{tcolorbox} 

\usepackage{xcolor} 
\usepackage[hidelinks]{hyperref} 

\usepackage[permil]{overpic}
\usepackage{pict2e} 
\usepackage{listings}

\newtheorem*{theorem*}{Theorem}

\theoremstyle{notation}

\numberwithin{equation}{section}
\theoremstyle{plain}
\newtheorem{definition}{Definition}[section]
\newtheorem{theorem}[equation]{Theorem}
\newtheorem{claim}[equation]{Theorem}
\newtheorem{remark}[equation]{Remark}
\newtheorem{corollary}[equation]{Corollary}
\newtheorem{lemma}[equation]{Lemma}
\newtheorem{conj}[equation]{Conjecture}
\newtheorem{proposition}[equation]{Proposition} 
\providecommand{\keywords}[1]
{
  \small	
  \textbf{\textit{Keywords---}} #1
}

\author{Eran Igra}

\title{Knots and Chaos in the Rössler System}
\address{Technion - Israel Institute of Technology}
\email{eranigra@simis.cn}
\begin{document}

\begin{abstract}
The Rössler System is one of the best known chaotic dynamical systems, exhibiting a plethora of complex phenomena - and yet, only a few studies tackled its complexity analytically. Inspired by recent numerical studies of the Rössler System, we introduce an idealized model for the Rössler System, prove its chaoticity, and study its bifurcations. 
\end{abstract}

\maketitle
\keywords{\textbf{Keywords} - The Rössler Attractor, Chaos Theory, Heteroclinic bifurcations, Topological Dynamics}
\section{Introduction}

In 1976, Otto E. Rössler introduced the following system of Ordinary Differential Equations, depending on parameters $A,B,C\in\mathbf{R}^{3}$ \cite{Ross76}:

\begin{equation} \label{Vect}
\begin{cases}
\dot{X} = -Y-Z \\
 \dot{Y} = X+AY\\
 \dot{Z}=B+Z(X-C)
\end{cases}
\end{equation}

Inspired by the Lorenz attractor (see \cite{Lo}), Otto E. Rössler attempted to find the simplest non-linear flow exhibiting chaotic dynamics. This is realized by the vector field above, as it has precisely one non-linearity, $XZ$ in the $\dot{Z}$ component. In more detail, the flow generated by the vector field above models the stretch-and-fold operation of a taffy machine (see \cite{Ros}) - and as observed by Rössler, the vector field appears to generate a chaotic attractor for $(A,B,C)=(0.2,0.2,5.7)$. In more detail, at these parameter values Rössler observed the first return map of the flow has the shape of a horseshoe (i.e., numerically), which is known to be chaotic (see \cite{S}).

\begin{figure}[h]
\centering
\begin{overpic}[width=0.5\textwidth]{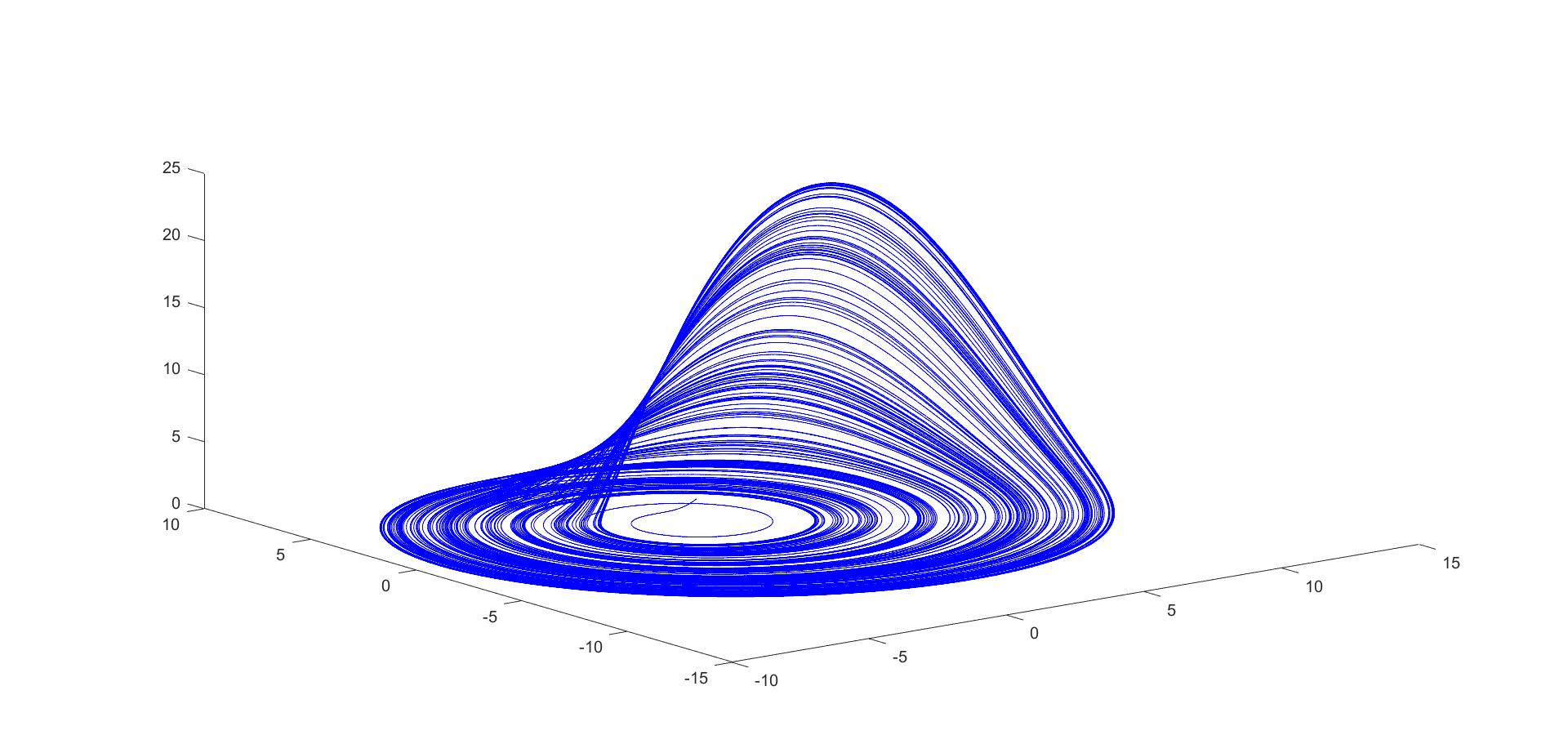}
\end{overpic}
\caption[Fig1]{The Rössler attractor at $(A,B,C)=(0.2,0.2,5.7)$}
\end{figure}

Since its introduction in 1976, the Rössler system was the focus of many numerical studies - despite the simplicity of the vector field, the flow gives rise to many non-linear phenomena (see, for example: \cite{MBKPS}, \cite{BBS}, \cite{G}, \cite{BB2}, \cite{SR}). One particular feature is that varying the parameters $A,B,C$ often leads to a change in the complexity of the system: that is, as the parameters are varied, more and more symbols appear in the first-return map of the flow (for more details, see \cite{MBKPS}, \cite{BBS},\cite{Le},\cite{RO}). In a topological context, several numerical studies noted this variation of parameters changes the topology of the attractor - see, for example, \cite{Le},\cite{RO}. \\

In contrast to the vast corpus of numerical studies, analytical results on the Rössler system are sparse. For example, in \cite{Pan}, the existence of periodic trajectories at some parameters was established; in \cite{CNV} the existence of an invariant Torus (and its breakdown) at some parameters was proven; and in \cite{LiLl} the dynamics of the flow at $\infty$ were analyzed. As far as chaotic dynamics go, their existence at $(A,B,C)=(0.2,0.2.5.7)$ was proven with rigorous numerical methods - see \cite{Zgli97},\cite{XSYS03}. To our knowledge, no studies on the Rössler system attempted to explain its nonlinear phenomena (and in particular, its chaotic dynamics) by purely analytical tools. It is precisely this gap that this paper aims to address. In this paper we analyze an idealized model of the Rössler, and prove analytically it is forced to behave chaotically - which we do based on the topology generated by the Rössler system.  \\ 

To state our main results, given parameter values $p=(A,B,C)$, denote by $F_p$ the corresponding vector field which generates the flow. We first prove the following result about the global dynamics of the  Rössler system in Section $2$ (see Lemma \ref{Clemma} and Th.\ref{th21}):
\begin{claim}
\label{th1} 
There exists an open set of parameters $O\subseteq\mathbf{R}^3$ s.t. $F_p$ can be extended to a vector field on $S^3$ with precisely three fixed points - two saddle-foci $P_{In},P_{Out}$ (of opposing indices) and a degenerate fixed point at $\infty$ of index $0$. Moreover, both $P_{In},P_{Out}$ admit heteroclinic trajectories connecting them to $\infty$.
\end{claim}
As can be seen, Th.\ref{th1} does not provide us with constrains on the flow which force the existence of complex dynamics. As such, in order to analytically describe the dynamics of the Rössler system we study the dynamics of an idealized model of the Rössler system modeled on very specific parameter values $p\in P$ - to which we refer as \textbf{trefoil parameters} (for the precise formulation, see Def.\ref{def32}).\\

\begin{figure}[h]
\centering
\begin{overpic}[width=0.3\textwidth]{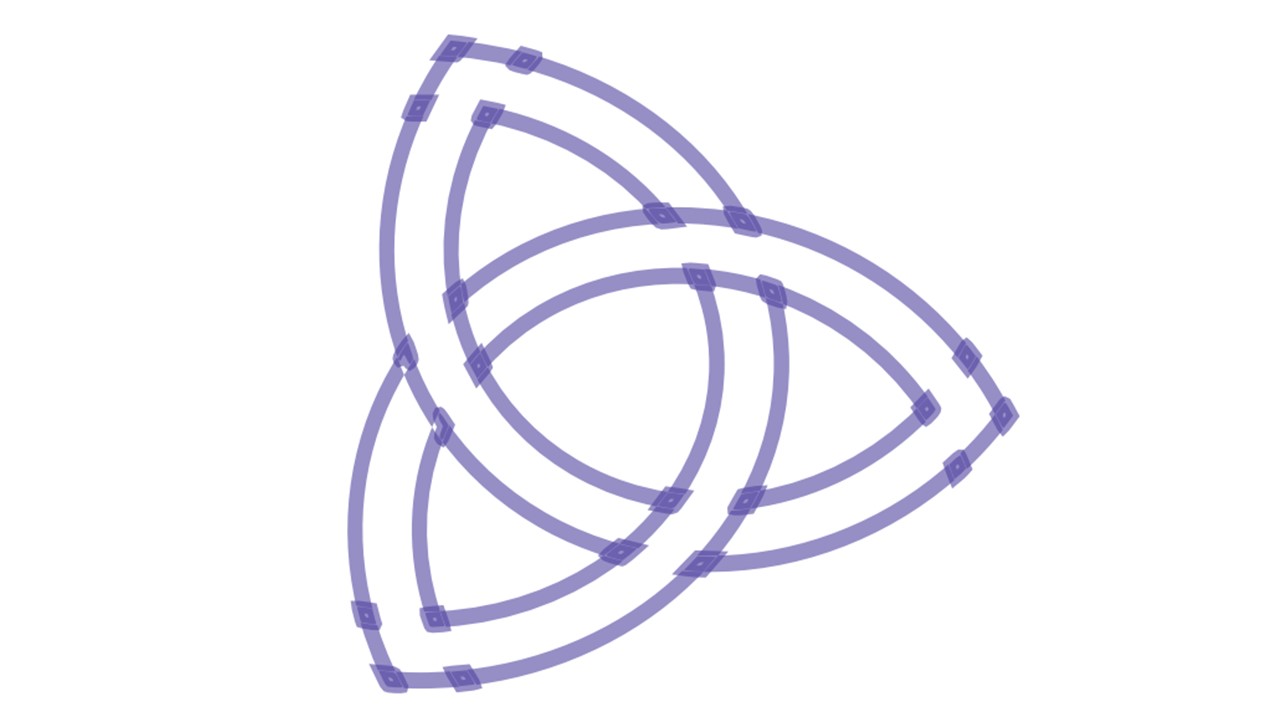}
\end{overpic}
\caption[A trefoil knot type.]{\textit{A trefoil knot type.}}
\label{trefoilknot}
\end{figure}

Roughly speaking, trefoil parameters are parameter values $(A,B,C)$ for which the Rössler system generates both a heteroclinic trefoil knot knot $\Lambda$ in the $3$-sphere $S^3$, whose knot type is that of the trefoil (see Fig.\ref{trefoilknot}), and bounded set of invariant dynamics. The existence of heteroclinic knots for the Rössler system was observed numerically (see Fig.5.B.1 in \cite{MBKPS} and the discussion at the beginning of Section 3). Motivated and inspired by both the topology of ${S}^3\setminus\Lambda$ and by Th.1 in \cite{KY}, using Th.\ref{th1} and the Betsvina-Handel Algorithm (see \cite{BeH}), we prove the following fact in Section $3$ (see Th.\ref{th31} in Sect.3):
\begin{claim}\label{th2}
The flow dynamics of trefoil parameters generate a bounded invariant set $Q$, which includes infinitely many periodic trajectories for the flow. Moreover, these periodic orbits are in surjective correspondence with those of the one-sided shift in $\{1,2\}^\mathbf{N}$. 
\end{claim}
Once Th.\ref{th2} is proven, we apply it to analyze the dynamical complexity of the Rössler system. Inspired by \cite{Zgli97} and using both Th.\ref{th2} and a Fixed-Point Index argument, we prove the following Theorem in Section $4$ (see Th.\ref{conti}):
\begin{claim}\label{th3}
Let $O$ be the parameter space from Th.\ref{th1}, and choose some $n>0$. Then, provided $v\in O$ corresponds to some Rössler system which is sufficiently $C^1-$close to the dynamics of a trefoil parameter, the said Rössler system generates at least $n$-distinct periodic trajectories.
\end{claim}
Theorem \ref{th3} has the following meaning - the $C^1$-closer the dynamics of the Rössler system corresponding to $v$ are to those of the idealized model of the Rössler system, the more complex they will be. With these ideas in mind, we remark that when the parameter space of the Rössler system was analyzed numerically, spiral bifurcation structures were observed (see, for example \cite{MBKPS}, \cite{BBS}, \cite{G} and \cite{BB2}). In all these studies, these spiral structures always accumulated at some point $p_0$, often referred to as a \textbf{periodicity hub}, which lies on the Shilnikov homoclinic curve (see \cite{LeS}). As was observed in\cite{SR} (see pg. 430), the dynamics around some periodicity hubs may, in fact, be heteroclinic. This suggests Th.\ref{th2} and \ref{th3} possibly have a part in explaining the emergence of such complex bifurcation phenomena.\\

Finally, we also stress that even though our results all relate to the Rössler system, in practice the arguments in Section $3$ and $4$ are mostly topological. As such, they exemplify the potential of topological methods in explaining the emergence of nonlinear phenomena. Moreover, our results also attest to the importance of bounded heteroclinic trajectories to the emergence of chaotic and complex dynamics in three-dimensional flows. In more detail, it is well-known that heteroclinic trajectories generate chaotic dynamics in the Lorenz system (see, for example, the results of \cite{Pi} and the references therein). As such, this study and \cite{Pi} lend further credence to the role of heteroclinic knots in the onset of complex dynamics. 
\subsubsection{Acknowledgments:}
The author would like to thank Tali Pinsky, Genadi Levin and Andrey Shilnikov for their helpful suggestions and support. Moreover, the author is particularly grateful to Noy Soffer-Aranov for her constant encouragement, and to the anonymous referee for the helpful comments and suggestions.
\section*{Preliminaries}
From now on, given $(a,b,c)\in\mathbf{R}^3$ we switch to the more convenient form of the Rössler system, defined by the following system of ODEs:
\begin{equation} \label{Field}
\begin{cases}
\dot{x} = -y-z \\
 \dot{y} = x+ay\\
 \dot{z}=bx+z(x-c)
\end{cases}
\end{equation}
Denote the vector field corresponding to $(a,b,c)\in\mathbf{R}^3$ by $F_{a,b,c}$ . This definition is slightly different from the one presented in Eq.\ref{Vect} - however, setting $p_1=\frac{-C+\sqrt{C^2-4AB}}{2A}$, it is easy to see that whenever $C^2-4AB>0$, $(X,Y,Z)=(x-ap_1,y+p_1,z-p_1)$ defines a change of coordinates between the vector fields in Eq.\ref{Vect} and Eq.\ref{Field}. \\

Before we continue, we introduce the following definition for chaotic dynamics which we will use throughout this paper. To do so, let $\sigma:\{1,2\}^\mathbf{N}\to\{1,2\}^\mathbf{N}$ denote the one-sided shift - we now define:

\begin{definition}{\textbf{Chaotic Dynamics}}\label{chaotic}
Let $F$ be a $C^\infty$ vector field on $\mathbf{R}^3$. We say $F$ is \textbf{chaotic} or \textbf{generates chaotic dynamics} provided there exists a bounded, invariant set $B\subseteq\mathbf{R}^3$ for the flow satisfying the following:
\begin{enumerate}
    \item $B$ includes infinitely many periodic orbits.
    \item There exists a cross-section $S\subseteq\mathbf{R}^3$ transverse to $B$ s.t. the following holds:
    \begin{itemize}
        \item The first-return map $f:B\cap S\to B\cap S$ is continuous.
        \item There exists some $\sigma-$invariant $\Sigma\subseteq\{1,2\}^\mathbf{N}$ which includes infinitely many periodic orbits for $\sigma$, and a continuous, surjective $\pi:B\cap S\to\Sigma$ s.t. $\pi\circ f=\sigma\circ\pi$.
        \item If $s\in\Sigma$ is periodic of minimal period $k$ w.r.t. $\sigma$, $\pi^{-1}(s)$ includes at least one periodic point of minimal period $k$ w.r.t. $f$.
    \end{itemize}
\end{enumerate}
\end{definition}

For example, any suspension of a Smale Horseshoe (see \cite{S}) generates chaotic dynamics. Moreover, as proven in Th.1.1 in \cite{Pi}, there exist parameter values for which the Lorenz system is chaotic.\\ 

Returning to the Rössler system, recall that given any parameter value $(a,b,c)=p$ we denote the corresponding vector field by $F_p$ (see Eq.\ref{Field}). Since this vector field depends on three parameters we now specify the region in the parameter space in which we prove our results. From now on unless stated explicitly otherwise, the parameter space $P\subseteq\mathbf{R}^3$ we consider throughout this paper is composed of parameters satisfying the following\label{eq:9}:

\begin{itemize}
\item \textbf{Assumption $1$} -  for every parameter $p\in P,p=(a,b,c)$ the parameters satisfy $a,b\in(0,1)$ and $c>1$. For every choice of such $p$, the vector field $F_p$ given by Eq.\ref{Field} always generates precisely two fixed points - $P_{In}=(0,0,0)$ and $P_{Out}=(c-ab,b-\frac{c}{a},\frac{c}{a}-b)$.
\item \textbf{Assumption 2 }- for every $p\in P$ the fixed points $P_{In},P_{Out}$ are both saddle-foci of opposing indices. In more detail, we always assume that $P_{In}$ has a one-dimensional stable manifold, $W^s_{In}$, and a two-dimensional unstable manifold, $W^u_{In}$. Conversely, we always assume $P_{Out}$ has a one-dimensional unstable manifold, $W^u_{Out}$, and a two-dimensional stable manifold, $W^s_{Out}$ (see the illustration in Fig.\ref{loci}). 
\item \textbf{Assumption 3 }- For every $p\in P$, let $\gamma_{In}<0$ and $\rho_{In}\pm i\psi_{In}$, $\rho_{In}>0$ denote the eigenvalues of $J_p(P_{In})$, the linearization of $F_p$ at $P_{In}$, and set $\nu_{In}=|\frac{\rho_{In}}{\gamma_{In}}|$. Conversely, let $\gamma_{Out}>0$, $\rho_{Out}\pm i\psi_{Out}$ s.t. $\rho_{Out}<0$ denote the eigenvalues of $J_p(P_{Out})$, the linearization at $P_{Out}$, and define $\nu_{Out}=|\frac{\rho_{Out}}{\gamma_{Out}}|$. We will refer to $\nu_{In},\nu_{Out}$ as the respective saddle indices at $P_{In},P_{Out}$, and we will always assume $(\nu_{In}<1)\lor(\nu_{Out}<1)$ - that is, for every $p\in P$ at least one of the fixed points satisfies the Shilnikov condition. 
\end{itemize}

It is easy to see the parameter space $P$ we are considering is open in the parameter space. Moreover, it includes the parameter space for the Rössler system considered in \cite{BBS},\cite{MBKPS} and \cite{G}, where many interesting bifurcation phenomena were observed.\\

\begin{figure}[h]
\centering
\begin{overpic}[width=0.4\textwidth]{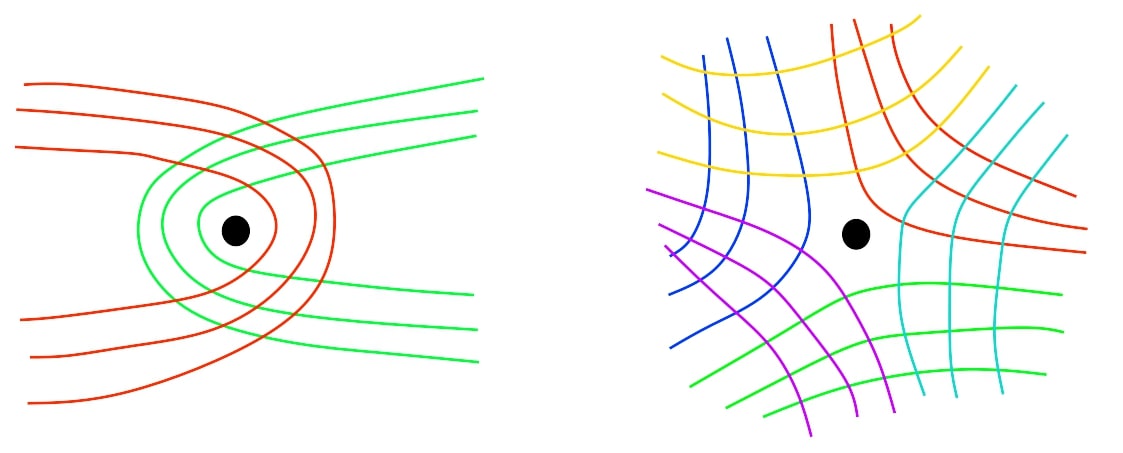}
\end{overpic}
\caption[Transverse folliations.]{\textit{Transverse folliations (with singularities) around the punctures of $S$ (i.e. the black discs). }}
\label{trannn}
\end{figure}

In addition, we will also need several results and notions from the theory of two-dimensional surface homeomorphisms. To begin, let $D$ denote a planar disc. Until the end of this section, $S$ would always denote a surface homeomorphic to $D\setminus\{x_1,...,x_n\}$, where $x_1,...,x_n$ are interior to $D$ and $n\geq2$. Additionally, $f:S\to S$ would always denote a homeomorphism which permutes the set $\{x_1,...,x_n\}$. We first introduce the following definition:

\begin{definition}
    \label{pseanosov} A homeomorphism $f:S\to S$ is \textbf{Pseudu-Anosov} provided there exist two folliations of $S$, $F^u$ and $F^s$, transverse to one another throughout $S$ (but not necessarily at the punctures $\{x_1,...,x_n\}$ - see the illustration in Fig.\ref{trannn}) and some $\lambda>0$ s.t. the following is satisfied:
    \begin{itemize}
        \item Both $F^s$ and $F^u$ are measured - i.e., if we move some leaf $L_1$ of $F^i$ to another leaf $L_2$ of $F^i$ by some isotopy of $S$ (where $i\in\{u,s\}$), the Borel measure on $L_2$ is the pushforward of the Borel measure on $L_1$.
        \item $f(F^u)=\lambda F^u$ while $f(F^s)=\frac{1}{\lambda}F^s$ - i.e., $f$ stretches uniformly the unstable folliation $F^u$ and squeezes uniformly the stable folliation $F^s$. We refer to $\lambda$ as the \textbf{expansion constant}.
    \end{itemize}
\end{definition}

\begin{figure}[h]
\centering
\begin{overpic}[width=0.4\textwidth]{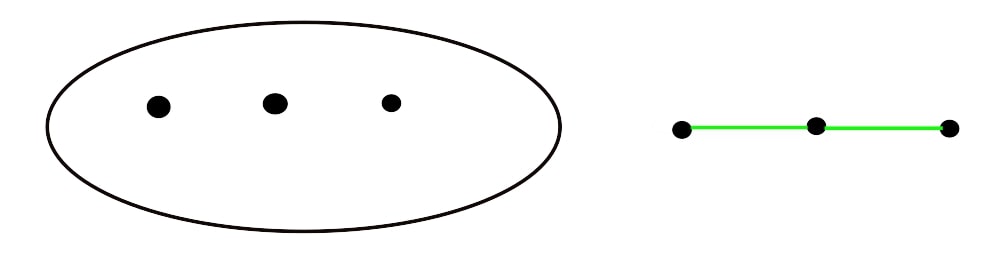}
\end{overpic}
\caption[A spine for a surface.]{\textit{On the left we have a surface $S$ homeomorphic to $D$ punctured at $3$ points, and on the right we have its spine, the graph $\Gamma$. As can be seen, $\Gamma$ has two edges and three vertices.} }
\label{spine}
\end{figure}

Pseudo-Anosov maps are essentially a generalized form of hyperbolic diffeomorphisms, like Smale's Horseshoe (see \cite{S}) - that is, they contract uniformly in one direction, and expand uniformly in another. The reason we are interested in Pseudo-Anosov maps is because as far as their dynamics are concerned, Pseudo-Anosov maps are dynamically minimal in the following sense - if $f:S\to S$ is a homeomorphism isotopic to a Pseudo-Anosov map $F:S\to S$, then the dynamics of $F$ are complex at least like those of $F$. More precisely, we have the following result proven as Th.2 and Remark 2 in \cite{Han}:

\begin{theorem}
    Let $F:S\to S$ be a Pseudo-Anosov map and let $f:S\to S$ be a homeomorphism isotopic to $F$. Then, there exists a closed set $Y\subseteq S$ and a continuous, surjective $\pi:Y\to S$ s.t. $\pi\circ f=F\circ \pi$.
\end{theorem}
We are now led to the following question - given a homeomorphism $f:S\to S$, when is it isotopic to a Pseudo-Anosov map? To introduce the answer to that question, we must first introduce several concepts. We begin with the notion of a spine:

\begin{definition}
    The \textbf{spine} of $S$ is a graph $\Gamma$ embedded in $\overline{S}$ with $n-1$ vertices and $n-2$ edges s.t. $\Gamma$ is a retract of $S$ (see the illustration in Fig.\ref{spine}).
\end{definition}

Now, let us consider a graph map $g:\Gamma\to \Gamma$, and let us denote by $S_1,...,S_{n-2}$ the edges of $\Gamma$ (see the illustration in Fig.\ref{spine}). We now define a matrix $A=\{a_{i,j}\}_{1\leq i,j\leq n-2}$ s.t. $a_{i,j}$ is the number of times $g(S_j)$ covers $S_i$. Finally, let us define the \textbf{spectral radius} of $A$ as the maximal eigenvalue for $A$ - and let us further note that since $A$ is a matrix with positive coefficients, its spectral radius is also positive. To continue, note that given a homeomorphism $f:S\to S$, by retracting $S$ to its spine $\Gamma$ with edges $S_1,...S_{n-2}$ we can reduce $f$ to a graph map $g:\Gamma\to\Gamma$, as described in Fig.\ref{cover2}. As a consequence, there exists a matrix $A$ generated by $g$ as described above, with a spectral radius $\gamma$. Then, we have the following fact, proven in Sections $3.4$ and $4.4$ of \cite{BeH}:

\begin{theorem}
    \label{betshan} With the notations above, whenever $\gamma>1$ $f$ is isotopic to a Pseudo-Anosov map $F:S\to S$. Moreover, if there exists an edges $S_i$ s.t. $g(S_i)$ covers itself at least twice (for some $1\leq i\leq n-2$), the dynamics of $f$ in the regions in $S$ collapsed to $S_i$ includes periodic orbits of all minimal periods.
\end{theorem}
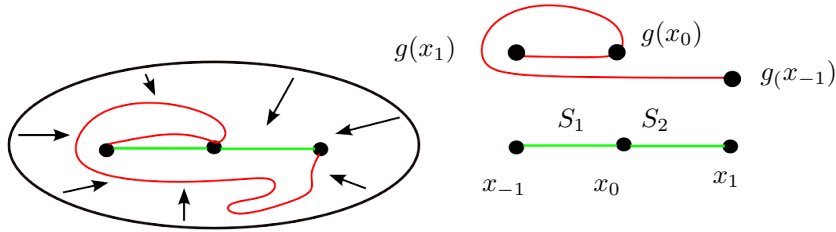
\begin{figure}[h]
\centering
\begin{overpic}[width=0.6\textwidth]{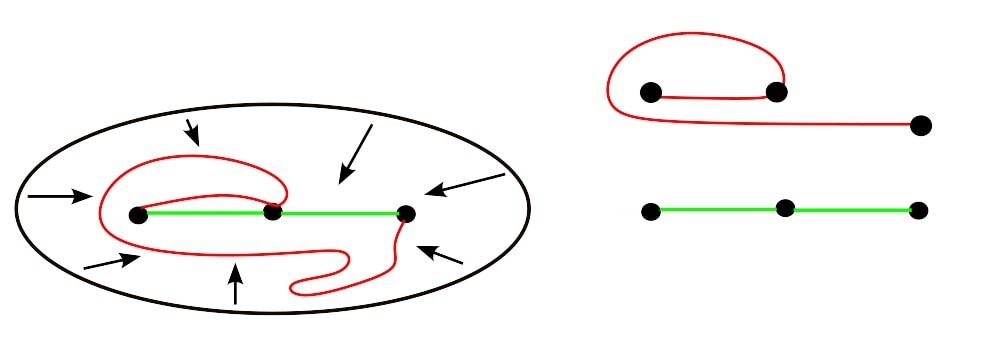}
\put(705,160){$S_1$}
\put(810,160){$S_2$}
\put(750,80){$x_0$}
\put(900,90){$x_1$}
\put(610,80){$x_{-1}$}
\put(960,220){$g_(x_{-1})$}
\put(810,270){$g(x_0)$}
\put(500,250){$g(x_1)$}
\end{overpic}
\caption[The graph map $g$.]{\textit{ $f$ distorts the graph $\Gamma$ inside $S$ into the red curve, s.t. $f(x_{-1})=x_1$, $f(x_1)=x_{-1}$ and $f(x_1)=x_1$ (while the outer circle remains fixed). Consequentially, $g(S_1)$ covers itself twice.}}\label{cover2}

\end{figure}

\section{The global dynamics of the Rössler system.}
\label{sect2}
In this section we perform basic qualitative analysis of the Rössler system for parameters $(a,b,c)$ in the parameter space $P$ - thus establishing certain key facts which will be used throughout this paper. This section is organized as follows: we first consider a cross-section for the flow given by $\{\dot{y}=0\}$, and study the local dynamics on it. Following that, we study the unbounded dynamics of the flow, and prove Th\ref{th21} (i.e., Th.\ref{th1} from the introduction), with which we conclude this section.\\

\begin{figure}[h]
    \centering
\begin{overpic}[width=0.5\textwidth]{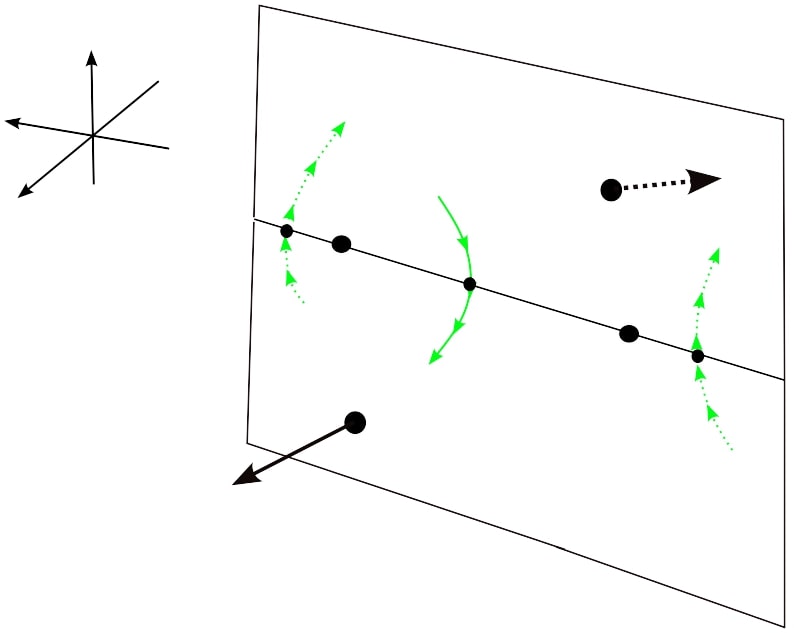}
\put(700,330){$P_{In}$}
\put(1050,200){$\{\dot{y}<0\}$}
\put(100,150){$\{\dot{y}>0\}$}
\put(390,420){$P_{Out}$}
\put(370,700){$U_p$}
\put(700,150){$L_p$}
\put(-10,630){$x$}
\put(-20,540){$y$}
\put(115,730){$z$}
\end{overpic}
\caption[The plane $\{\dot{y}=0\}$.]{\textit{The plane $\{\dot{y}=0\}$ and the cross sections $U_p$ and $L_p$ (and the directions of $F_p$ on them), separated by $l_p$. The green arcs represent flow lines tangent to $\{\dot{y}=0\}$ at different arcs on $l_p$ - see Lemma \ref{lem23}.}}
\label{planes}
\end{figure}
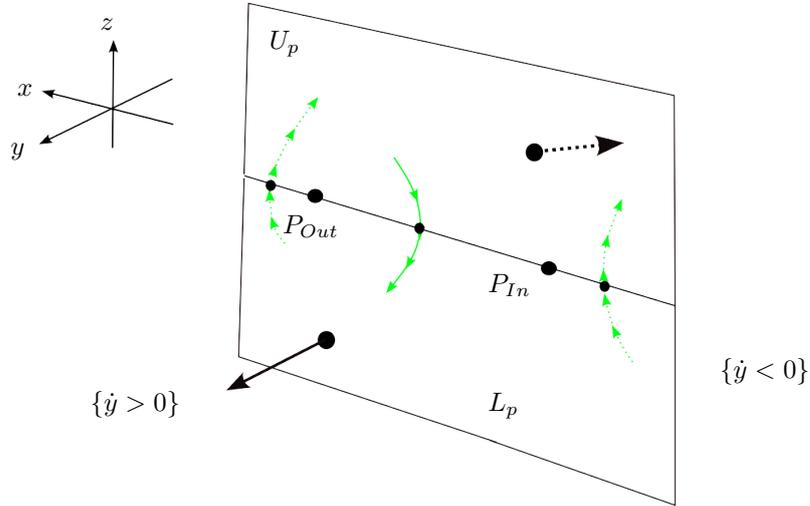

To begin, fix some parameter $p=(a,b,c)\in P$ and recall we denote the vector field generating the corresponding Rössler system by $F_p$ (see Eq.\ref{Field}). Now, consider the plane $\{(x,-\frac{x}{a},z)|x,z\in\mathbf{R}\}=\{\dot{y}=0\}$ \label{ydot}, and let $N_p=(1,a,0)$ denote the normal vector to $\{\dot{y}=0\}$ (with the velocity $\dot{y}$ taken w.r.t. $F_p$ - see Eq.\ref{Field}). By direct computation of the product $F_p(x,-\frac{x}{a},z)\bullet N_p$ we see $F_p$ is tangent to $\{\dot{y}=0\}$ precisely at the straight line $l_p=\{(t,-\frac{t}{t},\frac{x}{a})|t\in\mathbf{R}\}$. Hence, $\{\dot{y}=0\}\setminus l_p$ is consists of two components, both half-planes - let $U_p=\{(x,-\frac{x}{a},z)|x,z\in\mathbf{R},-z+\frac{x}{a}<0\}$ denote the upper half of $\{\dot{y}=0\}\setminus l_p$, and denote by $L_p=\{(x,-\frac{x}{a},z)|x,z\in\mathbf{R},-z+\frac{x}{a}>0\}$ the lower half (see the illustration in Fig.\ref{planes}).\label{eq:1} \\

By the definition of the regions $\{\dot{y}>0\}=\{(x,y,z)|x+ay>0\}$, $\{\dot{y}<0\}=\{(x,y,z)|x+ay<0\}$ correspond to the regions in front and below $\{\dot{y}=0\}$ (respectively - see Fig.\ref{planes}) - and by the sign of $F_p(x,-\frac{x}{a},z)\bullet N_p$ on $U_p$ and $L_p$ we immediately conclude:

\begin{itemize}
    \item On $U_p$, the vector field $F_p$ points inside $\{\dot{y}<0\}$ (see the illustration in Fig.\ref{planes}).
    \item On $L_p$, the vector field $F_p$ points inside $\{\dot{y}>0\}$ (see the illustration in Fig.\ref{planes}).
\end{itemize}

To continue, given an initial condition $s$, from now on we always denote its trajectory by $\gamma_s$ - parameterized s.t. $\gamma_s(0)=s$. As $l_p$, the tangency set of $F_p$ to $\{\dot{y}=0\}$ is a straight line, we immediately conclude the following, useful fact:

\begin{lemma}\label{obs}
Let $s\in\mathbf{R}^3$ be an initial condition that is not a fixed point, whose forward trajectory $\gamma_s$ is bounded and does does not limit to a fixed point - then, $\gamma_s$ intersects $U_p$ transversely infinitely many times. In particular, every periodic trajectory intersects $U_p$ transversely at least once.
\end{lemma}

Lemma \ref{obs} motivates us to introduce the following notations:

\begin{definition}
    \label{FRM} From now on throughout the remainder of this dissertation, we always denote by $f_p:\overline{U_p}\to\overline{U_p}$ the first-return map - wherever defined in $\overline{U_p}$.
\end{definition}

\begin{figure}[h]
\centering
\begin{overpic}[width=0.65\textwidth]{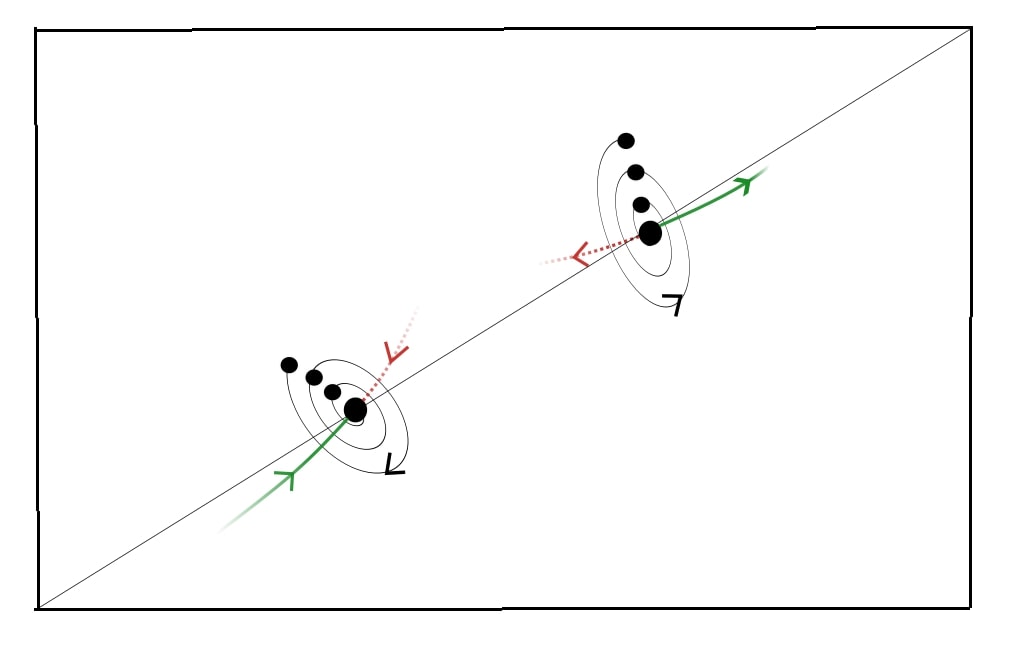}
\put(630,280){$W^s_{Out}$}
\put(750,430){$W^u_{Out}$}
\put(690,390){$P_{Out}$}
\put(390,130){$W^u_{In}$}
\put(410,220){$P_{In}$}
\put(440,420){$U_p$}
\put(570,180){$L_p$}
\put(250,110){$W^s_{In}$}
\end{overpic}
\caption[The local dynamics around the fixed points.]{\textit{The local dynamics around the fixed points. The green and red flow lines are the components of $W^s_{In},W^u_{Out}$.}}\label{loci}
\end{figure}
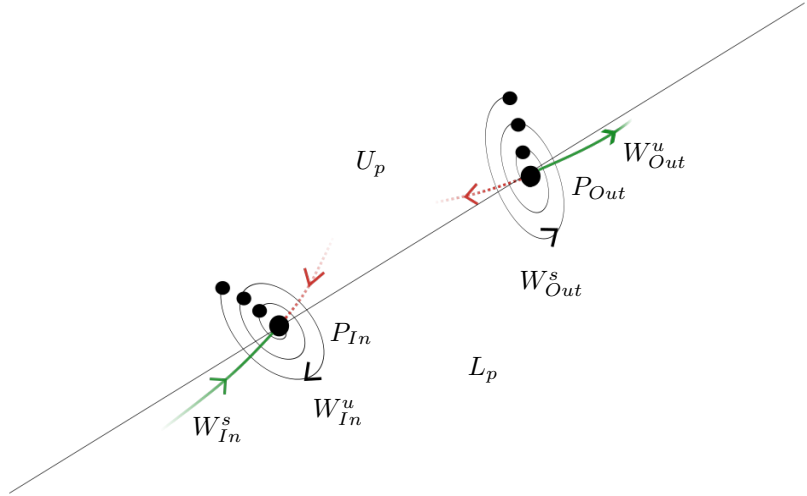

By Lemma \ref{obs}, $f_p$ is defined at least at every initial condition $s\in U_p$ whose forward trajectory is both bounded and does not limit to a fixed point. Continuing in our analysis of the plane $\{\dot{y}=0\}$, we now prove the following fact:

\begin{lemma}\label{cor211}
Let $p=(a,b,c)$ be some parameter value in $P$ - then, the corresponding two dimensional, invariant manifolds $W^u_{In},W^s_{Out}$ are transverse to $U_p,L_p$ at both $P_{In},P_{Out}$ (see the illustration in Fig.\ref{loci}). 
\end{lemma}
\begin{proof}
    Let us first consider the tangent vector to the straight line $l_p$, the tangency curve of $F_p$ to the plane $\{\dot{y}=0\}$ - as we parameterize $l_p$ by $l_p(x)=(x,-\frac{x}{a},\frac{x}{a})$, $x\in\mathbf{R}$, the tangent vector to $l_p(x)$ is given $v=(1,-\frac{1}{a},\frac{1}{a})$. Now, let us consider $J_p(P_{In})$, the Jacobian matrix of $F_p$ at the fixed-point $P_{In}=(0,0,0)$, given by the following matrix (see Eq.\ref{Field}):

\begin{equation}
    \begin{pmatrix}
    0 & -1 & -1\\
    1 & a & 0\\
    b & 0 & -c
\end{pmatrix}
\end{equation}

By direct computation, $J_p(P_{In})v=(0,0,b-\frac{c}{a})$ - as $a,b\in(0,1)$ and $c>1$, it follows we always have $ab\ne c$, i.e., $b-\frac{c}{a}\ne0$, which implies $v$ is not an eigenvector for $J_p(P_{In})$. Let us now remark the plane $\{\dot{y}=0\}$ is spanned by $(1,-\frac{1}{a},\frac{1}{a})=v$ and $(0,0,1)=\mu$, and that similarly to the computation above, $\mu$ is also not an eigenvalue for $J_p(P_{In})$.\\

Generalizing these arguments, we now prove that given any non-zero $\nu\in\{\dot{y}=0\}$, $\nu$ does not span an invariant subspace for $J_p(P_{In})$. To do so, recall that every $\nu\in\{\dot{y}=0\}$ is parameterized by $\nu=(x,-\frac{x}{a},z)$ (for some $x,z\in\mathbf{R}$). By computation, $J_p(P_{In})\nu=(-z+\frac{x}{a},0,bx+z(x-c))$ - which implies that whenever $x\ne 0$, $\nu$ is not an invariant direction for $J_p(P_{In})$. When $x=0$ and $z\ne0$ we have $\nu=(0,0,z)$, and consequentially $J_p(P_{In})\nu=(-z,0,-cz)$ - therefore, again, $\nu$ is not an invariant direction for $J_p(P_{In})$. Therefore, all in all, the plane $\{\dot{y}=0\}$ does not include any invariant directions for the Jacobian matrix $J_p(P_{In})$.\\

Now, let $W$ denote the two-dimensional unstable invariant subspace for $J_p(P_{In})$ in $\mathbf{R}^3$. Recall that since $P_{In}$ is a saddle-focus, $W$ is spanned by the vectors $r_1=\frac{\omega+\overline{\omega}}{2}$, $r_2=\frac{\omega-\overline{\omega}}{2i}$, for some $\omega\in\mathbf{C}^3$ which is an eigenvector of $J_p(P_{In})$ (where $\overline{\omega}$ is the complex-conjugate vector to $\omega$). It is easy to see both $r_1,r_2$ span eigenspaces for $J_p(P_{In})$, which, by the discussion above, implies $r_1,r_2\not\in\{\dot{y}=0\}$. Consequentially, as $W$ is spanned by $r_1$ and $r_2$, we conclude $W$ and $\{\dot{y}=0\}$ are transverse to one another. Consequentially, as the two-dimensional unstable manifold $W^u_{In}$ is tangent to $W$ at $P_{In}$ (see Th.2.7.1 in \cite{Per}) it follows $W^u_{In}$ must be transverse to $\{\dot{y}=0\}$ at $P_{In}$.\\

Similarly, by considering $J_p(P_{Out})$, the linearization of the saddle focus $P_{Out}=(c-ab,\frac{ab-c}{a},\frac{c-ab}{a})$, using a similar argument it again follows $\{\dot{y}=0\}$ does not include any invariant directions under $J_p(P_{Out})$ - therefore, a similar argument prove $W^s_{Out}$, the two-dimensional stable manifold of $P_{Out}$, is transverse to $\{\dot{y}=0\}$ at $P_{Out}$ and Lemma \ref{cor211} now follows (see the illustration in Fig.\ref{loci}).
\end{proof}

Having proven the invariant manifolds $W^u_{In}$ and $W^s_{Out}$ are transverse to $\{\dot{y}=0\}$ at the fixed points, we now study the local dynamics of $F_p$ on the tangency line $l_p$. To do so, first recall that given any $s\in\mathbf{R}^3$, we denote by $\gamma_s$ its trajectory, parameterized s.t. $\gamma_s(0)=s$ - and furthermore, recall we parameterize $l_p$ by $l_p(x)=(x,-\frac{x}{a},\frac{x}{a})$, where $x\in\mathbf{R}$. We now prove:
\begin{lemma}\label{lem23}
We have the equality $l_p=\{\dot{x}=0\}\cap\{\dot{y}=0\}$. Additionally, given $s\in l_p$, we also have the following:
\begin{itemize}
    \item For $x<0$ and $x>c-ab$, there exists some $\epsilon>0$ s.t. $\gamma_s(t)\in\{\dot{y}<0\}$ for all $t\ne 0,t\in(-\epsilon,\epsilon)$ (see the illustration in Fig.\ref{planes}). 
    \item Conversely, when $x\in(0,c-ab)$ there exists an $\epsilon>0$ s.t. for $t\in(-\epsilon,\epsilon),t\ne0$ we have $\gamma_s(t)\in\{\dot{y}>0\}$ (see the illustration in Fig.\ref{planes}). 
\end{itemize}

\end{lemma}
\begin{proof}
Note that given $p=(a,b,c)$ by Eq.\ref{Field} we have $\dot{x}=-y-z$, $\dot{y}=x+ay$ (where $x,y,z\in\mathbf{R}$) - which implies $\{\dot{x}=0\}\cap\{\dot{y}=0\}$ is parameterized by the set $\{(x,-\frac{x}{a},\frac{x}{a})|x\in\mathbf{R}\}$, i.e., the intersection is given by the line $l_p$. By computation, we have $F_p(l_p(x))=(0,0,bx+\frac{x^2}{a}-\frac{cx}{a})$. Now, let us note $bx+\frac{x^2}{a}-\frac{cx}{a}$ is a quadratic polynomial vanishing precisely at $P_{In},P_{Out}$ - therefore $F_p(l_p(x))$ points in the positive $z$-direction precisely when $x<0$ or $x>c-ab$ (conversely, it points in the negative $z$ direction when $x\in(0,c-ab)$).\\

On $U_p$ the vector field $F_p$ points into $\{\dot{y}>0\}$ (see the discussion preceding Lemma \ref{obs} and the illustration in Fig.\ref{planes}). Since $U_p$ is the upper hald-plane, it follows that for $x<0$ or $x>c-ab$ the trajectory of $s=l_p(x)$ can only arrive at $l_p(x)$ from $\{\dot{y}<0\}$ (the region behind $\{\dot{y}-0\}$ in Fig.\ref{planes}). Therefore, by the tangency of $F_p$ to the plane $\{\dot{y}=0\}$ at $s$ it follows the forward-trajectory of $s$ enters $\{\dot{y}<0\}$ upon leaving $s$ (as illustrated in Fig.\ref{planes}) - and it follows that whenever $x<0$ or $x>c-ab$, there exists an $\epsilon>0$ (depending on $s$) s.t. $\gamma_s(t)\in\{\dot{y}<0\}$ for all $t\ne 0,t\in(-\epsilon,\epsilon)$.\\

When $x\in(0,c-ab)$, $s=(x,-\frac{x}{a},\frac{x}{a})$, by the discussion above $F_p(s)$ points in the negative $z$-direction - which, using a similar argument, implies the backwards trajectory of $s$ arrives from $\{\dot{y}>0\}$, the region in front of $\{\dot{y}=0\}$ in Fig.\ref{planes}. By the tangency of $F_p$ to $\{\dot{y}=0\}$ at $s$, a similar argument to the one used in the previous paragraph now proves there exists some $\epsilon>0$ s.t. for $t\ne 0$, $t\in(-\epsilon,\epsilon)$ we have $\gamma_s(t)\in\{\dot{y}>0\}$ and Lemma \ref{lem23} now follows.
\end{proof}
Having concluded our study of the plane $\{\dot{y}=0\}$, we now begin studying the unbounded dynamics of the Rössler system in $\mathbf{R}^3$. We first recall Th.1 in \cite{LiLl}, where it was proven that for all $p\in P$, the behavior of the vector field $F_p$ on $\{(x,y,z)|x^2+y^2+z^2>r\}$ is independent of $r>0$ (for any sufficiently large $r$) - which proves one can extend $F_p$ to $S^3$ by adding $\infty$ as a fixed point for the flow. In this spirit, we conclude from Th.1 in \cite{LiLl}:

\begin{corollary}
\label{fixedinf}    For every parameter $p\in P$, the Rössler system can be extended to a flow on $S^3$ - and moreover, $\infty$ becomes a fixed-point for the flow.
\end{corollary}

In the same spirit, we now prove the following fact about the local dynamics around $\infty$:

\begin{lemma}
    \label{fixedinf2}
For every $p\in P$, the index of the fixed-point at $\infty$ w.r.t. the vector field $F_p$ is $0$.
\end{lemma}

\begin{proof}
Recall $l_p=\{\dot{x}=0\}\cap\{\dot{y}=0\}$ (see Lemma \ref{lem23}). Given $r>0$ consider $S_r=\{|w|=r,w\in\mathbf{R}^3\}$ - by definition, $F_p$ can point on $S_r$ at the $(0,0,-d),d>0$ direction\textbf{ }only on the intersection $S_r\cap l_p$. Further recall we earlier parameterized $l_p$ by $\{(x,-\frac{x}{a},\frac{x}{a}),x\in\mathbf{R}\}$, and that by Eq.\ref{Field} we have $F_p(l_p(x))=(0,0,bx+\frac{x^2-xc}{a})$ (see also the proof of Lemma \ref{lem23}). We now note that for any sufficiently large $|x|$ the polynomial $bx+\frac{x^2-xc}{a}$ is positive - therefore, given any sufficiently large $r>0$ the vector field $F_p$ cannot point at the $(0,0,-d),d>0$ direction on $S_r\cap l_p$. Consequentially, we conclude that provided $r>0$ is sufficiently large the function $\frac{F_p}{||F_p||}:S_r\to S^2$ is not surjective - which implies the index of the vector field $F_p$ on the fixed-point at $\infty$ can only be $0$. All in all it follows that given any parameter $p\in P$, as we extend the vector field $F_p$ to $\infty$ as described above $\infty$ is added as a degenerate fixed point for the Rössler system and the assertion follows. 
\end{proof}

Lemma \ref{fixedinf2} allows us to treat $F_p$ as a vector field on the $3$-sphere $S^3$ with precisely three fixed points - two saddle foci $P_{In},P_{Out}$, and a degenerate fixed point at $\infty$ of index $0$. Since the parameter $p$ lies inside the parameter space $P$, we already know both $P_{In}$ and $P_{Out}$ are saddle-foci, hence hyperbolic - therefore, due to the Hartman-Grobman Theorem we can easily describe the local dynamics around both by reducing them to their linearizations. However, unlike $P_{In}$ and $P_{Out}$, $\infty$ is not a hyperbolic fixed-point - which implies the local dynamics around it are potentially much more complicated. In order to study them, we first prove:

\begin{lemma}\label{Clemma}
    For any parameter $p=(a,b,c)\in P$, the vector field $F_p$ is not $C^1$ at $\infty$. In particular, the Rössler system is not smooth at the fixed point at $\infty$.
\end{lemma}
\begin{proof}
For any $(x,y,z)=v\in\mathbf{R}^3$ consider the linearized vector field $J_p(v)$ at $v$ - namely, the Jacobian matrix of $F_p$ at $v\in\mathbf{R}^3$ (see Eq.\ref{Field}). Now, choose some $\lambda\in\mathbf{R}$ and note the equation $det(J_p(v))=\lambda$ (where $det$ denotes the determinant) can be rewritten as $z=\frac{\lambda+c-x-ab}{a}$. This proves the set $S_\lambda=\{det(J_p(v))=\lambda\}$ is a plane in $\mathbf{R}^3$ - which implies that when we extend $\mathbf{R}^3$ to $S^3$ we have $\infty\in\partial S_\lambda$ (with $\partial S_\lambda$ taken in $S^3$). Consequentially for any $\lambda_1\ne\lambda_2$ we have $\{\infty\}=\partial S_{\lambda_1}\cap \partial S_{\lambda_2}$ - which proves $det(J_p(v))$ takes infinitely many values at any neighborhood of $\infty$. As such $F_p$ cannot be $C^1$ at $\infty$ and the Lemma follows.
\end{proof}
Lemma \ref{Clemma} proves that we cannot extend the Rössler system to a smooth flow on $S^3$ - which, at least on the surface, could complicate matters. In order to overcome this difficulty, will now prove throughout our parameter space the vector field $F_p$ always generates two unbounded heteroclinic trajectories: $\Gamma_{In}$ and $\Gamma_{Out}$, which connect $P_{In}$ and $P_{Out}$ to $\infty$ (respectively). As we will see below, the existence of such unbounded heteroclinic trajectories would allow us to perturb $F_p$ around $\infty$ to a smooth vector field, with an arbitrarily small loss of dynamical data.\\

To begin, recall we denote by $W^s_{In}$ the one-dimensional, stable manifold of the saddle-focus $P_{In}$, and that $W^u_{Out}$ denotes the one-dimensional unstable manifold of $P_{Out}$ (see the illustration in Fig.\ref{loci}). We now prove the following result, with which we conclude this section:

\begin{theorem}\label{th21}
For every parameter $p\in P$ the corresponding Rössler system generates two heteroclinic trajectories:
\begin{itemize}
    \item $\Gamma_{In}\subseteq W^s_{In}$, which connects $P_{In},\infty$ in $S^3$.
    \item $\Gamma_{Out}\subseteq W^u_{Out}$, which connects $P_{Out}$, $\infty$ in $S^3$.
\end{itemize}
As a consequence, for every sufficiently large $r>0$, there exists a smooth vector vector field on $S^3$, $R_p$, s.t.:
\begin{itemize}
    \item $R_p$ coincides with $F_p$ on the open ball $B_r(P_{In})$.
    \item $R_p$ has precisely two fixed points in $S^3$ - namely, the saddle foci $P_{In}$ and $P_{Out}$.
    \item $R_p$ generates a heteroclinic trajectory which connects $P_{In},P_{Out}$ and passes through $\infty$.
\end{itemize}
\end{theorem}
\begin{proof}
We prove Th.\ref{th21} in three stages. The sketch of the proof is as follows:
\begin{itemize}
    \item In Stage $I$ we prove the existence of an invariant manifold $\Gamma_{In}$, a component of $W^s_{In}$, which forms a heteroclinic connection between $P_{In},\infty$. 
    \item In Stage $II$ we prove the existence of $\Gamma_{Out}$ - an analogous component of $W^u_{Out}$ which forms a heteroclinic connection between $P_{Out},\infty$.
    \item Finally, in Stage $III$ we tie these results together, and conclude the proof of Th.\ref{th21} by constructing the vector field $R_p$, as described above. 
\end{itemize}

\subsection{Stage $I$ - the existence of $\Gamma_{In}$.}
\label{stagei}
In this subsection we prove the existence of a component $\Gamma_{In}\subseteq W^s_{In}$ which forms an unbounded heteroclinic trajectory, connecting $P_{In}$ and $\infty$. Begin by recalling the cross-section $\{\dot{y}=0\}$, and that it is composed of three sets: the half-planes $U_p$ and $L_p$, separated by the line $l_p$ (see the discussion before Lemma \ref{obs}). Additionally, recall that on $L_p$ the vector field $F_p$ points into $\{\dot{y}>0\}$, while on $U_p$ it points into $\{\dot{y}<0\}$ (see Fig.\ref{planes} and Fig.\ref{H1}) - while $F_p$ is tangent to $\{\dot{y}=0\}$ on $l_p$. Recall $l_p(x)=(x,-\frac{x}{a},\frac{x}{a})$, $x\in\mathbf{R}$, and set $l_1=\{l_p(x)|x<0\}$. By $l_p(0)=P_{In}=(0,0,0)$, it follows $l_1$ is an unbounded arc in $\mathbf{R}^3$ which connects the saddle focus $P_{In}$ and $\infty$ (in $S^3$). \\

Now, consider the half-plane $H_1=\{(x,0,z)|x<0,z\in\mathbf{R}\}$ - since the velocity $\dot{y}$ is given by $x-ay$ where the parameter $a$ is positive (see Eq.\ref{Field} and the discussion at page \pageref{eq:9}), by definition we have $H_1\subseteq\{\dot{y}<0\}=\{(x,y,z)|x<-ay\}$. Since the normal vector to $H_1$ is $(0,1,0)$, by computation it immediately follows that for $v\in H_1$, the dot product $F_p(v)\bullet(0,1,0)$ is negative. Or, in other words, on $H_1$ the vector field $F_p$ points inside the open region $\{(x,y,z)|y<0\}$ (see the illustration in Fig.\ref{H1}).\\ 
\begin{figure}[h]
\centering
\begin{overpic}[width=0.6\textwidth]{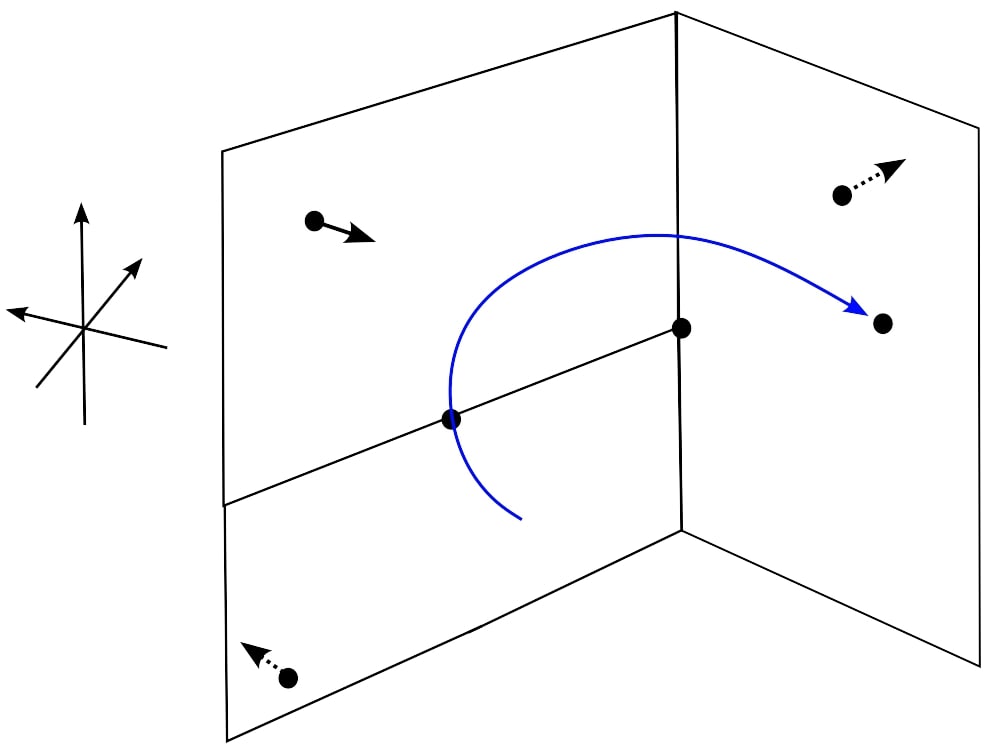}
\put(250,300){$l_1$}
\put(400,340){$s$}
\put(750,305){$H_1$}
\put(790,380){$\gamma_s(t(s))$}
\put(70,560){$z$}
\put(-15,445){$y$}
\put(150,495){$x$}
\put(155,690){$\{\dot{y}>0\}$}
\put(550,100){$\{\dot{y}<0\}$}
\put(450,570){$U_p$}
\put(350,190){$L_p$}
\put(570,440){$P_{In}$}
\put(620,530){$d_s$}
\end{overpic}
\caption[The quadrant $Q_1$.]{\textit{The quadrant $Q_1$, bounded by $H_1,U_p$ and $L_p$ (and $F_p's$ directions on them). The vector field is tangent to $l_1$, while $F_p$ points into $\{\dot{y}<0\}$ on $U_p$ and into $\{\dot{y}>0\}$ on $L_p$, and the flow line from $s\in l_1$ (denoted by $d_s$) connects to $\gamma_s(t(s))$.}}\label{H1}

\end{figure}
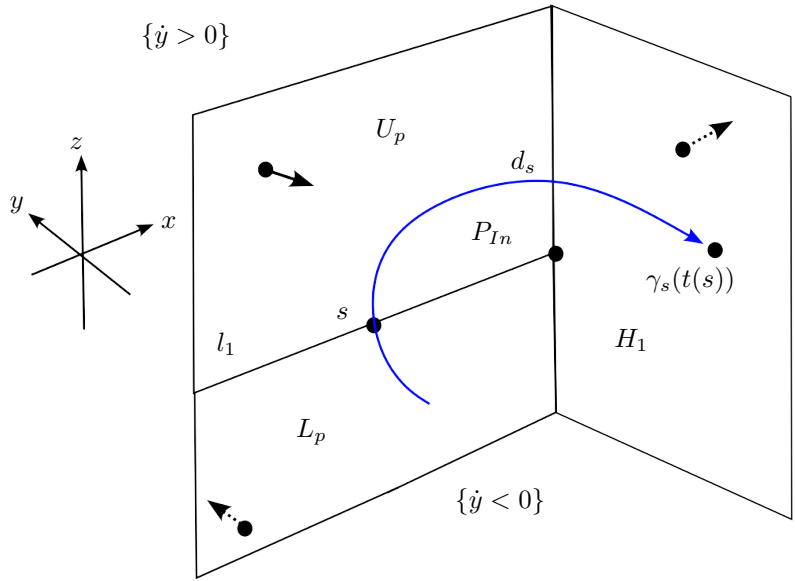

Let us note the two-dimensional set $H_1\cup\{\dot{y}=0\}$ traps a quadrant, $Q_1=\{\dot{y}\leq0\}\cap\{(x,y,z)|y>0\}$ (see the illustration in Fig.\ref{H1}) - by $l_1=\{l_p(x)|x<0\}$ and by $P_{Out}=(c-ab,\frac{ab-c}{a},\frac{c-ab}{a})$, $c-ab>0$ (see the discussion in page \pageref{eq:9}), we conclude both $l_1\subseteq Q_1$ and $P_{Out}\not\in\overline{Q_1}$.\\

Now, consider some $s\in l_1$ and its trajectory, $\gamma_s(t),t\in\mathbf{R}$ (recall we parameterize $\gamma_s$ s.t. $\gamma_s(0)=s$) - by Cor.\ref{cor211} the forward trajectory of $s$ enters the region  $\{\dot{y}<0\}$ immediately upon leaving $s$ (see the illustration in Fig.\ref{H1}). Moreover, let us remark that by definition, given $s\in l_1$ there exists some $x<0$ s.t. $s=(x,-\frac{x}{a},\frac{x}{a})$ - which, as the parameter $a$ is positive implies the $y-$coordinate of $s$ is also positive. Consequentially, it follows the forward trajectory of $s$ cannot be trapped in $Q_1$ forever - it either escapes $Q_1$ it by hitting $L_p$ transversely (i.e., entering the region $\{\dot{y}>0\}$), or by hitting $H_1$ transversely (i.e., entering the region $\{(x,y,z)|y<0\}$. Therefore, given $s\in l_1$, we conclude there exists a positive time $t(s)>0$ s.t. $\gamma_s(t(s))$ is the first intersection point between $H_1\cup{L_p}\cup l_1$ and the flow line $d_s$ connecting $\gamma_s(0)$ and $\gamma_s(t(s))$ (see the illustrations in Fig.\ref{H1} and Fig.\ref{branch}).\\

This motivates us to define $f:l_1\to H_1\cup L_p\cup l_1$ by $f(s)=\gamma_s(t(s))$ - by the discussion above, $f$ is well-defined. It is easy to see every component of $f(l_1)$ is a curve on the set $H_1\cup L_p\cup l_1$ (see Fig.\ref{branch} for an illustration). Now, set $V$ as the collection of flow lines $d_s$, $s\in l_1$, connecting $f(l_1)$ and $l_1$ - with these ideas in mind, we are now ready to give an overview for the proof of existence of $\Gamma_{In}$. In an ideal scenario, our tactic to prove the existence of the invariant manifold $\Gamma_{In}$ would be as follows:

\begin{itemize}
    \item First, we would like to define $C_{In}$ as the region trapped between $H_1,L_p$ and $V$ (see the illustration in Fig.\ref{branch}). Provided $C_{In}$ exists, it is easy to see it would form a topological cone with a tip at $P_{In}$.
    \item  Second, since by definition $\partial C_{In}\subseteq V\cup H_1\cup L_p$, for every $s\in\partial C_{In}$ it would follow $F_p(s)$ is either tangent to $\partial C_{In}$ (when $s\in V$) or points into either $\{\dot{y}>0\}$ or $\{(x,y,z)|y<0\}$ (when $s\in L_p$ or $H_1$, respectively). Since $Q_1$ is, by definition, a subset of $\{(x,y,z)|\dot{y}<0,y<0\}$ it would follow by $C_{In}\subseteq Q_1$ that whenever $s\in\partial C_{In}\setminus V$, $F_p(s)$ points outside of $Q_1$. As such, no trajectory can enter $C_{In}$ under the flow.
    \item Third, since by construction $C_{In}$ is a topological cone with a tip at $P_{In}$, by considering the linearization of $F_p$ around the saddle-focus $P_{In}$ we conclude it must include some invariant manifold $\Gamma_{In}$ - as such, since $C_{In}\subseteq Q_1\subseteq\{\dot{y}\leq0\}$ (and because no trajectory can escape $C_{In}$ under the inverse flow), by considering the backwards trajectory of $s\in\Gamma_{In}$, it would follow $\Gamma_{In}$ extends from $P_{In}$ to $\infty$.
\end{itemize}

\begin{figure}[h]
\centering
\begin{overpic}[width=0.3\textwidth]{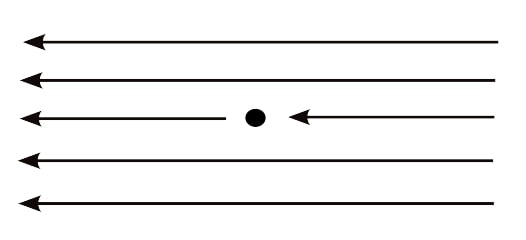}
\put(470,260){$\infty$}
\end{overpic}
\caption[An (almost) tubular flow at $\infty$.]{\textit{A flow around $\infty$ for the vector fields $F_r$ s.t. $\infty$ is a $0$-index fixed point.}}\label{LINF}

\end{figure}

In practice, the argument above has one major flaw - due to Lemma \ref{Clemma}, the dynamics around $\infty$ are hard to analyze, i.e., there is no reason to a-priori assume $f(l_1)$ includes an unbounded curve. As such, in order to overcome this difficulty, we first prove:
\begin{proposition}
        \label{CIN}
        For any $p\in P$ and any sufficiently large $r>0$ we can always deform the vector field $F_p$ to a vector field $F_r$ s.t.:

        \begin{itemize}
            \item $F_p$ and $F_r$ coincide on the ball of radius $r$ around the origin, $B_r(0)$.
            \item $F_r$ has precisely two fixed points in $\mathbf{R}^3$, the saddle-foci $P_{In}$ and $P_{Out}$ (in particular, $F_p$ and $F_r$ coincide around these fixed points).
            \item $F_r$ generates an invariant manifold $\Gamma_r$ for $P_{In}$, and moreover, $\Gamma_r$ is trapped inside $Q_1$ under the inverse flow generated by $-F_r$. 
        \end{itemize}
    \end{proposition}

In other words, Prop.\ref{CIN} implies that whatever the local dynamics of $F_p$ around the fixed point $\infty$ may be, in practice, we can always smoothly deform $F_p$ to another vector field, $F_r$ which satisfies the ideal scenario presented above (i.e., $F_r$ does generate a cone with a tip at $P_{In}$). In particular, the larger $r$ is, the more similar the dynamics of $F_r$ to those of $F_p$ (at least on compact sets in $\mathbf{R}^3$).
    
\begin{proof}
To begin, consider $r>0$ sufficiently large s.t. both fixed points $P_{In}$ and $P_{Out}$ are interior to the ball $B_r(0)$ - additionally, let us recall that by Cor.\ref{fixedinf2} the index of $F_p$ at its fixed point at $\infty$ is $0$ - hence the flow on $\mathbf{R}^3\setminus B_r(0)$ can be smoothly deformed to a tubular flow around a Fixed Point of Index $0$ at $\infty$ as in Fig.\ref{LINF}, without creating any new fixed points. Recall $l_1=\{l_p(x)|x<0\}$ (where $l_p(x)=(x,-\frac{x}{a},\frac{x}{a})$), i.e., it is a curve connecting the fixed point $P_{In}$ and $\infty$ through the quadrant $\overline{Q_1}$ - similarly, $\infty$ also lies on the closure of both $H_1$ and $L_p$ (where the closure in all cases is considered in the $3$-sphere $S^3$). We now define $F_r$ as the following vector field:

    \begin{itemize}
        \item $F_r$ coincides with $F_p$ on $B_n(0)$. Moreover, $F_r$ is a smooth vector field on $S^3$.
        \item $\infty$ is a fixed point of index $0$ for $F_r$, and the local dynamics around it are as in Fig.\ref{Qn} and Fig.\ref{LINF}.
        \item The sets $H_1,L_p, U_p$ and $l_1$ remain unchanged for $F_r$ - i.e., the local dynamics on $F_r$ on either one of these sets are orbitally equivalent to the local dynamics of $F_p$ on these sets. 
        \item Additionally, the regions $\{\dot{y}<0\}$ and $\{\dot{y}>0\}$ are the same for both $F_p$ and $F_r$ - that is, the quadrant $Q_1$ is trapped inside $\{\dot{y}\leq0\}$ w.r.t. both $F_p$ and $F_r$.
    \end{itemize}

It is easy to see the vector field $F_r$ satisfies the first two assertions of Prop.\ref{CIN} - therefore, to conclude the proof we must prove it also generates $\Gamma_r$, an invariant manifold trapped forever in $Q_1$ under the inverse flow. To do so, set $f_r:l_1\to H_1\cup L_p\cup l_1$ as the analogue (w.r.t. the vector field $F_r$) of the function $f$ defined above - i.e., $f_r$ is the first-hit map for initial conditions on $s\in l_1$ on $H_1\cup L_p\cup l_1$ w.r.t. $F_r$. When we smoothen $F_p$ around $\infty$ to create $F_r$ (as described in the proof of Cor.\ref{fixedinf2}), using Hopf's Theorem we deform $F_p$ to $F_r$ s.t. the flow moves around $\infty$ in a specific direction as described in Fig.\ref{Qn}.\\
 
 The constraints imposed on the flow by the behavior of $F_r$ on the cross-sections $H_1,L_p$ and $U_p$ now imply the trajectories of initial conditions $s\in l_1$ sufficiently close to $\infty$ collide with $H_1\cup L_p\cup l_1$ - which allows us to construct $F_r$ s.t. $\lim_{s\to\infty}f_r(s)=\infty$ (see the illustration in Fig.\ref{Qn}). Moreover, the same constraints imply that by slightly bumping flow lines emanating from $l_1$ outside $B_r(0)$ (if necessary) we ensure that $f_r$ is continuous at some neighborhood of $\infty$ in $l_1$ (see the illustration in Fig.\ref{Qn}).\\
    
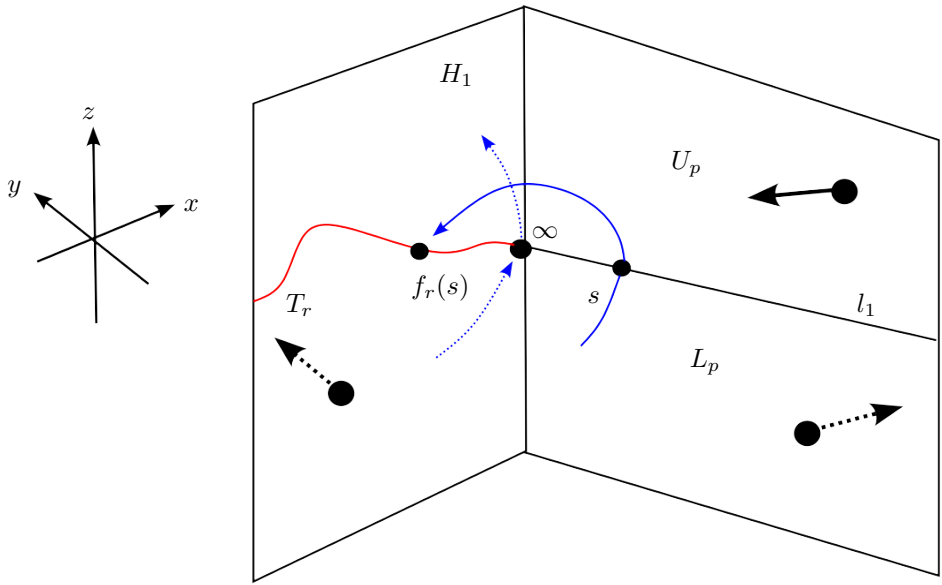
\begin{figure}[h]
\centering
\begin{overpic}[width=0.7\textwidth]{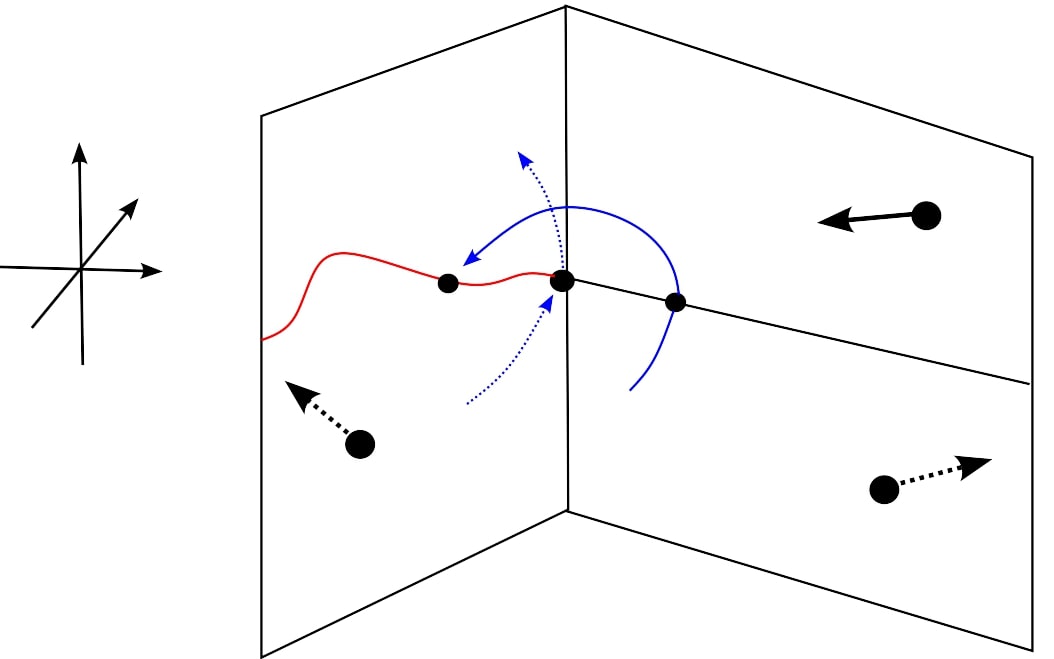}
\put(420,320){$f_r(s)$}
\put(700,455){$U_p$}
\put(285,300){$T_r$}
\put(65,495){$z$}
\put(130,445){$y$}
\put(160,360){$x$}
\put(450,550){$H_1$}
\put(610,310){$s$}
\put(720,240){$L_p$}
\put(900,300){$l_1$}
\put(550,380){$\infty$}
\end{overpic}
\caption[The curve $T_r$.]{\textit{The curve $T_r=f_r(l_1)$ and the directions of $F_r$ on $H_1,L_p$ and $U_p$ around $\infty$ (in $S^3$). As can be seen, $f_r$ is continuous around $\infty$. }}\label{Qn}

\end{figure}

Now, set $T_r=f_r(l_1)$ - it is easy to see every component of $T_r$ is a curve $H_1\cup L_p\cup l_1$ - and by the construction of $F_r$ above, as $\lim_{s\to\infty}f_r(s)=\infty$ at least one component of $T_r$ is a curve which stretches to $\infty$. Moreover, let us note that as $P_{In}$ is a saddle-focus with a two-dimensional unstable manifold $W^u_{In}$ transverse to $L_p$ at $P_{In}$ (see Lemma \ref{cor211}), it is easy to see $f_r$ is also continuous around $P_{In}$ and satisfies $\lim_{s\to P_{In}}f_r(s)=P_{In}$ (see the illustration in Fig.\ref{branch}). Now, set $V_r$ as the collection of flow-lines connecting $T_r$ and $l_1$ (see the illustration in Fig.\ref{branch}) - we claim:

\begin{lemma}
    \label{body}
    For every sufficiently large $r$, $F_r$ generates a three-dimensional body $C_{In}$ trapped between $V_r$, $L_p$ and $H_1$. Moreover, $P_{In}\in\partial C_{In}$.
\end{lemma}
\begin{proof}
Consider $f_r:l_1\to H_p\cup L_p\cup l_1$ - either $f_1$ is continuous on $l_1$ or it is not. If $f_1$ is continuous, as it is a first-hit map it follows by the Existence and Uniqueness Theorem, $T_r=f_r(l_1)$ is a simple curve in $H_p\cup L_p\cup l_1$, connecting $P_{In}$ and $\infty$. Whenever this is the case, it immediately follows there exists a three dimensional body $C_{In}$ trapped between $V_r, L_p$ and $H_1$ (see the illustration in Fig.\ref{branch}). Moreover, it is easy to see $P_{In}\in\partial C_{In}$.\\

    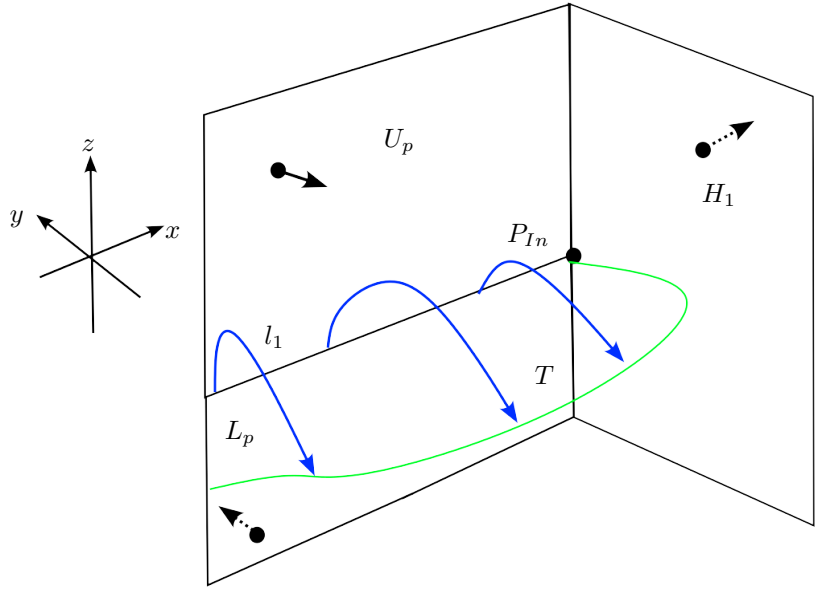
\begin{figure}[h]
\centering
\begin{overpic}[width=0.6\textwidth]{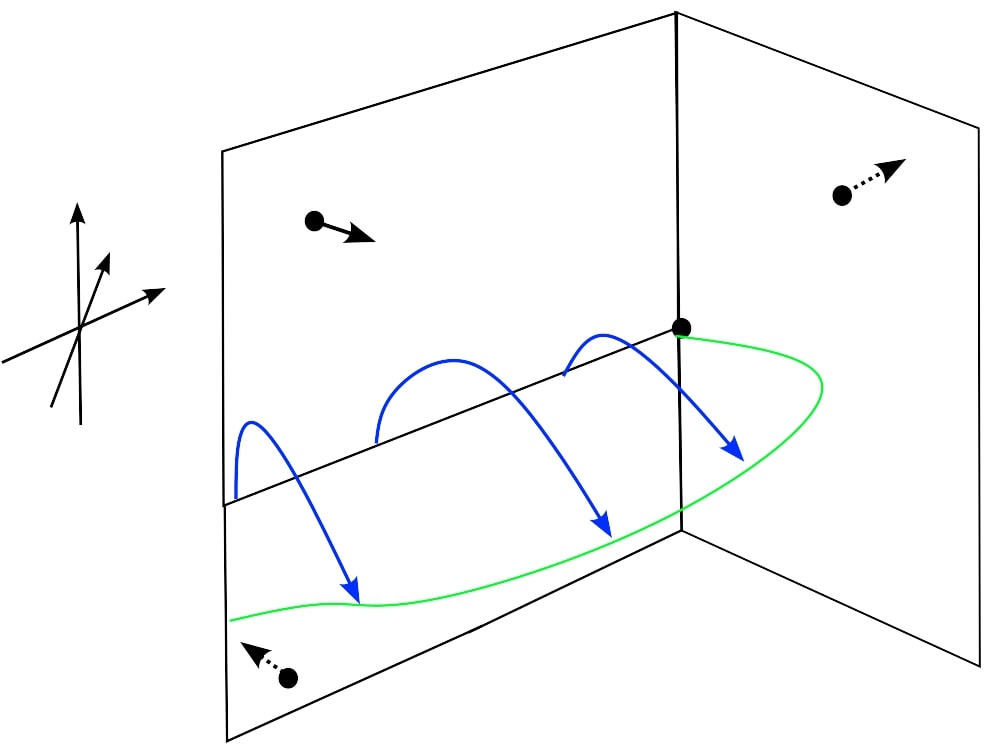}
\put(300,320){$l_1$}
\put(65,570){$z$}
\put(115,510){$y$}
\put(175,465){$x$}
\put(450,600){$U_p$}
\put(850,500){$H_1$}
\put(640,270){$T$}
\put(250,200){$L_p$}
\put(605,450){$P_{In}$}
\end{overpic}
\caption[The curve $T$, and the set $V$.]{\textit{The two-dimensional set $V$ (or $V_n$) is generated by the flow lines connecting $T$ (or $T_n$) and $l_1$. In this case, $T$ (or $T_n$) is a simple curve, connecting $P_{In}$ and $\infty$}.}\label{branch}

\end{figure}
    
Otherwise, $f_r$ is discontinuous at some $s\in l_1$ - which implies $T_r$ must include at least two components (as illustrated in Fig.\ref{Q112}). Therefore,  let $T^1$ denote the component of $T_r$ beginning at $P_{In}$ - since $F_r$ is transverse to both $H_1$ and $L_p$, it terminates at some $s_2\in l_1$. Or, in other words, there exists a sub-arc $J_1\subseteq l_1$, beginning at $P_{In}$ and terminating at some finite $s_1\in l_1$ s.t. $f_r(J_1)=T^1$ and $f_r(s_1)=s_2$ (as illustrated in Fig.\ref{branch}). This implies the flow lines connecting $J_1$ to $T^1=f_r(J_1)$ perform a wave-like movement, as sketched in Fig.\ref{branch} - which they perform s.t. $s_3=f_r(s_2)$ is strictly interior to $T^1$, as illustrated in Fig.\ref{Q11}. In particular, whenever $V_r$ is tangent to $l_1$ at $s_2$, the trajectories of initial conditions $s\in J_2$ (where $J_2=l_1\setminus J_1$) sufficiently close to $s_1$ slide below the surface connecting $J_1$ and $T^1$, as illustrated in Fig,\ref{Q11} - thus hitting $H_1\cup L_p$ inside the region enclosed between $T^1$ and $l_1$. Consequentially, given initial conditions $s\in J_2$ sufficiently close to $s_1$, $f_r(s)$ lies in the regions on $H_1\cup L_p$ trapped between $J_1$ and $T^1=f_r(J_1)$. Moreover, whenever this is the case $T_r$ is branched (as illustrated in Fig.\ref{Q11}). \\

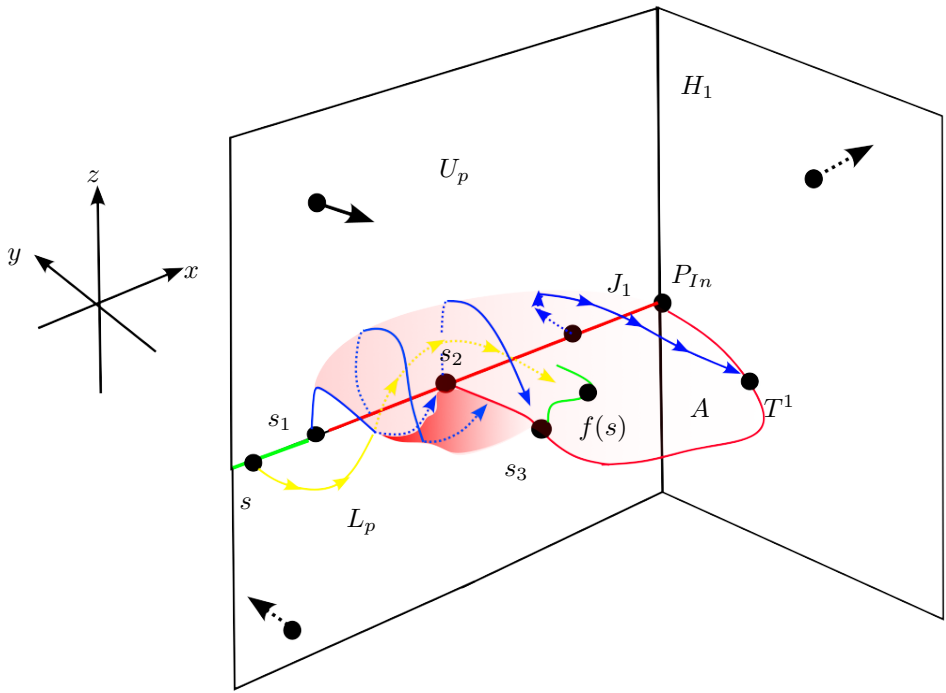
\begin{figure}[h]
\centering
\begin{overpic}[width=0.7\textwidth]{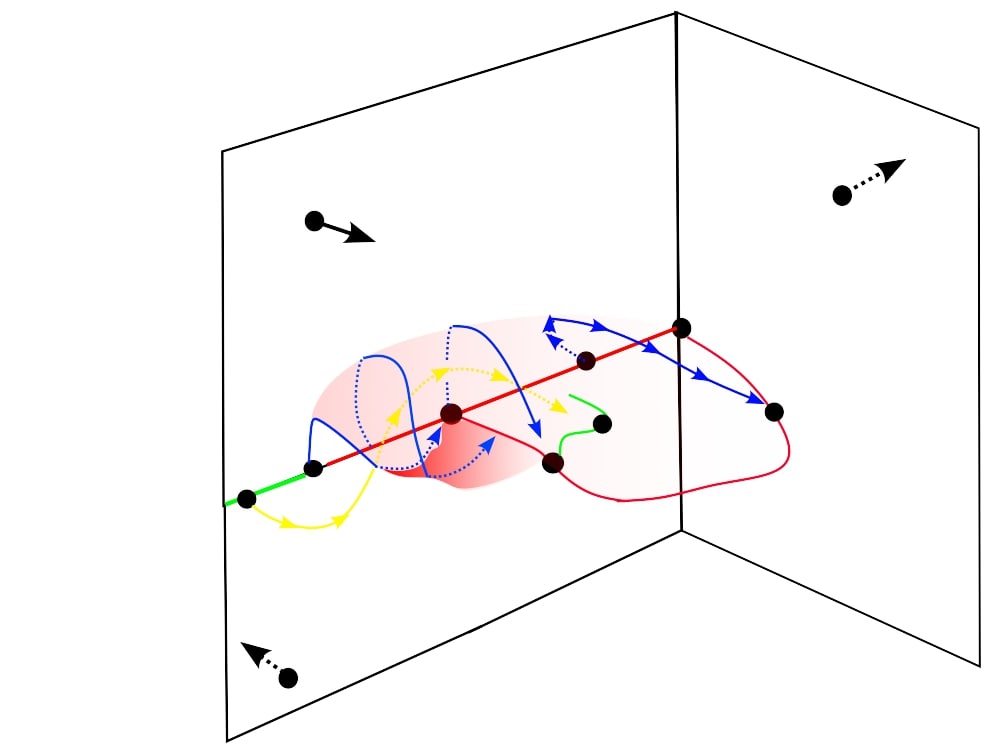}
\put(450,370){$s_2$}
\put(520,245){$s_3$}
\put(700,455){$P_{In}$}
\put(630,445){$J_1$}
\put(710,660){$H_1$}
\put(265,300){$s_1$}
\put(235,210){$s$}
\put(450,570){$U_p$}
\put(350,190){$L_p$}
\put(600,290){$f(s)$}
\put(800,310){$T^1$}
\put(720,310){$A$}
\end{overpic}
\caption[The surface $V_n$ in case of tangency.]{\textit{The red surface denotes the collection of flow lines connecting $J_1$ and $T^1$. Initial conditions $s$ on $J_2$ (the green arc) sufficiently close to $s_1$ flow along the flow line of $s_1$ below the red surface, thus hitting $A$ (as indicated by the yellow flow line, where $s_2=f_r(s_1)$, $s_3=f_r(s_2)$). $A$ denotes the two-dimensional region enclosed by $T^1$ and $l_1$ on $H_1\cup L_p$.}}\label{Q11}

\end{figure}

Now, set $A$ as the collection of two-dimensional domains on $L_p\cup H_1$ trapped between $T_r$ and $l_1$ (see the illustration in Fig.\ref{Q11} and \ref{Q112}). Using a similar argument to the one above it follows that given an initial condition $s\in l_1$ s.t. $f_r(s)\not\in\partial A$, then $f_r(s)$ is interior to $ A$. Hence, given an open arc of initial conditions $J$ on $l_1$ s.t. $f_r(J)\not\subseteq\partial A$, $f_r(J)$ is a collection of arcs in $A$, each with two endpoints on $\partial A$. To continue, recall that by our construction of $F_r$, $J$ is always bounded in $\mathbf{R}^3$ - therefore, the region trapped between $T_r$ and the flow lines connecting $J,f_r(J)$ is a topological tube with one opening inside $Q_1$ and another opening inside $A$ (as sketched in Fig.\ref{Q112}) - moreover, the interior of the said tube is a component of $Q_1\setminus V_r$.\\

As such, it follows $Q_1\setminus V_r$ includes a component $C_{In}$, satisfying $P_{In}\in\partial C_{In}$, which lies away from all such tubes - in particular, $C_{In}$ is trapped between $V_r,L_p$ and $H_1$. Moreover, since $f_r$ is continuous on initial conditions on $l_1$ that are sufficiently large, it follows $C_{In}$ is also unbounded. All in all, regardless of whether $f_r$ is continuous on $T_r$ or not, there always exists a connected three-dimensional region, $C_{In}$, trapped between $V_r$ and $L_p\cup H_1\cup l_1$ s.t. $P_{In}\in\partial C_{In}$. The proof of Lemma \ref{body} now follows. 
\end{proof}

\begin{figure}[h]
\centering
\begin{overpic}[width=0.7\textwidth]{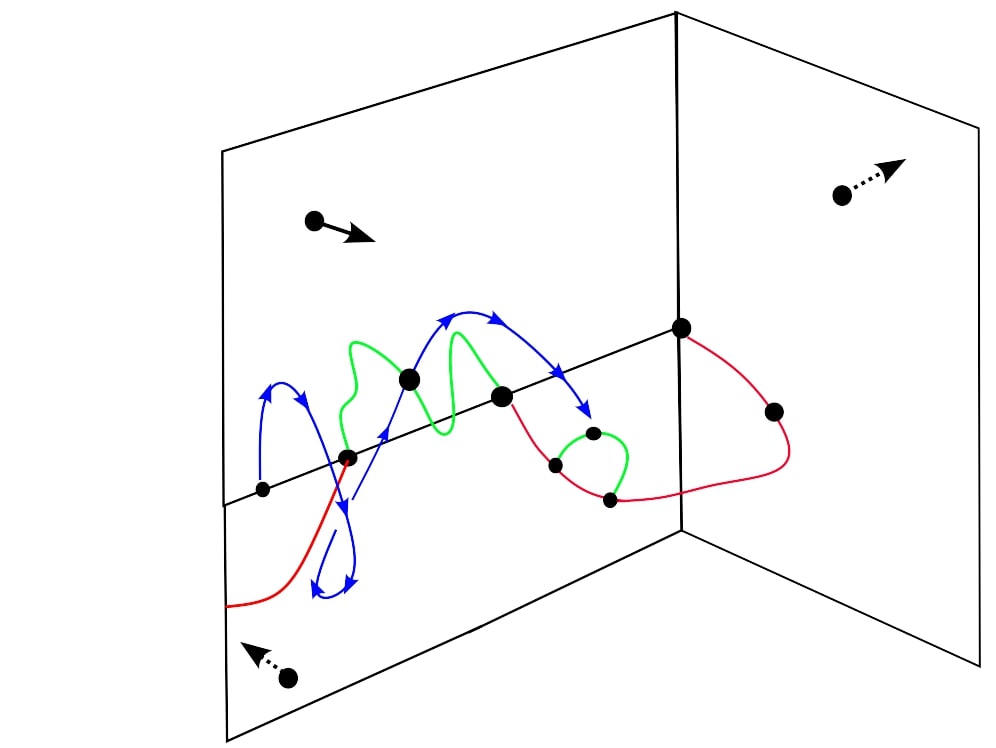}
\put(700,455){$P_{In}$}
\put(710,660){$H_1$}
\put(265,300){$s$}
\put(235,210){$A$}
\put(450,570){$U_p$}
\put(380,190){$L_p$}
\put(600,350){$f(s)$}
\put(800,310){$T^1$}
\put(720,310){$A$}
\put(620,430){$J_1$}
\end{overpic}
\caption[The case when $T_r$ has several components.]{\textit{In this scenario $T_r$ has two components, where the regions $A$ are trapped between $l_1$ and the components of $T_r$. Due to a tangency at the green curve, the points which flow to the green arc eventually flow into $A$, thus generating a tube connecting the two green arcs (and causing a branching at $T_r$ - as exemplified by the trajectory of $s$). }}\label{Q112}

\end{figure}
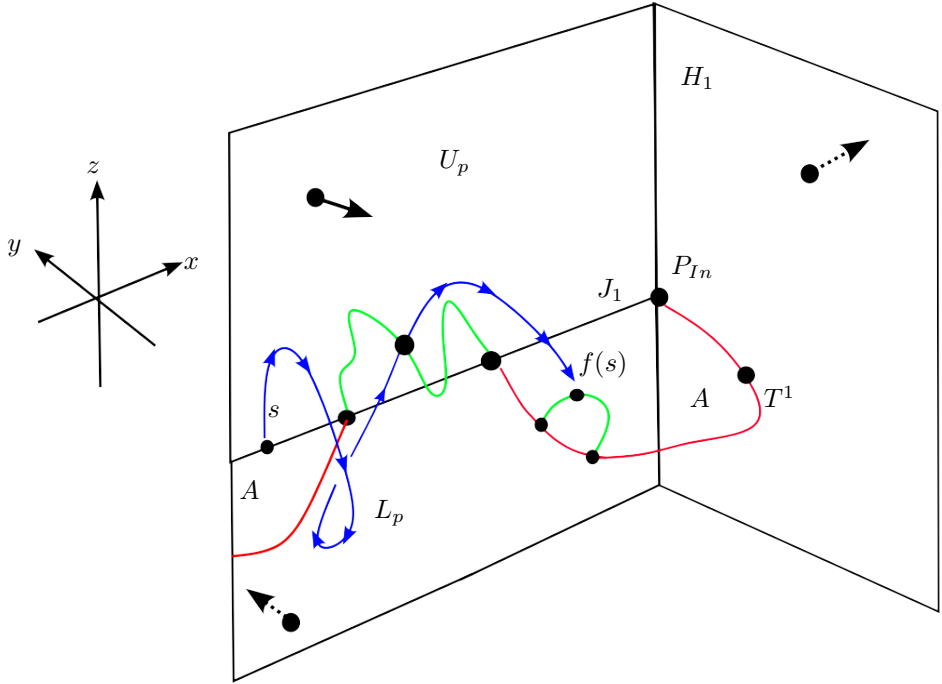

We are now ready to conclude the proof of Prop.\ref{CIN}. To do so, first note that given any $F_r$ (as defined at the beginning of the proof), by construction it follows the the three-dimensional body $C_{In}$ given by Lemma \ref{body} is trapped inside the quadrant $Q_1\subseteq\{(x,y,z)|y\leq0\}$. As a consequence, by $P_{Out}=(c-ab,\frac{ab-c}{a},\frac{c-ab}{a})$, $c>1,a,b\in(0,1)$ it follows $P_{Out}\not\in C_{In}$. Furthermore by construction $\partial C_{In}$ is a subset of the union $V_r\cup H_1\cup L_p$. Now, note that by definition, $F_r$ is tangent to $V_r$, and recall $F_r$ and $F_p$ are orbitally equivalent around $H_1$ and $L_p\subseteq\{\dot{y}=0\}$ - as such, since on $H_1$ and $L_p$ the vector field $F_p$ points into $\{(x,y,z)|y<0\}$ and $\{\dot{y}>0\}$ (respectively), the same is true for $F_r$ on $L_p$ and $H_1$. Therefore, since $Q_1\subseteq \{\dot{y}\leq0\}\cap\{(x,y,z)|y<0\}$ and because $\partial C_{In}\setminus V_r\subseteq L_p\cup H_1$, we conclude that for $s\in\partial C_{In}\setminus V_r$, $F_r(s)$ points out of $Q_1$ - and hence, outside of $C_{In}$ as well.\\

Therefore, since $\partial C_{In}\subseteq V_r\cup L_p\cup H_1$ it follows that for any $s\in\partial C_{In}$ $F_r(s)$ is either tangent to $\partial C_{In}$ (when $s\in V_r$), or points outside of $C_{In}$ (when $s\not\in V_r$). We may summarize our findings as follows:

\begin{itemize}
    \item $P_{Out}\not\in\overline{C_{In}}$ while $P_{In}\in\partial C_{In}$.
    \item $C_{In}\subseteq\{\dot{y}\leq0\}$, where $\dot{y}$ is taken w.r.t. $F_r$.
    \item Since $F_r$ is either tangent to $\partial C_{In}$ or points out of $C_{In}$ throughout $\partial C_{In}$, no trajectory can escape $\overline{C_{In}}$ w.r.t. $-F_{n}$.
\end{itemize}

It follows $C_{In}$ forms a topological cone with a tip at $P_{In}$ - therefore, by considering the linearization of $F_r$ it follows there must exist an invariant manifold $\Gamma_r$ for $P_{In}$ inside $C_{In}$ emanating from $P_{In}$ (w.r.t. the vector field $F_r$) - see Fig.\ref{cone1}. Moreover, if $\phi^r_t$ denotes the flow generated by $F_r$, by the discussion above we conclude that given any $s\in\Gamma_r$ and every $t<0$, $\phi^r_t(s)\in C_{In}\subseteq Q_1$. By this discussion we conclude $\Gamma_r$ is trapped in $Q_1$ under the inverse flow generated by $-F_r$, and the proof of Prop.\ref{CIN} is now complete.
\end{proof}

\begin{figure}[h]
\centering
\begin{overpic}[width=0.3\textwidth]{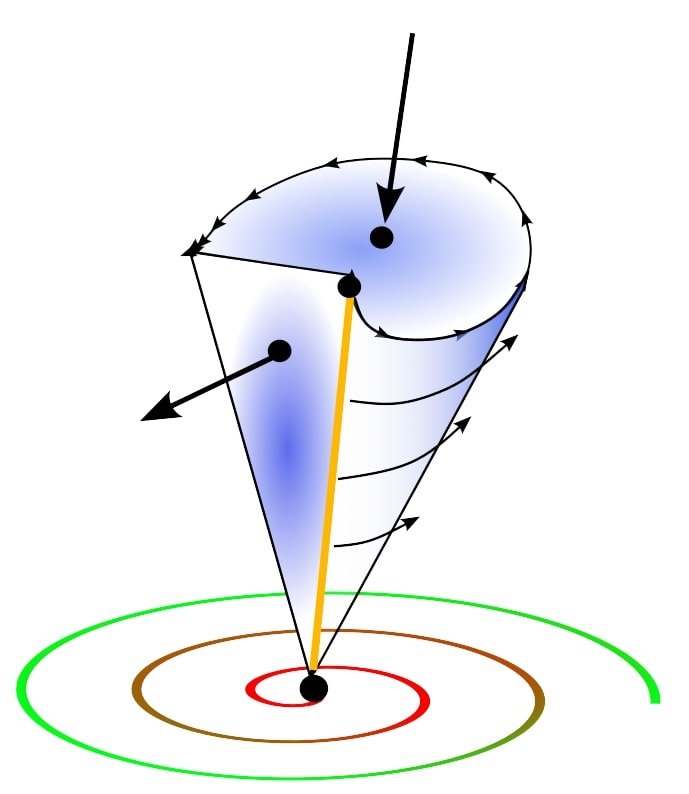}

\put(-50,40){$W^u_{In}$}
\put(650,450){$C$}
\put(540,920){$\Gamma_{n}$}
\put(425,60){$P_{In}$}
\end{overpic}
\caption[Generating a cone around $P_{In}$.]{\textit{By the Hartman-Grobman Theorem, by suspending a $l_1$ with the flow (the orange curve) we generate a topological cone into which no trajectory can enter. As such, it includes an invariant manifold of $P_{In}$ (w.r.t. $F_r$).}}
\label{cone1}
\end{figure}
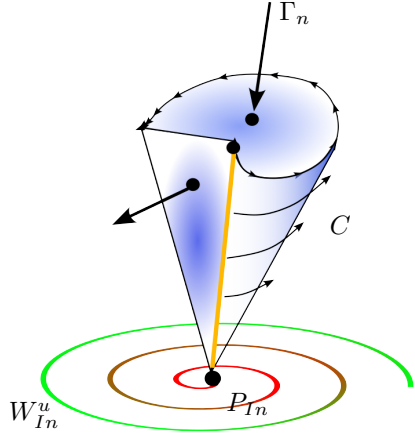

We now use Prop.\ref{CIN} to conclude the proof of Stage $I$ - namely, now prove the existence of an invariant, unbounded, one dimensional manifold $\Gamma_{In}$ for $P_{In}$ which we do in the following corollary:

\begin{corollary}\label{cor214}
    There exists a component $\Gamma_{In}\subseteq W^s_{In}$ - where $W^s_{In}$ is the one-dimensional, stable manifold for the saddle focus $P_{In}$ - s.t. $\Gamma_{In}$ is a heteroclinic trajectory for $F_p$ in $S^3$, connecting the fixed-points $P_{In}$ and $\infty$.
\end{corollary}
\begin{proof}
  Choose a sufficiently large $r>0$ s.t. Prop.\ref{CIN} holds, and consider an initial condition $s\in\Gamma_r\cap B_r(0)$. As $F_p$ and $F_r$ coincide on $B_r(0)$ and because $P_{In}\in B_r(0)$, provided $s$ is sufficiently close to $P_{In}$, $s$ is in some invariant manifold $\Gamma_{In}$ of $P_{In}$ w.r.t. to the original vector field $F_p$. Now, denote the flow generated by $F_p$ by $\phi^p_t$, and consider $\phi^p_t(s)$, for any $t<0$. Since $r$ can be chosen to be arbitrarily large and because $\Gamma_r\subseteq Q_1\subseteq\{\dot{y}\leq0\}$ for all $r$, we conclude that for every $t<0$ $\phi^p_t(s)\in\{\dot{y}\leq0\}$ - which implies $\Gamma_{In}$ is forever trapped in both $Q_1$ and $\{\dot{y}\leq0\}$. Recalling the two-dimensional unstable manifold $W^u_{In}$ is transverse to $\{\dot{y}=0\}$ at $P_{In}$ (see Cor.\ref{cor211}), it follows the backwards trajectory of every initial condition on $W^u_{In}$ fluctuates infinitely many times between $\{\dot{y}>0\}$ and $\{\dot{y}<0\}$ - therefore, $\Gamma_{In}\cap W^u_{In}=\emptyset$. Consequentially, $\Gamma_{In}$ can only be a component of the one-dimensional invariant manifold $W^s_{In}$.\\
    
Furthermore, since $P_{Out}=(c-ab,\frac{ab-c}{a},\frac{c-ab}{a})$, $c-ab>0$ (see the definition of the parameter space $P$ at page \pageref{eq:9}), by $Q_1\subseteq\{(x,y,z)|x\leq0\}$ it follows $P_{Out}\not\in\overline{Q_1}$. Therefore, given any $s\in\Gamma_{In}$, $P_{Out}\ne\lim_{t\to-\infty}\phi^p_t(s)$, i.e., $\Gamma_{In}$ is not a heteroclinic trajectory connecting $P_{In}$ and $P_{Out}$. Finally, since for every $s\in\Gamma_{In}$ and every $t<0$ we have $\phi^p_t(s)\in\{\dot{y}\leq0\}$, as $\Gamma_{In}$ is not a heteroclinic trajectory it follows $\lim_{t\to-\infty}\phi^p_t(s)=\infty$ and Cor.\ref{cor214} now follows.
\end{proof}
\subsection{Stage $II$ - the existence of $\Gamma_{Out}$.}
\label{stageii}
Having proven the existence of $\Gamma_{In}$, we now prove the analogous result for $P_{Out}$ - namely, we now prove the existence of a heteroclinic trajectory $\Gamma_{Out}$, a component of the one-dimensional unstable manifold $W^u_{Out}$, which connects $P_{Out},\infty$. Despite some slight differences, in practice the proof will be almost symmetric to the one in Stage $I$.\\

To begin, recall $P_{Out}=(c-ab,\frac{ab-c}{a},\frac{c-ab}{a})$, where $c-ab>0$ (see page \pageref{eq:9}) and that $l_p(c-ab)=P_{Out}$ - where $l_p$ is the tangency curve of the vector field $F_p$ to the plane $\{\dot{y}=0\}$, parameterized by $l_p(x)=(x,-\frac{x}{a},\frac{x}{a})$, $x\in\mathbf{R}$ (see the illustration in Fig.\ref{planes} and Fig.\ref{loci}). Now, set $l_2=\{l_p(x)|x>c-ab\}$ (see the illustration in Fig.\ref{H2}) - we prove the existence of $\Gamma_{Out}$ by making an analogous argument to the one used to prove the existence of $\Gamma_{In}$. To begin, consider the half-plane $H_2=\{(x,\frac{ab-c}{a},z)|x<c-ab,z\in\mathbf{R}\}$ (see the illustration in Fig.\ref{H2}). Again, by $\dot{y}=x+ay$ it follows, $H_2\subseteq\{\dot{y}<0\}$ - and since the normal vector to $H_2$ is $(0,1,0)$ by considering $F_p(x,\frac{ab-c}{a},z)\bullet(0,1,0)=x+ab-c$ it immediately follows that on $H_2$ the vector field $F_p$ points into the region $\{(x,y,z)|y<c-ab\}$ (see the illustration in Fig.\ref{H2}).\\

\begin{figure}[h]
\centering
\begin{overpic}[width=0.6\textwidth]{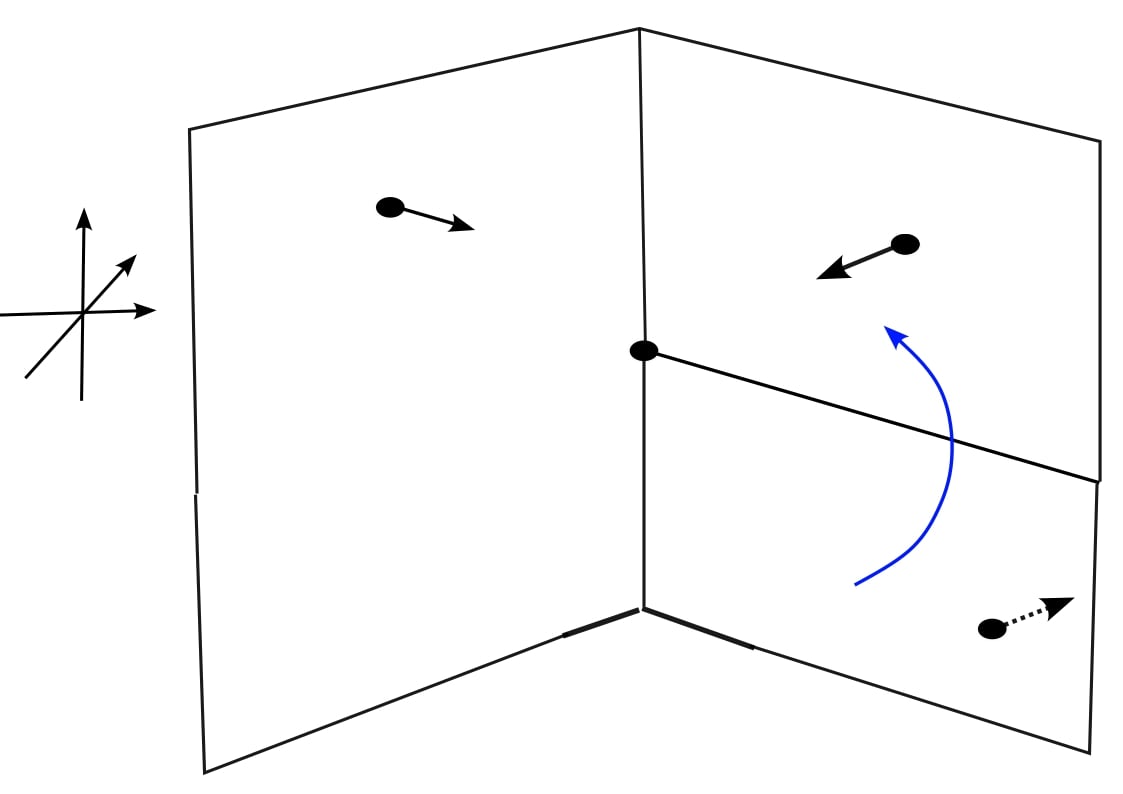}
\put(900,315){$l_2$}
\put(230,440){$H_2$}
\put(780,145){$L_p$}
\put(65,530){$z$}
\put(145,420){$y$}
\put(125,480){$x$}
\put(650,500){$U_p$}

\put(600,410){$P_{Out}$}
\end{overpic}
\caption[The quadrant $Q_2$.]{\textit{The quadrant $Q_2$ along with a flow line tangent to $l_2$ (and the directions of $F_p$ on $U_p,L_p$ and $H_2$). }}\label{H2}

\end{figure}
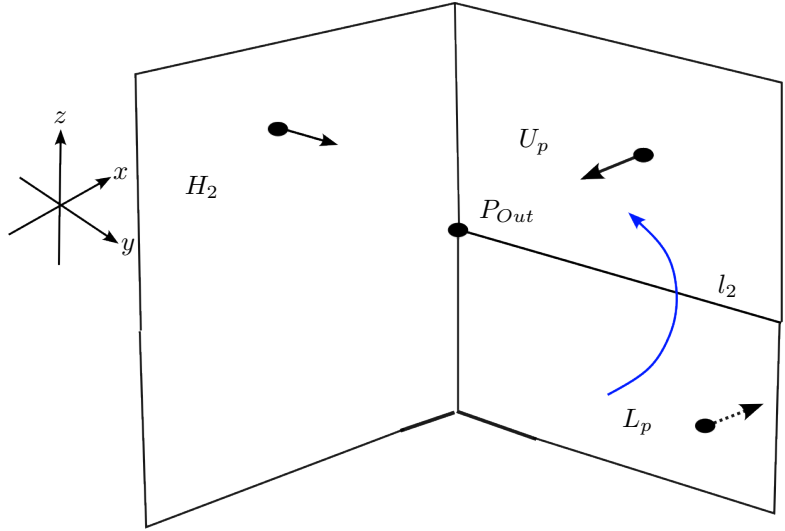

Consequentially, $H_2\cup\{\dot{y}=0\}$ traps a quadrant $Q_2=\{\dot{y}\leq0\}\cap\{(x,y,z)|y<\frac{ab-c}{a}\}$ (see the illustration in Fig.\ref{H2}) - by $l_2=\{l_p(x)|x>c-ab\}$, we conclude $l_2\subseteq \partial Q_2$.  Now, consider some $s\in l_2$ and its backwards trajectory. Recall that by Lemma \ref{cor211} the backwards trajectory of $s$ enters $\{\dot{y}<0\}$ immediately upon leaving $s$ (see the illustration in Fig.\ref{H2}) - moreover, let us remark that by definition, given $s\in l_2$, it is given by the coordinates $s=(x,-\frac{x}{a},\frac{x}{a})$ (for some $x>c-ab$) - hence, its $y-$coordinate is greater than $\frac{ab-c}{a}$. As such, the backwards trajectory of $s$ cannot be trapped in $Q_2$ forever - i.e., it either escapes by hitting the half-plane $U_p$ transversely (i.e., entering $\{\dot{y}>0\}$ - see the illustration in Fig.\ref{H2}), or by hitting $H_2$ transversely (i.e., entering $\{(x,y,z)|y>\frac{ab-c}{a}\}$).\\

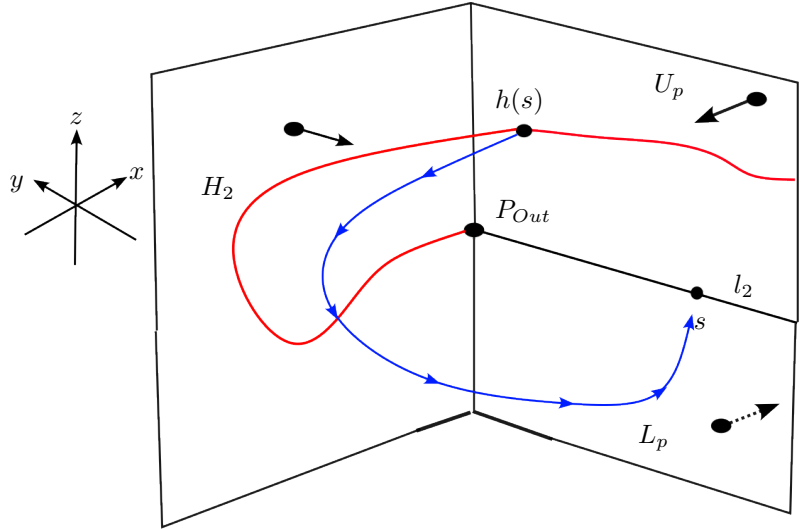
\begin{figure}[h]
\centering
\begin{overpic}[width=0.6\textwidth]{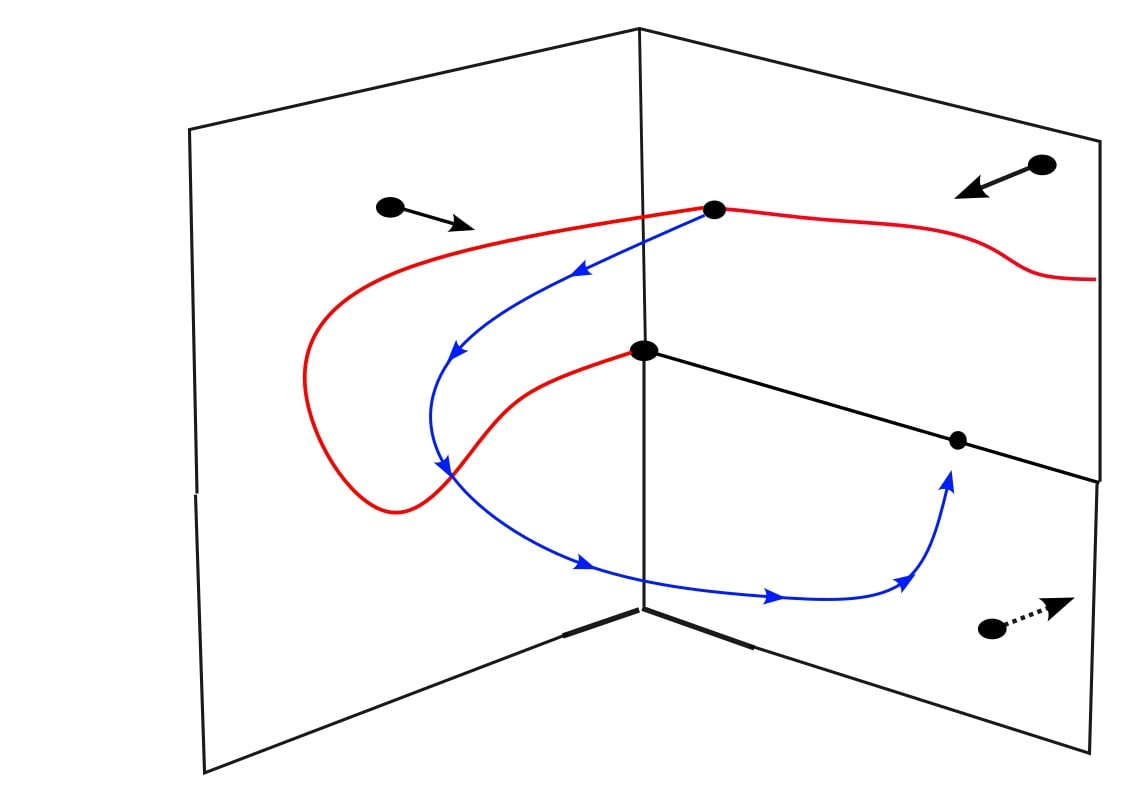}
\put(900,315){$l_2$}
\put(230,440){$H_2$}
\put(780,125){$L_p$}
\put(600,550){$h(s)$}
\put(800,570){$U_p$}
\put(850,270){$s$}
\put(600,410){$P_{Out}$}
\end{overpic}
\caption[The curve $J$.]{\textit{The quadrant $Q_2$ along with a flow line connecting a point $h(s)=\gamma_s(j(s))$ (or conversely, $h_n(s)$) in $J=h(l_2)$, the red curve, and $s\in l_2$. In this scenario $J$ is homeomorphic to an unbounded straight line. }}\label{Q2}

\end{figure}

    Therefore, given $s\in l_2$, there exists a strictly negative time $j(s)<0$ s.t. $\gamma_s(j(s))$ is the first intersection point between the two-dimensional set $H_2\cup U_p\cup l_2$ and the backwards trajectory of $s$ (see the illustration in Fig,\ref{H2}). Now, define $h:l_2\to H_2\cup U_p\cup l_2$ by $h(s)=\gamma_s(j(s))$ - that is, $h$ is the first-hit map w.r.t. the inverse flow for initial conditions in $l_2$ in $H_2\cup U_p\cup l_1$. Therefore, replacing the flow generated by $F_p$ with the inverse flow generated by $-F_p$, and replacing $f$ and $f(l_1)$ with $h$ and $h(l_2)$, using similar arguments to those used to prove Prop.\ref{CIN} we conclude :
\begin{proposition}
        \label{COUT}
        For any $p\in P$ and every sufficiently large $r>0$ we can always deform the vector field $-F_p$ to a vector field $G_r$ s.t.:
        \begin{itemize}
            \item $-F_p$ and $G_r$ coincide on $B_r(0)$, the ball of radius $r$ around the origin.
            \item $G_r$ has precisely two fixed points in $\mathbf{R}^3$ - the saddle-foci $P_{In}$ and $P_{Out}$ as saddle-foci (in particular, $-F_p$ and $G_r$ coincide around the fixed points).
            \item  $G_r$ generates an invariant manifold $\Gamma_r$ for $P_{Out}$ which is trapped in $Q_2$ under the flow generated by the vector field $-G_r$.
        \end{itemize}
    \end{proposition}
    \begin{figure}[h]
\centering
\begin{overpic}[width=0.3\textwidth]{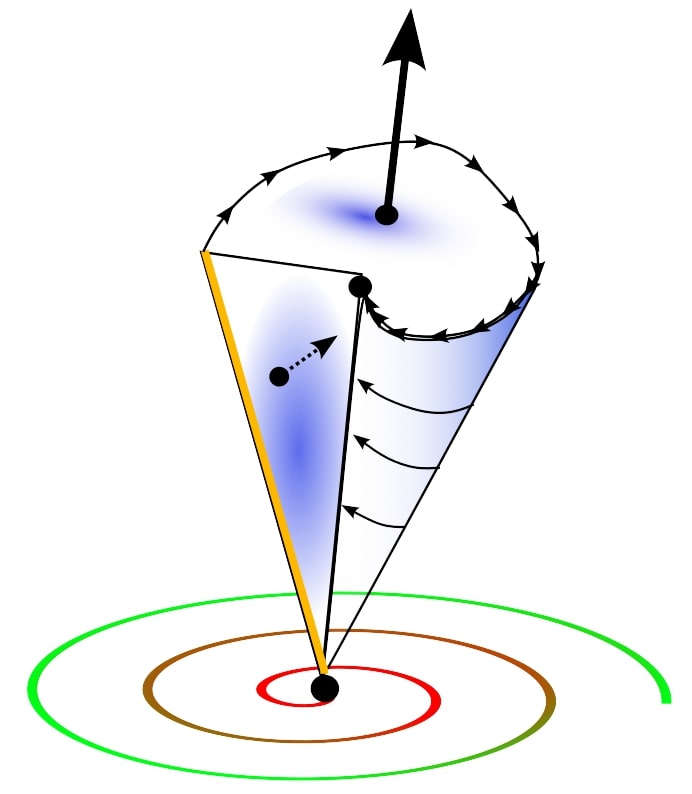}

\put(-50,40){$W^s_{Out}$}
\put(540,920){$\Gamma_n$}
\put(425,60){$P_{Out}$}
\end{overpic}
\caption[Generating a cone around $P_{Out}$.]{\textit{Similarly to the argument used to prove Prop.\ref{CIN}, by suspending a $l_2$ (the orange curve) with $G_r$ and considering the Hartman-Grobman Theorem we generate a topological cone with a tip at $P_{Out}$ which traps an invariant manifold of $P_{Out}$.}}
\label{cone2}
\end{figure}
We are now ready to conclude Stage $II$ of the proof - i.e., the prove the existence of $\Gamma_{Out}$, which we do using a symmetric argument to the one used to prove Cor.\ref{cor211}. Namely, replacing $F_p$ with $-F_p$ and $F_r$ with $G_r$ in the proof of Cor.\ref{cor211}, we conclude:
\begin{corollary}\label{cor215}
    There exists a component $\Gamma_{Out}$ of the one-dimensional, unstable invariant manifold of $P_{Out}$, $W^u_{Out}$, s.t. $\Gamma_{Out}$ is a heteroclinic trajectory connecting $P_{Out},\infty$.
\end{corollary}

\subsection{Stage $III$ - constructing $R_p$ and concluding the proof.}
\label{stageiii}
Having proven the existence of $\Gamma_{In},\Gamma_{Out}$, invariant, unbounded one-dimensional manifolds in $W^s_{In},W^u_{Out}$ (respectively), we can now conclude the proof of Th.\ref{th21}. That is, we will now use the fact ${\Gamma_{Out}}$ and ${\Gamma_{In}}$ connect at $\infty$ to prove the existence of a smooth vector field on $S^3$, $R_p$, s.t. the following conditions are satisfied:
\begin{itemize}
    \item For every sufficiently large $r>0$, we can construct $R_p$ s.t. it coincides on $\{(x,y,z)|x^2+y^2+z^2<r\}$ - where $F_p$ is the original vector field corresponding to the Rössler system at $p$ (see Eq.\ref{Field}). 
    \item $R_p$ has precisely two fixed points in $S^3$ - $P_{In},P_{Out}$.
    \item $R_p$ generates an unbounded heteroclinic trajectory, connecting $P_{In},P_{Out}$.
\end{itemize}

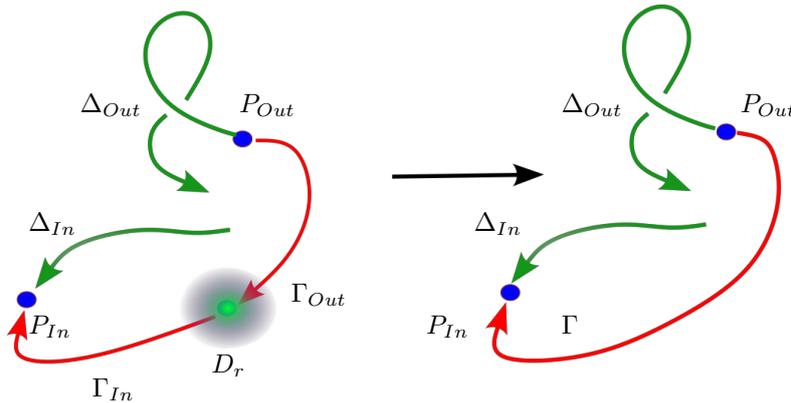
\begin{figure}[h]
\centering
\begin{overpic}[width=0.6\textwidth]{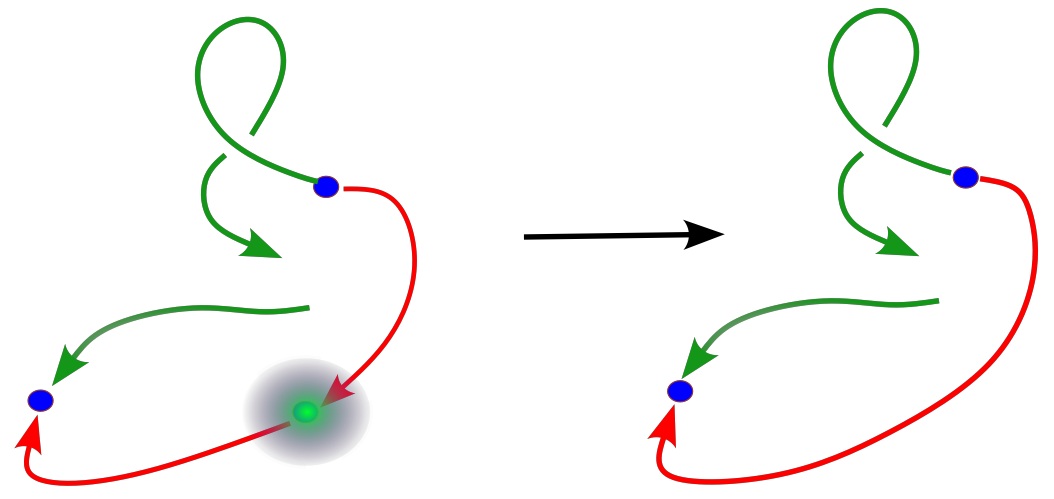}
\put(120,-30){$\Gamma_{In}$}
\put(270,0){$D_r$}
\put(370,90){$\Gamma_{Out}$}
\put(540,50){$P_{In}$}
\put(710,50){$\Gamma$}
\put(710,330){$\Delta_{Out}$}
\put(305,330){$P_{Out}$}
\put(105,330){$\Delta_{Out}$}
\put(935,330){$P_{Out}$}
\put(40,50){$P_{In}$}
\put(40,180){$\Delta_{In}$}
\put(600,180){$\Delta_{In}$}
\end{overpic}
\caption[The deformation of $F_p$ to $R_p$.]{\textit{The smooth deformation of $F_p$ to $R_{p}$ - $F_p$ remains unchanged in $\mathbf{R}^3\setminus D_r$. $\Delta_{In},\Delta_{Out}$ denote the respective components of $W^s_{In}\setminus\Gamma_{In},W^u_{Out}\setminus\Gamma_{Out}$. $D_r$ is the region $\{(x,y,z)|x^2+y^2+z^2\geq r^2\}$.}}
\label{deform21}
\end{figure}

To begin, first recall that by Lemma \ref{fixedinf2} the index of the vector field $F_p$ at $\infty$ is $0$ - hence, by Hopf's Theorem (see pg. 51 in \cite{Mil}) for every sufficiently large $r$, $F_p$ is homotopic on the set $D_r=\{w\in\mathbf{R}^3|||w||\geq r\}$ to a vector field which generates a tubular flow (that is, a flow given by the vector field $F(x,y,z)=(1,0,0)$). Therefore, applying Hopf's Theorem, for any sufficiently large $r$ we now smoothly deform $F_{p}$ inside $D_r$ by some smooth homotopy which removes the fixed point at $\infty$ - as the index at $\infty$ is $0$ we can do so s.t. no new fixed points are generated in $D_r$. More importantly, because $\infty\in\overline{\Gamma_{In}}\cup \overline{\Gamma_{Out}}$ (where the closure is taken at $S^3$) we construct this perturbation s.t. $\Gamma_{In},\Gamma_{Out}$ connect to a heteroclinic connection passing through $\infty$ - the curve $\Gamma$ (see the illustration in Fig.\ref{deform21}).\\

Now, denote the vector field constructed above by $R_p$ - by construction it is a smooth vector field in $S^3$, and by the paragraphs above it satisfies the following:
\begin{itemize}
    \item Provided our initial $r$ is sufficiently large, $R_p$ coincides with the vector field $F_p$ (i.e., with the Rössler system) on $\{(x,y,z)|x^2+y^2+z^2>r\}$. 
    \item $R_p$ has precisely two fixed points in $S^3$ - $P_{In}$ and $P_{Out}$, both saddle foci (of opposing indices).
    \item $P_{In}$ and $P_{Out}$ are connected by $\Gamma$, a heteroclinic trajectory passing through $\infty$.
\end{itemize}

All in all, Th.\ref{th21} is now proven.
\end{proof}
\begin{remark}
In \cite{MBKPS} it was observed a component of $W^{u}_{Out}$ always tends to infinity. Th.\ref{th21} (and in particular, Cor.\ref{cor215}) is an analytic proof of this numerical observation.
\end{remark}

\section{Chaotic dynamics and Trefoil Knots}

\begin{figure}[h]
\centering
\begin{overpic}[width=0.6\textwidth]{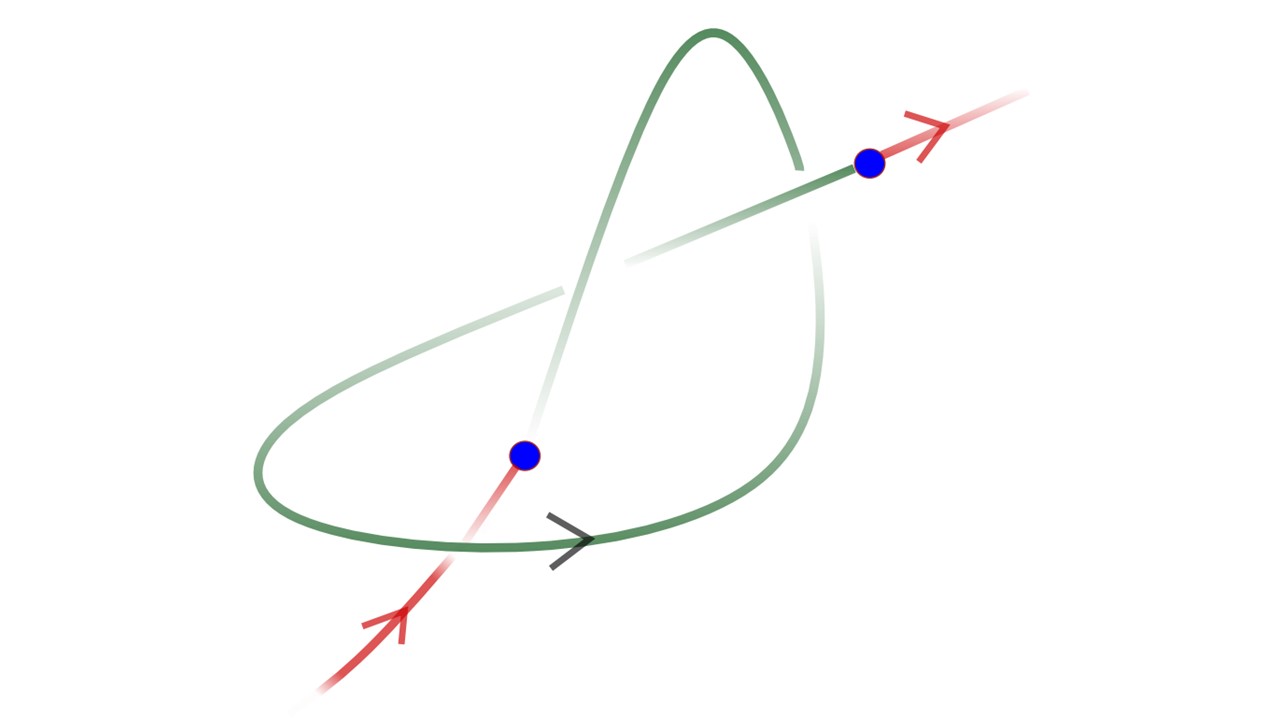}
\put(440,210){$P_{In}$}
\put(400,350){$\Theta$}
\put(660,370){$P_{Out}$}
\end{overpic}
\caption[A heteroclinic trefoil knot.]{\textit{A heteroclinic trefoil knot for the Rössler system (see Def.\ref{def32}). $\Theta$ denotes the bounded heteroclinic connection.}}
\label{fig7}
\end{figure}

Having characterized the global dynamics of the flow in Th.\ref{th21}, we now turn our gaze to prove a sufficient condition for the chaoticity of the flow (per Def.\ref{chaotic}). In more detail, in this section we prove that under some idealized assumptions (to which we refer as trefoil parameters), the Rössler system behaves chaotically. To begin, recall that given a parameter $p=(a,b,c)\in P$ we denote by $W^s_{In}$ and by $W^u_{Out}$ the respective one-dimensional invariant manifolds for the saddle-foci $P_{In},P_{Out}$ (see the discussion in page \pageref{eq:9}). We first define:

\begin{definition}\label{def31}
Consider $p\in P$ for which there exists a bounded heteroclinic trajectory for the corresponding Rössler system in $W^s_{In}\cap W^u_{Out}$, connecting $P_{In},P_{Out}$. In that case we say $p$ is a \textbf{heteroclinic parameter} - see the illustration at Fig.\ref{fig7}. 
\end{definition}

Having defined heteroclinic parameters for the Rössler system we now aim to introduce trefoil parameters as a geometric model that forms idealized version of the Rössler system - and in order to do so in a meaningful way, we first recall the results of the numerical studies. As observed numerically in \cite{MBKPS}, there exists a partition of the parameter space $P$ to three subsets: $S_1,S_2$, and $S_3$, on which the dynamics is characterized as follows:

\begin{itemize}
    \item For parameters $p\in S_1$, the corresponding Rössler system generates an attractor (not necessarily chaotic). Moreover, it was observed throughout $S_1$ that the two-dimensional $W^s_{Out}$ appears to shield trajectories from escaping to $\infty$ - while a bounded component of the one-dimensional $W^u_{Out}$ repels them towards the attractor (see the discussion in Section $V$ in \cite{MBKPS}). Additionally, the attractor intersects transversely with some cross section, on which the first-return map is well-defined. 
    \item For parameters $p\in S_2$, the dominant behavior is unbounded. In more detail, the two-dimensional $W^u_{Out}$ no longer shields trajectories from diverging to $\infty$ - and consequentially most trajectories appear to diverge to $\infty$. 
    \item For parameters $p\in S_3$, the attractor coincides with the two-dimensional manifold $W^s_{Out}$. In addition, $S_3$ forms a bifurcation set between $S_1$ and $S_2$.
\end{itemize}

\begin{figure}[h]
\centering
\begin{overpic}[width=0.6\textwidth]{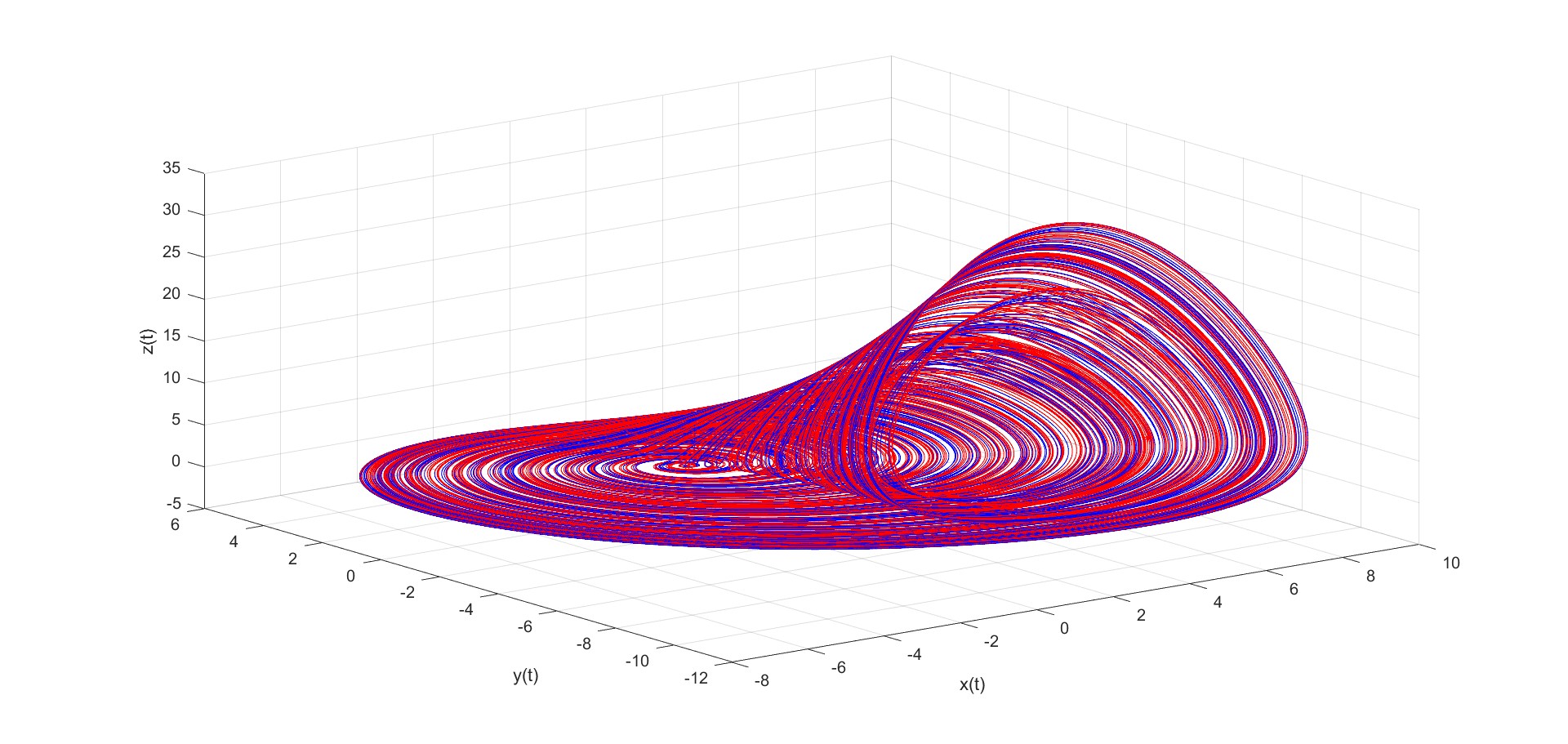}
\end{overpic}
\caption[A heteroclinic attractor.]{\textit{The attractor at parameter values $(a,b,c)\approx(0.468,0.3,4.615)$, sketched in Matlab. In these parameter values the vector field was observed numerically to generate a heteroclinic trefoil knot as in Fig.\ref{fig7} (see Fig.B.1 \cite{MBKPS}).}}
\label{attractortref}
\end{figure}

These numerical observations are too general to work with analytically - therefore, in order to rigorously study the dynamics of the Rössler system we now introduce another, more specific topological constrain on the dynamics which sharpens certain aspects of these observations.\\

To do so, recall that as observed numerically in \cite{MBKPS} there exists a parameter $p=(a,b,c)\in S_1$, $(a,b,c)\approx(0.468,0.3,4.615)$ at which the Rössler system generates a bounded heteroclinic trajectory $\Theta$ as in Fig.\ref{fig7} (see Fig.5.B.1 in \cite{MBKPS}). By numerical approximation, the attractor at the said parameters is inseparable from the fixed-points $P_{In},P_{Out}$, and appears to be wrapped around $\Theta$ (see the illustration in Fig.\ref{attractortref}). Now, recall the invariant one-dimensional manifolds $\Gamma_{In}$ and $\Gamma_{Out}$ given by Th.\ref{th21}. Under ideal circumstances, $\Theta$ is not linked with either $\Gamma_{In}$ or $\Gamma_{Out}$ - in which case the set $\Lambda=\Gamma_{In}\cup\Gamma_{Out}\cup\Theta\cup\{P_{In},P_{Out},\infty\}$ forms a heteroclinic trefoil knot in $S^3$ (see the illustration in Fig.\ref{fig7}).\\

\begin{figure}[h]
\centering
\begin{overpic}[width=0.5\textwidth]{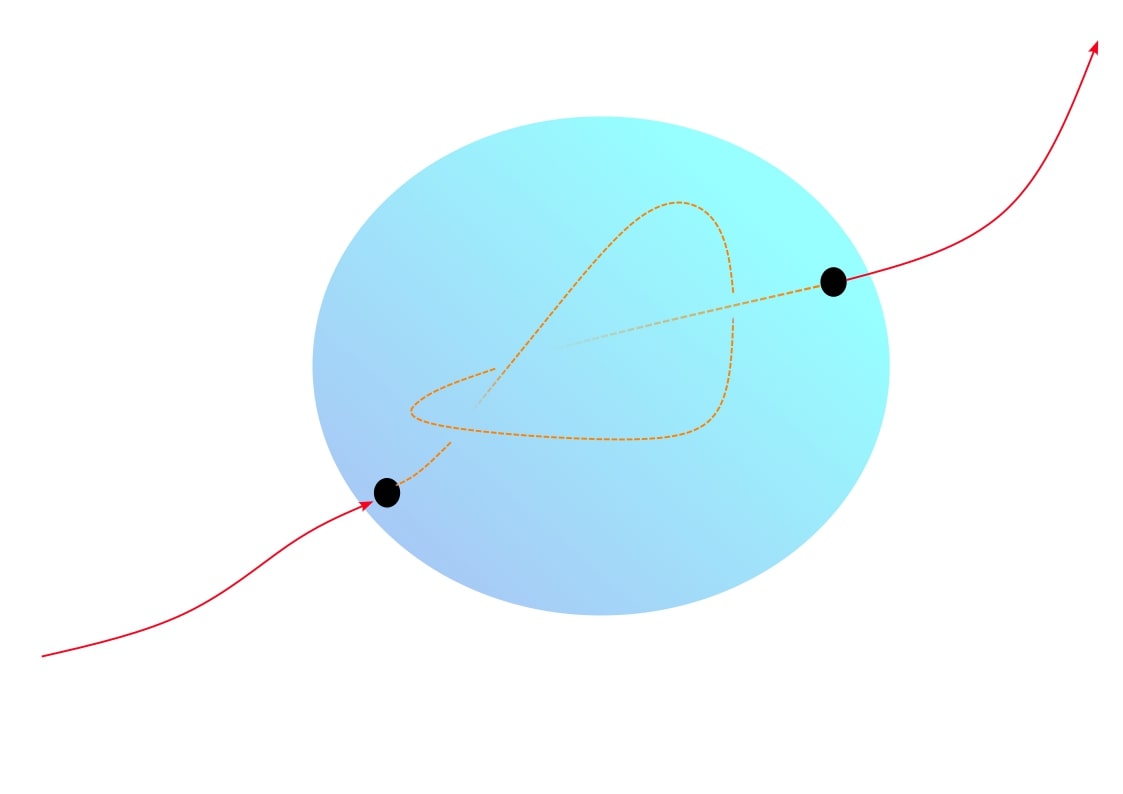}
\put(460,250){$\Theta$}
\put(270,190){$P_{In}$}
\put(160,80){$\Gamma_{In}$}
\put(400,480){$B_\alpha$}
\put(680,390){$P_{Out}$}
\put(860,420){$\Gamma_{Out}$}
\end{overpic}
\caption[The ball $B_\alpha$.]{\textit{The heteroclinic trajectory $\Theta$ trapped inside the topological ball $B_\alpha$, sketched as a the blue sphere.}}
\label{boundedtref}
\end{figure}
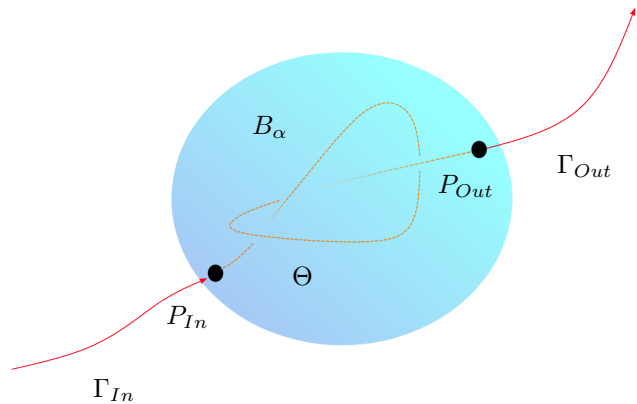

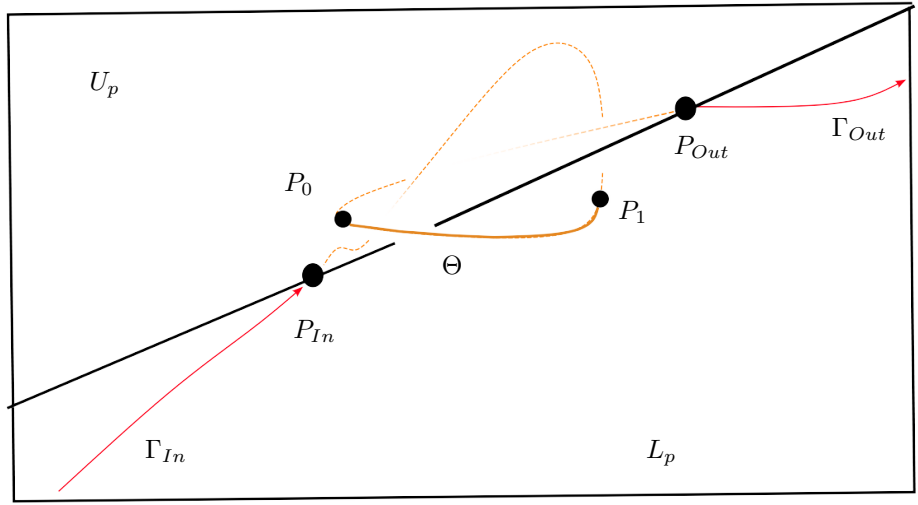
\begin{figure}[h]
\centering
\begin{overpic}[width=0.7\textwidth]{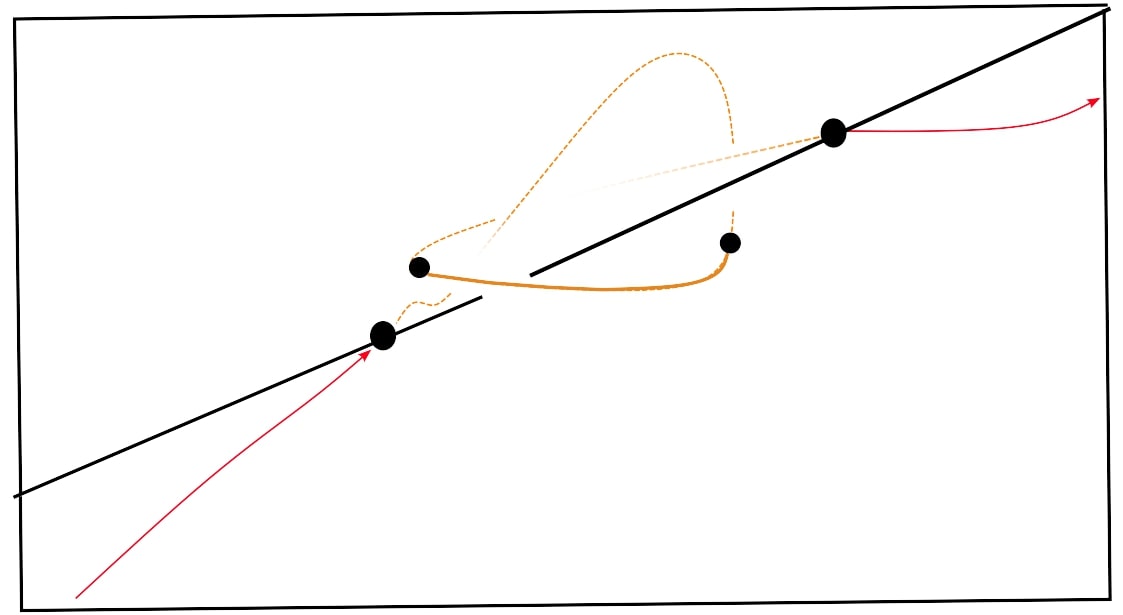}
\put(480,250){$\Theta$}
\put(320,180){$P_{In}$}
\put(160,50){$\Gamma_{In}$}
\put(310,340){$P_0$}
\put(670,310){$P_1$}
\put(100,450){$U_p$}
\put(700,50){$L_p$}
\put(730,380){$P_{Out}$}
\put(900,400){$\Gamma_{Out}$}
\end{overpic}
\caption[The intersection of the heteroclinic trefoil with $\{\dot{y}=0\}$.]{\textit{The heteroclinic trajectory $\Theta$ winds once around $P_{In}$ - hence it intersects the half-plane $U_p$ at $P_0$ and the half-plane $L_p$ at $P_1$.}}
\label{intersect}
\end{figure}

We now combine the existence of a heteroclinic trefoil knot with the general numerical observations above to propose a geometric model for the Rössler system defined as follows:

\begin{definition}\label{def32}
Let $p\in P$ be a heteroclinic parameter for the Rössler system which generates a bounded heteroclinic trajectory $\Theta$ as in Fig.\ref{fig7}. We say a smooth vector field $F_p$ of $\mathbf{R}^3$ is a \textbf{trefoil parameter which models the dynamics at $p$} (or in short, a \textbf{trefoil parameter}) provided it is a geometric model of the Rössler system at the parameter $p$, satisfying the following properties:
\begin{itemize}
    \item $F_p$ satisfies all the assumptions and conclusions of all the results proven in Sect.\ref{sect2} - i.e., it generates a geometric flow modeled after the Rössler system.
    \item The set $\Lambda$ (as defined above) forms a trefoil knot in $S^3$.
    \item The two-dimensional manifolds $W^u_{In}$ and $W^s_{Out}$ coincide. This condition implies $\overline{W^u_{In}}=\overline{W^s_{Out}}$ forms the boundary of an open topological ball - which, from now on, we always denote by $B_\alpha$. It is easy to see $\Theta\subseteq B_\alpha$, while $\Gamma_{In},\Gamma_{Out}\not\subseteq B_\alpha$ (see the illustration in Fig.\ref{boundedtref}).
    \item  $\Theta\cap \overline{U_p}=\{P_0\}$ is a point of transverse intersection (where $U_p$ is as in Lemma \ref{obs}) - see the illustration at Fig.\ref{intersect} and Fig.\ref{trefint}.
\end{itemize}
\end{definition}

We now claim the geometric model defined above captures key aspects of the numerical observations of the Rössler system. To see why, note that in addition to creating a heteroclinic trajectory as $\Theta$, in trefoil parameter the notion of the two-dimensional invariant manifold $W^s_{Out}$ shielding flow lines from escaping to $\infty$ is made precise. Later on, at the end of the next section, we will discuss how the results we will now prove to trefoil parameters can be used to study the non-idealized Rössler system.\\

Having defined trefoil parameters, let us give a brief overview of what lies ahead in this section. In subsection \ref{sect31} and \ref{sect32} we study the existence and basic properties of the first-return map for trefoil parameters. Finally, in subsection \ref{sect33} we prove Th.\ref{th31} - where we establish that at trefoil parameters the first-return map has infinitely many periodic orbits, thus implying chaotic behavior for the flow per Def.\ref{chaotic}. 

\subsection{The first-return map at trefoil parameters.}
\label{sect31}

\begin{figure}[h]
\centering
\begin{overpic}[width=0.6\textwidth]{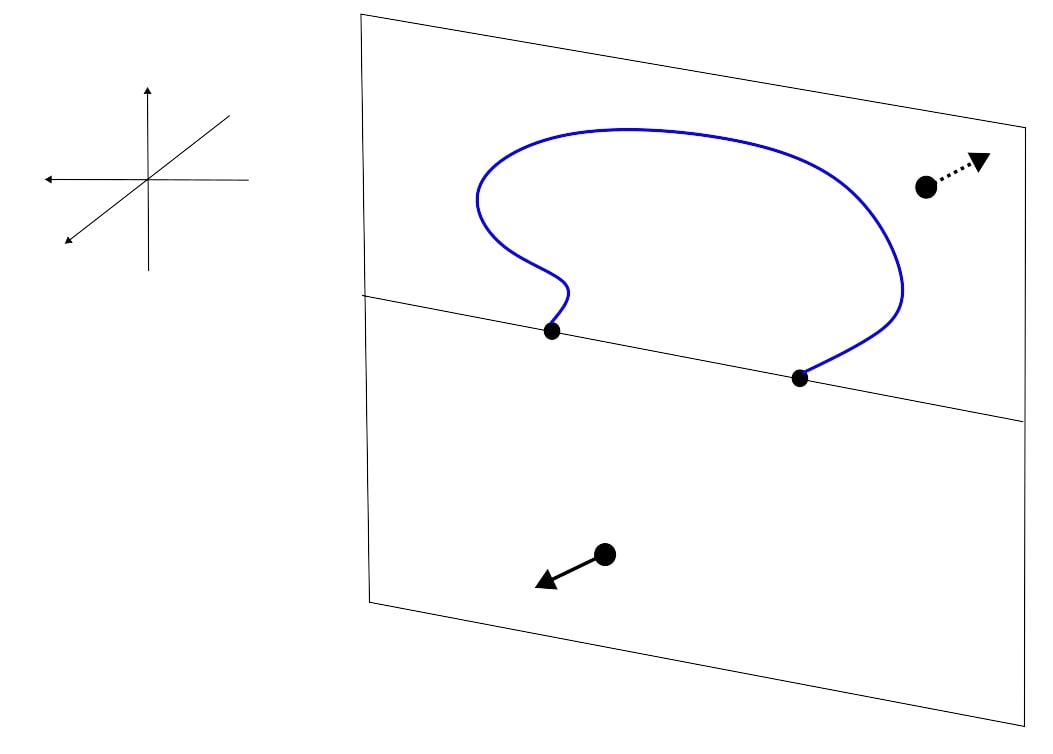}
\put(750,300){$P_{In}$}
\put(630,210){$L_p$}
\put(470,350){$P_{Out}$}
\put(850,340){$l_1$}
\put(125,640){$z$}
\put(30,480){$y$}
\put(15,540){$x$}
\put(650,600){$U_p$}
\put(550,530){$L$}
\put(650,380){$l$}
\put(400,430){$l_2$}
\put(20,205){$\{\dot{y}>0\}$}
\put(1090,270){$\{\dot{y}<0\}$}
\end{overpic}
\caption[$l_1$, $l_2$ and $l$.]{\textit{T}\textit{he curves $l_1,l_2$ and $l$ on $\{\dot{y}=0\}$ along with the directions of $F_p$ on the plane $\{\dot{y}=0\}$- these curves are all sub-arcs of $l_p$, the tangency set of the vector field $F_p$ to $\{\dot{y}=0\}$. In this illustration, we have $\beta_{In}=\beta_{Out}=L$, while $s_{In}=P_{Out}$, $s_{Out}=P_{In}$.}}\label{Qq}

\end{figure}
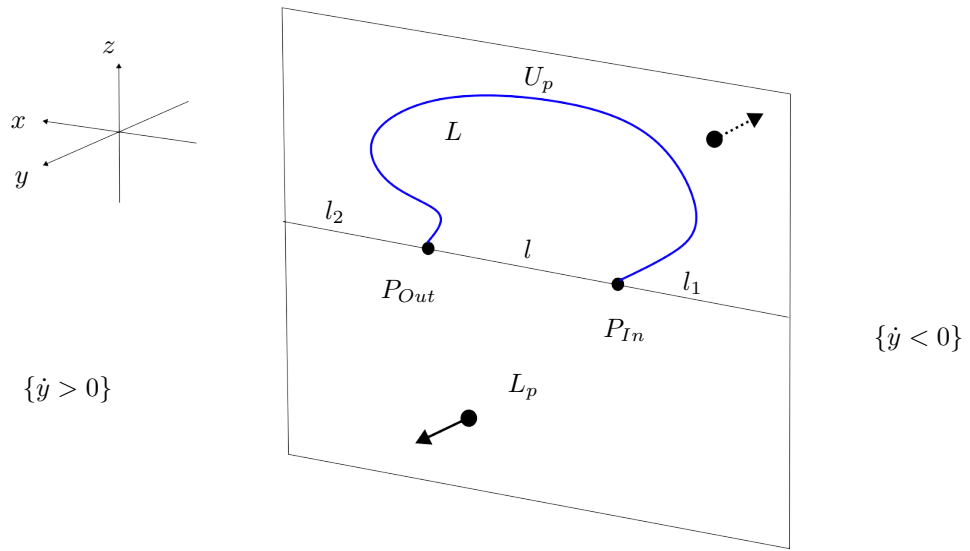

Let $F_p$ be a trefoil parameter for the Rössler system, and consider the open topological ball $B_\alpha$ given by Def.\ref{def32}. Recall the cross-section $U_p$ is an open half plane (see the discussion before Lemma \ref{obs}), and set $D_\alpha=B_\alpha\cap U_p$ (see the illustration in Fig.\ref{fig8}) - we first prove:
\begin{lemma}
    \label{noempty}
Let $F_p$ be a trefoil parameter - then, the set $D_\alpha$ is non-empty, $P_0$ is interior to $D_\alpha$, and $P_{In},P_{Out}\in\partial D_\alpha$.
\end{lemma}
\begin{proof}
 Recall Lemma \ref{cor211}, where we proved $W^u_{In},W^s_{Out}$ are both transverse to $U_p$ at the fixed-points $P_{In},P_{Out}$ (respectively). Since $F_p$ is a trefoil parameter, by $\partial B_\alpha=W^u_{In}\cup\{P_{In},P_{Out}\}=W^s_{Out}\cup\{P_{In},P_{Out}\}$, it follows $U_p\setminus W^u_{In}$ intersects both components of  $\mathbf{R}^3\setminus \partial B_\alpha$ - which proves $D_\alpha=U_p\cap B_\alpha\ne\emptyset$. We similarly conclude $P_{In},P_{Out}\in\partial B_\alpha\cap\partial U_p$, i.e., $P_{In},P_{Out}\in\partial D_\alpha$.\\
 
To conclude the proof it remains to show $P_0$ is interior to $D_\alpha$. Because $\Theta\subseteq B_\alpha$ and because by Def.\ref{def32} $P_0$ is interior to the cross-section $\overline{U_p}$ (see the illustrations in Fig.\ref{boundedtref} and Fig.\ref{intersect}), it follows $P_0$ is in $\overline{D_\alpha}$. Note that $\partial D_\alpha\subseteq\partial U_p\cap\partial B_\alpha$ - therefore, since $P_0$ is interior to both ${U_p}$ and $B_\alpha$ it follows $P_0$ is also be interior to $D_\alpha$ and Lemma \ref{noempty} now follows (see the illustration in Fig.\ref{fig8}). 
\end{proof}

Set $S=W^u_{In}=W^s_{Out}$ - by $\partial B_\alpha=\{P_{In},P_{Out}\}\cup S$ it follows $\partial B_\alpha$ is invariant under the flow, which implies $\overline{B_\alpha}$ is also invariant under the flow. Therefore, since $\overline{B_\alpha}$ is bounded, given any $s\in \overline{B_\alpha}$, its trajectory is bounded as well. Combining this fact with Lemma \ref{obs}, we now prove:
\begin{lemma}\label{firstret}
Let $F_p$ be a trefoil parameter for the Rössler system Then, the first-return map $f_p:\overline{D_\alpha}\setminus\{P_0\}\rightarrow \overline{D_\alpha}\setminus\{P_0\}$ is well-defined - and conversely, so is its inverse $f^{-1}_p:\overline{D_\alpha}\setminus\{P_0\}\rightarrow \overline{D_\alpha}\setminus\{P_0\}$. Moreover, whenever $s\in\overline{D_\alpha}\setminus\{P_0\}$ is not a fixed point, $f_p(s)$ is also not a fixed-point.
\end{lemma}
\begin{proof}
Recall the forward trajectory of an initial condition $s\in\overline{D_\alpha}$ tends to a fixed-point precisely when $s\in W^s_{In}\cup W^s_{Out}$ - and that $(W^s_{In}\cup W^s_{Out})\cap B_\alpha=W^s_{Out}\cup\Theta$, where $\Theta$ is the heteroclinic trajectory given by Def.\ref{def32}.\\

We now prove Lemma \ref{firstret}, in three short steps. To begin, consider $s\in W^s_{Out}$, the stable two-dimensional, invariant manifold of the saddle focus $P_{Out}$. By Lemma \ref{cor211}, $W^s_{Out}$ is transverse to $U_p$ at $P_{Out}$ - which implies the forward trajectory of $s$ hits the cross-section $\overline{U_p}$ transversely infinitely many times as it spirals to $P_{Out}$ (see the illustration in Fig.\ref{loci}). Consequentially, by $\partial B_\alpha=W^s_{Out}\cup\{P_{In},P_{Out}\}$ it follows that given $s\in\partial B_\alpha$ which is not a fixed-point, its forward trajectory hits $\overline{U_p}$ transversely infinitely many times - and each such intersection point between the trajectory of $s$ and $\overline{U_p}$ is not a fixed point .\\
 
Now, consider an initial condition $s\in B_\alpha\setminus\Theta$ - as $s$ is interior to $B_\alpha$, it is not a fixed point itself. Moreover, since $B_\alpha$ is an open topological ball in $\mathbf{R}^3$, by $\partial B_\alpha=W^s_{Out}\cup\{P_{In},P_{Out}\}$ we conclude $s\not\in W^s_{Out}$, hence its trajectory does not limit to $P_{Out}$ - and by $\Theta=W^s_{In}\cap \overline{B_\alpha}$ it follows the forward trajectory of $s$ does not limit to $P_{In}$ either. Therefore, since $B_\alpha$ is bounded and invariant under the flow, so is the forward trajectory of $s$: which, by Lemma \ref{obs}, proves the forward trajectory of $s$ hits $\overline{U_p}$ transversely infinitely many times. Moreover, because $s$ is not a fixed point (and doesn't get attracted to one through some stable manifold) every intersection point between its trajectory and $\overline{U_p}$ is not fixed point.\\

We conclude that for all $s\in\overline{B_\alpha}\setminus(\Theta\cup\{P_{In},P_{Out}\})$ its forward trajectory hits the cross-section $\overline{U_p}$ transversely infinitely many times - and moreover, every point of intersection between the said trajectory and $\overline{U_p}$ is not a fixed point. Recalling $D_\alpha=B_\alpha\cap U_p$, because $\overline{B_\alpha}$ is invariant under the flow the trajectory of any initial condition $s\in \overline{D_\alpha}\setminus\Theta$ returns to $\overline{U_p}$ precisely at $\overline{D_\alpha}\setminus\Theta$. Since $\overline{D_\alpha}\cap\Theta=\{P_0\}$, it now follows the first-return map $f_p:\overline{D_\alpha}\setminus\{P_0\}\to\overline{D_\alpha}\setminus\{P_0\}$ is well defined - and whenever $s\in\overline{D_\alpha}\setminus\{P_0\}$ is not a fixed point, neither is $f_p(s)$. Using a similar argument applied to the inverse flow, we conclude the inverse first-return map $f^{-1}_p:\overline{D_\alpha}\setminus\{P_0\}\to\overline{D_\alpha}\setminus\{P_0\}$ is also well-defined and Lemma \ref{firstret} now follows.   
\end{proof}
\begin{remark}
    \label{full}
    Using the invariance of $B_\alpha$ under the flow, the argument used to prove Lemma \ref{firstret} implies both $f_p(\overline{D_\alpha}\setminus\{P_0\})=\overline{D_\alpha}\setminus\{P_0\}$ and $f^{-1}_p(\overline{D_\alpha}\setminus\{P_0\})=\overline{D_\alpha}\setminus\{P_0\}$.
\end{remark}

We conclude this section by proving the set $D_\alpha$ is a topological disc in the half-plane $U_p$. To do so, first recall the curve $l_p$ is parameterized by $l_p(x)=(x,-\frac{x}{a},\frac{x}{a})$, $x\in\mathbf{R}$ (see Lemma \ref{lem23} and the illustrations in Fig.\ref{planes} and Fig.\ref{loci}). In particular, recall $l_p$ is the tangency curve of $F_p$ to the plane $\{\dot{y}=0\}$, and that $\{\dot{y}=0\}\setminus l_p$ is composed of two components - the half planes $U_p$ and $L_p$ (see the discussion preceding Lemma \ref{obs}, and the illustration in Fig.\ref{Qq}). As shown immediately before the proof of Lemma  \ref{obs}, on $U_p$ the vector field $F_p$ points into $\{\dot{y}<0\}$ (the region in front of $\{\dot{y}=0\}$), while on $L_p$, $F_p$ points into $\{\dot{y}>0\}$ (the region behind $\{\dot{y}=0\}$) - see the illustration in Fig.\ref{Qq}.\\

Additionally, we will also need to recall some notations from Section 2 - and introduce two new ones:

\begin{itemize}\label{Ll}
    \item From now on, we always denote by $l$ the set $\{l_p(x)|x\in(0,c-ab)\}$ - that is, $l$ is the sub-arc on $l_p$ connecting $P_{In},P_{Out}$ (see Fig.\ref{trefint} and Fig.\ref{Qq}).
    \item Similarly to what we did in Sect.2, we define $l_1=\{l_p(x)|x<0\}$, and $l_2=\{l_p(x)|x>c-ab\}$ (see the illustration in Fig.\ref{Qq}).
    \item Finally, we set $L=W^u_{In}\cap {U_p}=W^s_{Out}\cap{U_p}=\partial B_\alpha\cap U_p$ (see the illustration in Fig.\ref{Qq} and .\ref{trefint}). Since $W^u_{In}=W^s_{Out}=S\subseteq\partial B_\alpha$, by $D_\alpha=B_\alpha\cap U_p$ it follows $L\subseteq\partial D_\alpha$. 
\end{itemize}

With these ideas in mind, we now prove:

\begin{proposition}
\label{arccor}    Let $F_p$ be a trefoil parameter. Then, the set $L$ is a curve in $U_p$, with one endpoint at $P_{In}$ and another at $P_{Out}$, as illustrated in both Fig.\ref{trefint} and Fig.\ref{Qq}. Consequentially, $D_\alpha$ is a open topological disc on the half-plane $U_p$, satisfying $\partial D_\alpha=L\cup l\cup\{P_{In},P_{Out}\}$ - see the illustrations in Fig.\ref{trefint} and Fig.\ref{Qq}.
\end{proposition}

\begin{proof}
    As we assume $F_p$ is a trefoil parameter, by definition we have $W^u_{In}=W^s_{Out}$ - and therefore also $\partial B_\alpha=W^s_{Out}\cup\{P_{In},P_{Out}\}=W^u_{In}\cup\{P_{In},P_{Out}\}$ (we will use both these representations throughout the proof). The proof of Prop.\ref{arccor} is rather technical and based on directly analyzing the local dynamics of the vector field $F_p$ in $\mathbf{R}^3$.\\

To give an outline of our argument, note that that by Lemma \ref{cor211}, $\partial B_\alpha$ is transverse to $U_p$ at both fixed-points. This implies there exists a component $\beta_{In}\subseteq W^u_{In}\cap U_p$ which is a curve beginning at $P_{In}$ and terminates at some point $s_{In}\in l_p$ - similarly, there also exists a component $\beta_{Out}\subseteq W^s_{Out}\cap U_p$ which begins at $P_{Out}$ and terminates at some $s_{Out}\in l_p$ (see the illustrations in Fig.\ref{d1} and \ref{d2}, respectively). We now sketch the proof of Prop.\ref{arccor} - where $l,l_1$ and $l_2$ are the components of $l_p\setminus\{P_{In},P_{Out}\}$ as defined above: 

\begin{itemize}
  \item First, we prove $s_{In}\not\in l_1\cup l$ (see Lemma \ref{nol1}). 
  \item Second, we prove - using an almost symmetric argument - that $s_{Out}\not\in l_2\cup l$ (see Lemma \ref{nol}).
  \item Finally, we tie these results together and conclude $s_{In}=P_{Out}$, $s_{Out}=P_{In}$, from which Prop.\ref{arccor} would follow.

\end{itemize}

To begin, recall we denote by $f_p:\overline{D_\alpha}\setminus\{P_0\}\to\overline{D_\alpha}\setminus\{P_0\}$ the first-return map, and recall $f_p$ is defined throughout $\partial B_\alpha\cap \overline{U_p}\subseteq\partial D_\alpha$ (see Lemma \ref{firstret}). Per the outline above, we first prove:

\begin{lemma}
    \label{nol1}
With the notations and definitions above, whenever $F_p$ is a trefoil parameter, $s_{In}\not\in l_1\cup l$. 
\end{lemma}
\begin{proof}
    We prove Lemma \ref{nol1} by contradiction - we first rule out the possibility $s_{In}\in l_1$, after which we rule out $s_{In}\in l$, using similar arguments.\\
    
    To begin, assume $s_{In}\in l_1$, which implies $\beta_{In}$ is a curve in ${U_p}$ whose closure connects $P_{In}$ and $l_1$ as in Fig.\ref{d1}. Since by definition $\beta_{In}\subseteq W^u_{In}\cap U_p$, as $W^u_{In}$ is the unstable invariant manifold of $P_{In}$ it follows both $f_p(\beta_{In})\not\subseteq \beta_{In}$ and $\beta_{In}\subseteq f_p(\beta_{In})$. This proves there exists an initial condition $s_1$, strictly interior to $\beta_{In}$, s.t. $f_p(s_1)=s_{In}$ as illustrated in Fig.\ref{d1}. Since by Lemma \ref{lem23} the trajectory of $s_1$ flows to $s_{In}$ through $\{\dot{y}\leq0\}$, the $y-$coordinate of $s_1$ is strictly greater than that of $s_{In}$ (see the illustration in Fig.\ref{d1}).\\

\begin{figure}[h]
\centering
\begin{overpic}[width=0.6\textwidth]{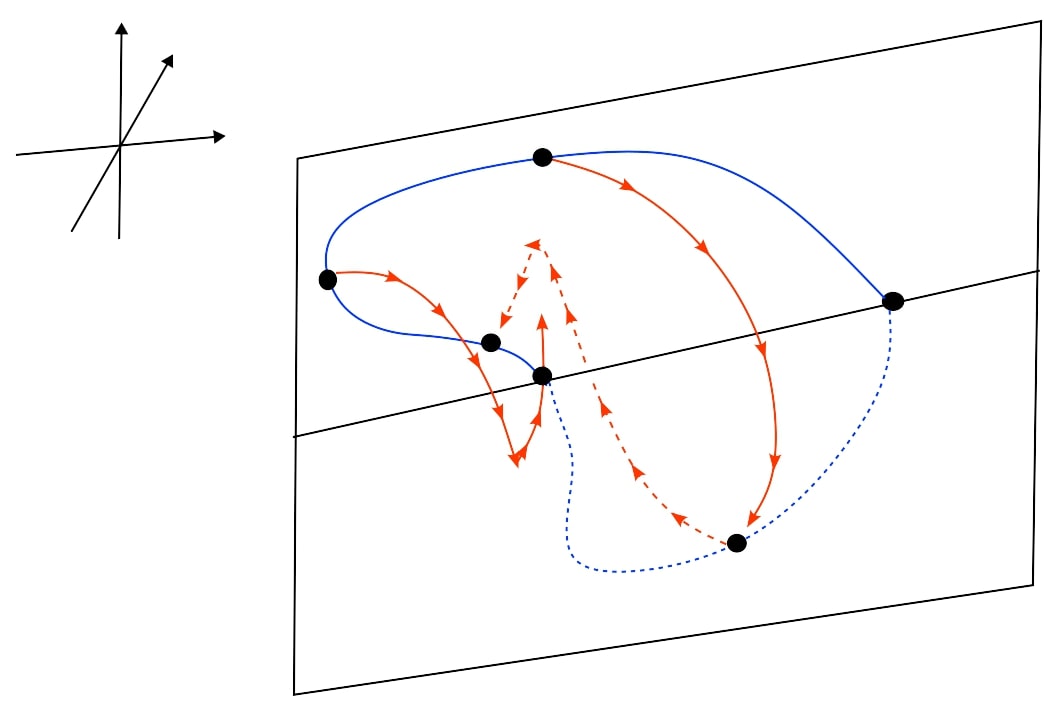}
\put(430,70){$L_p$}
\put(550,300){$s_{In}$}
\put(870,320){$P_{In}$}
\put(100,660){$z$}
\put(170,625){$x$}
\put(220,535){$y$}
\put(850,500){$U_p$}
\put(320,320){$\beta_{In}$}
\put(650,540){$\alpha_{In}$}
\put(380,220){$l_1$}
\put(330,430){$s_1$}
\put(20,205){$\{\dot{y}>0\}$}
\put(1090,270){$\{\dot{y}<0\}$}
\end{overpic}
\caption[The case when $s_{In}\in l_1$.]{\textit{The curve $\beta_{In}$ when $s_{In}\in l_1$. The sub-arc $\alpha_{In}$ flows to $\beta_{In}$ by hitting $L_p$, while the $y-$coordinate of $s_1$ is greater than that of $s_{In}$. This scenario will imply a contradiction.}}\label{d1}

\end{figure}
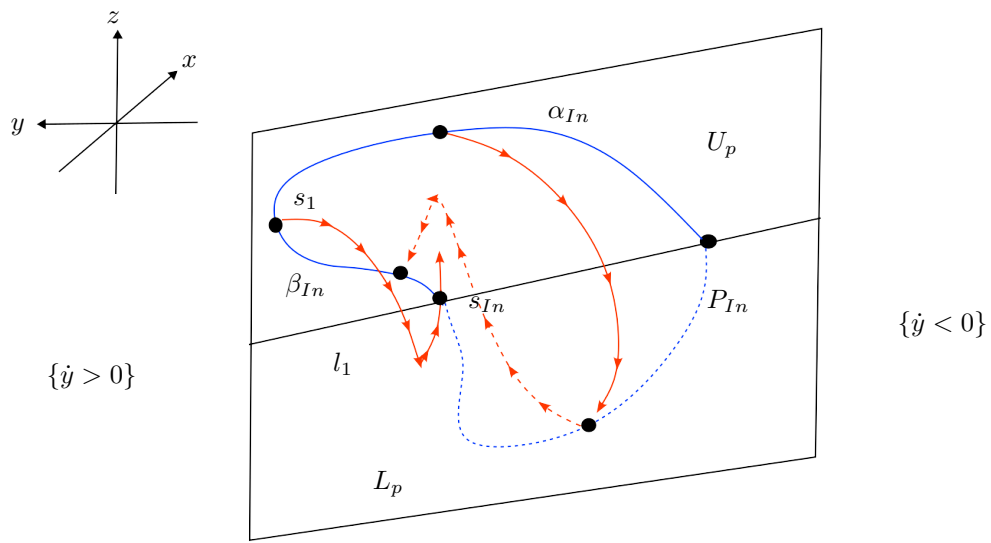

    Now, let $\alpha_{In}\subseteq\beta_{In}$ denote the sub-arc of $\beta_{In}$ connecting $P_{In}$ and $s_1$, and let us consider the two-dimensional set $V$ generated by the collection of flow-lines connecting $\alpha_{In}$ to $\beta_{In}$ (as illustrated in Fig.\ref{d1}) - around $P_{In}$, $V$ forms the boundary of a topological cone $C$ with a tip at $P_{In}$. Moreover, $V$ and the plane $\{\dot{y}=0\}$ trap between them a topological cap $\mu=C\cap\{\dot{y}\geq0\}$, as illustrated in Fig.\ref{d1}. As a consequence, one can enter the cone $C$ only through the (closed) two-dimensional region $\{\dot{y}\leq0\}$. Now, write $s_{In}=(x_{In},y_{In},z_{In})$ and consider the half-plane $H=\{(x,y_{In},z)|x+ay_1<0\}$ - by the discussion above, writing $s_1=(x_1,y_1,z_1)$ we know $y_1>y_{In}$, which, as $s_1$ flows to $s_{In}$ implies $V$ intersects $H\cup L_p$ in a collection of curves connecting $\alpha_{In}$ and $s_{In}$, as illustrated in Fig.\ref{DD}.\\

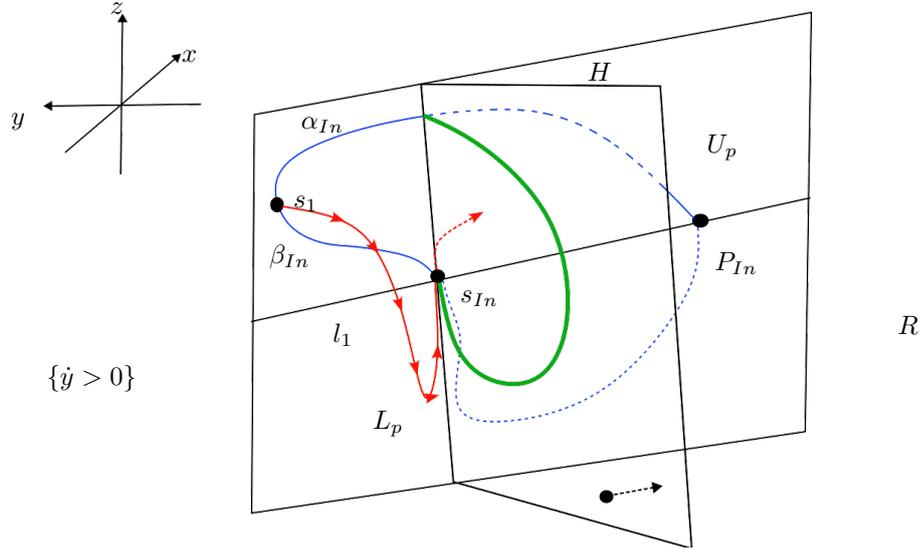
\begin{figure}[h]
\centering
\begin{overpic}[width=0.6\textwidth]{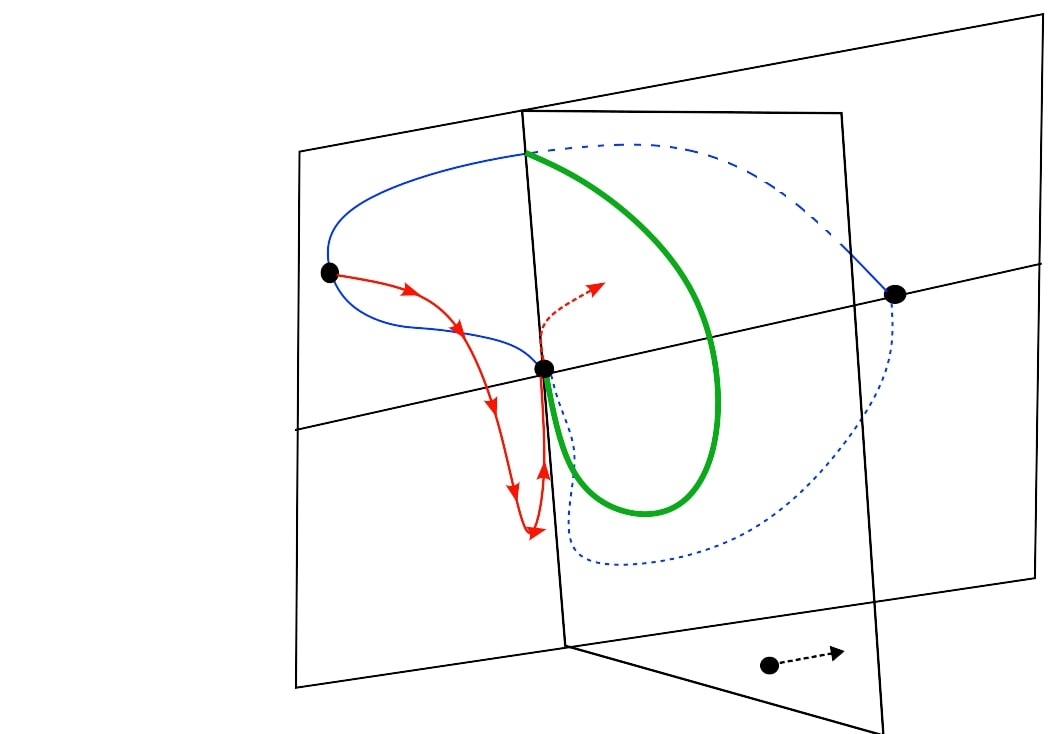}
\put(430,150){$L_p$}
\put(540,310){$s_{In}$}
\put(860,350){$P_{In}$}
\put(850,500){$U_p$}
\put(300,360){$\beta_{In}$}
\put(340,530){$\alpha_{In}$}
\put(700,590){$H$}
\put(380,260){$l_1$}
\put(330,430){$s_1$}
\put(20,205){$\{\dot{y}>0\}$}
\put(1090,270){$R$}
\end{overpic}
\caption[The intersection of $C$ with $H$.]{\textit{The intersection of the cone $C$ with $H$, sketched as the dark green curve, along with the directions of $F_p$ on $H$. The trajectory of $s_{In}$ enters the region $R=\{(x,y,z)|\dot{y}<0,y<y_{In}\}$ immediately after leaving $s_{In}$.}}\label{DD}

\end{figure}
    
The normal vector to $H$ is $(0,1,0)$, therefore, given $s\in H$, $s=(x,y_{In},z)$, by definition $F_p(s)\bullet(0,1,0)=\dot{y}=x+ay_{In}<0$ (see the illustration in Fig.\ref{DD}) - therefore, on $H$ the vector field $F_p$ points into the region $\{(x,y,z)|y<y_{In},\dot{y}<0\}$. Consequentially, because $s_{In}\in l_1$, we know $F_p(s_{In})$ is tangent to the plane $\{\dot{y}=0\}$ (see Lemma \ref{lem23}) - which implies the trajectory of $s_{In}$ enters the region $C\cap\{(x,y,z)|y<y_{In},\dot{y}<0\}$ upon leaving $s_{In}$ (as illustrated in Fig.\ref{DD}). Since $s_{In}$ lies on the trajectory of $s_1$, the trajectory of every $s\in W^u_{In}$ sufficiently close to $s_1$ also enters $C\cap\{y<y_{In},\dot{y}<0\}$ along the trajectory of $s_1$ (see the illustration in Fig.\ref{d2}).\\

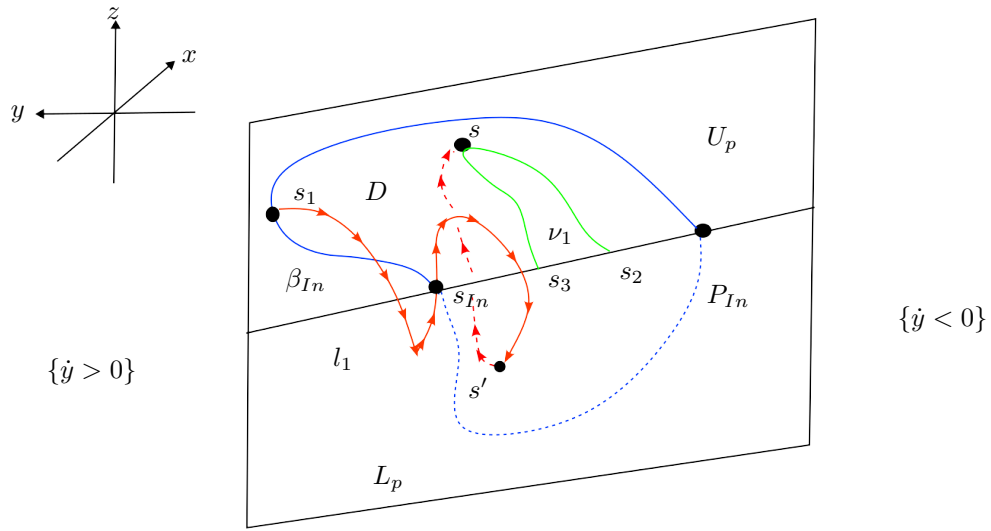
\begin{figure}[h]
\centering
\begin{overpic}[width=0.6\textwidth]{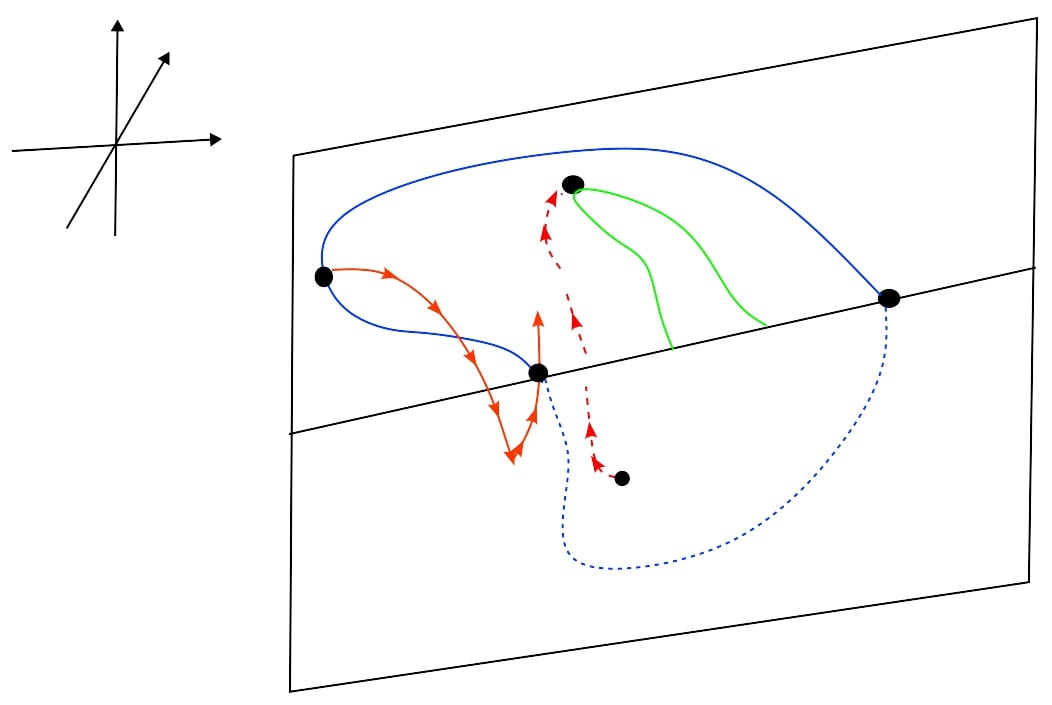}
\put(430,70){$L_p$}
\put(530,300){$s_{In}$}
\put(550,180){$s'$}
\put(550,510){$s$}
\put(850,300){$P_{In}$}
\put(100,660){$z$}
\put(170,625){$x$}
\put(220,535){$y$}
\put(850,500){$U_p$}
\put(320,320){$\beta_{In}$}
\put(650,380){$\nu_1$}
\put(740,330){$s_2$}
\put(650,320){$s_3$}
\put(380,220){$l_1$}
\put(330,430){$s_1$}
\put(20,205){$\{\dot{y}>0\}$}
\put(1090,270){$\{\dot{y}<0\}$}
\put(420,430){$D$}
\end{overpic}
\caption[The intersection of $W^u_{In}$ with $D$.]{\textit{The trajectory of $s_1$ flows to some $s'\in L_p$, after which it hits $D$ transversely at $s$. Consequentially, there exists a component $\nu_1$ of $W^u_{In}\cap U_p$ in $D$. We use this to generate a contradiction.}}\label{d2}

\end{figure}

    Now, set $D=C\cap U_p$ - by definition, $D$ is a Jordan domain on $U_p$, trapped between the curves $\beta_{In}$ and $l_1$ (as illustrated in Fig.\ref{d2}). As the trajectory of $s_1$ enters $C$ upon leaving $s_{In}$, by Lemma \ref{firstret} and by $f_p(s_1)=s_{In}$ there must exists some minimal $k>1$ s.t. $f^k_p(s_1)\in D$ - which implies there exists a component $\nu_{1}$ of $W^u_{In}\cap U_p$ in $D$ (see the illustration in Fig.\ref{d2}). Since $W^u_{In}$ is a surface and both $\nu_1$ and $\beta_{In}$ are components of the transverse intersection $W^u_{In}\cap U_p$, we have $\nu_{1}\cap\beta_{In}=\emptyset$ - therefore, we conclude $\nu_{1}$ is a curve in $D$ with two endpoints in $l_1\cap\partial D_1$, which we denote by $s_2$ and $s_3$ (see Fig.\ref{d2}). As there are no homoclinic trajectories in $W^u_{In}$ per assumption (as $p$ is a trefoil parameter), both $s_2$ and $s_3$ are not $P_{In}$, and are also distinct from one another (as illustrated in Fig.\ref{d3}).\\
    
   To continue, let us suspend $\nu_{1}$ with the flow - as $s_2,s_3\in l_1\cap\partial D$, using a similar argument to the one used to prove the trajectory of $s_{In}$ enters $C$, it follows the forward trajectories of $s_2$ and $s_3$ remain trapped in $C$ until hitting $D$ transversely. Therefore, $f(\nu_{1})=\nu_{2}$ is also a collection of curves inside $D$, which lie on components of $W^u_{In}\cap D$. Let us note every component of $W^u_{In}\cap D$ is a curve with two endpoints on $l_1\cap\partial D_1$ (as illustrated in Fig.\ref{d3}). Additionally, we have $\nu_1\cap \nu_2=\emptyset$, $\nu_2\cap\beta_{In}=\emptyset$.\\
   
Repeating this process, as every component of $W^u_{In}\cap D$ has two endpoints on $l_1\cap \partial D$, it follows $\nu_{2}$ flows through the cone $C$ to $\nu_{3}=f_p(\nu_2)$, another collection of curves trapped inside $D$ - and again, every component of $\nu_3$ lies on some component of $W^u_{In}\cap D$, hence the same is true for $\nu_4=f_p(\nu_3)$. Repeating this suspension over and over, we conclude the sequence $\{\nu_n\}_n$, $\nu_n=f_p(\nu_{n-1})$, $n>2$ satisfies the following:
    
    \begin{itemize}
        \item For every $n$, the trajectories of initial conditions in $\nu_{n-1}$ flow to $\nu_n$ through the cone $C$.
        \item For every $n$, we have $\nu_n\subseteq D=C\cap U_p$. 
    \end{itemize}

Since $D$ is a Jordan domain trapped between $\beta_{In}$ and $l_1$, it follows by $s_{In}\ne P_{Out}$ that $P_{Out}\not\in\overline{D}$ (see the illustration in Fig.\ref{d3}). This implies that given any initial condition $s\in\nu_1\subseteq W^u_{In}\cap U_p$, by $f^n_p(s)\in\nu_{n+1}$ we have $\lim_{n\to\infty}f^n_p(s)\ne P_{Out}$. However, because $F_p$ is a trefoil parameter, by Def.\eqref{def32} we must have $W^u_{In}=W^s_{Out}$ - which implies that for every $s\in W^u_{In}\cap U_p$ we have $\lim_{n\to\infty}f^n_p(s)=P_{Out}$. This is a contradiction, from which we conclude $s_{In}\not\in l_1$, i.e., $\overline{\beta_{In}}$ cannot connect $P_{In}$ with $l_1$ through $U_p$.\\

\begin{figure}[h]
\centering
\begin{overpic}[width=0.6\textwidth]{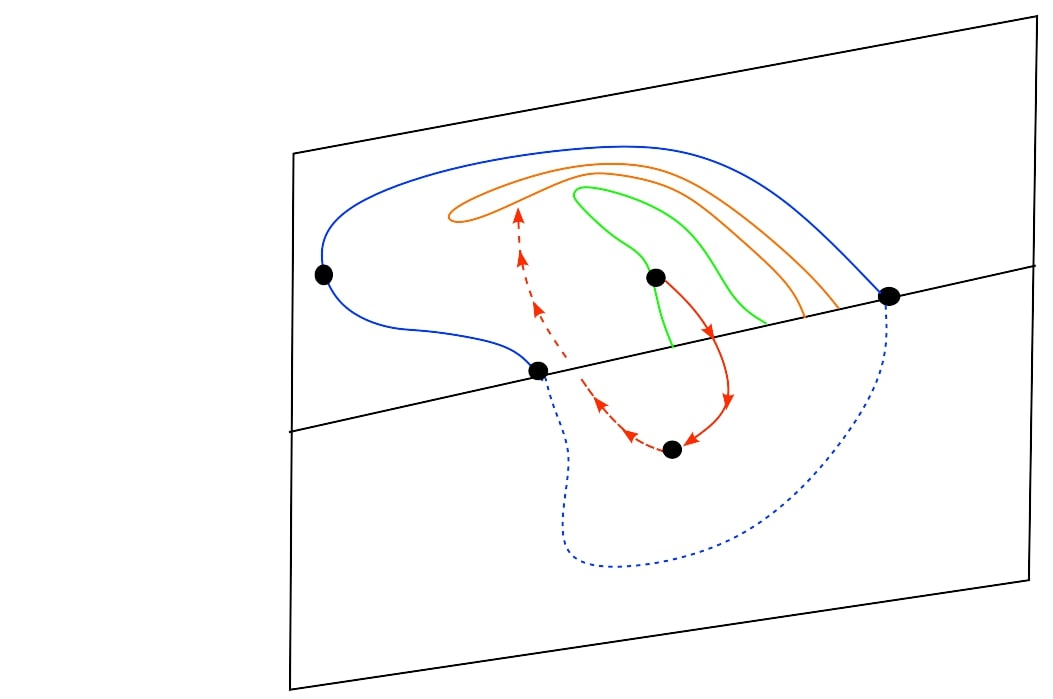}
\put(430,70){$L_p$}
\put(500,270){$s_{In}$}
\put(550,510){$\nu_2$}
\put(850,300){$P_{In}$}
\put(850,500){$U_p$}
\put(320,320){$\beta_{In}$}
\put(650,380){$\nu_1$}
\put(380,220){$l_1$}
\put(330,420){$s_1$}
\put(20,205){$\{\dot{y}>0\}$}
\put(1090,270){$\{\dot{y}<0\}$}
\put(420,400){$D$}
\end{overpic}
\caption[The sequence $\{\nu_n\}_n$.]{\textit{The sequence of curves $\nu_n$ which flow to one another and remain trapped in $D$. This yields a contradiction.}}\label{d3}

\end{figure}
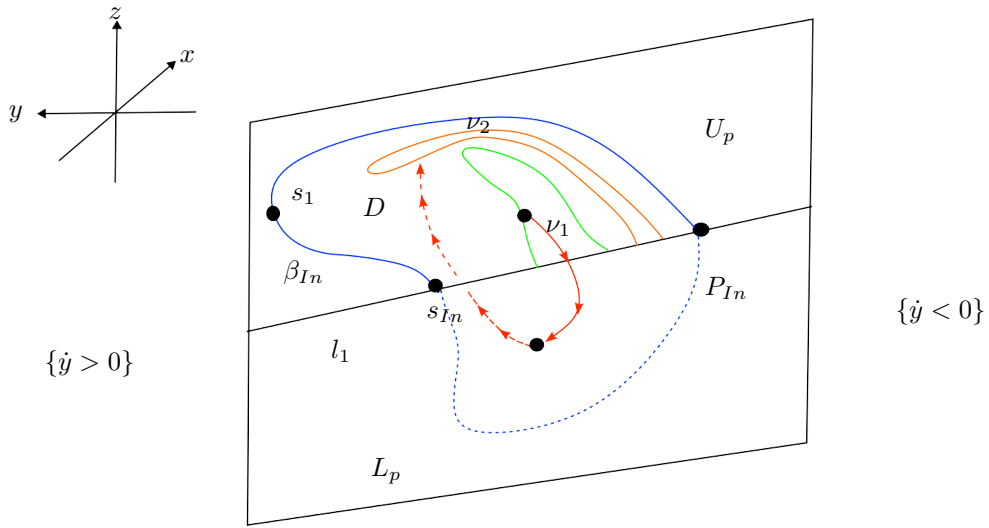

Having proven $s_{In}\not\in l_1$, we now prove $s_{In}\not\in l$ using a similar argument - thus concluding the proof of Lemma \ref{nol1}. To this end, assume by contradiction $s_{In}\in l$, as illustrated in Fig.\ref{d4}. Recalling Lemma \ref{cor211}, again we conclude there exists some interior $s_1\in\beta_{In}$ s.t. $f_p(s_1)=s_{In}$, and some sub-arc $\alpha_{In}\subseteq\beta_{In}$ connecting $P_{In}$ and $s_1$ (as illustrated in Fig.\ref{d4}) s.t. $f_p(\alpha_{In})=\beta_{In}$. Note that $V$, the collection of flow lines connecting $\alpha_{In}$ and $\beta_{In}$ traps a topological cone $C$ - which intersects $L_p$ in a curve $\rho_{In}$, as indicated in Fig.\ref{d4} (it is easy to see $P_{Out}\not\in \overline{C}$). By definition, the trajectory of $s_1$ flows through $\{\dot{y}<0\}$ to some $s\in\rho_{In}$ after which it flows through $\{\dot{y}>0\}$ to $s_{In}$. This implies the $y-$component of $s$ is smaller than that of $s_{In}$ (see the illustration in Fig.\ref{d4}).\\

Similar arguments to those used before prove the flow lines connecting $\rho_{In}$ and $\beta_{In}$ and the plane $\{\dot{y}=0\}$ trap between them a topological cap $\mu'\subseteq\{\dot{y}\leq0\}$ - moreover, setting $s_{In}=(x_{In},y_{In},z_{In})$ and considering the half-plane $H'=\{(x,y_{In},z)|x+ay_{In}>0\}$, using similar arguments we conclude the trajectory of $s_{In}$ enters the region $C\cap\{(x,y,z)|y>y_{In},\dot{y}>0\}$ upon leaving $s_{In}$. Again, the trajectories of initial conditions $s\in W^u_{In}$ sufficiently close to $s_1$ also enter $C$ - and similarly to the previous arguments, after entering $C$ the trajectories of such $s$ can hit $U_p$ transversely only in the region $C\cap U_p$ (i.e., the analogue of $D$ from before). Using a similar argument we again derive a contradiction from which we conclude $s_{In}\not\in l$ - and Lemma \ref{nol1} now follows.
\end{proof}
\begin{figure}[h]
\centering
\begin{overpic}[width=0.65\textwidth]{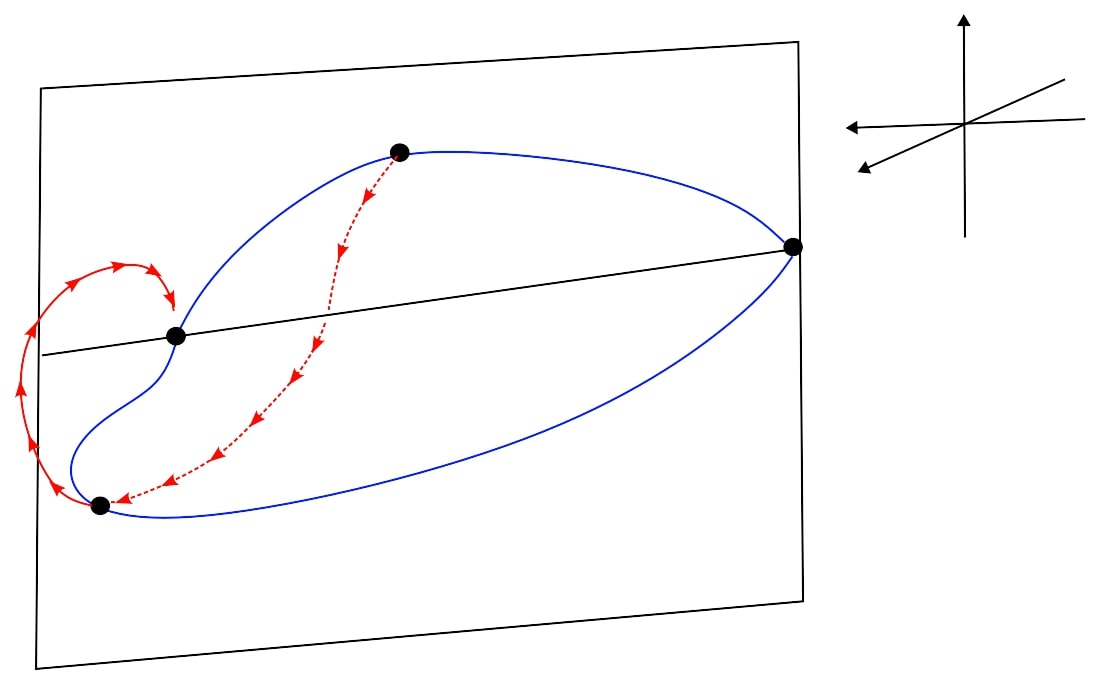}
\put(430,80){$L_p$}
\put(150,270){$s_{In}$}
\put(530,530){$U_p$}
\put(875,630){$z$}
\put(760,465){$y$}
\put(755,515){$x$}
\put(350,290){$l$}
\put(350,220){$\rho_{In}$}
\put(740,380){$P_{In}$}
\put(100,120){$s$}
\put(330,510){$s_1$}
\put(-250,205){$\{\dot{y}<0\}$}
\put(890,270){$\{\dot{y}>0\}$}
\put(440,420){$\beta_{In}$}
\end{overpic}
\caption[The case where $s_{In}\in l$.]{\textit{The case where $s_{In}\in l$ - as can be seen, there still exists some $s_1\in\beta_{In}$ which flows to $s\in L_p$, after which it flows to $s_{In}=f_p(s_1)$. The $y-$coordinate of $s$ is lesser than that of $s$ - which implies a contradiction.}}\label{d4}

\end{figure}
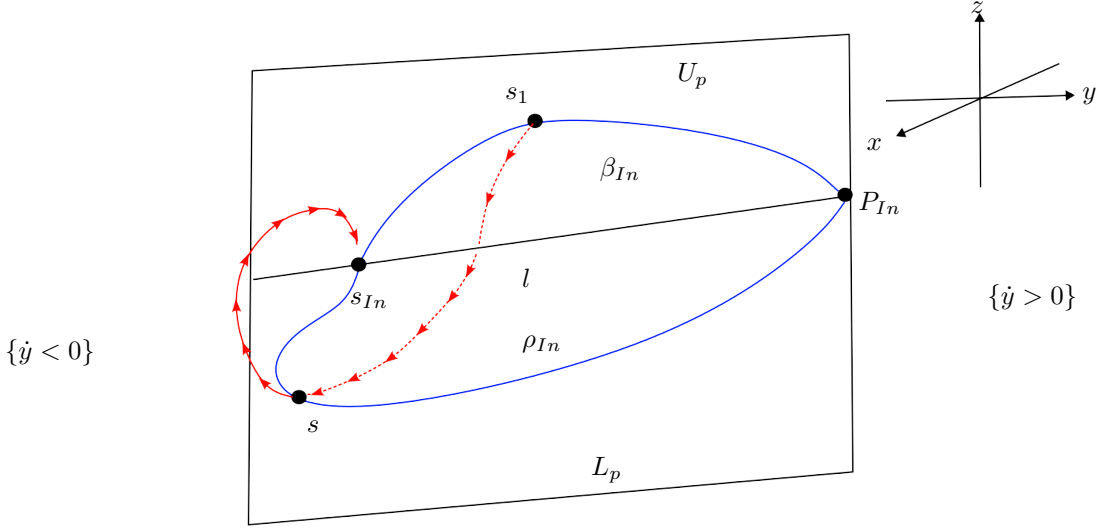

    Having proven Lemma \ref{nol1}, we now prove its analogue for $\beta_{Out}$ - the component of $W^s_{Out}\cap U_p$ which begins at $P_{Out}$ (see the illustration in Fig.\ref{Qq} and Fig.\ref{d5}). To do so, recall we denote by $s_{Out}$ the termination point of $\overline{B_{Out}}$ on $l_p$, s.t. $s_{Out}\ne P_{Out}$ (see the illustration in Fig.\ref{Qq} and Fig.\ref{d5}), and recall $l_2=\{l_p(x)|x>c-ab\}$ (see the discussion immediately before Prop.\ref{arccor}). We now prove:
    \begin{lemma}
        \label{nol}
 Whenever $F_p$ is a trefoil parameter, we have $s_{Out}\not\in l_2\cup l$.     \end{lemma}
    \begin{proof}
        We prove Lemma \ref{nol} using similar arguments to those used to prove Lemma \ref{nol1} - again, we first prove $s_{Out}\not\in l_2$, after which we prove $s_{Out}\not\in l$ (both by contradiction).\\
        
To begin, assume by contradiction $s_{Out}\in l_2$, as illustrated in Fig.\ref{d5}. Since $\beta_{Out}\subseteq W^s_{Out}\cap U_p$, because $W^s_{Out}$ is the two-dimensional stable manifold of $P_{Out}$ we conclude there exists some $s_1\in\beta_{Out}$ s.t. $f^{-1}_p(s_1)=s_{Out}$ - and by Lemma \ref{lem23} we know the backwards trajectory of $s_1$ arrives at $s_{Out}$ through $\{\dot{y}\leq0\}$ (passing first through the half-plane $L_p$, as indicated in Fig.\ref{d5}). Similarly to the proof of Lemma \ref{nol1}, let $\alpha_{Out}$ denote the sub-arc of $\beta_{Out}$ beginning at $P_{Out}$ and terminating at $s_1$ - again, using similar arguments, it follows the inverse flow lines connecting $\alpha_{Out}$ and $\beta_{Out}$ trap a topological cone $C$ (see the illustration in Fig.\ref{d5}). It is easy to see $P_{In}\not\in\overline{C}$.\\
        
Now, set $s_{Out}=(x_{Out},y_{Out},z_{Out})$ and $H=\{(x,y_{Out},z)|x+ay_{Out}<0\}$. Using a similar argument to the one used to prove Lemma \ref{nol1} (applied to the inverse flow) it follows the backward trajectories of initial conditions $s\in W^s_{Out}$ sufficiently close to $s_1$ enter $C\cap\{(x,y,z)|y>y_{Out}\}$ under the inverse flow, and never escape it. By $P_{In}\not\in\overline{C}$, we conclude that for $s\in W^s_{Out}$ sufficiently close to $s_1$ we have $\lim_{n\to\infty}f^{-n}_p(s)\ne P_{In}$ - and by $W^s_{Out}=W^u_{In}$, we again derive a contradiction which implies $s_{Out}\not\in l_2$.\\
\begin{figure}[h]
\centering
\begin{overpic}[width=0.6\textwidth]{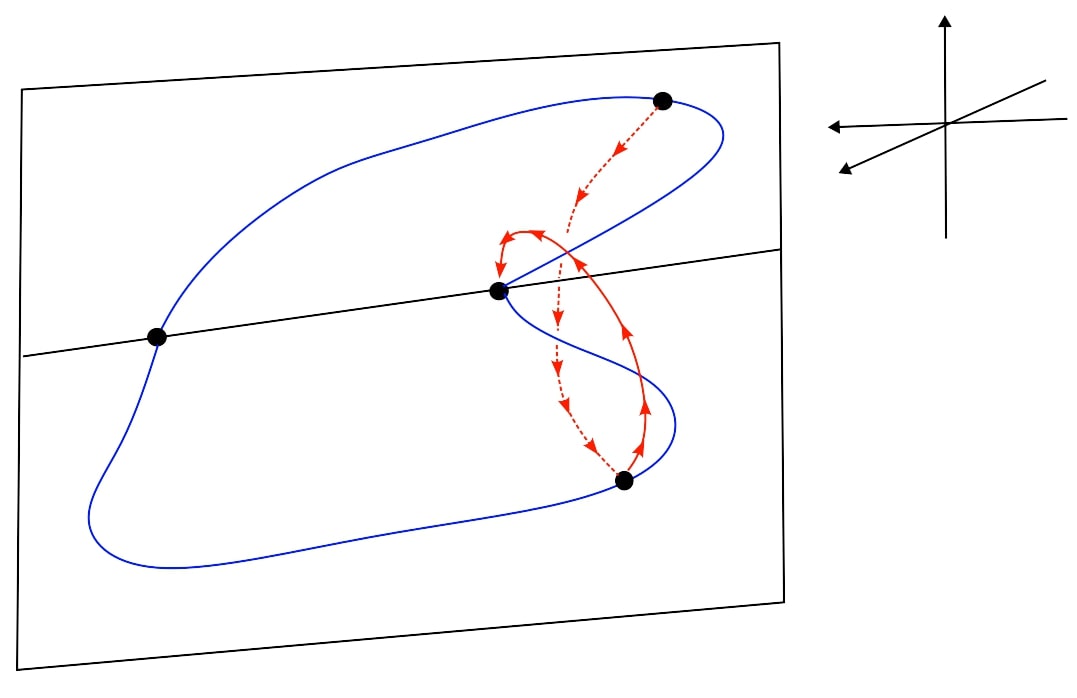}
\put(430,85){$L_p$}
\put(150,270){$P_{Out}$}
\put(600,570){$s_1$}
\put(865,640){$z$}
\put(750,460){$y$}
\put(740,525){$x$}
\put(740,380){$l_2$}
\put(550,140){$s$}
\put(130,510){$U_p$}
\put(-250,205){$\{\dot{y}>0\}$}
\put(890,270){$\{\dot{y}<0\}$}
\put(440,490){$\alpha_{Out}$}
\put(570,420){$\beta_{Out}$}
\put(400,320){$s_{Out}$}
\end{overpic}
\caption[The case where $s_{Out}\in l_2$.]{\textit{By $s_{Out}\in l_2$ there exists some $s_1\in\beta_{Out}$ which flows to $s\in L_p$, after which it flows (under the inverse flow) to $s_{Out}=f^{-1}_p(s_1)$ (the arc $\alpha_{Out}$ satisfies $f^{-1}_p(\alpha_{Out})=\beta_{Out}$). As $s$ flows to $s_{Out}$ through $\{\dot{y}<0\}$ under the inverse flow, the $y-$coordinate of $s_{Out}$ is greater than that of $s$ - from which, using similar arguments we again derive a contradiction.}}\label{d5}

\end{figure}
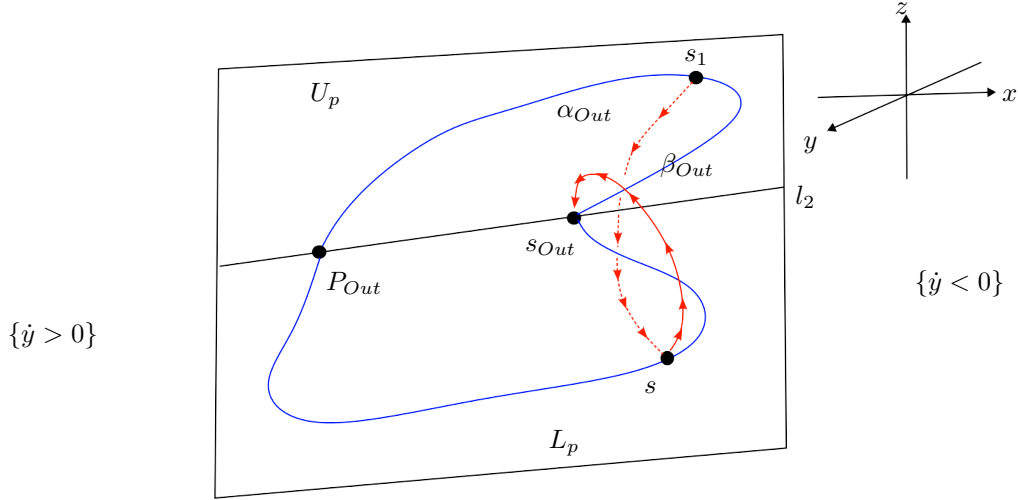

To conclude the proof of Lemma \ref{nol}, it remains to prove $s_{Out}\not\in l$. To do so, assume by contradiction we have $s_{Out}\in l$ - then, again, suspending $\beta_{Out}$ with the inverse flow (as indicated in Fig.\ref{d6}) we conclude there exists some $s_1\in \beta_{Out}$ s.t. $f^{-1}_p(s_1)=s_{Out}=(x_{Out},y_{Out},z_{Out})$. Defining $\alpha_{Out}\subseteq\beta_{Out}$ analogously, we again generate a topological cone $C$ with a tip at $P_{Out}$ - again, $P_{In}\not\in\overline{C}$ (see the illustration in Fig.\ref{d6}).\\

As the backwards trajectory connecting $s_1$ to $s_{Out}$ arrives from $\{\dot{y}\geq0\}$, using symmetric arguments, it follows the backward trajectories of initial conditions $s\in W^s_{Out}$ sufficiently close to $s_1$ eventually enter $C\cap\{(x,y,z)|\dot{y}>0,y<y_{Out}\}$  - and remain trapped in it forever. Again, by $W^s_{Out}=W^u_{In}$ we have a contradiction - which implies $s_{Out}\not\in l$ as well and Lemma \ref{nol} now follows.
    \end{proof}
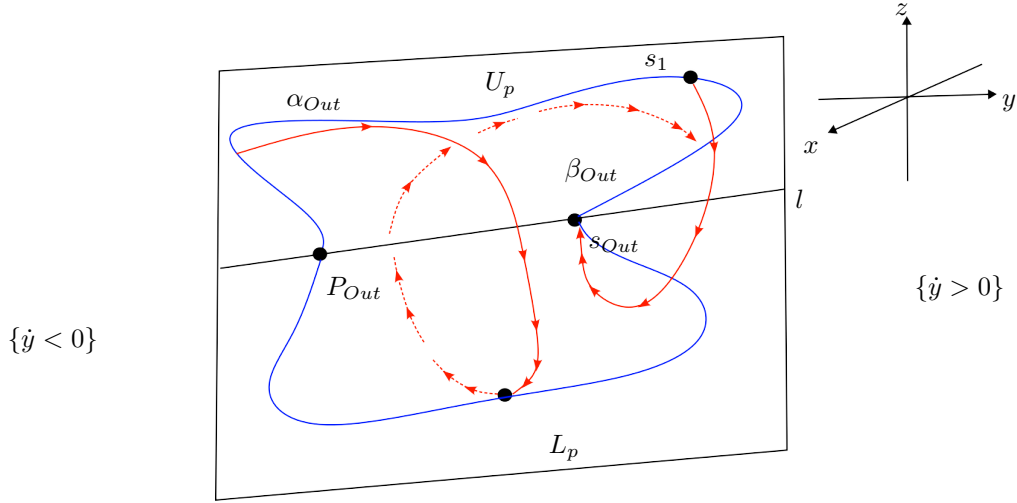
\begin{figure}[h]
\centering
\begin{overpic}[width=0.6\textwidth]{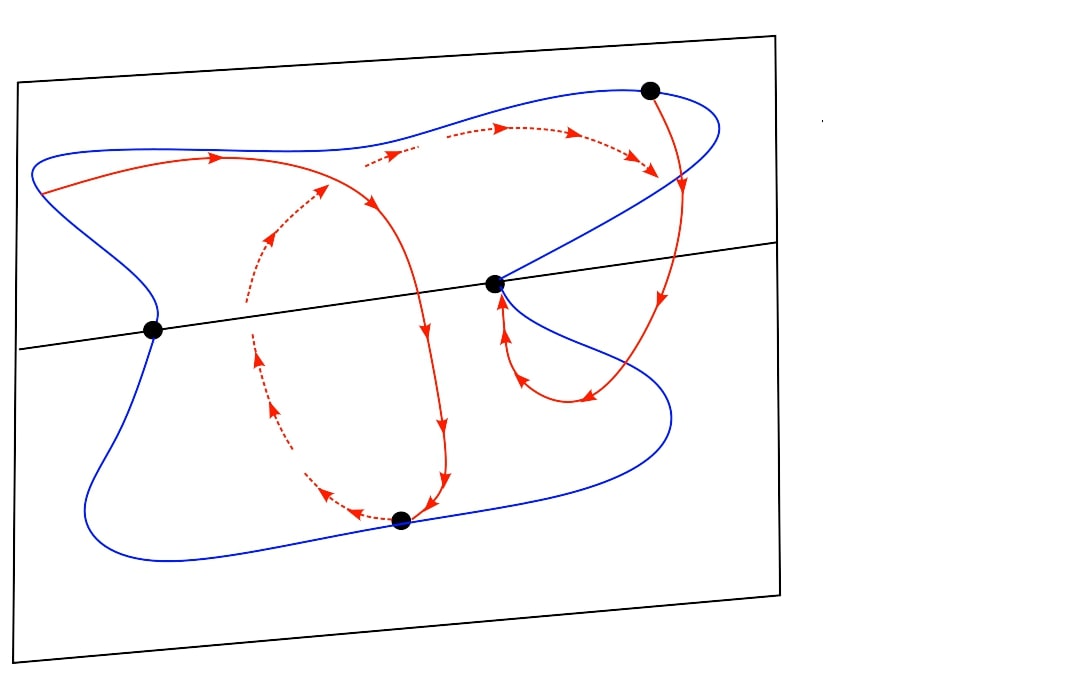}
\put(430,70){$L_p$}
\put(150,270){$P_{Out}$}
\put(350,530){$U_p$}
\put(740,380){$l$}
\put(100,510){$\alpha_{Out}$}
\put(-250,205){$\{\dot{y}<0\}$}
\put(890,270){$\{\dot{y}>0\}$}
\put(550,560){$s_1$}
\put(450,420){$\beta_{Out}$}
\put(480,330){$s_{Out}$}
\end{overpic}
\caption[The case where $s_{Out}\in l$.]{\textit{The case where $s_{Out}\in l$ - again, there must exist some $s_1\in\beta_{Out}$ which flows to $s_{Out}=f^{-1}_p(s_1)$ under the inverse flow (the arc $\alpha_{Out}$ satisfies $f^{-1}_p(\alpha_{Out})=\beta_{Out}$). The $y-$coordinate of $s_{Out}$ is greater than that of $s_1$ - which implies a contradiction.}}\label{d6}

\end{figure}

From Lemmas \ref{nol1} and \ref{nol}, we conclude the points $s_{In}$ and $s_{Out}$ lie outside of $l_1\cup l$ and $l_2\cup l$ (respectively). We can almost prove $\beta_{In}=\beta_{Out}$ - before doing so, we will also need the following fact, a corollary of both Lemma \ref{nol1} and Lemma \ref{nol}:
    \begin{corollary}
        \label{nol2}
        $s_{In}\not\in l_2$ and $s_{Out}\not\in l_1$.
    \end{corollary}
    \begin{proof}
    We prove Cor.\ref{nol2} by contradiction. If $s_{In}\in l_2$, it immediately follows $\beta_{In}$ separates $\beta_{Out}$ from $l_1$ (see the illustration in Fig.\ref{d7}). This implies $\overline{\beta_{Out}}$ connects $P_{Out}$ and some $s_1\in l\cup l_2$, or, in other words, $s_{Out}\in l\cup l_2$. Since this is impossible by Lemma \ref{nol}, we must have $s_{In}\not\in l_2$. Using a symmetric argument we conclude $s_{Out}\not\in l_1$ and the assertion follows.
    \end{proof}
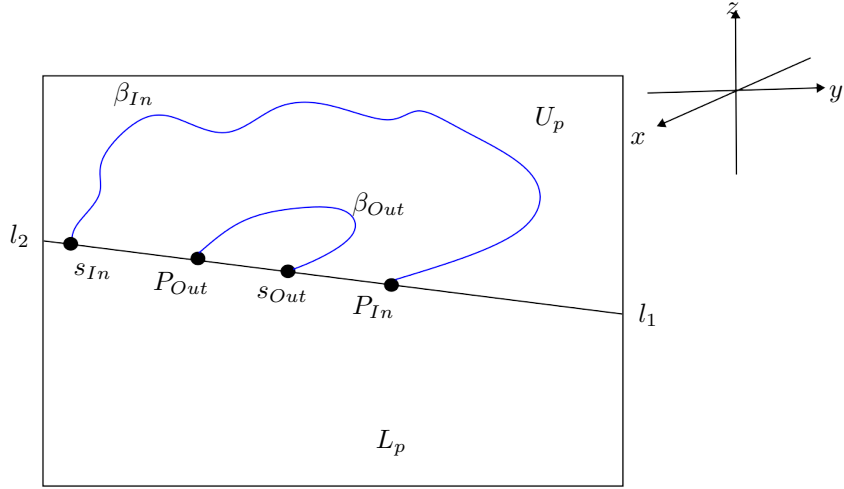
\begin{figure}[h]
\centering
\begin{overpic}[width=0.6\textwidth]{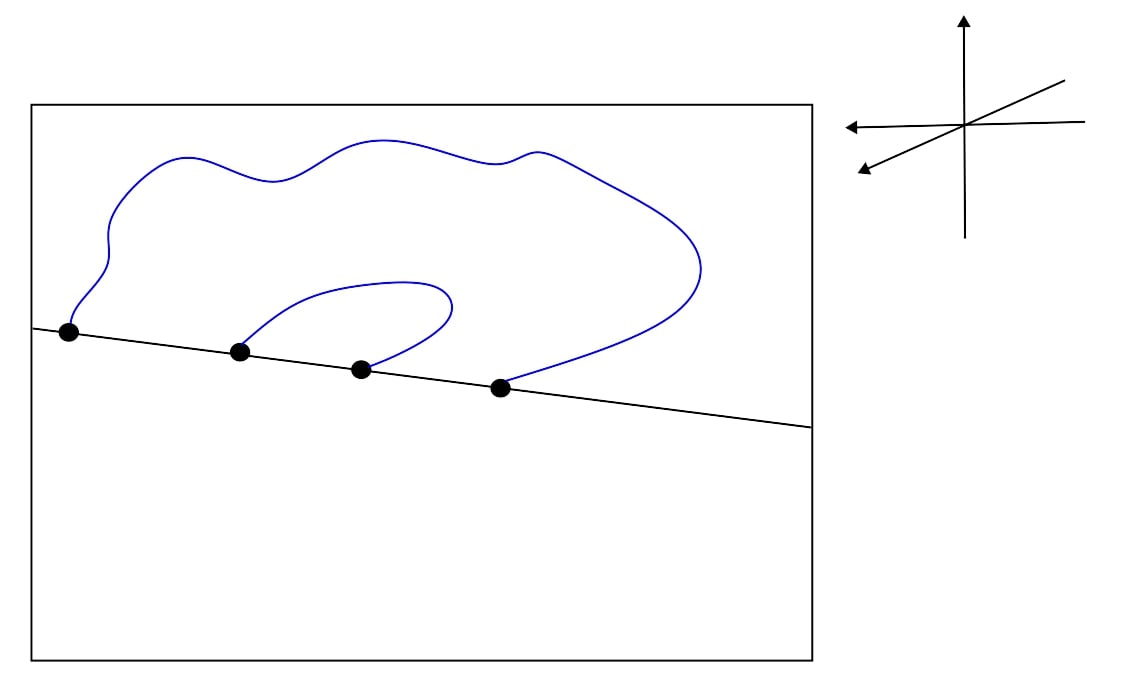}
\put(430,70){$L_p$}
\put(150,270){$P_{Out}$}
\put(850,600){$z$}
\put(730,485){$x$}
\put(740,450){$y$}
\put(760,230){$l_1$}
\put(-30,330){$l_2$}
\put(100,490){$\beta_{In}$}
\put(630,480){$U_p$}
\put(400,370){$\beta_{Out}$}
\put(400,240){$P_{In}$}
\put(280,260){$s_{Out}$}
\put(50,280){$s_{In}$}
\end{overpic}
\caption[The case where $s_{In}\in l_2$.]{\textit{The case where $s_{In}\in l_2$ - as $\beta_{In}$ would then separate $\beta_{Out}$ and $l_1$, it would force $s_{Out}$ to be in either $l_2$ or $l$ (where $l$ is the open arc connecting $P_{In}$ and $P_{Out}$). This contradicts Lemma \ref{nol}.}}\label{d7}

\end{figure}

    Having proven Lemmas \ref{nol1}, \ref{nol} and Cor.\ref{nol2}, we now conclude the proof of Prop.\ref{arccor}. To do so, let us first remark that since $F_p$ is a trefoil parameter there are no homoclinic trajectories - which implies $\beta_{In}$ cannot be a closed loop which begins and terminates at $P_{In}$. Consequentially, as $\overline{\beta_{In}}$ cannot connect $P_{In}$ through $U_p$ with either $l_1,l_2,l$ or $P_{In}$, the only possibility is that it terminates at $P_{Out}$, i.e., $s_{In}=P_{Out}$. This immediately implies $L=\beta_{In}=\beta_{Out}=W^u_{In}\cap U_p=W^s_{Out}\cap U_p$ is a simple curve in $U_p$, connecting $P_{In}$ and $P_{Out}$ through $U_p$, as illustrated in Fig.\ref{Qq} and Fig.\ref{trefint}.\\
    
    It now immediately follows the union $T=L\cup l\cup\{P_{In},P_{Out}\}$ is a Jordan curve on $\overline{U_p}$, which implies $U_p\setminus T$ consists of two components - a bounded and an unbounded set. As $B_\alpha$ is bounded we conclude the bounded set is $D_\alpha=B_\alpha\cap U_p$ which implies $D_\alpha$ is a topological disc on $U_p$ as illustrated in Fig.\ref{trefint}. The proof of Prop.\ref{arccor} is now complete.
\end{proof}

\subsection{The discontinuity properties of the first-return map}
\label{sect32}

\begin{figure}[h]
\centering
\begin{overpic}[width=0.6\textwidth]{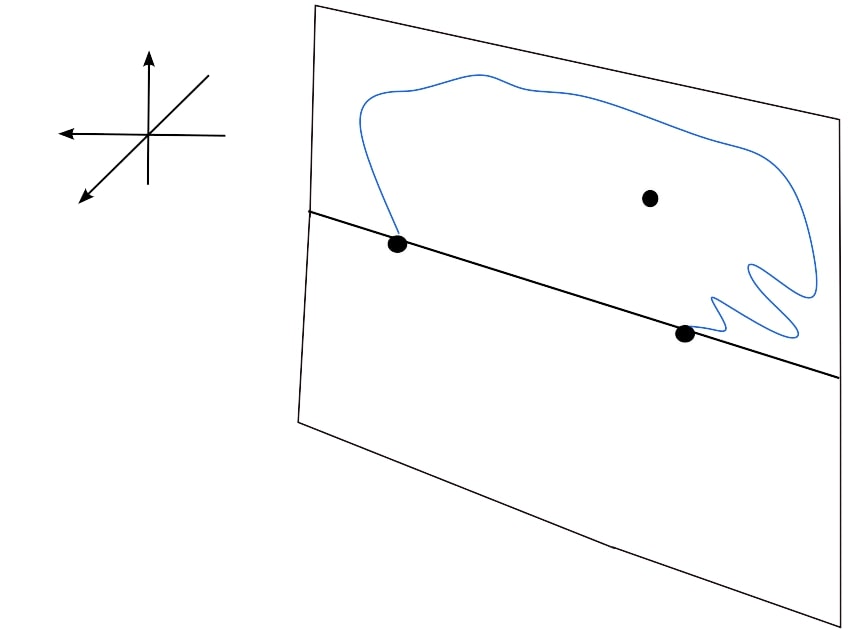}
\put(700,340){$P_{In}$}
\put(670,530){$P_0$}
\put(1050,200){$\{\dot{y}<0\}$}
\put(100,150){$\{\dot{y}>0\}$}
\put(360,440){$P_{Out}$}
\put(570,450){$l$}
\put(770,610){$L$}
\put(330,750){$U_p$}
\put(460,550){$D_\alpha$}
\put(700,200){$L_p$}
\put(45,580){$x$}
\put(60,490){$y$}
\put(160,690){$z$}
\end{overpic}
\caption[The disc $D_\alpha$.]{\textit{The Jordan domain $D_\alpha\subseteq U_p$, bounded by the blue curve $L=\partial B_\alpha\cap U_p$ (see Def.\ref{def32}), the fixed points, and the boundary arc $l$. $l$ lies on the line $l_p$, the tangency set of $F_p$ to $\{\dot{y}=0\}$ which separates $L_p$ from $U_p$.}}
\label{trefint}
\end{figure}
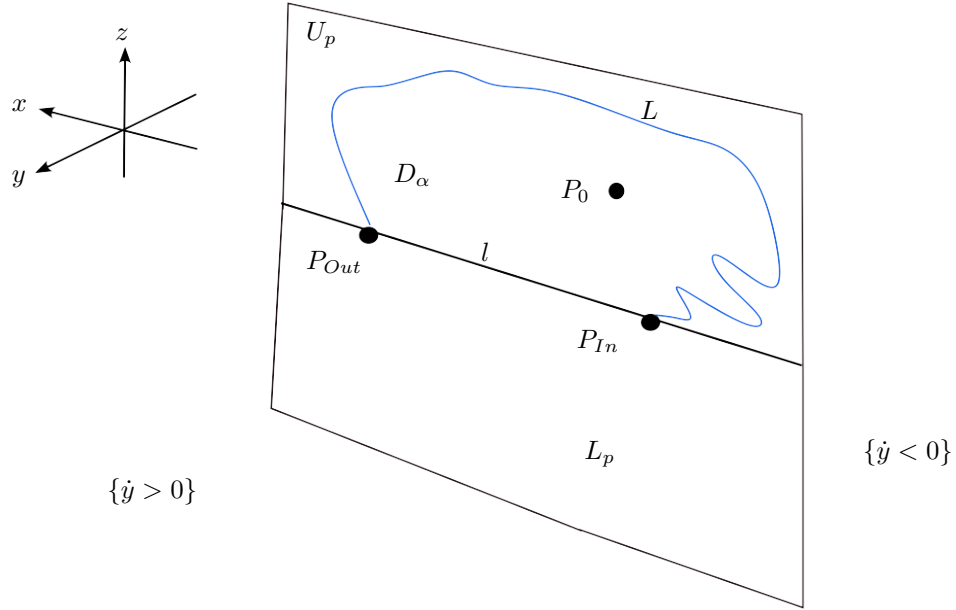

Having proven the existence of the first-return $f_p:\overline{D_\alpha}\setminus\{P_0\}\to\overline{D_\alpha}\setminus\{P_0\}$ map at trefoil parameters and that $D_\alpha$ is a topological disc, in this subsection we now study its continuity and discontinuity properties of $f_p$.\\

To begin, let $F_p$ be a trefoil parameter. Recall the heteroclinic trefoil knot intersects the cross-section $U_p$ at a single point, $P_0$, interior to the topological disc $D_\alpha$ (see the illustration in Fig.\ref{trefint}). We begin with the following technical result, where we characterize the discontinuity set of $f_p$ in the punctured disc $D_\alpha\setminus\{P_0\}$:

\begin{proposition}\label{lem33}
Let $F_p$ be a trefoil parameter for the Rössler system. Then, the first-return map $f_p:\overline{D_\alpha}\setminus\{P_0\}\to\overline{D_\alpha}\setminus\{P_0\}$ is continuous at the fixed-point $P_{In}$ and on the curve $L$ - moreover, the discontinuities of $f_p$ in the punctured topological disc $D_\alpha\setminus\{P_0\}$, denoted by $Dis(f_p)$, are given by $f^{-1}_p(l)\cap D_\alpha$. In addition, there exists a curve $\delta\subseteq Dis(f_p)$ satisfying:
\begin{itemize}
    \item $\delta$ is a component of $Dis(f_p)$, and it is homeomorphic to an open interval.
    \item $\delta$ has two endpoints: $P_0$ and some point $\delta_0\in \overline{l}$. Moreover, $\delta_0\ne P_{In}$.
    \item $P_{In}\in\overline{f_p(\delta)}$.
    \item If $\rho$ is a component of $Dis(f_p)$ s.t. $\rho\ne\delta$, it is a curve with (at least) two endpoints on $\overline{l}$. Moreover, $P_0\not\in\overline{\rho}$.
\end{itemize}
\end{proposition}
\begin{proof} 
First, recall that since $p$ is a trefoil parameter the two-dimensional invariant manifolds $W^u_{In},W^s_{Out}$ coincide. Moreover, recall that by Prop.\ref{arccor} $D_\alpha$ is an open topological disc bounded by the Jordan curve $l\cup L\cup\{P_{In},P_{Out}\}$ - and that $L=W^u_{In}\cap U_p$ is a curve connecting $P_{In}$ and $P_{Out}$ through $U_p$ (see the illustration in Fig.\ref{trefint}). By the invariance of $W^u_{In}$ under the flow we have $f_p(L)=L$.\\

We now prove $f_p$ is continuous on $L$. To see why, note that by $L\subseteq U_p$ it follows $L$ lies away from $l_p$, the tangency curve of the vector field $F_p$ to the plane $\{\dot{y}=0\}$ - see the illustration in Fig.\ref{trefint}. Therefore, for every $s\in L$ the flow line connecting $s$ to $f_p(s)\in L$ is transverse to the cross-section $U_p$ at $f_p(s)$. Consequentially, the first-return map $f_p$ is continuous throughout $L$.\\

Using a similar argument, we now prove $f_{p}$ is continuous around $P_{In}$ (in $\overline{D_\alpha}$) - to do so, recall that by Lemma \ref{cor211} $W^u_{In}$, is transverse to $\overline{U_p}$ at $P_{In}$. Therefore, since $P_{In}$ is a saddle-focus initial conditions on the arc $l_p$ sufficiently close to $P_{In}$ are mapped inside $U_p$ by $f_{p}$ (see the illustration in Fig.\ref{connn}). By the orientation-preserving property of the flow and because both $f_p(L)=L$ and $f_p(P_{In})=P_{In}$, it follows $f_p$ squeezes some sector in $\overline{D_\alpha}$ trapped between $L$ and $l$ inside $\overline{D_\alpha}$ (see the illustration in Fig.\ref{connn}). Consequentially, $f_p$ is continuous on some neighborhood of $P_{In}$ in $\overline{D_\alpha}$.\\

\begin{figure}[h]
\centering
\begin{overpic}[width=0.6\textwidth]{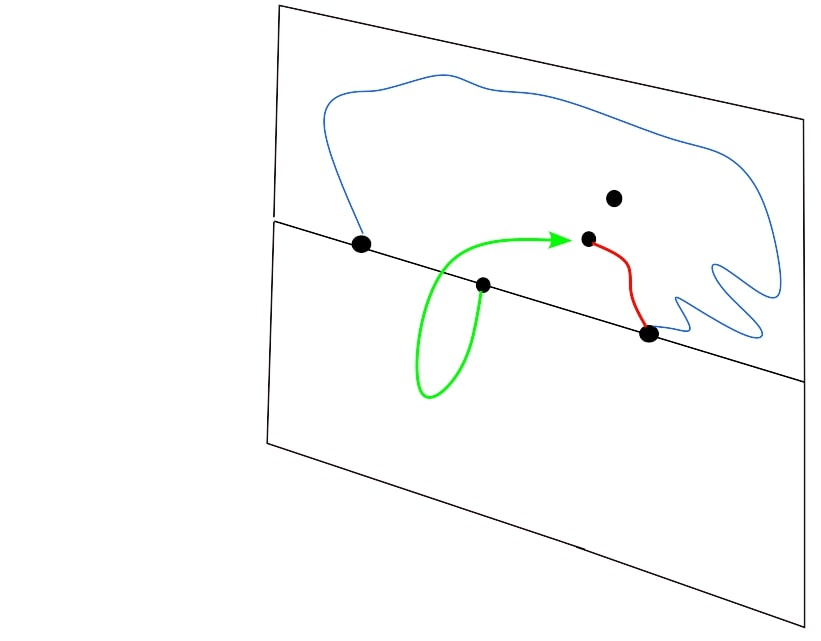}
\put(700,350){$P_{In}$}
\put(690,570){$P_0$}
\put(1050,200){$\{\dot{y}<0\}$}
\put(100,150){$\{\dot{y}>0\}$}
\put(540,400){$s_1$}
\put(680,510){$s_2$}
\put(340,450){$P_{Out}$}
\put(680,420){$l$}
\put(460,530){$D_\alpha$}
\put(370,720){$U_p$}
\put(700,150){$L_p$}
\put(920,570){$L$}
\end{overpic}
\caption[The trajectories of initial conditions close to $P_{In}$.]{\textit{The trajectory of an initial condition $s_1\in l$ which is sufficiently close to $P_{In}$ - because $P_{In}$ is a saddle-focus, the trajectory of $s_1$ hits the cross-section transversely. The red curve is the image under $f_p$ of the arc connecting $P_{In}$ and $s_1$.}}
\label{connn}
\end{figure}
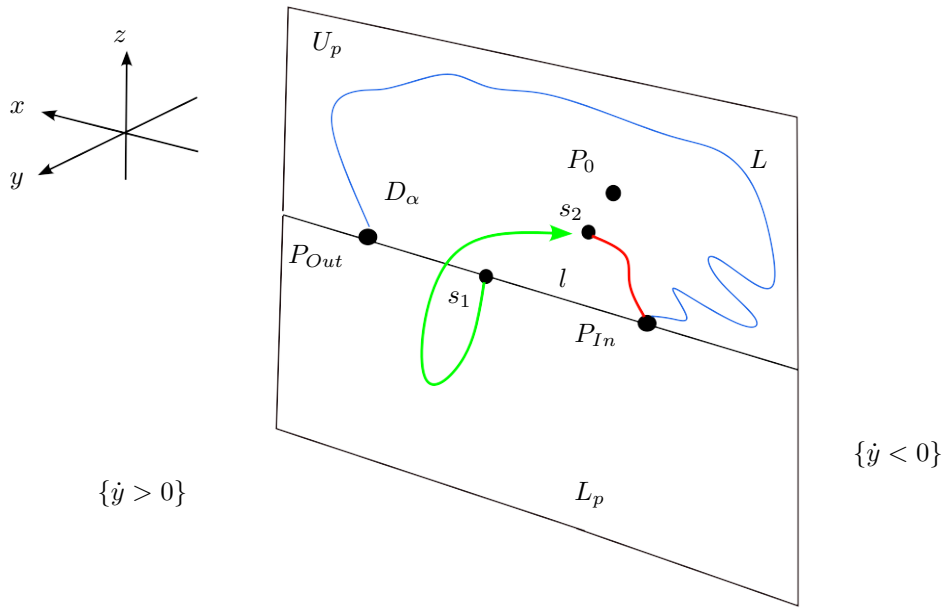

We now prove $Dis(f_p)=f^{-1}_p(l)\cap D_\alpha$. We first prove that for an initial condition $s\in D_\alpha\setminus\{P_0\}$ to be a discontinuity point for $f_p$, a necessary condition is that $f_p(s)$ lies in $l$. To see why this is so, recall the following facts:

\begin{itemize}
    \item By Lemma \ref{noempty}, the fixed points $P_{In}$ and $P_{Out}$ lie in $\partial D_\alpha$. As such, any $s\in D_\alpha\setminus\{P_0\}$ is not a fixed point for the flow - hence, by Lemma \ref{firstret} $f_p(s)\ne P_{In},P_{Out}$.
    \item   By the continuity of $f_p$ on $L\subseteq\partial D_\alpha$, there are no discontinuity points for $f_p$ on $L$.
    \item Since the vector field $F_p$ is transverse to $U_p$ it is also transverse to the (open) punctured disc $D_\alpha\setminus\{P_0\}\subseteq U_p$ - which proves that whenever $f_p(s)\in D_\alpha\setminus\{P_0\}$, $f_p$ is continuous at $s$.
    \end{itemize}

By $\overline{D_\alpha}\setminus\{P_0\}=L\cup l\cup\{P_{In},P_{Out}\}\cup(D_\alpha\setminus\{P_0\})$ we conclude that whenever $f_p$ is discontinuous at $s\in D_\alpha\setminus\{P_0\}$, the only possibility is $f_p(s)\in l$. Now, recall $l\subseteq l_p$ - where $l_p$ is the tangency set of $F_p$ to $\overline{U_p}$ (see the discussion immediately before Lemma \ref{obs}). This implies that given $s\in D_\alpha\setminus\{P_0\}$ whose trajectory is tangent to $\overline{D_\alpha}$ at $f_p(s)\in l$, there exists a two-dimensional disc on $D_\alpha$, centered at $s$, whose image under $f_p$ is torn in two as illustrated in Fig.\ref{Disc}. In other words, whenever $s\in D_\alpha\setminus\{P_0\}$ satisfies $f_p(s)\in l$, $f_p$ is discontinuous at $s$.\\

Therefore, since the condition $f_p(s)\in l$ is both sufficient and necessary for an initial condition $s\in D_\alpha\setminus\{P_0\}$ to be a discontinuity for $f_p$, we conclude the discontinuity set of $f_p$ in $D_\alpha\setminus\{P_0\}$ is given by $f^{-1}_p(l)\cap D_\alpha=Dis(f_p)$. Furthermore, by $l\cap L=\emptyset$ and by the invariance of $L$ it is also immediate that $L$ does not intersect the closure of $Dis(f_p)$.\\

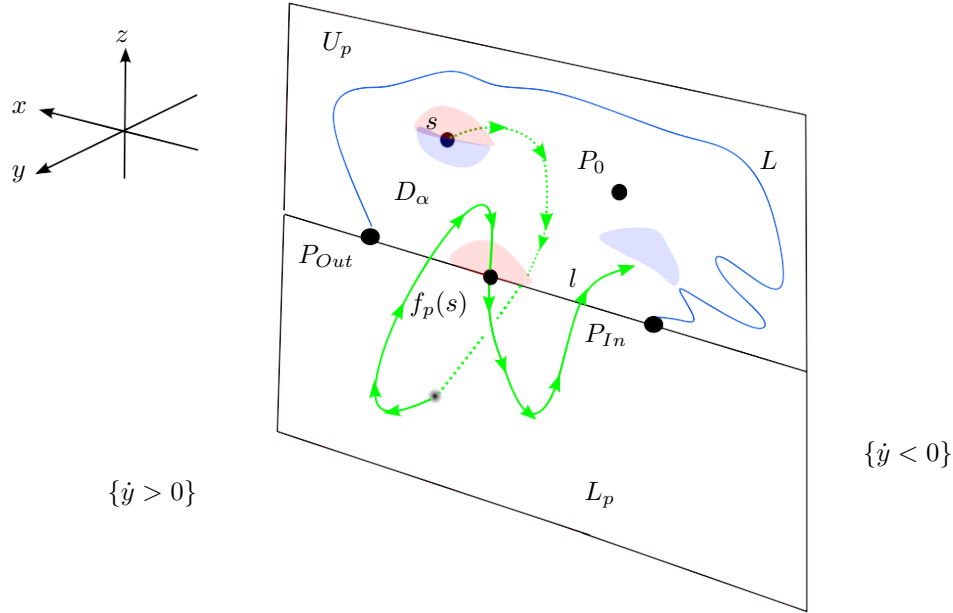
\begin{figure}[h]
\centering
\begin{overpic}[width=0.6\textwidth]{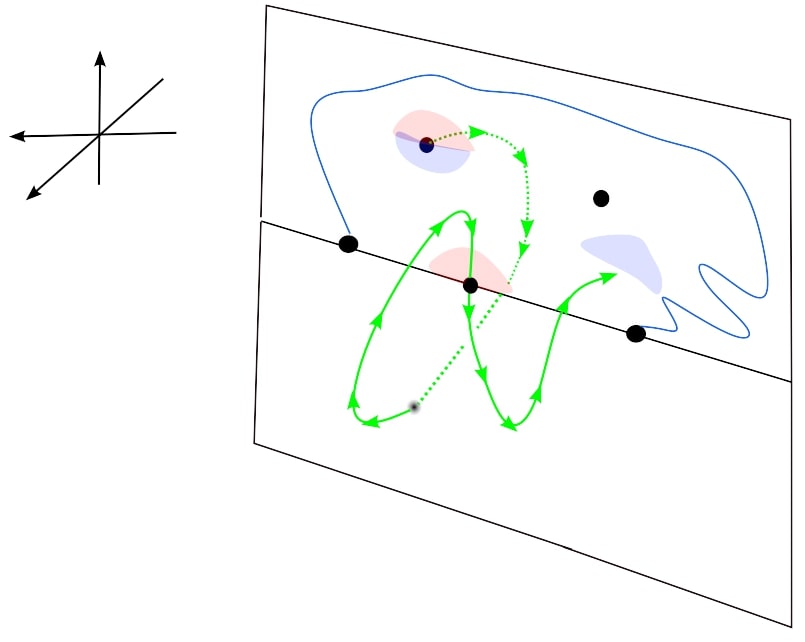}
\put(700,350){$P_{In}$}
\put(690,570){$P_0$}
\put(1050,200){$\{\dot{y}<0\}$}
\put(100,150){$\{\dot{y}>0\}$}
\put(480,390){$f_p(s)$}
\put(340,450){$P_{Out}$}
\put(680,420){$l$}
\put(460,530){$D_\alpha$}
\put(370,720){$U_p$}
\put(500,620){$s$}
\put(700,150){$L_p$}
\put(-20,620){$x$}
\put(0,510){$y$}
\put(115,740){$z$}
\put(920,570){$L$}
\end{overpic}
\caption[Discontinuities for the first-return map $f_p$.]{\textit{The trajectory of an interior point $s\in D_\alpha$ s.t. $f_p(s)\in l$. Typically, this would mean there are nearby trajectories that hit $D_\alpha$ in its interior when the trajectory of $s$ hits the boundary at $f_p(s)$, and other trajectories which continue with the trajectory of $s$ and hit the interior of $D_\alpha$ around $f^2_p(s)$.}}
\label{Disc}
\end{figure}

We now prove the existence of the discontinuity curve $\delta$ as posited above. To do so, consider some closed loop $\zeta$ in $D_\alpha$ surrounding $P_0$. Since $p$ is a trefoil parameter, the trajectory of $P_0$ flows to $P_{In}$ in infinite time without ever hitting $\overline{D_\alpha}$ along the way (see Def.\ref{def32}) - hence, $f_{p}$ cannot be continuous on $\zeta$: for if it were continuous on $\zeta$, $f_{p}(\zeta)$ would be a closed curve in $\overline{D_\alpha}$ surrounding $P_{In}$. Because this is impossible by $P_{In}\in\partial D_\alpha$, from the discussion above we conclude there must exist a discontinuity point $\zeta_0\in\zeta$ - i.e., there exists some $\zeta_0\in\zeta$  s.t. $f_{p}(\zeta_0)\in l$.\\

Varying the loop $\zeta$ continuously in $D_\alpha$ we construct an open arc $\delta$ of discontinuity points, whose closure $\overline{\delta}$ is a curve connecting $P_0$ and some unique point $\delta_0$ in $\overline{l}$ (see the illustration in Fig.\ref{fig8} and Fig.\ref{imagedelta}). As $P_0$ tends to $P_{In}$ in infinite time, it follows $f_p(\delta)$ is an arc on $l$ s.t. $P_{In}\in\overline{f_p(\delta)}$ - similarly, we also have $f^j_p(\delta_0)\in\overline{f_p(\delta)}\subseteq \overline{l}$ for some $j\geq1$ (see the illustration in Fig.\ref{imagedelta}). It is easy to see that because $f_p$ is continuous around $P_{In}$, it follows $\delta_0\ne P_{In}$ - and moreover, by $Dis(f_p)\cap L=\emptyset$ we have $L\cap\delta=\emptyset$. \\

\begin{figure}[h]
\centering
\begin{overpic}[width=0.8\textwidth]{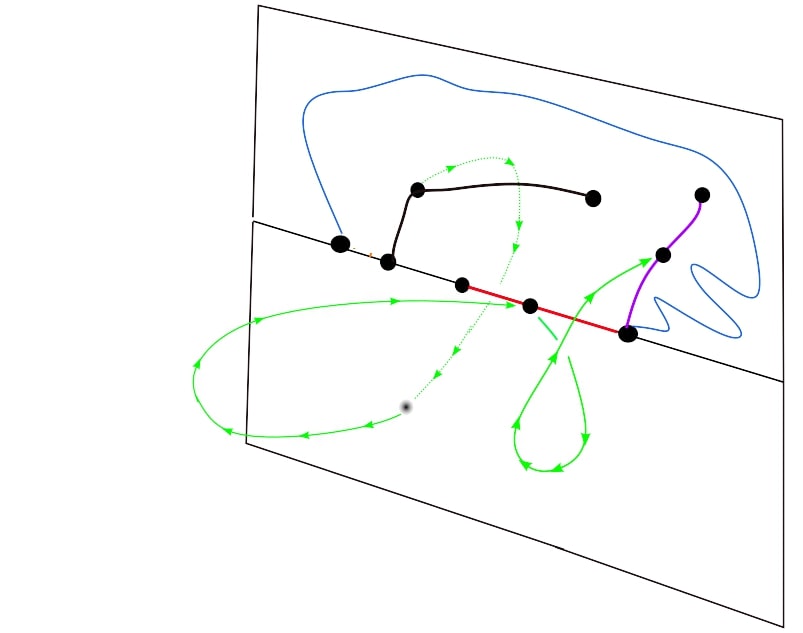}
\put(720,350){$P_{In}$}
\put(690,500){$P_0$}
\put(1000,200){$\{\dot{y}<0\}$}
\put(100,150){$\{\dot{y}>0\}$}
\put(580,370){$f_p(s)$}
\put(340,460){$P_{Out}$}
\put(450,430){$\delta_0$}
\put(600,570){$\delta$}
\put(680,420){$f_p(\delta)$}
\put(540,460){$f_p(\delta_0)$}
\put(870,480){$f^2_p(s)$}
\put(900,550){$f^2_p(\delta_0)$}
\put(490,550){$s$}
\put(370,730){$U_p$}
\put(460,630){$D_\alpha$}
\put(560,695){$L$}
\put(700,150){$L_p$}
\end{overpic}
\caption[Discontinuities for the FRM, 1]{\textit{The trajectory of some $s\in \delta$, $s\ne\delta_0$. The flow line connecting $s,f_p(s)$ and $f^2_p(s)$ lies in $\{\dot{y}\geq0\}$, hence $f^2_p(s)\not\in l$ (where $l$ is the arc connecting $P_{In},P_{Out}$ - moreover, $f_p(\delta)$ is a sub-arc of $l$, and $f^2_p(\delta)$, the purple curve, lies away from $l$). }}
\label{imagedelta}
\end{figure}

We now prove that $f^{-2}_p(l)\cap \delta=\emptyset$, which we will later use to prove $\delta$ is a component of $Dis(f_p)$. To do so, first note that for any $s\in l$ the flow line connecting $s,f_p(s)$ lies strictly inside $\{\dot{y}\geq0\}$ (see the illustration in Fig.\ref{imagedelta}). This proves that for $s\in l$, the $y$-coordinate of $f_p(s)$ is strictly greater than that of $s$ - and by $f_p(\delta)\subseteq l$ we conclude that for $s\in \delta$, the $y$-coordinate of $f^2_p(s)$ is greater than that of $f_p(s)$. Second, by $\delta\subseteq D_\alpha$ and $f_p(\delta)\subseteq\partial D_\alpha$, as $D_\alpha$ is an open topological disc we conclude $f_p(\delta)\cap\delta=\emptyset$ - which further yields $f^2_p(\delta)\cap f_p(\delta)=\emptyset$ (see the illustration in Fig.\ref{imagedelta}).\\

Now, recall we parameterize $l$ by $(x,-\frac{x}{a},\frac{x}{a})$, $x\in(0,c-ab)$ - by $P_{In}=(0,0,0)$ it follows the $y$-coordinate on $l$ increases monotonically as $s\in l$ tends to $P_{In}$. Therefore, because $\overline{f_p(\delta)}$ is an arc with one endpoint at $P_{In}$, given $v\in l\setminus f_p(\delta)$ and $\mu\in f_p(\delta)$ the $y-$coordinate of $v$ is smaller that that of $\mu$ - and since by previous paragraph for every $s\in\delta$ the $y$-coordinate of $f^2_p(s)$ is greater than that of $f_p(s)$, we conclude $f^2_p(\delta)\cap l\setminus f_p(\delta)=\emptyset$. Combined with $f_p(\delta)\cap f^2_p(\delta)=\emptyset$, it follows $f^2_p(s)\not\in l$, i.e., $\delta\cap f^{-2}_p(l)=\emptyset$ (see the illustration in Fig.\ref{imagedelta}). \\

\begin{figure}[h]
\centering
\begin{overpic}[width=0.7\textwidth]{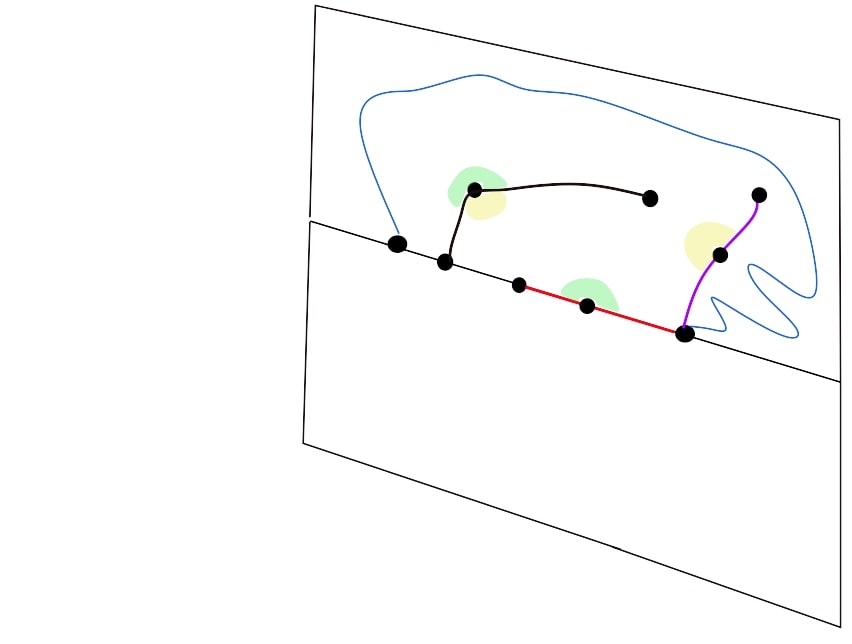}
\put(800,310){$P_{In}$}
\put(720,540){$P_0$}
\put(1050,200){$\{\dot{y}<0\}$}
\put(100,150){$\{\dot{y}>0\}$}
\put(630,350){$f_p(s)$}
\put(420,430){$P_{Out}$}
\put(510,400){$\delta_0$}
\put(600,540){$\delta$}
\put(450,590){$D_\alpha$}
\put(870,450){$f^2_p(s)$}
\put(900,540){$f^2_p(\delta_0)$}
\put(530,535){$s$}
\put(400,680){$U_p$}
\put(560,675){$L$}
\put(800,150){$L_p$}

\end{overpic}
\caption[The trajectories of initial conditions around $\delta$.]{\textit{The trajectories on initial conditions $s\in\delta$ - the green region hits the green half-disc above $f_p(\delta)$ (denoted by the red arc), while the yellow half-disc hits close to $f^2_p(s)$ (which lies on $f^2_p(\delta)$, denoted by the purple arc). }}
\label{imagedelta2}

\end{figure}
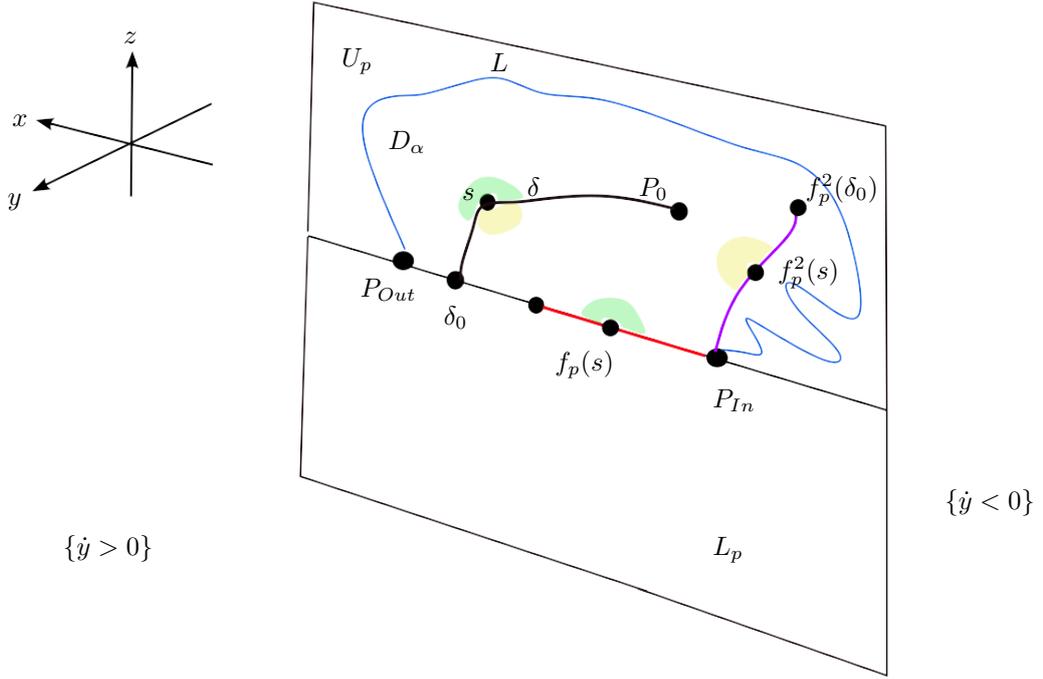

Consequentially, given $s\in\delta$, the trajectories of initial conditions $\omega\not\in\delta$ sufficiently close to $s$ hit $D_\alpha\setminus\{P_0\}$ either above the arc $f_p(\delta)$ and close to $f_p(s)$, or, alternatively, close to $f^2_p(s)$, hence away from $l$ (see the illustration in Fig.\ref{imagedelta2}). Therefore, for every $s\in\delta$ there exists a neighborhood $N_s$ s.t. $N_s\setminus\delta$ consists of precisely two components - and both components of $N_s\setminus\delta$ do not intersect $Dis(f_p)$, i.e., $f_p$ is continuous on both components (see the illustration in Fig.\ref{imagedelta2}). Therefore, by $Dis(f_p)=f^{-1}_p(l)\setminus l=f^{-1}_p(l)\cap D_\alpha\setminus\{P_0\}$, we conclude $\delta$ is a component of $Dis(f_p)$ homeomorphic to an open interval with endpoints $P_0$ and $\delta_0$ (see the illustration in Fig.\ref{imagedelta} and \ref{fig8}).\\

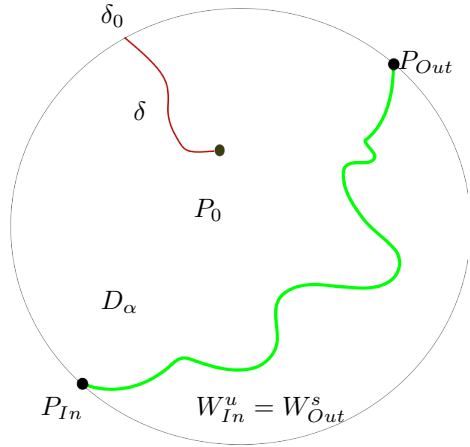
\begin{figure}[h]
\centering
\begin{overpic}[width=0.35\textwidth]{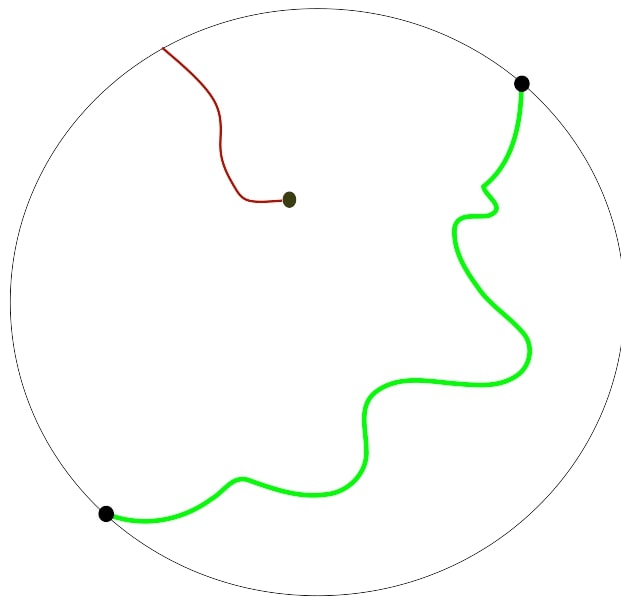}
\put(80,70){$P_{In}$}
\put(395,70){ $W^u_{In}= W^s_{Out}$}
\put(410,500){$P_0$}
\put(210,300){$D_\alpha$}
\put(210,920){$\delta_0$}
\put(280,710){$\delta$}
\put(850,820){$P_{Out}$}
\end{overpic}
\caption[The geography of $U_p$ and $D_\alpha$.]{\textit{The geography of the cross-section $U_p$ for a trefoil parameter $F_p$, sketched as a disc for simplicity (with the curve $\delta$ sketched in red). The green arc denotes $W^u_{In}\cap\overline{U_p}=W^s_{Out}\cap\overline{U_p}$. The set $D_\alpha$ corresponds to $B_\alpha\cap \overline{U_p}$.}}
\label{fig8}
\end{figure}

To conclude the proof of Lemma \ref{lem33} it remains to prove that given any component $\rho\subseteq Dis(f_p)$ s.t. $\rho\ne\delta$, then $\overline{\rho}$ is a curve with (at least) two endpoints on $l$, and $P_0\not\in\overline{\rho}$ - as illustrated in Fig.\ref{fig15}. We will do so by first proving that if $\overline{\rho}$ has an endpoint strictly interior to $D_\alpha$, then $\rho=\delta$. To begin, note that interior points in $D_\alpha\setminus\{P_0\}$ cannot form components in $Dis(f_{p})$: since given an interior point $s\in D_\alpha$ s.t. $f_p(s)\in l$, because $l$ is an open arc the continuity of the flow implies there exists a curve $\gamma\subseteq D_\alpha$ s.t, $f_{p}(\gamma)$ is a sub-arc on $l$ - and $s$ is interior to $\gamma$ (see the illustration in Fig.\ref{dosc}). Therefore, every component in $Dis(f_p)$ is a curve.\\

\begin{figure}[h]
\centering
\begin{overpic}[width=0.7\textwidth]{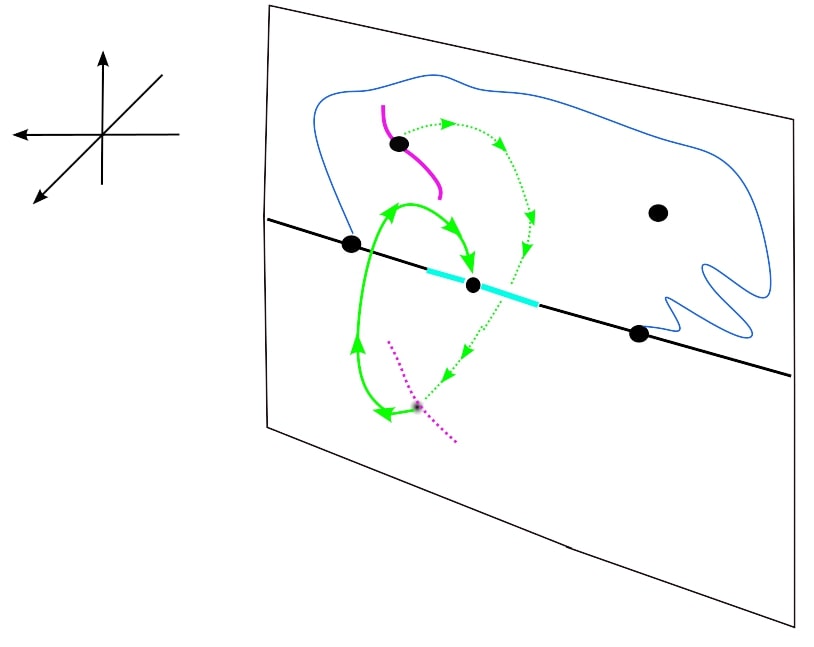}
\put(720,350){$P_{In}$}
\put(770,490){$P_0$}
\put(1050,200){$\{\dot{y}<0\}$}
\put(100,150){$\{\dot{y}>0\}$}
\put(500,400){$f_p(s)$}
\put(340,460){$P_{Out}$}
\put(480,570){$s$}
\put(370,730){$U_p$}
\put(700,580){$D_\alpha$}
\put(560,705){$L$}
\put(680,150){$L_p$}
\put(-5,625){$x$}
\put(10,530){$y$}
\put(115,740){$z$}
\end{overpic}
\caption[Interior points are not components of $Dis(f_p)$.]{\textit{$\gamma\subseteq Dis(f_p)$, the  pink curve, flows towards $l$, passes through $L_p$ at the dashed pink curve, after which it hits $l$ at $\theta=f_p(\gamma)$, the cyan arc. As can be seen, if $f_p(s)$ is interior to $\theta$, $s$ is interior to $\gamma$ (and vice versa).}}
\label{dosc}
\end{figure}
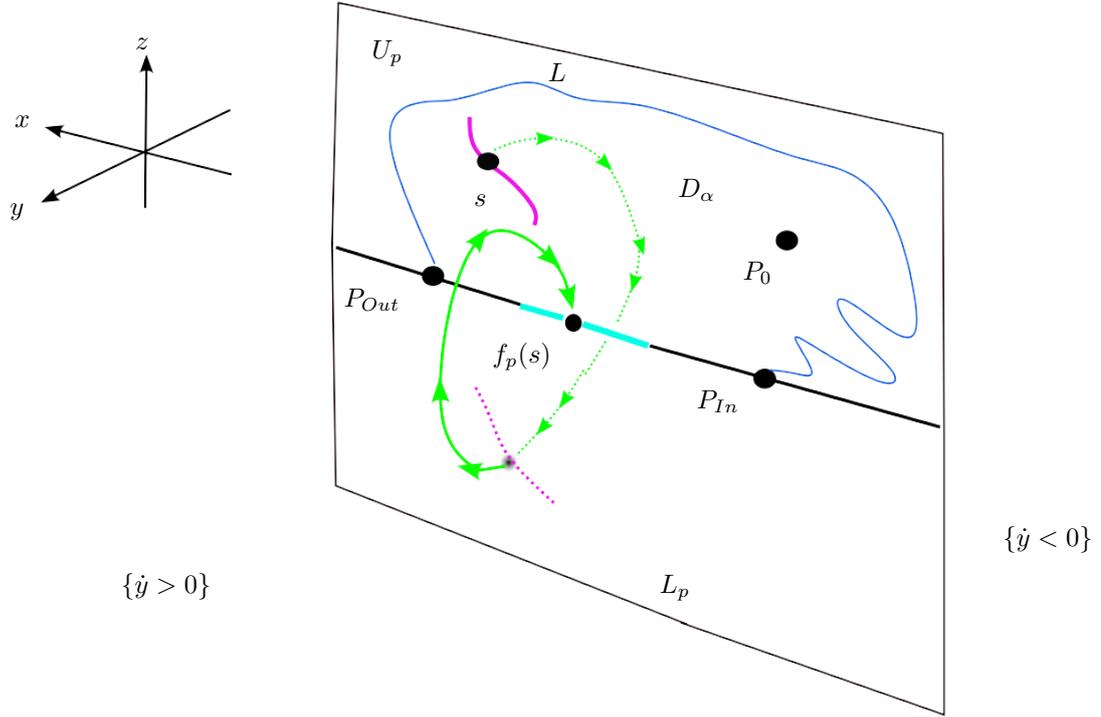

Now, let $\rho$ be a component of $Dis(f_p)$ - by previous paragraph, we already know $\rho$ is a curve. We now prove that if the closure of $\rho$ has an endpoint interior to $D_\alpha$, that endpoint is $P_0$. To do so, parameterize the closure of $\rho$, denoted by $\overline{\rho}$, by $[0,1]$ s.t. $\overline{\rho}(1)$ is the endpoint interior to $D_\alpha$ - since the fixed points are on the boundary of $D_\alpha$, $\overline{\rho}(1)$ is not a fixed point. Now, assume by contradiction that  $\overline{\rho}(1)\ne P_0$  - therefore, the trajectory of $\overline{\rho}(1)$ does not limit to $P_{In}$. Therefore, by Lemma \ref{firstret} $f_p(\overline{\rho}(1))$ is defined and not a fixed point - which, by the continuity of the flow and by $f_p(\rho)\subseteq l$, implies $f_p(\overline{\rho}(1))\in \overline{l}$ - and since $\overline{l}=l\cup\{P_{In},P_{Out}\}$ (see the illustration in Fig.\ref{trefint}), we conclude $f_p(\overline{\rho}(1))$ is interior to $l$.\\

Since $l$ is an open arc and because $\overline{\rho}(1)$ is interior to $D_\alpha$ per assumption, by the continuity of the inverse flow $f^{-1}_{p}$ is continuous on a neighborhood of $f_{p}(\overline{\rho}(1))$ inside the arc $l$ (see Fig.\ref{dosc}) - or, in other words, $\overline{\rho}(1)$ is interior to some component of $Dis(f_p)$. Because this contradicts the maximality of $\rho$ as a component of $Dis(f_p)$ it follows $\overline{\rho}(1)=P_0$ - that is, we have just proven that if $\rho$ is a component of $Dis(f_p)$ with an endpoint strictly interior to $D_\alpha$, that endpoint can only be $P_0$.\\

To conclude the proof, recall that since both $\rho$ and $\delta$ are components of $Dis(f_p)$ we have $f_{p}(\rho),f_{p}(\delta)\subseteq l$ - recalling $l$ is an open arc on $\partial D_\alpha$ whose endpoints are $P_{In},P_{Out}$, since $P_0$ lies in $\overline{\rho}\cap\overline{\delta}$ and because $P_0$ tends to $P_{In}$ in infinite time, we conclude $P_{In}\in\overline{f_{p}(\delta)}\cap\overline{f_{p}(\rho)}$. As both $f_{p}(\delta),f_{p}(\rho)$ are arcs on  $l$ with an endpoint at $P_{In}$, it follows $f_{p}(\rho)\cap f_{p}(\delta)\ne\emptyset$ - therefore, by the Existence and Uniqueness Theorem $\rho=\delta$. Consequentially, if $\rho$ is a component of $Dis(f_p)$ s.t. $\rho\ne\delta$, then $\rho$ has two endpoints in $l$. A similar argument proves that whenever $\rho\ne\delta$ we have $P_0\not\in\rho$ and Prop.\ref{lem33} follows (see Fig.\ref{fig15} for an illustration).
\end{proof}
\begin{remark}
    \label{extend}
    Given $s\in l$, the same arguments used above also imply that whenever $f_p(s)$ is interior to $D_\alpha$, $f_p$ is continuous on a neighborhood of $s$ in $\overline{D_\alpha}$.
\end{remark}
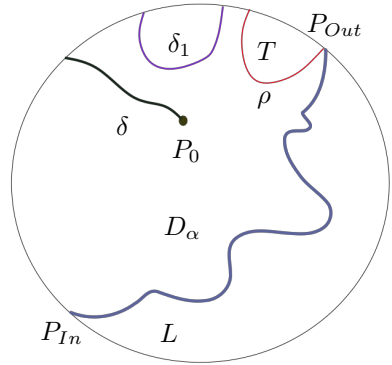
\begin{figure}[h]
\centering
\begin{overpic}[width=0.30\textwidth]{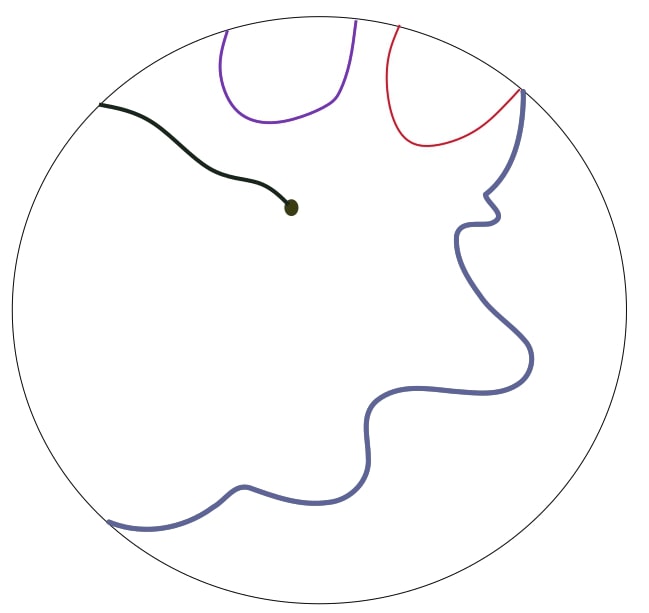}
\put(90,70){$P_{In}$}
\put(390,65){$L$}
\put(280,590){$\delta$}
\put(410,800){$\delta_1$}
\put(400,335){$D_\alpha$}
\put(420,530){$P_0$}
\put(760,850){$P_{Out}$}
\put(635,670){$\rho$}
\put(635,780){$T$}
\end{overpic}
\caption[More discontinuities for the first-return map.]{\textit{The discontinuity curves for $f_{p}$ in $D_\alpha$ - $\delta,\rho$, and $\delta_1$, another possible discontinuity curve. In this sketch, $U_p$ is drawn for simplicity as a disc. $T$ is the sector trapped between $l$ and $\rho$.}}
\label{fig15}
\end{figure}

The main takeaway from Prop.\ref{lem33} is that the point $P_0$ generates a discontinuity curve $\delta$ which, due to the existence of the heteroclinic trefoil, is associated with the fixed point $P_{In}$. However, this does not imply the first-return map is continuous at the second fixed point, $P_{Out}$. Using the same notations and ideas used to prove Prop.\ref{lem33} we now prove:
\begin{corollary}\label{disconcurve}
 Let $F_p$ be a trefoil parameter. Then, there exists a unique curve $\rho\subseteq Dis(f_{p})$ s.t. $P_{Out}\in\overline\rho$. Moreover, there exists a bounded sector $T$, trapped between $\rho$ and $l$ with a tip at $P_{Out}$, s.t. $P_0\in\partial f_p(T)$ (see Fig.\ref{fig15} and Fig.\ref{fig1555} for illustrations).
\end{corollary}
\begin{proof}
Let us recall the two-dimensional stable manifold for $P_{Out}$, $W^s_{Out}$, is transverse to the cross-section $U_{p}$ at $P_{Out}$ (see Lemma \ref{cor211} and the illustration in Fig.\ref{loci}). This implies the backwards trajectories of initial conditions on the arc $l$ sufficiently close to $P_{Out}$ all hit $U_{p}$ transversely - or, in other words, there exists a curve $\rho\subseteq Dis(f_{p})$, $P_{Out}\in\overline\rho$, s.t. $f_{p}(\rho)$ is an arc on $l$ with an endpoint at $P_{Out}$ (see the illustration in Fig.\ref{fig1555}) - that is, the flow pushes $\rho$ over the arc $f_p(\rho)$ - it is easy to see that by the Existence and Uniqueness Theorem, $\rho$ is unique. Now, let $T$ denote the open sector trapped between $\rho$ and $l$ (see the illustration in Fig.\ref{fig155} and Fig.\ref{fig1555}), and let us consider a sequence $\{x_n\}_n\subseteq T$ s.t. $x_n\to P_{Out}$.\\

\begin{figure}[h]
\centering
\begin{overpic}[width=0.30\textwidth]{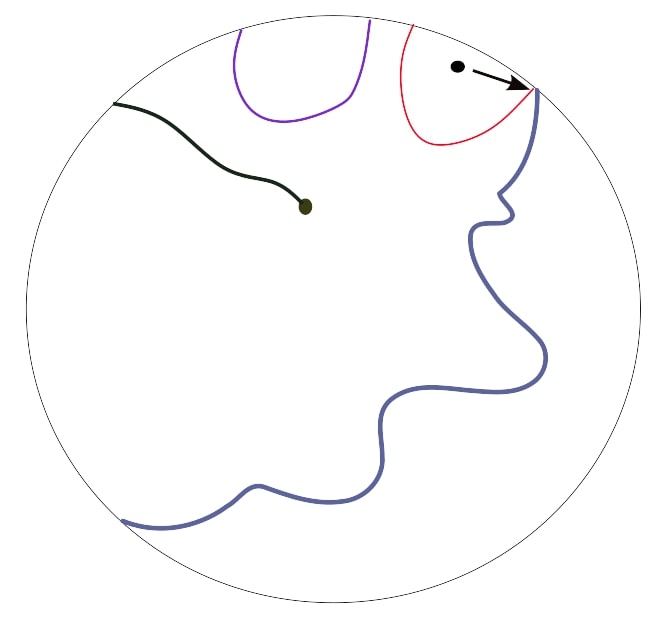}
\put(100,70){$P_{In}$}
\put(390,65){$L$}
\put(280,590){$\delta$}
\put(410,800){$\delta_1$}
\put(400,335){$D_\alpha$}
\put(420,530){$P_0$}
\put(760,850){$P_{Out}$}
\put(635,670){$\rho$}
\put(635,780){$x_n$}
\end{overpic}
\caption[The sequence $\{x_n\}_n$.]{\textit{$x_n$ tends to $P_{Out}$ from inside the sector $T$.}}
\label{fig155}
\end{figure}
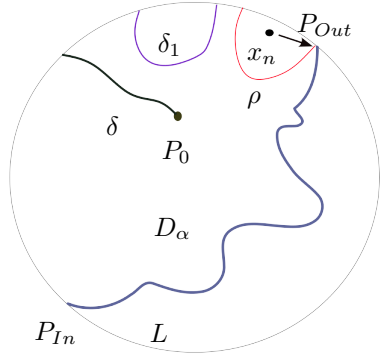

Now, recall the heteroclinic trajectory $\Theta$ given by Def.\ref{def32} flows from $P_{Out}$ to $P_{In}$, while hitting $D_\alpha$ transversely at $P_0$ along the way. Additionally, note that since the first-return map pushes $\rho$ down on $f_p(\rho)$, it follows $f_p(T)\cap T=\emptyset$, hence $f_p(x_n)\not\to P_{Out}$. Therefore, because $\{x_n\}_n$ tends to $P_{Out}$ it follows that for $n$ sufficiently large $x_n$ leaves $D_\alpha$ and travels along the bounded heteroclinic trajectory $\Theta$ until hitting $D_\alpha$ (see the illustration in Fig.\ref{fig1555}). Hence, since $x_n\to P_{Out}$ and because $\Theta$ hits $D_\alpha$ transversely at $P_0$ we have $f_p(x_n)\to P_{0}$ and Cor.\ref{disconcurve} now follows (see the illustration in Fig.\ref{fig155} and Fig.\ref{fig1555}). 
\end{proof}
\begin{figure}[h]
\centering
\begin{overpic}[width=0.7\textwidth]{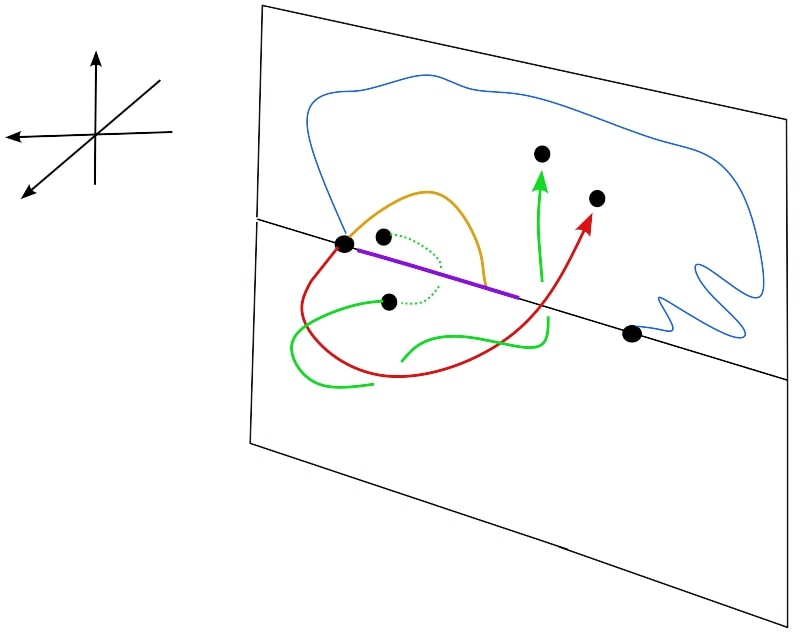}
\put(720,380){$P_{In}$}
\put(710,590){$P_0$}
\put(340,480){$P_{Out}$}
\put(550,610){$f_p(x_n)$}
\put(560,400){$f_p(\rho)$}
\put(490,570){$\rho$}
\put(500,520){$x_n$}
\put(370,730){$U_p$}
\put(460,650){$D_\alpha$}
\put(560,710){$L$}
\put(700,160){$L_p$}
\put(-15,620){$x$}
\put(0,535){$y$}
\put(110,740){$z$}
\end{overpic}
\caption[The trajectories of $\{x_n\}_n$.]{\textit{The flow line connecting $x_n$ and $f_p(x_n)$, spiralling along the heteroclinic trajectory (the red curve). The purple arc denotes $f_p(\rho)$ - it is easy to see that since the flow pushes $\rho$ on $l$, $\rho$ is the unique curve of $Dis(f_p)$ s.t. $P_{Out}\in\overline{\rho}$.}}
\label{fig1555}
\end{figure}
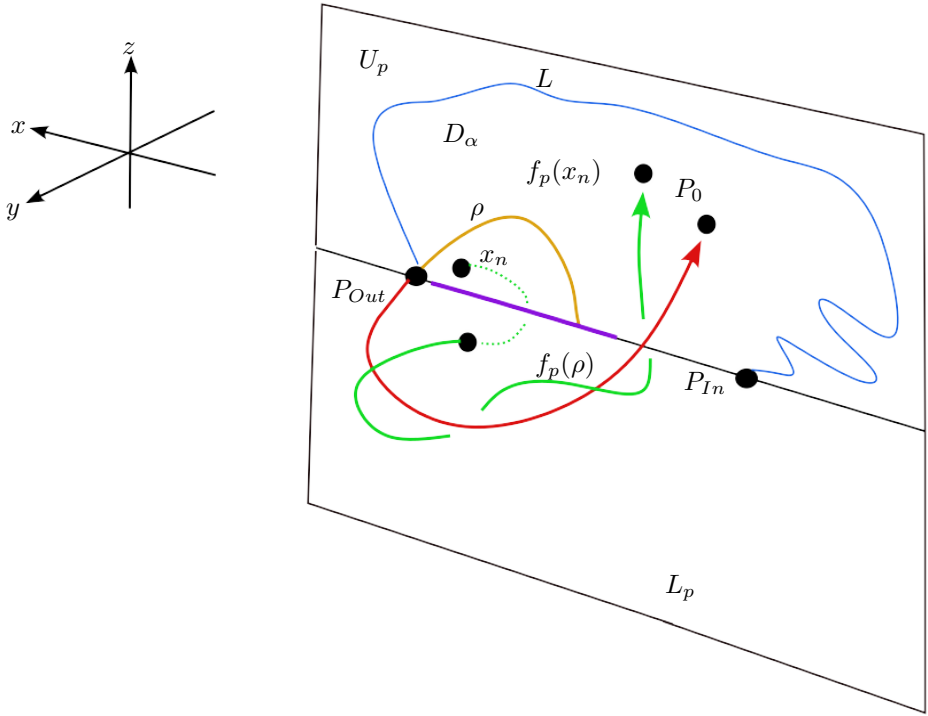

To continue, having studied the discontinuity properties of $f_p$ in Prop.\ref{lem33} and Cor.\ref{disconcurve}, we now study the continuity properties of the first-return map $f_p$. To do so, consider the set $I=(\overline{D_\alpha}\setminus\{P_0\})\setminus \cup_{n>0}f^{-n}_p(l)$ - that is, $I$ is the (maximal) invariant set of $f_p$ in $(\overline{D_\alpha}\setminus\{P_0\})\setminus l$. We now prove the following corollary of Prop.\ref{lem33}:

\begin{corollary}
    \label{nosep}
    Let $C$ be a component of $I$. Then, for every $n>0$, $f^n_p$ is continuous on $C$ - and moreover, $I$ is dense in $\overline{D_\alpha}$. 
\end{corollary}
\begin{proof}
By Prop.\ref{lem33} it immediately follows that given $C$ as above and any $n>0$, $f^n_p$ is continuous on $C$. Therefore, it remains to prove $I$ is dense in $\overline{D_\alpha}$. To do so, recall $f_p$ is a first-return map generated by a smooth flow - hence given any $n\in \mathbf{N}$, $f^{-n}_p(l)$ is a collection of smooth curves (some of which are possibly singletons - see the illustration in Fig.\ref{BAIRE}). As such, for every $n\geq 0$, $I_n=\overline{D_\alpha}\setminus (\cup_{0\leq k\leq n}\overline{f^{-n}_p(l)})$ is dense in $\overline{D_\alpha}$ and open in $\overline{D_\alpha}$. $\overline{D_\alpha}$ is a complete and separable metric space, hence, by the Baire Category Theorem it follows $I=\cap_{n\geq1} I_n$ is dense in $\overline{D_\alpha}$ and Cor.\ref{nosep} follows. 
\end{proof}
 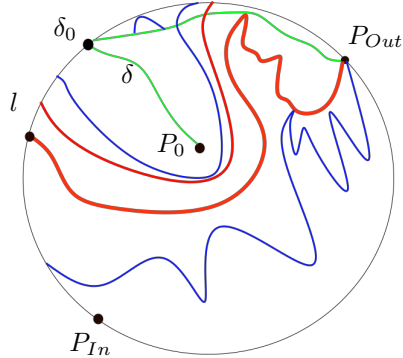
\begin{figure}[h]
\centering
\begin{overpic}[width=0.3\textwidth]{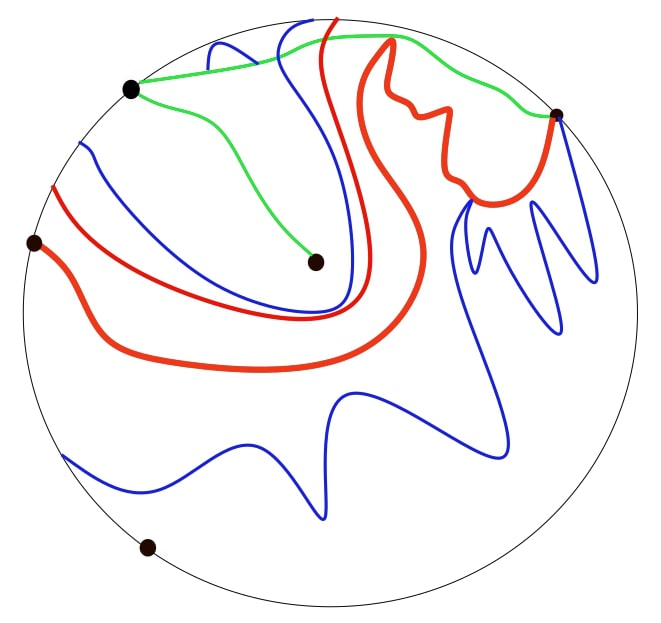}
\put(150,30){$P_{In}$}

\put(280,710){$\delta$}
\put(110,830){$\delta_0$}
\put(0,640){$l$}
\put(370,535){$P_0$}
\put(850,810){$P_{Out}$}

\end{overpic}
\caption[The set $I_3$.]{\textit{The set $I_3$ inside the disc $D_\alpha$- the green curves correspond to $f^{-1}_p(l)$, the red to $f^{-2}_p(l)$, and the blue to $f^{-3}_p(l)$.}}
\label{BAIRE}
\end{figure}
Having proven Prop.\ref{lem33}, Cor.\ref{disconcurve} and Cor.\ref{nosep}, we now study the question of how the curves in $f^{-1}_p(l)$ (and by extension, $Dis(f_p)$) are configured inside the cross-section $\overline{D_\alpha}$. The reason we are interested in this question is because the configuration of $f^{-1}_p(l)$ dictates how the first-return map $f_p$ can (and cannot) behave. As we will now soon prove, there are precisely three ways in which $f^{-1}_p(l)$ can be configured inside the cross-section $D_\alpha$. To do so, we first prove the following technical Lemma:

\begin{lemma}
\label{lemhp}    Let $p$ be a trefoil parameter, let $\delta$ be the curve given by Prop.\ref{lem33}, and let $\rho$ be as in Cor.\ref{disconcurve}. Then there exists a cross section $H_p$, a component of $(\overline{D_\alpha}\setminus\{P_0\})\setminus f^{-1}_p(l)$, satisfying:
\begin{itemize}
    \item $H_p$ is a topological disc, and $f_p$ is continuous on $H_p$.
    \item $P_{In}$, $L$, $\rho$ and $P_{Out}$ all lie in $\partial H_p$ (see the illustration in Fig.\ref{DISC1}).
    \item Either $\delta\subseteq \partial H_p$ (in which case $H_p$ is homeomorphic to a slit disc), or $\overline{\delta}\cap\partial H_p\subseteq\{\delta_0\}$ (see the illustrations in Fig.\ref{OPT}).
\end{itemize}
\end{lemma}
\begin{proof}
By Prop.\ref{lem33}, every component of $f^{-1}_p(l)\cap D_\alpha$ is either the curve $\delta$, or alternatively, a curve with at least two endpoints on $l$. Consequentially, as ${D_\alpha}\setminus\{P_0\}$ is homeomorphic to an open punctured disc and because $P_0\in\overline{\delta}$, every component of $({D_\alpha}\setminus\{P_0\})\setminus f^{-1}_p(l)$ is an open topological disc. Recalling that by Prop.\ref{lem33} the discontinuities of $f_p$ in $D_\alpha\setminus\{P_0\}$ are given by $f^{-1}_p(l)\cap D_\alpha$, it follows the the components of $({D_\alpha}\setminus\{P_0\})\setminus f^{-1}_p(l)$ are the continuity sets of $f_p$ inside ${D_\alpha}\setminus\{P_0\}$.\\

Since by Prop.\ref{lem33} $f_p$ is continuous on a neighborhood of $P_{In}$ in $D_\alpha$, there exists a neighborhood of $P_{In}$ in $\overline{D_\alpha}$ which lies away from $f^{-1}_p(l)$ (see the illustration in Fig.\ref{connn}). Consequentially, there exists a unique component of $(D_\alpha\setminus\{P_0\})\setminus f^{-1}_p(l)$, which we denote by $H_p$, s.t. $P_{In}\in\partial H_p$. By the discussion above, $f_p$ is continuous on $H_p$, and $H_p$ is a topological disc (see the illustration in Fig.\ref{DISC1}).\\
    
We now prove that in addition to the fixed-point $P_{In}$, the sets $L$, $P_{Out}$ and $\rho$ also lie in $\partial H_p$. To do so, recall the curve $L$ is the intersection of the cross-section $U_p$ with the two-dimensional, unstable, invariant manifold $W^u_{In}$ (see the illustration in Fig.\ref{DISC1}) - as proven in Prop.\ref{arccor}, $L$ is curve in $\partial D_\alpha\cap U_p$ with endpoints $P_{In}$, $P_{Out}$. Since by Prop.\ref{lem33} $f_p$ is continuous at both $L$ and a neighborhood of $P_{In}$ in $\overline{D_\alpha}$, we conclude $\overline{L}\subseteq\partial H_p$. Therefore, by $P_{Out}\in \overline{L}$ it follows $P_{Out}$ is also in $\partial H_p$ (see the illustration in Fig.\ref{optA} or \ref{optB}). Finally, since by Cor.\ref{disconcurve} $\rho\subseteq f^{-1}_p(l)$ is the unique curve in $f^{-1}_p(l)\cap D_\alpha$ with an endpoint at $P_{Out}$, it follows $\rho\subseteq\partial H_p$ as well (see Fig.\ref{DISC1}). That is, we have just proven $L,P_{Out}$ and $\rho$ are all subsets of $\partial H_p$.\\ 

\begin{figure}[h]
\centering
\begin{overpic}[width=0.7\textwidth]{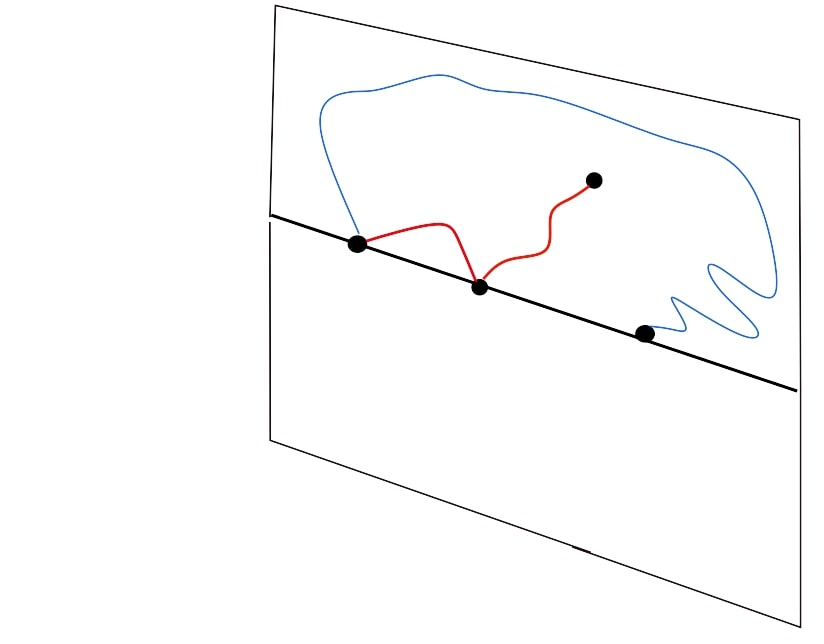}
\put(750,320){$P_{In}$}
\put(380,450){$P_{Out}$}
\put(640,520){$\delta$}
\put(580,370){$\delta_0$}
\put(800,630){$L$}
\put(350,730){$U_p$}
\put(460,620){$H_p$}
\put(480,520){$\rho$}
\put(700,580){$P_0$}
\put(350,265){$L_p$}
\end{overpic}
\caption[The cross-section $H_p$.]{\textit{The cross-section $H_p\subseteq D_\alpha$ (where $D_\alpha$ is the Jordan domain trapped between $L$ and $l$), trapped between $f^{-1}_p(l)$, $l$, and $L$ - where $l$ is the interval connecting $P_{In}$ and $P_{Out}$, and $f^{-1}_p(l)$ is the red curve (in this illustration, $\rho$ is the red arc connecting $P_{Out}$ and $\delta_0$). It follows $H_p$ is an open topological disc.}}
\label{DISC1}
\end{figure}
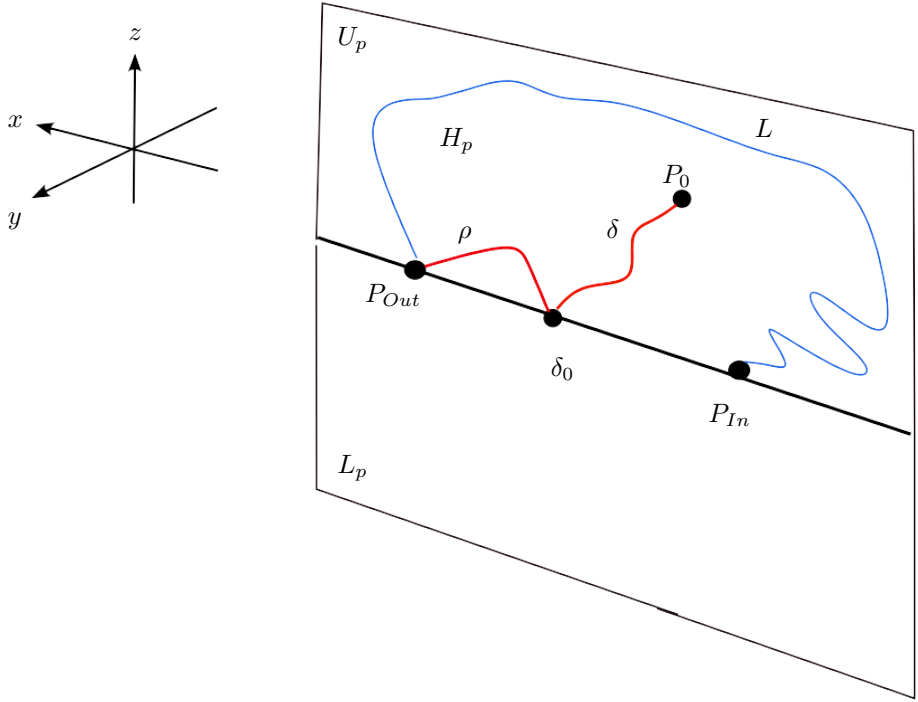

To prove Lemma \ref{lemhp} it remains to prove either $\delta\subseteq \partial H_p$ or $\overline{\delta}\cap \partial H_p\subseteq\{\delta_0\}$ - where $\delta\subseteq f^{-1}_p(l)\cap D_\alpha$ is the curve given by Prop.\ref{lem33}, and $\delta_0$ is the endpoint of $\delta$ on $l$ (see the illustration in Fig.\ref{DISC1}). As $H_p$ is a component of $D_\alpha\setminus f^{-1}_p(l)$, by $\delta\subseteq f^{-1}_p(l)$ we already know $\delta\cap H_p=\emptyset$ - which proves $\delta$ intersects $\overline{H_p}$ at most in $\partial H_p$. Now, let us recall $\partial H_p$ is composed of arcs in $\partial D_\alpha$ (where the boundary is taken in $\overline{U_p}$) and in $f^{-1}_p(l)\cap D_\alpha$ - therefore, since every $s\in\delta$ is interior to $D_\alpha$, and because $\delta$ is also a component of $f^{-1}_p(l)\cap D_\alpha$ homeomorphic to an open interval, it follows that if $\delta\cap\partial H_p\ne\emptyset$ the only possibility is $\delta\subseteq\partial H_p$. Consequentially, whenever $\delta\not\subseteq\partial H_p$, the only possibility is $\partial H_p\cap\overline{\delta}\subseteq\{\delta_0\}$ (see the illustration in Fig.\ref{OPT} or \ref{optA}). The proof of Lemma \ref{lemhp} is now complete.
\end{proof}

We are now ready to prove that given a trefoil parameter $F_p$, there are precisely three options as to how the collection of curves $f^{-1}_p(l)$ can be arranged inside $D_\alpha$ - and we do so by considering the possible topology of $\partial H_p$ (as dictated by Lemma \ref{lemhp}). First, whatever the case, by Lemma \ref{lemhp} the fixed points $P_{In}$ and $P_{Out}$ and the arcs $L$ and $\rho$ are inside $\partial H_p$. Additionally, by Prop.\ref{lem33} either $\delta_0=P_{Out}$ or $\delta_0$ is interior to $l$. Therefore, given any trefoil parameter $F_p$ there are precisely three possibilities as to how $\delta,\delta_0$ and the cross-section $H_p$ can be arranged inside $\overline{D_\alpha}$ (see Fig.\ref{OPT}):

\begin{itemize}
    \item The first possibility is $P_{Out}=\delta_0$. Since $\delta_0$ is the unique intersection point of $\overline{\delta}$ with $l$ (see Prop.\ref{lem33}), whenever $\delta_0=P_{Out}$ the curve $\delta$ is an arc in the open disc $D_\alpha$ - with one endpoint at $P_0$ and another at $P_{Out}$ (see the illustration in Fig.\ref{OPT} and \ref{optc}). By $f_p(P_{Out})=P_{Out}$ and because the trajectory of $P_0$ tends to $P_{In}$ we have $f_p(\delta)=l$ - i.e. $f^{-1}_p(l)=\delta=\rho$. This implies $H_p=D_\alpha\setminus\delta$, hence $H_p$ is homeomorphic to a slit disc (see the illustration in Fig.\ref{OPT} and Fig.\ref{optc}). We refer to this scenario as \textbf{Case $C$}.
    \item   The second case we consider is the scenario where $\delta\subseteq\partial H_p$ and $\delta_0\ne P_{Out}$ (see the illustration in Fig.\ref{OPT} and Fig.\ref{optB}). As $\delta$ is a component of $f^{-1}_p(l)\cap D_\alpha$, it again follows $H_p$ is homeomorphic to a slit disc. We refer to this scenario as \textbf{Case} $B$.
    \item The third (and final) possibility is that $\overline{\delta}\cap\partial H_p\subseteq\{\delta_0\}$ (as illustrated in both Fig.\ref{OPT} and Fig.\ref{optA}). In light of the discussion above, whenever $\overline{\delta}\cap\partial H_p\subseteq\{\delta_0\}$ we have $\delta_0\ne P_{Out}$ - we refer to this scenario as \textbf{Case} $A$. 
\end{itemize}

It is easy to see that by Lemma \ref{lemhp}, there are no other possible configurations of $\partial H_p$. Summarizing our results, we obtain:

\begin{corollary}
    \label{uniquness} Let $F_p$ be a trefoil parameter - then, the configuration of $\partial H_p$ in $\overline{D_\alpha}$ falls into either Case $A$, Case $B$ or Case $C$.
\end{corollary}

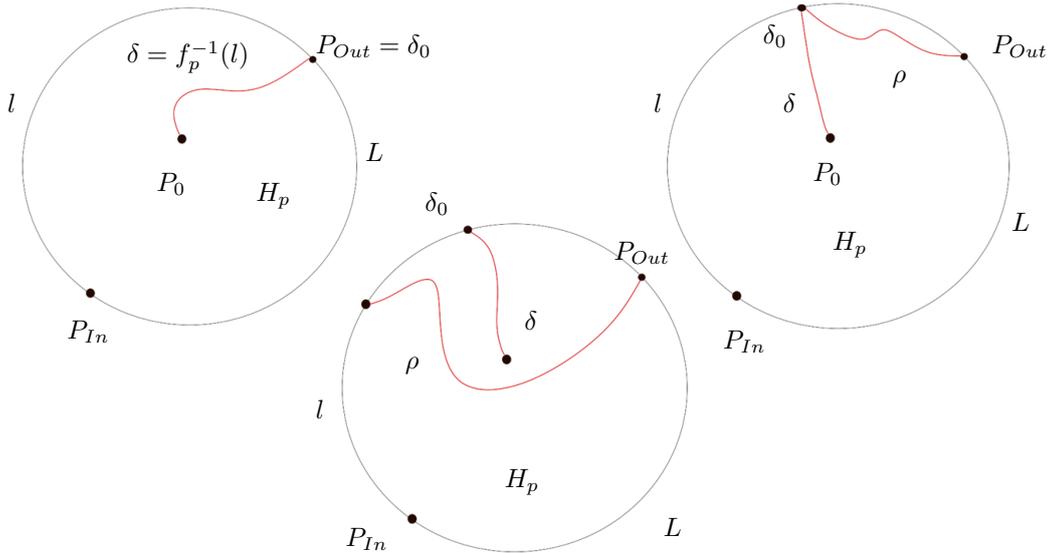
\begin{figure}[h]
\centering
\begin{overpic}[width=0.75\textwidth]{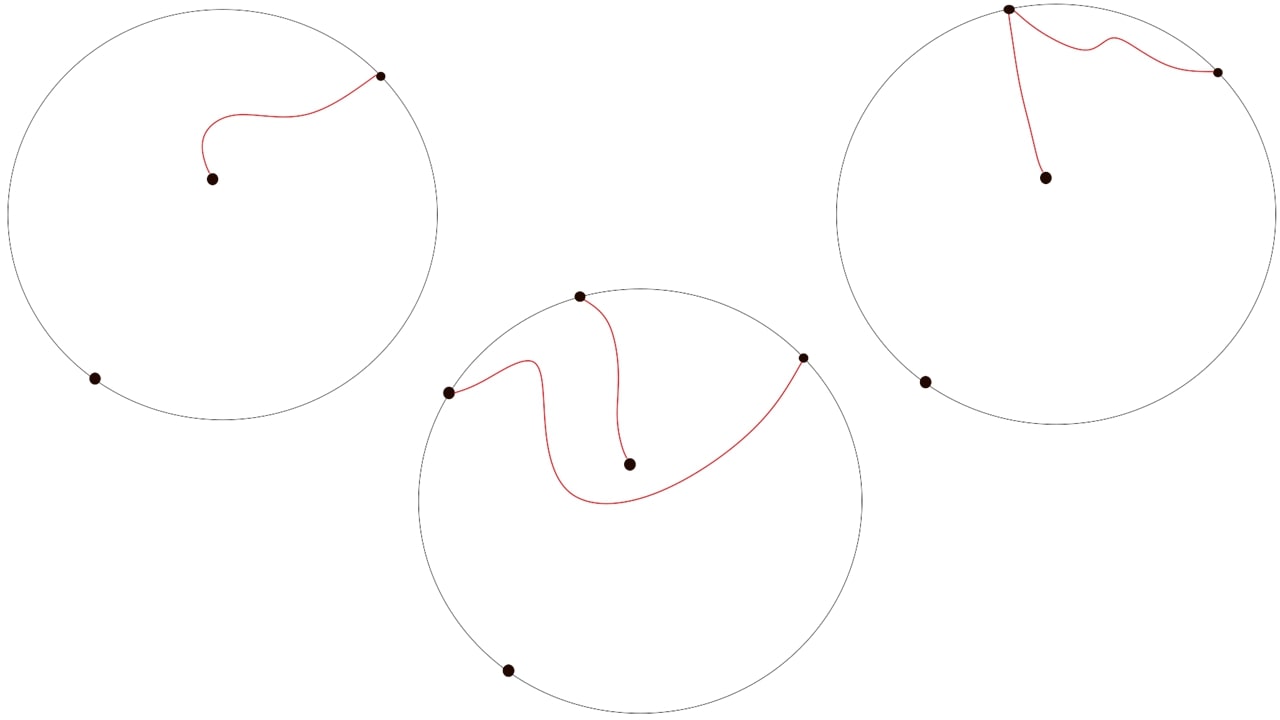}%
\put(980,510){$P_{Out}$}
\put(1000,330){$L$}
\put(410,350){$\delta_0$}
\put(140,370){$P_0$}
\put(770,450){$\delta$}
\put(750,520){$\delta_0$}
\put(330,10){$P_{In}$}
\put(490,70){$H_p$}
\put(50,220){$P_{In}$}
\put(600,300){$P_{Out}$}
\put(240,360){$H_p$}
\put(110,500){$\delta=f^{-1}_p(l)$}
\put(350,400){$L$}
\put(300,510){$P_{Out}=\delta_0$}
\put(800,380){$P_0$}
\put(650,20){$L$}
\put(390,190){$\rho$}
\put(510,230){$\delta$}
\put(710,210){$P_{In}$}
\put(820,310){$H_p$}
\put(880,480){$\rho$}
\put(-10,450){$l$}
\put(640,450){$l$}
\put(300,140){$l$}
\end{overpic}
\caption[The three possibilities for $H_p$.]{\textit{The three possibilities for the topology of $\partial H_p$ inside $D_\alpha$ (for convenience, we sketch $D_\alpha$ as a topological disc). In Case $C$ (upper left), $\delta=f^{-1}_p(l)=\rho$, and $P_{Out}=\delta_0$. In Case $B$ (upper right) we have $\delta_0\ne P_{Out}$ and $\delta\subseteq \partial H_p$, while in Case $A$ (lower-center) we have $\overline{\delta}\cap\partial H_p\subseteq\{\delta_0\}$. }}\label{OPT}
\end{figure}

At this point we remark that despite its technical nature, Cor.\ref{uniquness} will be of extreme importance to the proof of chaoticity for the Rössler system at trefoil parameters in the next section - primarily because each scenario forces the first-return map $f_p:H_p\to\overline{D_\alpha}$ to behave somewhat differently.\\

Moreover, before we conclude this subsection and move to prove Th.\ref{th31}, we further remark all the results proven so far study only the bounded dynamics of the flow - therefore, before concluding this section, for the sake of completeness we study the unbounded dynamics at trefoil parameters. To do so, recall that given a parameter $F_p$, we always denote by $F_p$ the corresponding vector field. Inspired by Th.\ref{th21} from and Shilnikov's Homoclinic Bifurcation Theorem, we now prove:
\begin{proposition}\label{cor216}
Let $F_p$ be a trefoil parameter. Then, there exists a sequence of smooth vector fields $\{K_n\}_n$ on $S^3$ s.t.:
\begin{itemize}
    \item Given any compact $K\subseteq\mathbf{R}^3$, $K_n\to F_p$ on $K$ (w.r.t. the $C^\infty$ metric).
    \item Every $K_n$ generates an unbounded set of periodic trajectories. 
\end{itemize}

\end{proposition}
\begin{proof}
Let $F_p$ be a trefoil parameter. By Th.\ref{th21}, we can always perturb $F_p$ into a smooth vector field on $S^3$, $R_p$ s.t. the following is satisfied:
\begin{itemize}
    \item $R_p$ coincides with $F_p$ around the fixed points.
    \item Given any compact $K\subseteq\mathbf{R}^3$, we can choose $R_p$ s.t. it coincides with $F_p$ on $K$.
    \item $R_p$ generates an unbounded heteroclinic trajectory $\Gamma$ which flows from $P_{Out}$ towards $P_{In}$ (in infinite time). 
\end{itemize}

Since $p=(a,b,c)$ lies in the parameter space $P$, per our assumption on the parameter space, at least one of the fixed points $P_{In}$ and $P_{Out}$ always satisfies the Shilnikov Condition (see the discussion in page \pageref{eq:9}). Assume first that $P_{Out}$ satisfies the Shilnikov condition - since $R_p$ is a smooth vector field on $S^3$, we can smoothly deform $R_p$ around the fixed point $P_{In}$ by connecting $\Gamma$ and some trajectory on $S=W^{u}_{In}=W^{s}_{Out}$ to generate $\zeta$, a homoclinic trajectory to $P_{Out}$ - in particular, we choose this perturbation s.t. $R_p$ and $K_n$ coincide around $P_{Out}$. Since $\Gamma$ is unbounded, so is $\zeta$.\\

Denote by $K_n$ the vector field described above - it is easy to see $K_n$ is a smooth vector field on $S^3$, and that we can choose it s.t. its $C^\infty$ distance from $R_p$ is arbitrarily small. Therefore, given any compact $K\subseteq\mathbf{R}^3$, since $R_p$ can be chosen to coincide with $F_p$ on $K$, $K_n$ can be chosen to be arbitrarily close to $F_p$ on $K$ (in the $C^\infty$ metric).\\

Since $P_{Out}$ satisfies the Shilnikov Condition w.r.t. $F_p$ and $R_p$, by our construction of $K_n$ it also satisfies it w.r.t. $K_n$. Therefore, as $P_{Out}$ satisfies the Shilnikov Condition and because $K_n$ generates a homoclinic trajectory $\zeta$ to $P_{Out}$, from Shilnikov's Theorem $K_n$ generates $\Omega$ - a collection of periodic trajectories which are dense around $\zeta$  (see \cite{LeS} or Th.13.8 in \cite{SSTC}). As such, since $\zeta$ is unbounded so is $\Omega$. When $P_{In}$ satisfies the Shilnikov Condition, similar arguments (when applied to the inverse flow) imply the same result and Prop.\ref{cor216} now follows. 
\end{proof}

\subsection{Chaotic dynamics at trefoil parameters}
\label{sect33}
In this subsection we prove the Rössler system at trefoil parameter generates complex dynamics (which we do in Th.\ref{th31}). In more detail, we will prove that given a trefoil parameter, its dynamics include infinitely many periodic trajectories - which manifest as periodic orbits for the first-return map $f_p$, of all minimal periods.\\

\begin{figure}[h]
\centering
\begin{overpic}[width=0.6\textwidth]{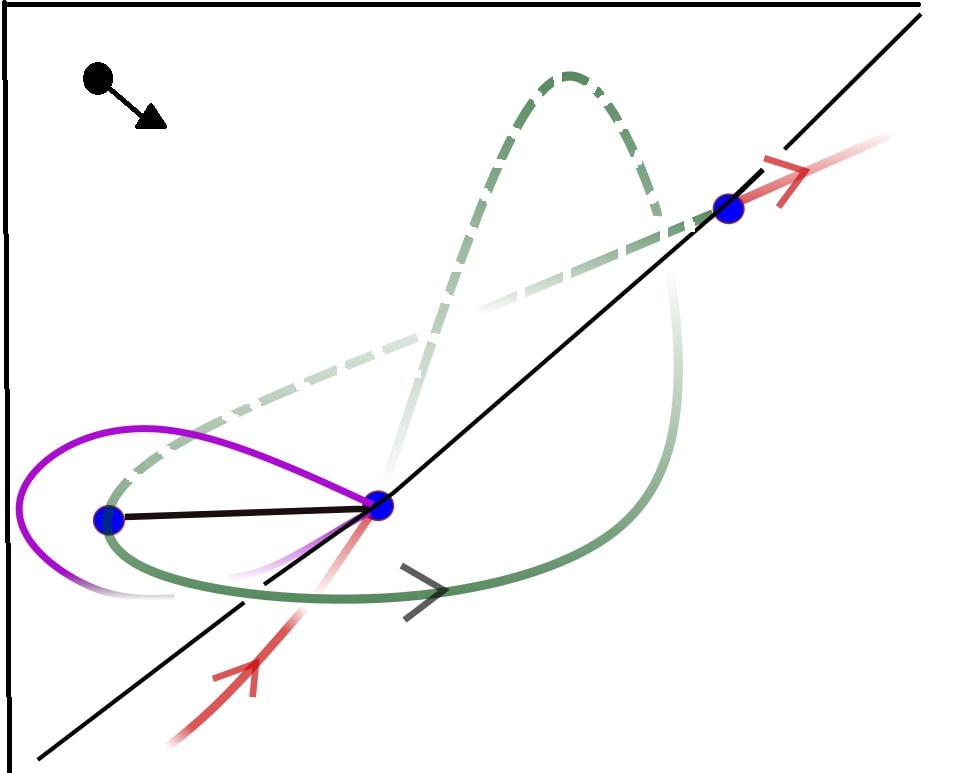}
\put(420,250){$P_{In}$}
\put(250,290){$\gamma$}
\put(80,300){$P_0$}
\put(100,390){$f_p(\gamma)$}
\put(800,550){$P_{Out}$}
\put(560,470){$\Theta$}
\end{overpic}
\caption[Suspending a curve along the trefoil knot.]{\textit{Flowing $\gamma$ along the trefoil. $\Theta$ is the green trajectory.}}
\label{fig16}
\end{figure}
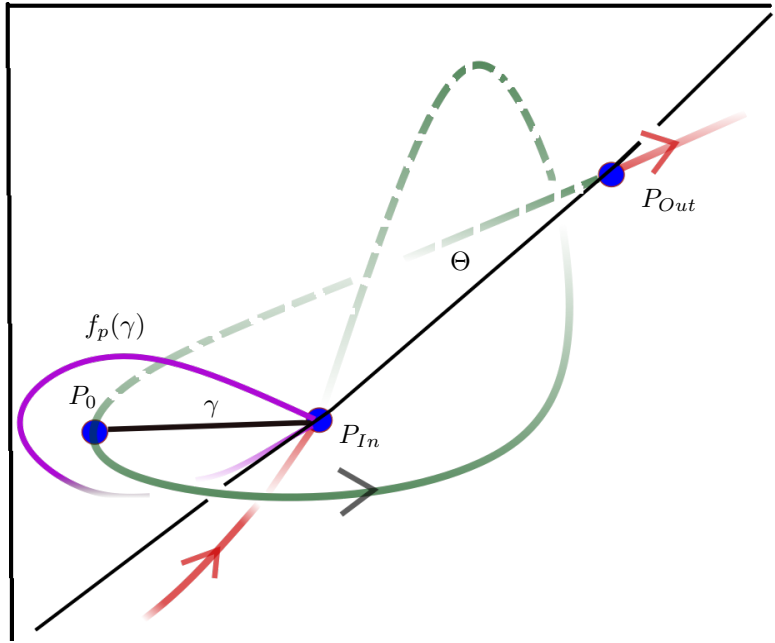

The motivation behind the proof is relatively simple. To introduce it, consider a scenario as described in Fig.\ref{fig16}. That is, given a trefoil parameter $F_p$ and the corresponding bounded heteroclinic trajectory $\Theta$, choose some curve $\gamma\subseteq U_p$ s.t. $\gamma$ connects the fixed point $P_{In}$ and $P_0$ on the cross-section $U_p$.\\

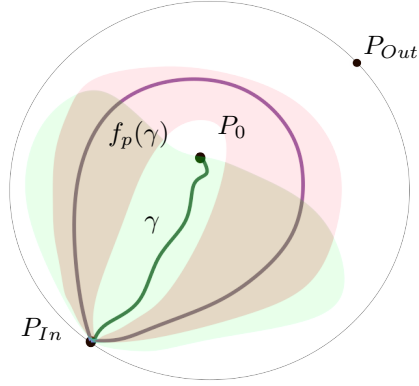
\begin{figure}[h]
\centering
\begin{overpic}[width=0.55\textwidth]{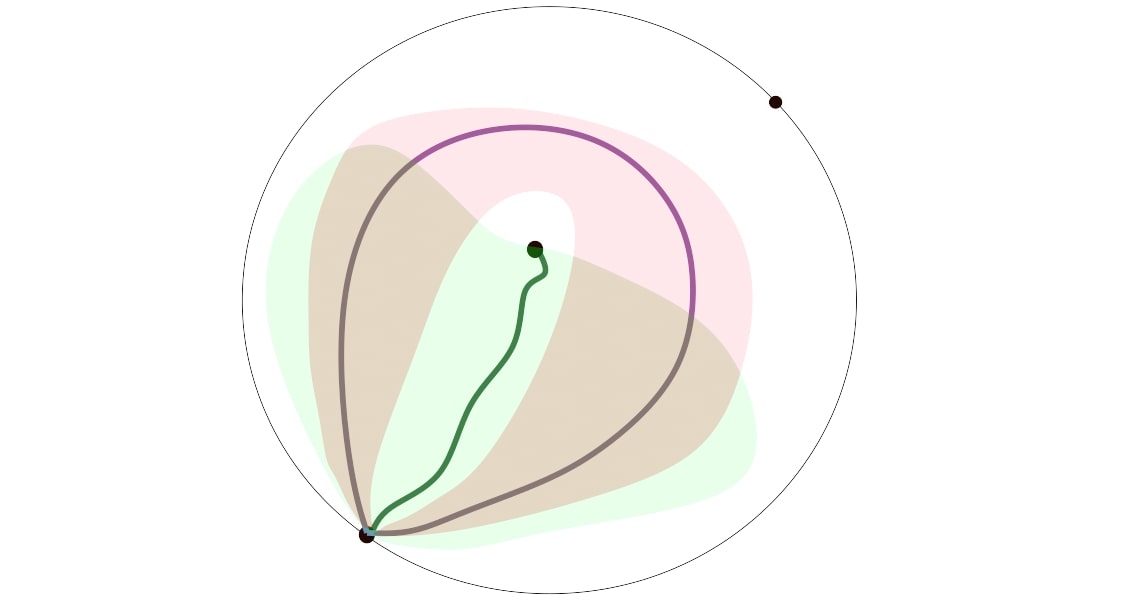}
\put(230,70){$P_{In}$}
\put(500,350){$P_{0}$}
\put(400,220){$\gamma$}
\put(350,340){$f_p(\gamma)$}
\put(700,460){$P_{Out}$}
\end{overpic}
\caption[The first-return map as a horseshoe.]{\textit{Imposing $f_p(\gamma)$ on $U_p$ - since $f_p(\gamma)$ is a closed loop, heuristically, there exists a degenerate rectangle $R$ (the green region) with sides $AB,CD$ corresponding to $P_{In},P_0$ (respectively) s.t. $f_p(R)$ (the red region) looks like a singular Smale Horseshoe.}}
\label{fig17}
\end{figure}

Now, let us suspend the curve $\gamma$ along $\Theta$, the bounded component of the heteroclinic trefoil knot (see the illustration in Fig.\ref{fig16}). By the topology of the heteroclinic trefoil knot, we conclude it constrains $\gamma$ to return to $U_p$ as in Fig.\ref{fig16} - that is, heuristically, $f_p(\gamma)$, the image under the first-return map is a closed loop on the cross-section $U_p$ which begins and terminates at the fixed point $P_{In}$. Therefore, we can impose $f_p(\gamma)$ on $U_p$ roughly as in Fig.\ref{fig17} - which, heuristically, implies the existence of a rectangle $R\subseteq U_p$ s.t. $f_p$ is Smale Horseshoe map on $R$. As such, it should follow the dynamics of $F_p$ must include infinitely many periodic trajectories in $R$.\\

In practice, due to Prop.\ref{lem33} and Cor.\ref{disconcurve} the heuristic described above is far from obvious. Namely, the assumption we can choose a curve $\gamma\subseteq U_p$ s.t. $f_{p}(\gamma)$ is a closed loop is far from trivial: recalling both Prop.\ref{lem33} and Cor.\ref{disconcurve}, it is easy to see one cannot assume $f_p(\gamma)$ is even connected in $\overline{D_\alpha}$. However, as we will see below, provided we are sufficiently careful, this heuristic becomes a rigorous proof.
\begin{theorem}
    \label{th31}
    Let $(a,b,c)=p\in P$ be a trefoil parameter, and denote by $f_p:H_p\to \overline{D_\alpha}$ the first-return map as defined in the previous section. Then, there exists an invariant set $Q\subseteq H_p$ s.t. the first-return map $f_p:Q\to Q$ is continuous and chaotic per Def.\ref{chaotic}. 
\end{theorem}

\begin{proof}
To begin, recall the bounded cross-section $D_\alpha\subseteq U_p$, $D_\alpha=B_\alpha\cap U_p$ (where $B_\alpha$ is the topological ball generated by the equality $W^u_{In}=W^s_{Out}$ - see Def.\ref{def32}). Additionally, recall the first-return map $f_p:\overline{D_\alpha}\setminus\{P_0\}\to\overline{D_\alpha}\setminus\{P_0\}$ - as proven in Prop. \ref{arccor} $D_\alpha$ is a bounded Jordan domain on the half-plane $U_p$ (see the illustration in Fig.\ref{trefint}), and the first-return map $f_p$ is well-defined throughout the punctured (topological) disc $\overline{D_\alpha}\setminus\{P_0\}$ (see Lemma \ref{firstret}). Additionally, recall we denote by $H_p$ is the component of $\overline{D_\alpha}\setminus f^{-1}_p(l)$ s.t. $P_{In}\in\partial H_p$ (see Lemma \ref{lemhp}). Our strategy of proof is as follows - we prove the existence of a set $Q$ inside the cross-section $\overline{H_p}$ s.t. the following is satisfied:

    \begin{itemize}
        \item $f_p$ is continuous on $Q$.
        \item There exists a continuous $\pi:Q\to\{1,2\}^\mathbf{N}$ s.t. $\pi\circ f_p=\sigma\circ\pi$. 
        \item  $\pi(Q)$ includes any periodic $s\in\{1,2\}^\mathbf{N}$ - and if the minimal period of $s$ is $k$, $\pi^{-1}(s)$ includes a periodic point for $f_p$ of minimal period $k$.
    \end{itemize}

 Now, set $Per$ as the set of periodic trajectories for the flow intersecting $Q$. It is easy to see that since $Q\subseteq\overline{H_p}\subseteq \overline{D_\alpha}\subseteq \overline{B_\alpha}$, because $\overline{B_\alpha}$ is bounded $Per$ is bounded as well - and that $\pi$ gives us the continuous and surjective map between $Per$ and the periodic orbits for $\sigma$ in $\{1,2\}^\mathbf{N}$. As such, to prove Th.\ref{th31} it would suffice to prove the assertions above.\\
 
 To do so, recall that per Cor.\ref{uniquness}, given a trefoil parameter $F_p$ it falls into either Case $A$, Case $B$ or Case $C$, s.t. in each case the topology of the cross-section $H_p$ is somewhat different. Therefore, we will prove Th.\ref{th31} for each of these cases separately - despite the subtle differences in the arguments, the proof in each case would be essentially the same. As we will see, in each case we extend the first-return map $f_p:\overline{H_p}\to\overline{D_\alpha}\setminus\{P_0\}$ to a homeomorphism of the disc which we isotope to a Pseudo-Anosov homeomorphism, $G$, which acts as a Smale Horseshoe map on some rectangle. The existence of $G$ would suffice to imply Th.\ref{th31} - because, as we will prove below, all the periodic orbits for $G$ persist when we continuously deform it back to $f_p$. As such, since we will prove $G$ acts as a Smale Horseshoe map on some rectangle, Th.\ref{th31} would follow. 
\subsubsection{Stage $I$ - proving Th.\ref{th31} for Case $A$:}

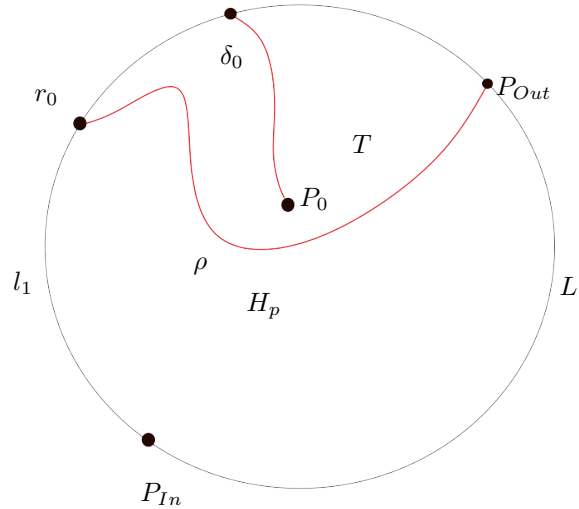
\begin{figure}[h]
\centering
\begin{overpic}[width=0.40\textwidth]{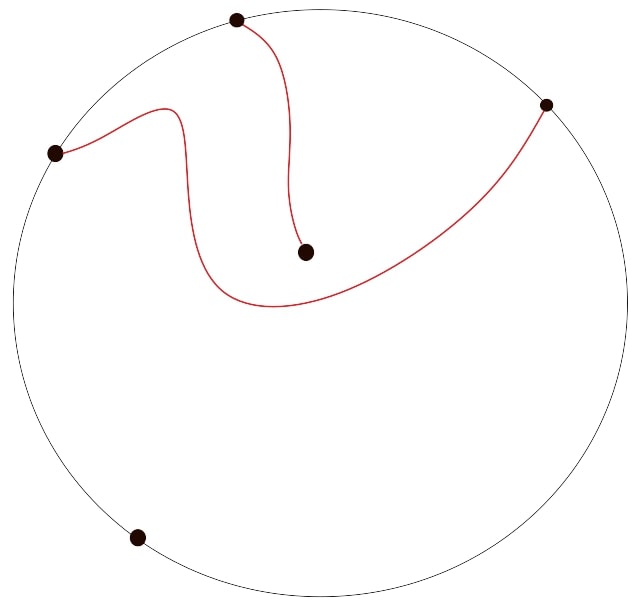}
\put(200,-10){$P_{In}$}
\put(500,550){$P_{0}$}
\put(600,650){$T$}
\put(400,350){$H_p$}
\put(300,430){$\rho$}
\put(0,750){$r_0$}
\put(-40,390){$l_1$}
\put(350,820){$\delta_0$}
\put(870,760){$P_{Out}$}
\put(990,380){$L$}
\end{overpic}
\caption[Case $A$.]{\textit{Case $A$ - $\delta$ lies away from $\overline{H_p}$ - in particular, $\delta_0\ne P_{Out}$}.
}
\label{optA}
\end{figure}

As stated above, we begin by proving Th.\ref{th31} for Case $A$ - that is, in Case $A$, we have $\partial H_p\cap\overline{\delta}\subseteq\{\delta_0\}$ . To begin, let us recall the first-return map $f_p:\overline{D_\alpha}\setminus\{P_0\}\to\overline{D_\alpha}\setminus\{P_0\}$ is continuous around $P_{In}$ (see Prop.\ref{lem33}) - which implies there exists a maximal arc $l_1\subseteq l$, beginning at $P_{In}$ and terminating at some $r_0\in l\cap f^{-1}_p(l)$ s.t. $l_1\subseteq\partial H_p$ (see the illustration in Fig.\ref{optA}) - consequentially, $f_p$ is continuous in $l_1$ (in particular,  $r_0$ is the closet point in $f^{-1}_p(l)\cap l$ to the fixed-point $P_{In}$). Therefore, for every $s\in l_1$ which is strictly interior to $l_1$, $s\not\in f^{-1}_p(l)$  - i.e., for such $s$, $f_p(s)$ is interior to $\overline{D_\alpha}$. Since by definition $r_0\in f^{-1}_p(l)\cap l$ it follows $f_p(l_1)$ is an arc in $D_\alpha$ with endpoints $r_1=f_p(r_0)\in l$ and $P_{In}$, which implies $\overline{f_p(l_1)}\cap\partial D_\alpha=\{P_{In},r_1\}$ (see the illustration in Fig.\ref{first1}).\\

Additionally, since the flow line connecting $r_0,r_1$ lies in $\{\dot{y}\geq0\}$ and since $r_0,r_1$ are both in $l$, recalling $l$ is parameterized by $l(x)=(x,-\frac{x}{a},\frac{x}{a})$, $x\in(0,c-ab)$ and that $P_{In}=(0,0,0)$, it follows $r_1$ is closer to $P_{In}$ than $r_0$ (as illustrated in Fig.\ref{first1}). Since $r_0\ne P_0$, it follows $r_1\ne P_{In}$ - that is, $\overline{f_p(H_p)}$ misses an arc on $l$ which begins at $P_{In}$ (see the illustration in Fig.\ref{first1}).\\

Now, recall that $P_0$ has no pre-images, as it lies on the heteroclinic trefoil. Now, recall the curve $\rho$ from Cor.\ref{disconcurve}, and that by Lemma \ref{lemhp} the curve $\rho$ also lies in $\partial H_p$, we conclude by Cor.\ref{disconcurve} there exists no $\{x_n\}\in\overline{H_p}$ s.t. $f_p(x_n)\to P_0$ (or, in other words, the sector $T$ given by Cor.\ref{disconcurve} lies outside of $H_p$ - see the illustration in Fig.\ref{optA}). Consequentially, we have $P_0\not\in f_p(\overline{D_\alpha}\setminus\{P_0\})$. Moreover, let us note that as $H_p$ is a component of $D_\alpha\setminus f^{-1}_p(l)$, by Prop.\ref{lem33} it follows the first-return map is continuous on $\overline{H_p}$, hence all in all we conclude $P_0\not\in\overline{f_p(H_p)}$. Now, consider the set $\Gamma$ of curves $\gamma$ in $D_\alpha$ which connect $P_{In}$ and $P_0$, we see $\Gamma\cap H_p$ folliates $H_p$. Finally, recall that for every $\gamma\in\Gamma$, as we suspend $\gamma$ with the flow, $\gamma$ loops around the heteroclinic trajectory - as depicted in Fig.\ref{fig16}. As such, since $\Gamma\cap H_p$ folliates $H_p$, all in all, as $P_0\not\in\overline{f_p(H_p)}$ and since $r_1$ is closer to $P_{In}$ then $r_0$, we conclude $f_p(H_p)$ can be imposed on $D_\alpha$ like a crescent-shaped figure with one tip at $P_{In}$ and another at $r_1$ (see the illustration in Fig.\ref{first1}).\\

\begin{figure}[h]
\centering
\begin{overpic}[width=0.40\textwidth]{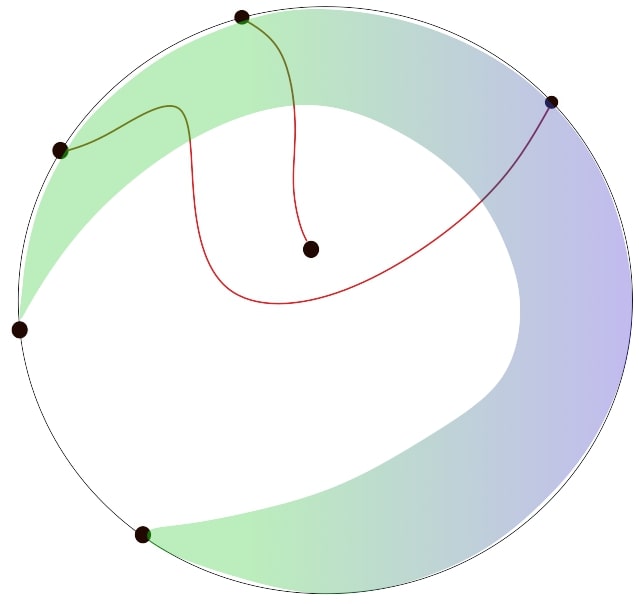}
\put(200,-10){$P_{In}$}
\put(500,550){$P_{0}$}
\put(400,350){$H_p$}
\put(300,430){$\rho$}
\put(600,350){$f_p(l_1)$}
\put(0,750){$r_0$}
\put(-40,390){$r_1$}
\put(350,820){$\delta_0$}
\put(590,850){$f_p(\rho)$}
\put(870,760){$P_{Out}$}
\put(990,380){$L=f_p(L)$}
\end{overpic}
\caption[The first-return map in Case $A$.]{\textit{Imposing $f_p(H_p)$ on $D_\alpha$ in Case $A$. $f_p(l_1)$ is the arc connecting $P_{In}$ and $r_1$.}}
\label{first1}
\end{figure}
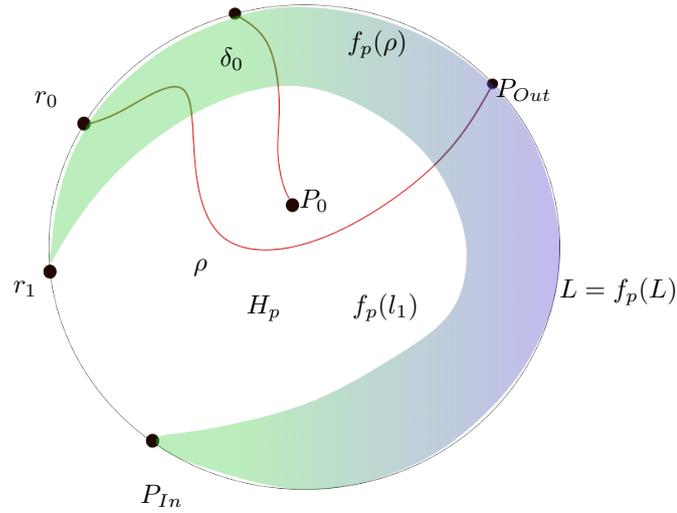

There is not much that can be said on the invariant set of $f_p$ in $\overline{H_p}$ from Fig.\ref{first1} alone - let alone factor its dynamics to those of the one-sided shift $\sigma:\{1,2\}^\mathbf{N}\to\{1,2\}^\mathbf{N}$ (let us remark, however, that it is easy to see it is non-empty - by definition, it includes both fixed points and the arc $L$). In order to say something smart about it, let us first extend it as follows - first, let us sketch $D_\alpha$ as a half-disc instead of a disc or a half plane, and let us glue to $D_\alpha$ another half-disc, $C_\alpha$, at the curve $l$, s.t. $V_\alpha=\overline{D_\alpha}\cup C_\alpha$ forms a disc (see the illustration in Fig.\ref{CALPHA1}).\\

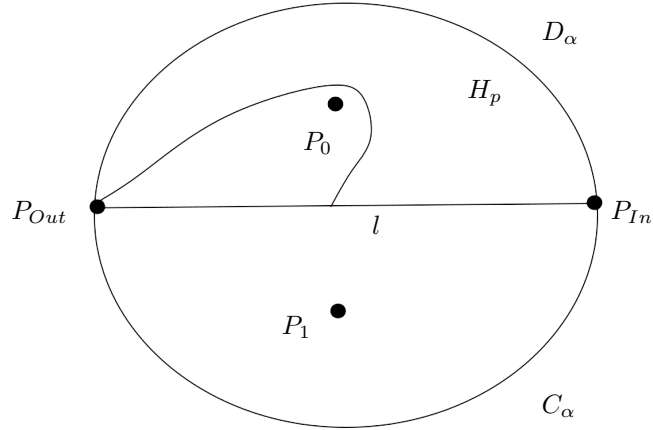
\begin{figure}[h]
\centering
\begin{overpic}[width=0.40\textwidth]{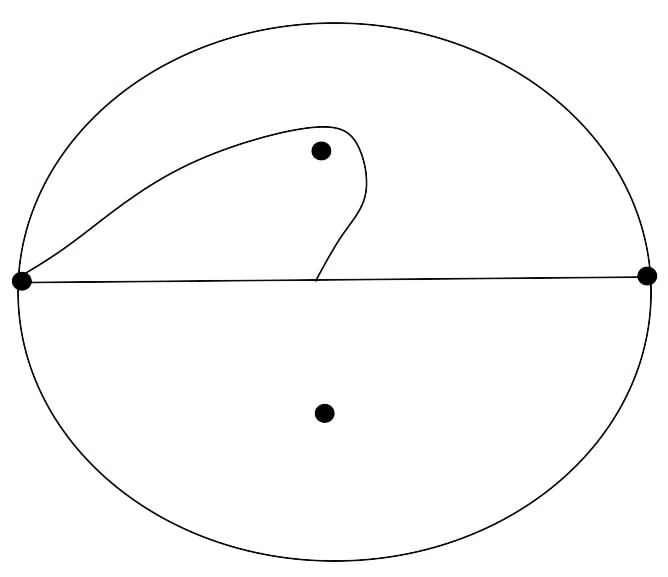}
\put(420,550){$P_{0}$}
\put(730,650){$H_p$}
\put(550,390){$l$}
\put(-130,420){$P_{Out}$}
\put(380,200){$P_1$}
\put(870,760){$D_\alpha$}
\put(870,50){$C_\alpha$}
\put(1000,420){$P_{In}$}
\end{overpic}
\caption[The disc $V_\alpha$.]{\textit{$D_\alpha$ is the upper half disc, glued to $C_\alpha$, the lower half disc at the curve $l$ - together, they constitute $V_\alpha$}.}
\label{CALPHA1}
\end{figure}

Now, choose some $P_1$, an interior point to $C_\alpha$ (as illustrated in Fig.\ref{CALPHA1}), and let $F:V_\alpha\to V_\alpha$ be some homeomorphism s.t. the following is satisfied (see the illustration in Fig.\ref{F1}):

\begin{itemize}
    \item $F|_{\overline{H_p}}=f_p$, i.e., $F$ coincides with $f_p$ on $\overline{H_p}$ (we allow $F$ to differ from $f_p$ on $\overline{D_\alpha}\setminus\overline{H_p}$).
    \item $F(D_\alpha\setminus H_p)$ includes an region in $V_\alpha$, connecting $P_1$ and $l$.
    \item $F(\partial V_\alpha)=\partial V_\alpha$, and in particular - $F(P_{In})=P_{In}$, while $F(P_0)=P_1$ and $F(P_1)=P_0$.
\end{itemize}

Our strategy of proof now is as follows: we first consider the surface $S=V_\alpha\setminus\{P_0,P_1\}$, as illustrated in Fig.\ref{F1}, after which we study the isotopy class of $F$ in $S$. That is, we consider the isotopy class of all maps $F'$ isotopic to $F$ in $S$, s.t. $F'(P_{In})=P_{In}$, $F'(P_0)=P_1$, and $F'(P_1)=P_0$. As we will now prove, this isotopy class includes a Pseudo-Anosov map whose periodic orbits are all isotopy stable - from which the conclusion would follow.\\ 

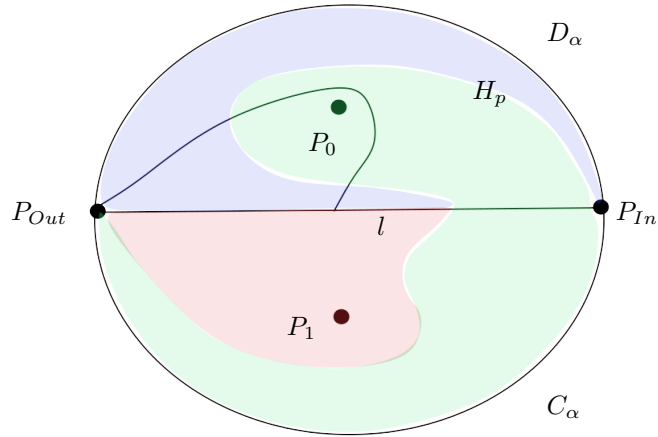
\begin{figure}[h]
\centering
\begin{overpic}[width=0.40\textwidth]{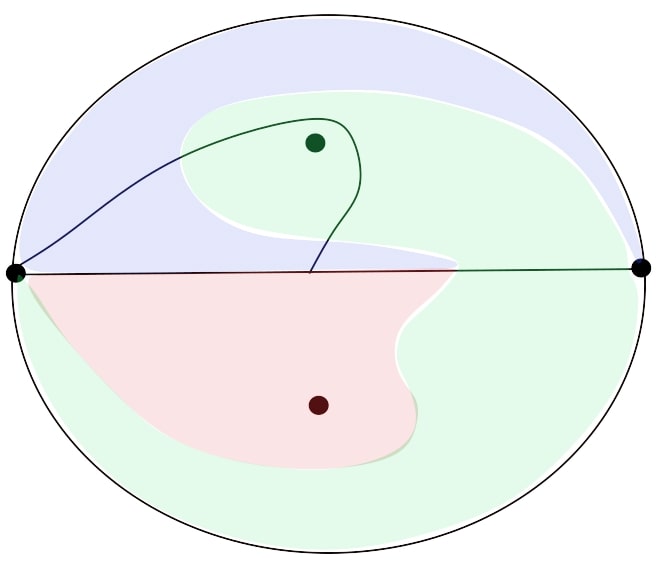}
\put(420,550){$P_{0}$}
\put(730,650){$H_p$}
\put(550,390){$l$}
\put(-140,420){$P_{Out}$}
\put(380,200){$P_1$}
\put(870,760){$D_\alpha$}
\put(870,50){$C_\alpha$}
\put(1000,420){$P_{In}$}
\end{overpic}
\caption[The homeomorphism $F$.]{\textit{The action of $F$ on $V_\alpha$ - the blue region corresponds to $F(H_p)$, the red to $F(D_\alpha\setminus H_p)$, and the green to $F(C_\alpha)$}.}
\label{F1}
\end{figure}

To begin, per the outline in Sect.1 of the Betsvina-Handel Algorithm (for more details, see \cite{BeH} or \cite{Bo}), let us choose a spine for $S$, i.e., a graph embedded in $S$, which is a retract of $S$ - in particular, we choose a graph $T\subseteq S$, with vertices at $P_{In},P_0$ and $P_1$, with two edges: $T_1$ connecting $P_0$ to $P_{In}$ and $T_2$, connecting $P_1$ to $P_{In}$ (see the illustration in Fig.\ref{T13}). Now, consider $F(T)$ - recalling $f_p(H_p)$ lies away $P_0$ and that $F|_{\overline{H_p}}=f_p$, by $F(P_{In})=P_{In}$, $F(P_0)=P_1$ and $F(P_1)=P_0$ we conclude $F(T)$ looks as in Fig.\ref{T1}.\\

\begin{figure}[h]
\centering
\begin{overpic}[width=0.40\textwidth]{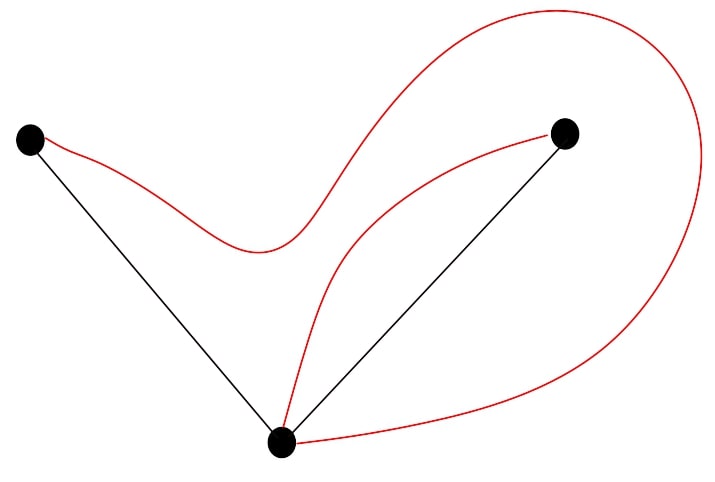}
\put(380,550){$F(T_1)$}
\put(250,250){$F(T_2)$}
\put(280,20){$P_{In}$}
\put(-90,420){$P_1$}
\put(70,300){$T_2$}
\put(650,250){$T_1$}
\put(800,420){$P_0$}
\end{overpic}
\caption[The graph $T$ and its image under $F$.]{\textit{The graph $T$ with vertices $P_{In},P_0$ and $P_1$ (in black), and its image under $F$ (in red).}}
\label{T1}
\end{figure}
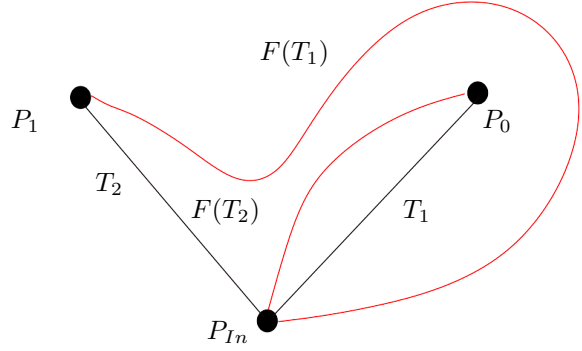

Consequentially, $F(T_1)$ begins at $P_{In}$, surrounds $T_1$ from both sides, and then stretches along $T_2$ until terminating at $P_1$ - conversely, $F(T_2)$ is an edge connecting $P_{In}$ to $P_0$, stretching along $T_1$ (see the illustration in Fig.\ref{T1}). Now, let us retract $F$ as described in Fig.\ref{T13} to a graph map $g:T\to T$, and consider the matrix $A=\{a_{i,j}\}_{1\leq i,j\leq 2}$, where $a_{i,j}$ denotes the times the curve $g(T_i)$ covers the side $T_j$ - it is easy to see $A$ is simply the matrix:
\begin{equation*}
\begin{pmatrix}
    n&1\\
    1&0
\end{pmatrix}
\end{equation*}
Where $n\geq2$. That is - $g(T_1)$ covers $T_1$ at least twice, while covering $T_2$ only once (conversely, $g(T_2)$ only covers $T_1$ once) . By computation, the spectral radius of $A$ is $\frac{\sqrt{n^2+4}+n}{2}$ (i.e., its maximum positive eigenvalue), with a corresponding eigenvector $(\frac{\sqrt{n^2+4}+n}{2},1)$. As $n\geq2$, by $\sqrt{n^2+4}>2-n$ it follows $\frac{\sqrt{n^2+4}+n}{2}>1$ which proves the spectral radius is greater than $1$. \\

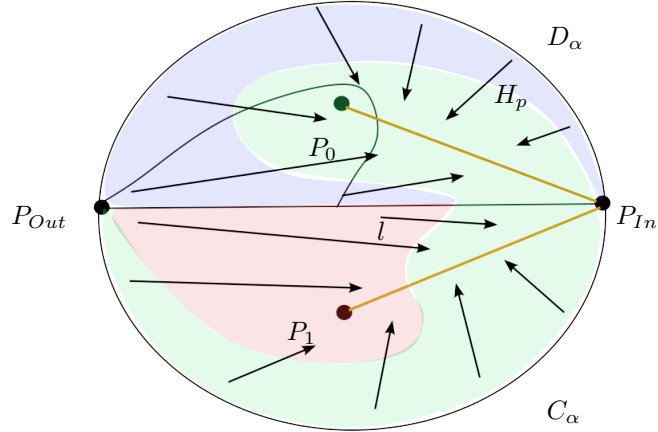
\begin{figure}[h]
\centering
\begin{overpic}[width=0.40\textwidth]{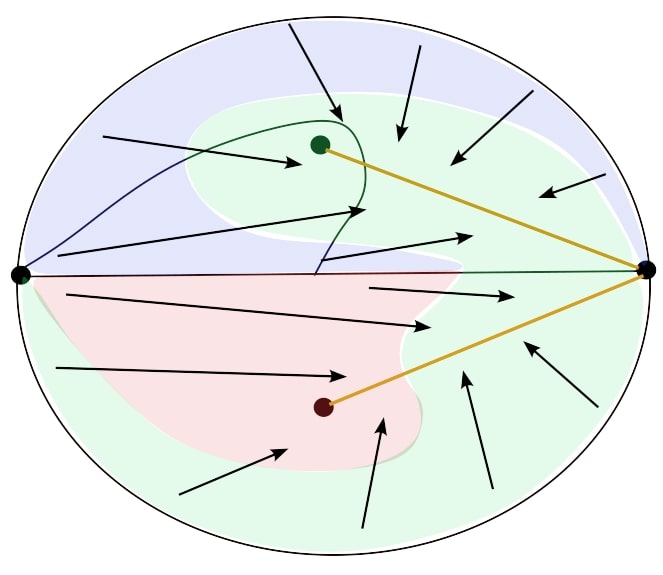}
\put(420,550){$P_{0}$}
\put(770,650){$H_p$}
\put(550,390){$l$}
\put(-140,420){$P_{Out}$}
\put(380,200){$P_1$}
\put(870,760){$D_\alpha$}
\put(870,50){$C_\alpha$}
\put(1000,420){$P_{In}$}
\end{overpic}
\caption[Collapsing $F$ to a graph map.]{\textit{Collapsing $F$ to a graph map - as can be seen, $H_p$ and $f_p(H_p)$ are collapsed on $T_1$, the edge connecting $P_{In}$ and $P_0$.}}
\label{T13}
\end{figure}

Now, let us begin applying the Betsvina-Handel Algorithm to $F$ -  namely, going over steps $1-7$ in \cite{BeH} we deform $F$ isotopically s.t. the graph map $g$ is "straightened" - i.e., we deform $F$ isotopically in $S$ which ensures the resulting graph map $g$ is injective on every edge in $T$ (see the illustration in Fig.\ref{T15}). In particular, let us remark that all the dynamics of $F$ in $H_p$ (and consequentially, those of the first-return map $f_p$ in $H_p$) are collapsed by the deformation precisely to $g(T_1)\cap T_1$ - and moreover, any point in the invariant set of $g$ in $T_1$ corresponds under the deformation to some point in the invariant set of $F$ in $\overline{H_p}$, i.e,, to that of $f_p$ in $\overline{H_p}$. Consequentially, provided we prove $F$ has infinitely many periodic orbits in $\overline{H_p}$, we're done.\\

\begin{figure}[h]
\centering
\begin{overpic}[width=0.40\textwidth]{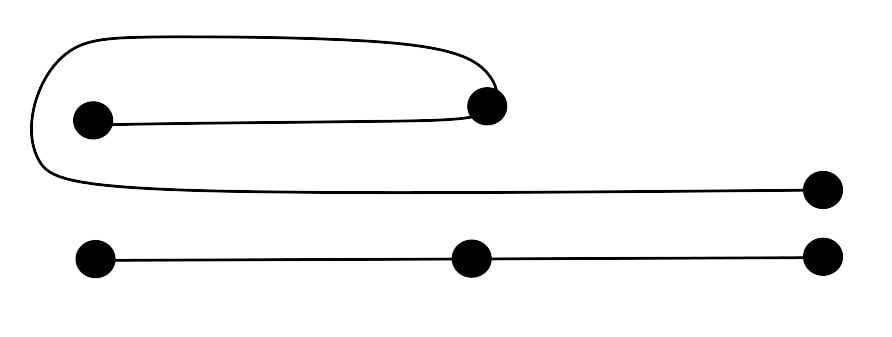}
\put(450,0){$P_{In}$}
\put(600,250){$g(P_{In})$}
\put(0,50){$P_{0}$}
\put(-150,230){$g(P_1)$}
\put(950,50){$P_1$}
\put(970,140){$g(P_0)$}
\end{overpic}
\caption[The induced graph map, $g$.]{\textit{As can be seen, $g$ acts on $T$ as follows - $g(T_1)$ covers $T_1$ twice and stretched towards $P_1$, while $g(T_2)$ covers $T_1$.}}
\label{T14}
\end{figure}
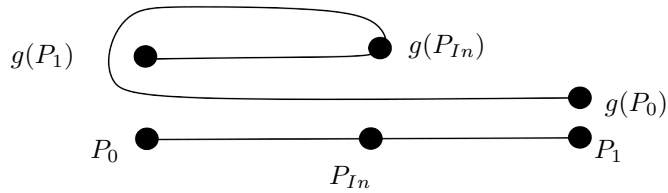

To do so, following Ch.$3.4$ and Ch.$4.4$ in \cite{BeH}, as the spectral radius of the matrix $A$ is $\frac{\sqrt{n^2+4}+n}{2}>1$, it follows $F:S\to S$ can be isotopically deformed to a Pseudo-Anosov homeomorphism $G:S\to S$. Consequentially, since $g(T_1)$ covers $T_1$ twice, it follows by the argument given in Ch$.3.4$ and Ch.4.4 of \cite{BeH} (see also Th.\ref{betshan}) that there exists a $G$-invariant set, $I\subseteq S$ which includes infinitely many periodic orbits - moreover, by Th.\ref{betshan}, the set $I$ lies precisely in the regions of $S$ which are retracted to $T_1$, that is, in $\overline{H_p}$. In fact, it is easy to see that on $\overline{H_p}$, $G$ acts as a Smale Horseshoe map (see the illustration in Fig.\ref{T15}). Namely, we have just proven:

\begin{corollary}
    \label{hors1} $g$ has periodic orbits of all minimal periods in $T_1$ - consequentially, so does $G$ in $\overline{H_p}$, which from now on we denote by $ABCD$. In particular, $G:ABCD\to ABCD$ is a Smale Horseshoe map.
\end{corollary}

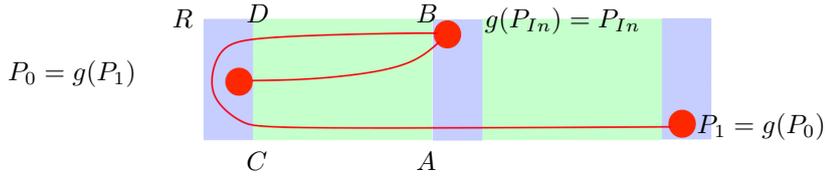
\begin{figure}[h]
\centering
\begin{overpic}[width=0.40\textwidth]{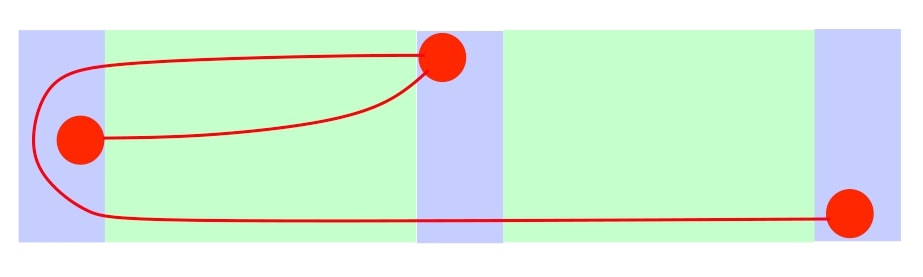}
\put(420,-20){$A$}
\put(420,260){$B$}
\put(550,250){$g(P_{In})=P_{In}$}
\put(100,-20){$C$}
\put(100,260){$D$}
\put(-40,250){$R$}
\put(-350,150){$P_0=g(P_1)$}
\put(950,50){$P_1=g(P_0)$}

\end{overpic}
\caption[The action of $G$ on $ABCD$.]{\textit{The red curve is $g(T)$ embedded in $S$ - due to these properties, $G$ acts on $ABCD$ as a Smale Horseshoe map (which is continuously deformed to a subset of $\overline{H_p}$ by the isotopy). The region $R$ lies in $S$ outside $ABCD$.}}
\label{T15}
\end{figure}

Now, let us recall that by definition, $G$ is a homeomorphism of $\overline{S}$ which permutes the points $\{P_0,P_{In},P_1\}$ as follows: $G(P_{In})=P_{In}$, $G(P_1)=P_0$, and $G(P_0)=P_1$- additionally, recall the following result given as Th.1, Th.2 and Remark $2.4$ in \cite{Han}:

\begin{claim}
\label{stability}    Let $S$ be a surface (possibly a with boundary) and let $\phi:S\to S$ be Pseudo-Anosov homeomorphism - then, we have the following:
\begin{itemize}
    \item If $x_1$ is periodic of minimal period $k$ for $\phi$, then given any $\psi:S\to S$ which is isotopic to $\phi$, when $\phi$ is isotoped to $\psi$ the point $x_1$ is continuously deformed to $x_2$, a periodic point for $\psi$ of minimal period $k$. That is, the periodic dynamics of $\phi$ are \textbf{unremovable} - i.e., they are not destroyed by an isotopy.
    \item For every homeomorphism $\psi:S\to S$ isotopic to $\phi$, there exists a closed $Y\subseteq S$ and a continuous surjection $\pi:Y\to S$ s.t. $\pi\circ \psi=\phi\circ\pi$. In particular, if $x\in S$ is periodic of minimal period $n$ for $\phi$, then $\pi^{-1}(x)$ includes a periodic orbit of minimal period $n$ for $\psi$ - i.e., the periodic dynamics are \textbf{uncollapsible}. In particular, when we isotope $\phi$ to $\psi$, no two periodic orbits of $\phi$ collapse into one another.
\end{itemize}
\end{claim}

As a consequence of Th.\ref{stability} and Th.\ref{betshan}, every periodic orbit generated by $G$  also persists as a periodic orbit for $F$  - and moreover, since by the isotopy it follows $H_p$ is continuously deformed to the rectangle $ABCD$, we conclude that when we return from $G$ to $F$, the periodic orbits of $G$ in $ABCD$ are all deformed to periodic orbits for $F$ in $\overline{H_p}$ - and moreover, no two periodic orbits are destroyed, collapsed into one another, or change their minimal period. Consequentially, $F$ generates infinitely many periodic orbits in its invariant set in $\overline{H_p}$ - and all in all, since $F|_{\overline{H_p}}=f_p$, the first-return map at trefoil parameter in Case $A$, we conclude:
\begin{corollary}
    \label{isotret1}
    Let $F_p$ be a trefoil parameter s.t. Case $A$ is satisfied - then, the first-return map $f_p:\overline{H_p}\to\overline{D_\alpha}$ generates infinitely many periodic orbits, of all minimal periods.
\end{corollary}

Having proven Cor.\ref{isotret1}, it follows that in order to conclude the proof of Th.\ref{th31} in Case $A$, we now need to define symbolic dynamics on some $f_p$-invariant $Q\subseteq H_p$, prove $f_p$ is continuous on $Q$ - and most importantly, show that the periodic orbits given by Cor.\ref{isotret1} all lie in $Q$. We first prove we can define symbolic dynamics on a set $I$, some invariant set of $f_p$ in $\overline{H_p}$. We do so as follows:

\begin{lemma}
\label{sym1}    There exists a curve $\Delta\subseteq\overline{H_p}$, s.t. $\overline{H_p}\setminus\Delta$ includes two components - $D_1$ and $D_2$. Consequentially, let $I'$ denote the maximal invariant set in $\overline{H_p}\setminus\Delta$ - then, there exists a continuous $\pi':I'\to\{1,2\}^\mathbf{N}$ s.t. $\pi'\circ f_p=\sigma\circ\pi'$. Moreover, all the periodic orbits given by Cor.\ref{isotret1} are in $I'$.
\end{lemma}
\begin{proof}
  Recall $G(ABCD)$ is a Smale Horseshoe, which implies there is a region $R$ in the surface $S$ (as indicated in Fig.\ref{T15}) s.t. $G^{-1}(R)\cap ABCD=V_1$ is a sub-rectangle of $ABCD$, connecting the $AC$ and $BD$ sides. Moreover, the different periodic orbits of $G$ in $ABCD$ are separated by the pre-images of $G$. As we deform isotopically $G$ back to $F$, because $ABCD$ is deformed to $\overline{H_p}$, the rectangle $V_1$ is deformed to $\Delta'$, some subset of $\overline{H_p}$ - and as $F$ coincides with $f_p$ on $\overline{H_p}$, it follows the different periodic orbits given by Cor.\ref{isotret1} are separated from one another by the pre-images of $\Delta'$ (w.r.t $f_p$). This proves there exists a curve $\Delta\subseteq\Delta'$ s.t. $\overline{H_p}\setminus\Delta$ is composed of two components: $D_1$ and $D_2$ (see the illustration in Fig.\ref{part1}).\\ 
  
Note that $I'$, the maximal invariant set of $f_p$ in $\overline{H_p}\setminus\Delta$, is a strict subset of the invariant set of $f_p$ in $\overline{H_p}$ - therefore, because $f_p$ is continuous on $\overline{H_p}$ it is also continuous on $I'$ as well. Consequentially, there exists a continuous $\pi':I'\to\{1,2\}^\mathbf{N}$ s.t. both $\pi'\circ f_p=\sigma\circ\pi'$, and $\pi'(x)=\{s_n\}_{n\geq0}$, where $s_i$ is $1$ whenever $f^i_p(x)\in D_1$ and $2$ otherwise. Moreover, it follows by the discussion above that all the periodic orbits given by Cor.\ref{isotret1} are in $I'$ and Lemma \ref{sym1} now follows.
\end{proof}
\begin{figure}[h]
\centering
\begin{overpic}[width=0.40\textwidth]{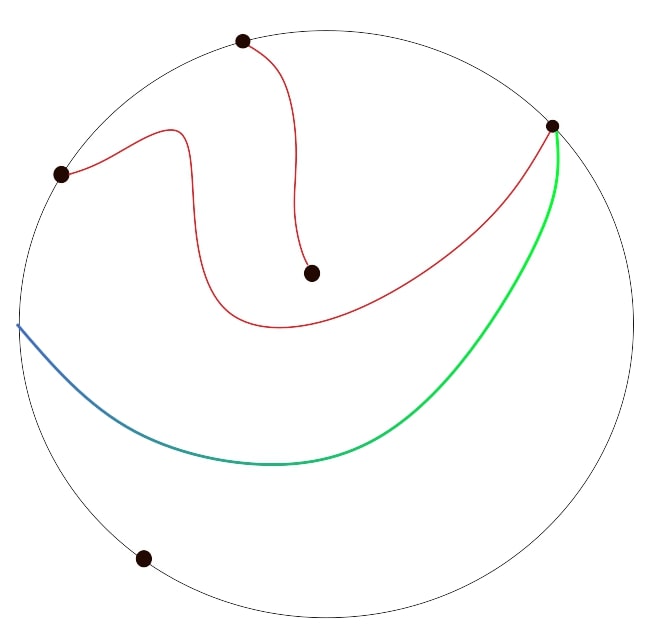}
\put(200,-10){$P_{In}$}
\put(500,550){$P_{0}$}
\put(400,350){$D_2$}
\put(500,150){$D_1$}
\put(370,750){$\delta$}
\put(-40,420){$\Delta$}
\put(350,940){$\delta_0$}
\put(870,760){$P_{Out}$}
\put(990,380){$L$}
\end{overpic}
\caption[The curve $\Delta$ in Case $A$.]{\textit{The curve $\Delta\subseteq f^{-1}_p(\delta)$ in Case $A$.}}
\label{part1}
\end{figure}
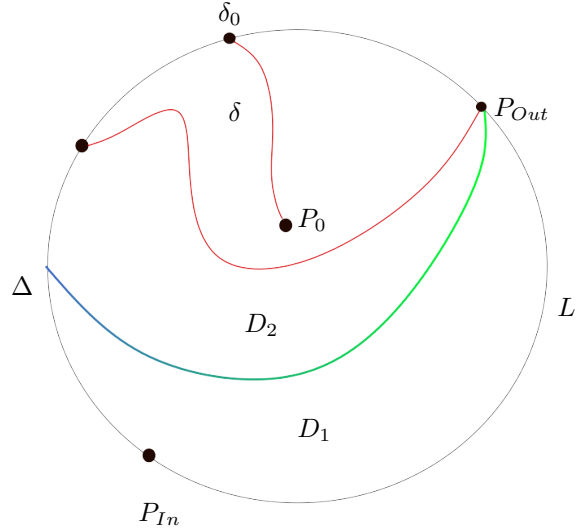
    
Having proven Lemma \ref{sym1}, we now prove:

\begin{lemma}
    \label{continu1} Let $F_p$ be a trefoil parameter falling into case $A$. Then, given any periodic $s\in\{1,2\}^\mathbf{N}$, we have $s\in\pi'(I')$.
\end{lemma}
\begin{proof}
The proof is similar to that of Lemma \ref{sym1}. As we isotopically deform $F:\overline{H_p}\to S$ to $G:ABCD\to S$, the set $\Delta'$ introduced in the proof of Lemma \ref{sym1} is continuously deformed to the sub-rectangle $V_1$ inside $ABCD$, defined by $G(V_1)\cap ABCD=\emptyset$ (see the illustration in Fig.\ref{T16}). Recall the symbolic dynamics of $G$ on its invariant set in $ABCD$ are defined by $L_1$ and $L_2$, the sub-rectangles composing $ABCD\setminus V_1$ -  and similarly to how $\Delta'$ is deformed to $V_1$, the isotopy also deforms $D_1$ to $L_1$ and $D_2$ to $L_2$. Consequentially, it follows that if the periodic orbit $\{x_1,...,x_k\}$ for $G$ corresponds to some symbol in $\{1,2\}^\mathbf{N}$, as it is deformed to $\{y_1,...,y_k\}$, a periodic orbit for $f_p=F|_{\overline{H_p}}$, its symbolic dynamics do not change - i.e., if $\{x_1,...,x_k\}$ corresponds to some periodic $s\in\{1,2\}^\mathbf{N}$, the same is true for $\{y_1,...,y_k\}$. Lemma \ref{continu1} now follows.
\end{proof}

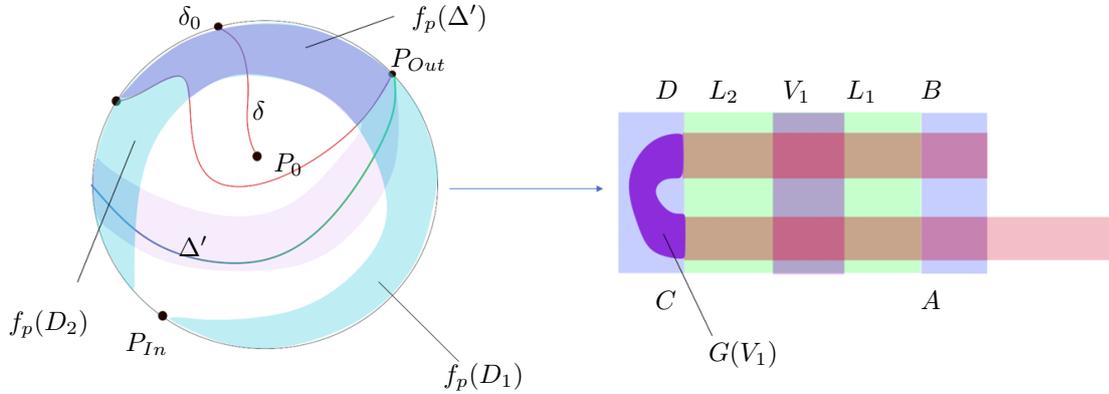
\begin{figure}[h]
\centering
\begin{overpic}[width=0.8\textwidth]{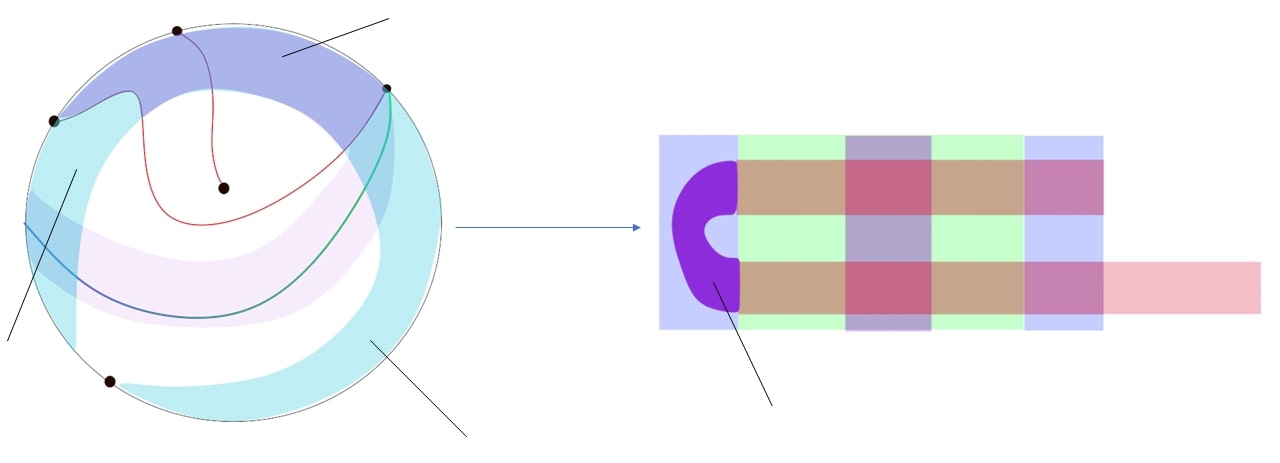}
\put(50,30){$P_{In}$}
\put(190,200){$P_{0}$}
\put(320,340){$f_p(\Delta')$}
\put(350,0){$f_p(D_1)$}
\put(170,250){$\delta$}
\put(100,120){$\Delta'$}
\put(100,340){$\delta_0$}
\put(300,300){$P_{Out}$}
\put(-60,50){$f_p(D_2)$}
\put(600,20){$G(V_1)$}
\put(550,70){$C$}
\put(550,270){$D$}
\put(670,270){$V_1$}
\put(800,270){$B$}
\put(800,70){$A$}
\put(730,270){$L_1$}
\put(600,270){$L_2$}
\end{overpic}
\caption[Deforming $\Delta$ to $L_3$.]{\textit{The isotopy of $F|_{\overline{H_p}}=f_p$ on $H_p$ to $G:ABCD\to S$ deforms the purple region $\Gamma$ to the sub-rectangle $L_3$.}}
\label{T16}
\end{figure}

We now conclude the proof of Th.\ref{th31} for Case $A$. To do so, recall the factor map $\pi'$ given by Lemma \ref{sym1}, and for any periodic $s\in\{1,2\}^\mathbf{N}$ of minimal period $k$ which is not the constant $\{1,1,1...\}$, let us denote by $D_s$ the component of $\pi'^{-1}(s)$ s.t. the following holds:

\begin{itemize}
    \item $D_s$ contains a periodic orbit of minimal period $k$, $\{x_1,...,x_k\}$.
    \item  $\{x_1,...x_k\}$ is deformed isotopically to a periodic orbit for $G$, of the same minimal period and the same symbol.
\end{itemize}

Such a $D_s$ exists by both Th.\ref{stability} and Cor.\ref{continu1}. For the constant $s=\{1,1,1...\}$, set $D_s=\{P_{In}\}\cup L$. Finally, let $Q$ denote the collection of such $D_s$ - by definition, $Q\subseteq I'$, where $I'$ is as in Cor.\ref{sym1}. Setting $\pi=\pi'|_Q$ and summarizing our results, we conclude the following:

\begin{itemize}
    \item By Cor.\ref{sym1}, the first-return map $f_p$ is continuous on $Q$.
    \item Additionally, by Cor.\ref{sym1}, there exists a continuous $\pi:Q\to\{1,2\}^\mathbf{N}$ s.t. $\pi\circ f_p=\sigma\circ\pi$ (where $\sigma:\{1,2\}^\mathbf{N}\to\{1,2\}^\mathbf{N}$ denotes the one-sided shift).
    \item By Lemma \ref{continu1}, $\pi(Q)$ includes every periodic symbol in $\{1,2\}^\mathbf{N}$.
    \item Moreover, by our construction of $D_s$ above, if $s$ is periodic of minimal period $k$, $\pi^{-1}(s)$ includes a periodic orbit of minimal period $k$ for $f_p$.
\end{itemize}

All in all, we have completed the proof of Th.\ref{th31} for Case $A$.

\subsubsection{Stage $II$ - proving Th.\ref{th31} for Case $B$:}

\begin{figure}[h]
\centering
\begin{overpic}[width=0.7\textwidth]{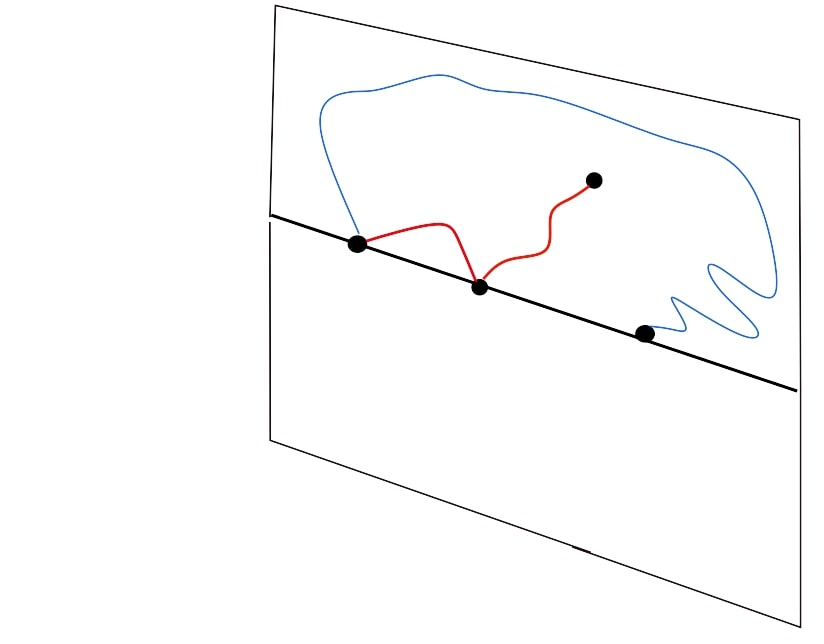}
\put(750,320){$P_{In}$}
\put(380,450){$P_{Out}$}
\put(640,520){$\delta$}
\put(580,370){$\delta_0$}
\put(800,630){$L$}
\put(350,730){$U_p$}
\put(460,620){$H_p$}
\put(500,470){$T$}
\put(700,580){$P_0$}
\put(350,265){$L_p$}
\end{overpic}
\caption[Case $B$.]{\textit{The cross-section $H_p$, in case $B$, with $f^{-1}_p(l)$ sketched as a curve.}}
\label{optB}
\end{figure}
Having proven Th.\ref{th31} for trefoil parameters $F_p$ which fall into Case $A$, we now prove Th.\ref{th31} for trefoil parameters which fall into Case $B$ - that is, we now assume the curve $\delta$ and the cross-section $H_p$ satisfy the following:

\begin{itemize}
    \item Both $\delta$ and $P_{In}$ lie in $\partial H_p$ (see the illustration in Fig.\ref{optB}).
    \item Additionally, we assume $\delta_0$, the endpoint on $\delta$ in $l$, satisfies $\delta_0\ne P_{Out}$.
\end{itemize}

Even though the heart of the proof will be essentially the same as in Case $A$, we will need to slightly adjust several of our arguments - if only because the first-return map $f_p$ in Case $B$ behaves somewhat differently. We begin with the following Lemma, which studies just how the first-return map folds $H_p$ inside $D_\alpha$:

\begin{lemma}
    \label{corloop}
    Let $F_p$ be a trefoil parameter falling into Case $B$. Then, there exists a curve $\gamma\subseteq H_p$ satisfying:
    \begin{itemize}
        \item $\overline{\gamma}\cap\partial H_p=\{P_{In},P_0\}$.
        \item $\overline{f_p(\gamma)}$ is a closed loop in $D_\alpha$, s.t. $\overline{f_p(\gamma)}\cap\partial D_\alpha=\{P_{In}\}$ - i.e., $f_p(\gamma)$ begins and terminates at $P_{In}$.
        \item $\overline{f_p(\gamma)}$ winds once around $P_0$ - in particular, it separates $P_0$ from $\partial D_\alpha\setminus\{P_{In}\}$ (see the illustration in Fig.\ref{DALPHA133}).
    \end{itemize}

\end{lemma}
\begin{proof}
By Prop.\ref{lem33}, the components of $f^{-1}_p(l)\cap D_\alpha$ are either $\delta$, or curves with at least two endpoints on $l$. Consequentially, since $H_p$ is a component of ${D_\alpha}\setminus f^{-1}_p(l)$, by Prop.\ref{lem33} it follows $H_p$ is a simply connected two-dimensional set, bounded by curves: i.e., a topological disc (see the illustration in Fig.\ref{optB} - in Case $B$ $H_p$ is not a Jordan domain as $\partial H_p$ is not homeomorphic to $S^1$). As such, since by Assumption $2$ both $P_{In},\delta$ lie on $\partial H_p$, because $H_p$ is a topological disc and because $P_0\in\delta$ there exists a curve $\gamma\subseteq H_p$ s.t. $\overline\gamma\cap\partial H_p=\{P_{In},P_{0}\}$ - i.e., $\gamma$ satisfies the following (see the illustration in Fig.\ref{DALPHA133}):
\begin{itemize}
    \item $\gamma$ is homeomorphic to an open interval.
    \item The endpoints of $\gamma$ are simply $P_{In}$ and $P_0$.
    \item Every point on $\gamma$ is strictly interior to $H_p$.
\end{itemize}

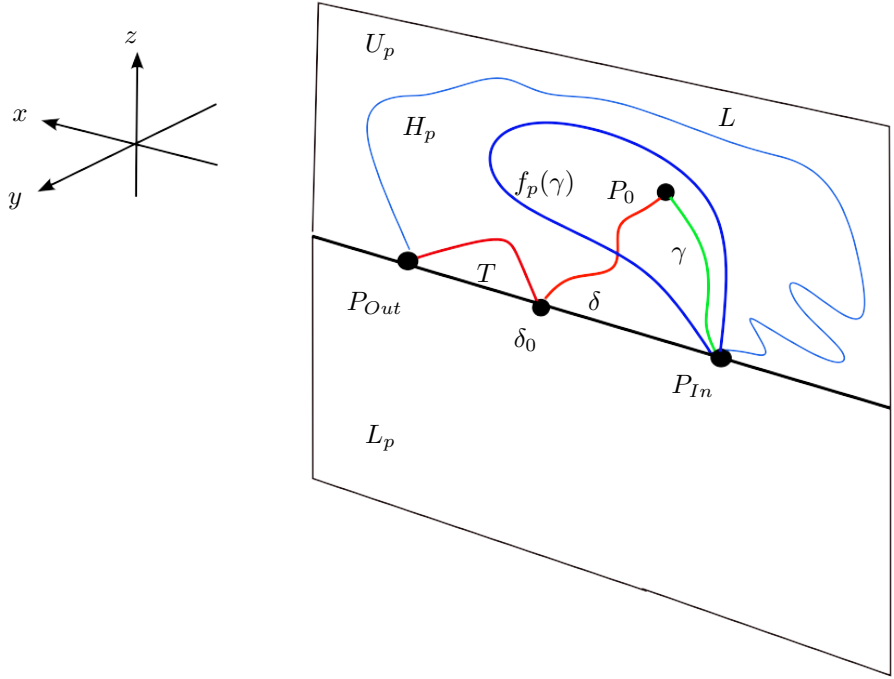
\begin{figure}[h]
\centering
\begin{overpic}[width=0.7\textwidth]{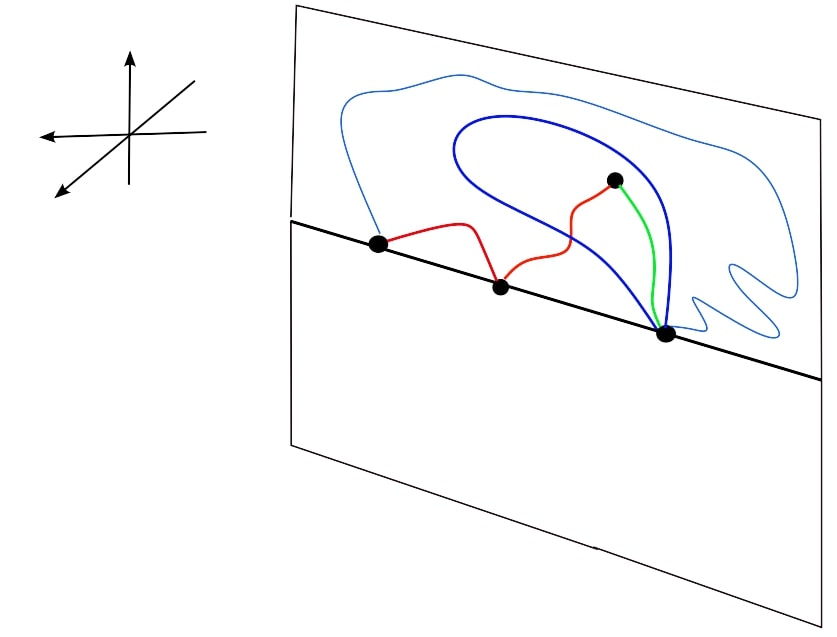}

\put(750,465){$\gamma$}
\put(580,540){$f_p(\gamma)$}
\put(750,320){$P_{In}$}
\put(400,410){$P_{Out}$}
\put(660,410){$\delta$}
\put(580,370){$\delta_0$}
\put(800,610){$L$}
\put(420,690){$U_p$}
\put(460,600){$H_p$}
\put(680,530){$P_0$}
\put(420,265){$L_p$}
\put(540,440){$T$}
\put(25,590){$x$}
\put(40,515){$y$}
\put(150,710){$z$}
\end{overpic}
\caption[$\gamma$ and $f_p(\gamma)$.]{\textit{$\gamma$ is the green curve inside $H_p$, which returns to the cross-section as $f_p(\gamma)$ - the closed blue loop, which begins and terminates at $P_{In}$ and winds once around $P_0$. Again, for simplicity we sketch $f^{-1}_p(l)$ as a curve.}}
\label{DALPHA133}
\end{figure}

As a consequence, since $H_p$ is a component of ${D_\alpha}\setminus f^{-1}_p(l)$, it trivially follows $\gamma\cap f^{-1}_p(l)=\emptyset$ - therefore, by Prop.\ref{lem33} $f_p$ is continuous on $\gamma$. Now, recall the curve $\rho$ given by Cor.\ref{disconcurve} lies on $\partial H_p$ (see Lemma \ref{lemhp}) - and that the sector $T$ given by Cor.\ref{disconcurve} lies outside of $H_p$ (see the illustration in Fig.\ref{optB}). Therefore, there is no sequence $\{x_n\}_n\subseteq\overline{H_p}$ s.t. $f_p(x_n)\to P_0$, which implies that for every sequence $\{s_n\}_n\subseteq\gamma$, $f_p(s_n)\not\to P_0$ - hence, $P_0\not\in\overline{f_p(\gamma)}$.\\

Now, since the trajectory of $P_0$ tends to the fixed-point $P_{In}$ (in infinite time), we conclude $f_p(\gamma)$ is a closed loop in $D_\alpha$ which begins and terminates at $P_{In}$. Moreover, by the existence of the heteroclinic trefoil knot (see Def.\ref{def32}), it follows we can now apply the argument given at the beginning of this section and conclude $f_p(\gamma)$ separates $P_0,\partial D_\alpha$ - i.e., $f_p(\gamma)$ is a closed loop in $D_\alpha$ as in Fig.\ref{fig16}, Fig.\ref{fig17} and Fig.\ref{DALPHA133}. In particular, $\overline{f_p(\gamma)}$ is a curve which begins and terminates at the fixed-point $P_{In}$, winding once around $P_0$. Lemma \ref{corloop} is now proven.
\end{proof}
\begin{figure}[h]
\centering
\begin{overpic}[width=0.35\textwidth]{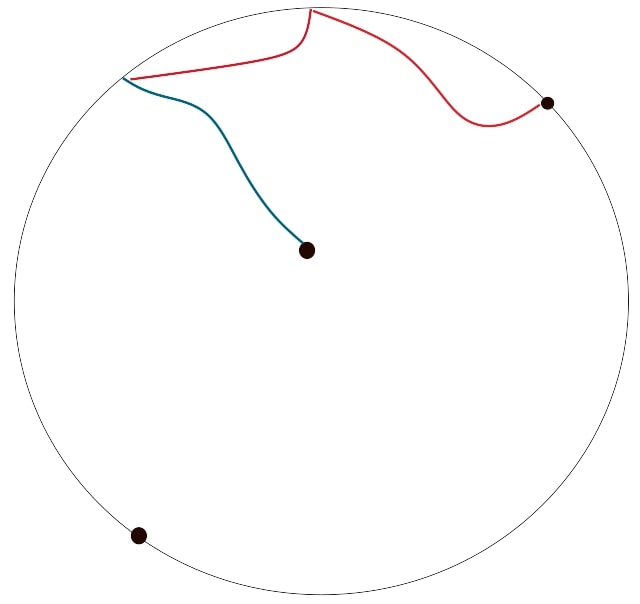}
\put(100,70){$P_{In}$}
\put(570,85){$L$}
\put(290,700){$\delta$}
\put(110,830){$\delta_0$}
\put(690,790){$T$}
\put(-20,680){$L_1$}
\put(370,535){$P_0$}
\put(670,335){$H_p$}
\put(850,810){$P_{Out}$}
\end{overpic}
\caption[$D_\alpha$ and $H_p$ as topological discs.]{\textit{The cross section $D_\alpha$, sketched as a topological disc (rather then a subset of a half-plane, as in Fig.\ref{DALPHA133}). The sub-arc $L_1$ denotes $[P_{In},\delta_0)$, and $L$ denotes $\partial B_\alpha\cap U_p$ - the boundary arc connecting $P_{In},P_{Out}$ (see Def.\ref{def32}.)}}
\label{DALPHA1222}
\end{figure}
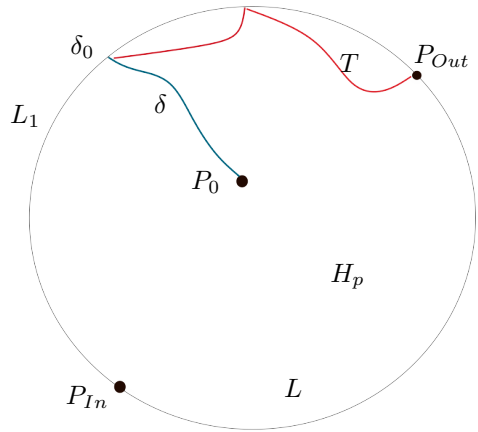

We are now ready to sketch $f_p(H_p)$ inside $D_\alpha$. To do so, let us recall $f_p(L)=L$ (where $L=\partial B_\alpha\cap U_p$ is a curve connecting $P_{In},P_{Out}$ - see Def.\ref{def32}). Now, note the following facts:

\begin{itemize}
    \item Replicating the argument in Cor.\ref{corloop} above, as the sector $T$ lies outside $H_p$ (and as $T$ and $H_p$ are separated by the curve $\rho$), for every $\{x_n\}\in H_p$, $f_p(x_n)\not\to P_0$ - that is, $P_0\not\in\overline{f_p(H_p)}$.
    \item Recall $f_p(\delta)$ is an arc on $l$ connecting $P_{In}$ and $f_p(\delta_0)$, and that $f^2_p(\delta)$ lies strictly inside $D_\alpha$ (see the proof of Prop.\ref{lem33}) - i.e., $f^2_p(\delta)$ is an arc interior to $D_\alpha$, beginning at $P_{In}$ and terminating at $f^2_p(\delta_0)$.
    \item  Finally, recall that in Case $B$ the cross-section $H_p$ is homeomorphic to a slit topological disc, with the slit corresponding to $\delta$ (see the discussion immediately after Lemma \ref{lemhp}). As such, for every $s\in\delta$ there exists some $r>0$ s.t. $B_r(s)\setminus\delta$ includes two components, $C_1$ and $C_2$ s.t. $C_i\subseteq H_p$, $i=1,2$, and for $\{s_n\}_n\subseteq C_i$, $s_n\to s$, $f_p(s_n)\to f^i_p(s), i=1,2$ (where $B_r(s)$ denotes a planar disc on $D_\alpha$).
\end{itemize}

As such, it follows $f_p(H_p)$ is as in Fig.\ref{firstpic} - that is, the domain winds around $P_0$, and encloses a domain $D_0\subseteq D_\alpha$ s.t. $P_0\in D_0$. Moreover, $f_p(\delta)\in\partial f_p(H_p)$ and $f^2_p(\delta)$ is at the interior of $\overline{f_p(H_p)}$ (see the illustration in Fig.\ref{firstpic}). In particular, the first-return map $f_p$ tears $H_p$ along $\delta$, in accordance with the two-sided limits above.\\

\begin{figure}[h]
\centering
\begin{overpic}[width=0.5\textwidth]{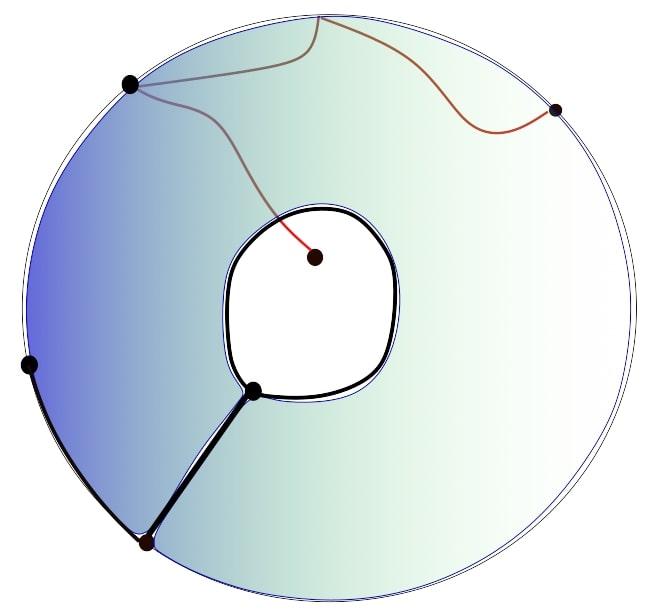}
\put(150,70){$P_{In}$}
\put(620,85){$L$}
\put(-50,250){$f_p(\delta)$}
\put(350,235){$f^2_p(\delta)$}
\put(280,710){$\delta$}
\put(140,850){$\delta_0$}
\put(390,500){$P_0$}
\put(450,450){$D_0$}
\put(850,810){$P_{Out}$}
\put(700,800){$T$}
\put(-110,380){$f_p(\delta_0)$}
\put(390,370){$f^2_p(\delta_0)$}
\end{overpic}
\caption[The image of the first return map.]{\textit{$f_p(H_p)$, imposed on the cross-section $D_\alpha$. The red curve denotes $f^{-1}_p(l)$ - for simplicity, it is drawn as a curve.}}
\label{firstpic}
\end{figure}
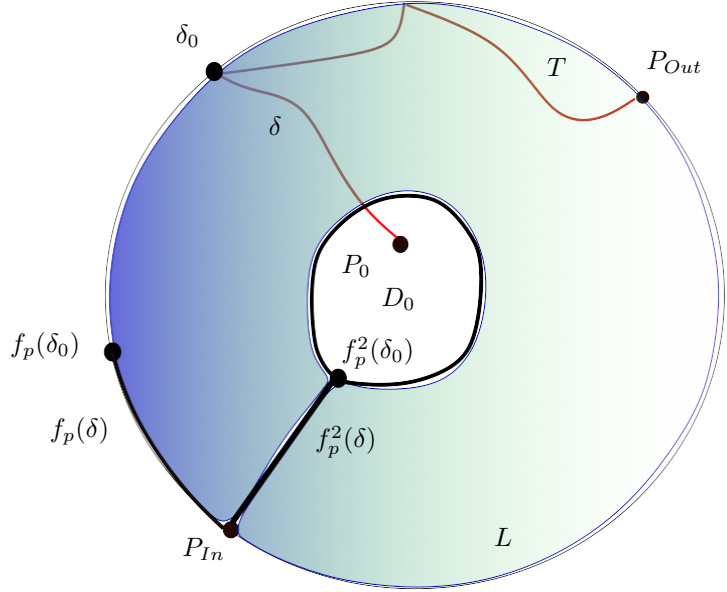
Having sketched the first-return map for $f_p$, we see we cannot immediately apply the same arguments used to prove Case $A$ - with the reason being that the fact $f_p(H_p)$ is glued to itself at $f^2_p(\delta_0)$ simply does not allow us to extend $f_p$ to a disc homeomorphism as we did in Stage $I$. In order to overcome this difficulty, we begin by expanding the fixed point $P_{In}$ to a directions sphere, with a fixed point $P'_{In}$ on its boundary. Moreover, we open up the curve $\delta$ to an open sector $\delta'$ by splitting $\delta\cap\partial H_p$ to two curves, $\delta_1$ and $\delta_2$, which enclose the said sector (see the illustration in Fig.\ref{DEF1}).\\

\begin{figure}[h]
\centering
\begin{overpic}[width=0.40\textwidth]{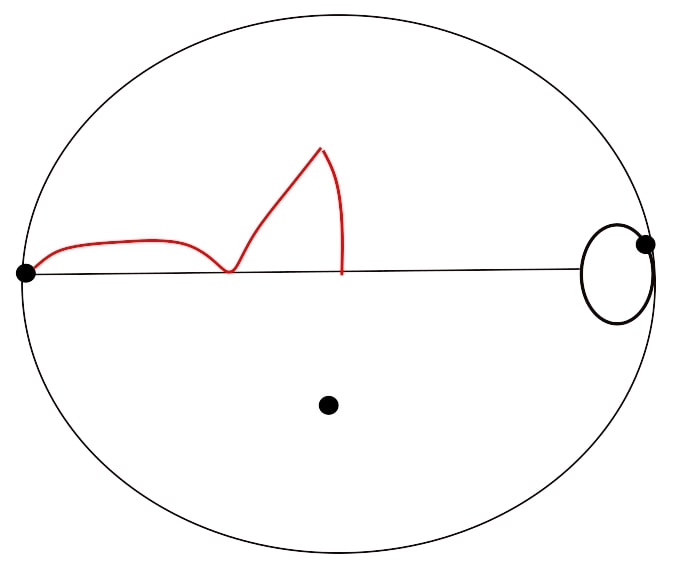}
\put(350,550){$\delta_1$}
\put(530,550){$\delta_2$}
\put(730,650){$H_p$}
\put(550,390){$l$}
\put(-130,420){$P_{Out}$}
\put(380,200){$P'_1$}
\put(870,760){$D_\alpha$}
\put(870,50){$C_\alpha$}
\put(970,450){$P'_{In}$}
\end{overpic}
\caption[Amending $H_p$.]{\textit{Amending $H_p$ - $P_{In}$ is opened to a sphere of direction with a unique fixed point on its boundary, $P'_{In}$, while $\delta$ is opened to the sector trapped between $\delta_1,\delta_2$ and $l$. Again, we glue $C_\alpha$ to $D_\alpha$ at $l$.}}
\label{DEF1}
\end{figure}
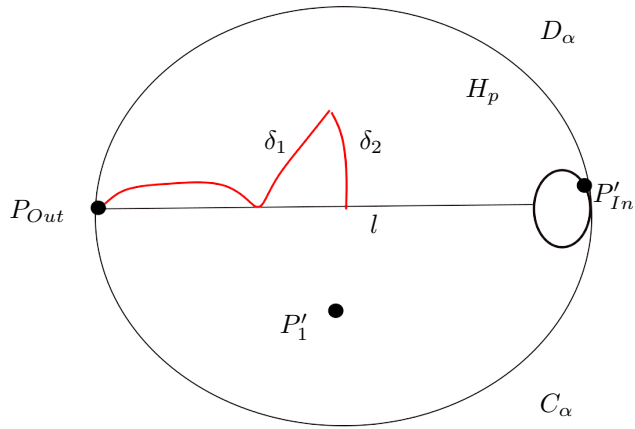

Now, let $f_0:\overline{H_p}\to\overline{D_\alpha}$ denote a continuous, injective function which is conjugate to $f_p$ away from the open sector $\delta'$, the set $\cup_{n\geq0}f^{-n}_0(\delta')$, and the directions sphere $P_{In}$ - in particular, we choose $f_0$ s.t. for $x\in H_p$, $x\to\delta_i$, we have $f_0(x)\to f_0(\delta_i)$, $i=1,2$. As such, we can sketch the graph of $f_0$ as in Fig.\ref{DEF2}. In other words, we obtain $f_0$ from $f_p$ by opening the components of $\cup_{n\geq0}f^{-n}_p(\delta)$ to open sets - and after this deformation, we are now able to extend $f_0$ to a homeomorphism of $S$, a twice-punctured disc, similarly to what we did in Case $A$.\\

To do so, choose some point, $P'_0$, inside the sector $\delta'$ which lies away from $f_0(\overline{H_p})\cap\delta'$, and make it a periodic point of period $2$ (see the illustrations in Fig.\ref{DEF2} and Fig.\ref{DEF1}) - i.e., similarly to Case $A$ we extend $f_0$ to $F$, some homeomorphism of the disc, s.t. the following is satisfied:
\begin{itemize}
    \item Similarly to Case $A$, we glue $C_\alpha$ to $D_\alpha$ at $l$, and generate a disc $V_\alpha=C_\alpha\cup D_\alpha$.
    \item  We choose the extension $F$ s.t. $F(P'_{In})=P'_{In}$, and $F(P'_0)=P'_1$, $F(P'_1)=P'_0$.
\end{itemize}

\begin{figure}[h]
\centering
\begin{overpic}[width=0.40\textwidth]{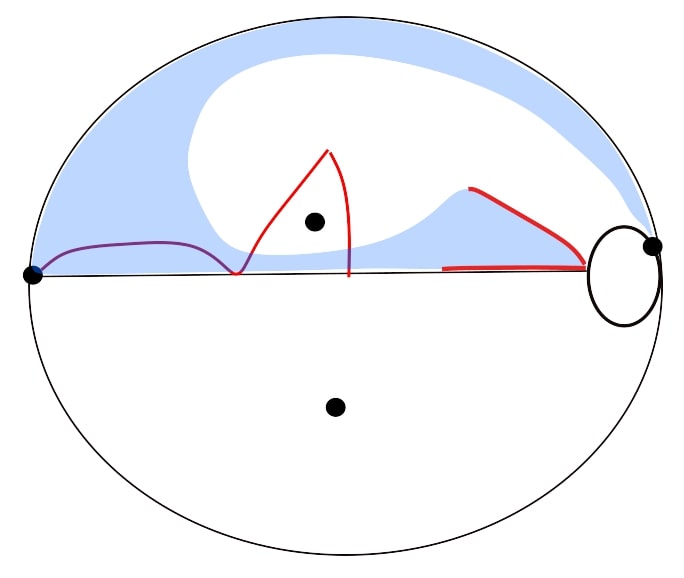}
\put(350,550){$\delta_1$}
\put(530,550){$\delta_2$}
\put(730,650){$H_p$}
\put(550,390){$l$}
\put(650,360){$f_0(\delta)$}
\put(-130,420){$P_{Out}$}
\put(380,200){$P_1$}
\put(870,760){$D_\alpha$}
\put(870,50){$C_\alpha$}
\put(970,450){$P'_{In}$}
\end{overpic}
\caption[$f_0(H_p)$ imposed on the amended cross-section.]{\textit{The new map $f_0$ - as can be seen, we can apply the logic of Case $A$ to analyze it (the blue region denotes $f_0(H_p)$). Moreover, $f_0$ is conjugate to $f_p$ away from $\delta_1,\delta_2$ and $P'_{In}$.}}
\label{DEF2}
\end{figure}
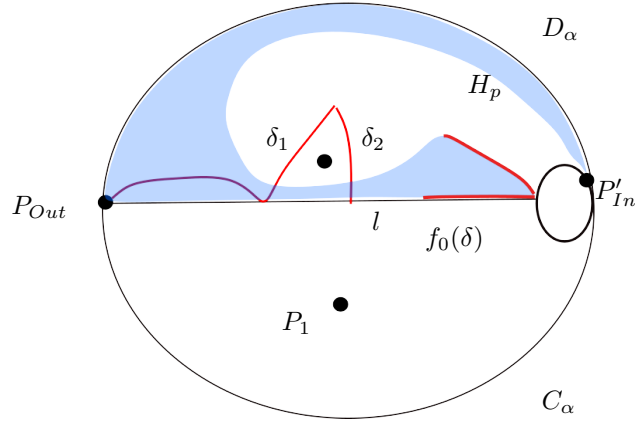

Consequentially, we reduced $f_0$ to the same situation as in Case $A$ - as such, using the same arguments (and in particular, Th.\ref{stability}), we immediately conclude:

\begin{corollary}
\label{periodic2}    With the notations above, $f_0$ has infinitely many periodic orbits in $\overline{H_p}$ - of all minimal periods. Additionally, the dynamics of $f_0$ on its invariant are isotopically deformed to that of a Smale Horseshoe map $G:ABCD\to V_\alpha$ (where $ABCD$ is some topological rectangle in $V_\alpha$) - in particular, we have the following:

    \begin{itemize}
        \item As we isotopically deform $F$ to $G$, $\overline{H_p}$ is isotoped to the rectangle $ABCD$.
        \item Every periodic orbit of minimal period $k$ for $G$ is isotopically deformed to a periodic orbit for $f_0$ in $\overline{H_p}$, of the same minimal period.
        \item If two periodic orbits $\{x_1,...,x_k\}$ and $\{y_1,...,y_j\}$ for $G$ in $ABCD$ are deformed to the same periodic orbit $\{z_1,...,z_t\}$ for $f_0$, then $k=j=t$ and $x_i=y_i$, $1\leq i\leq k$.

    \end{itemize}
\end{corollary}

We would now like to prove the analog of Lemmas \ref{sym1} and \ref{continu1} for the first-return map $f_p$ -  however, since $f_0$ is not exactly the same function as $f_p$ we need to be a little bit more careful. We begin with the following fact:

\begin{lemma}
 \label{continu2}   Let $F_p$ be a trefoil parameter falling into Case $B$. Then, $f_p$ generates infinitely many periodic orbits, of every minimal period. Moreover, $f_p$ is continuous at the said periodic orbits.
\end{lemma}
\begin{proof}
    We first claim that save perhaps for one periodic orbit, every periodic orbit for $G$ of minimal period $k$, is deformed to a periodic orbit $\{x_1,...,x_k\}$ for $f_0$ which lies away from the sector $\delta'$ and the direction sphere $P_{In}$. To do so, let $G:ABCD\to ABCD$ denote the said Smale Horseshoe map from Cor.\ref{periodic2}, and consider a periodic orbit for $G$ which is not deformed to $P'_{In}$, the fixed-point for $f_0$ - consequentially, by the the said orbit must lie away from the directions sphere $P_{In}$.\\

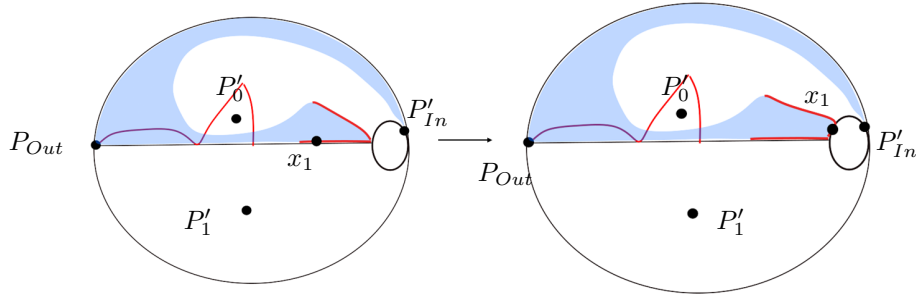
\begin{figure}[h]
\centering
\begin{overpic}[width=0.7\textwidth]{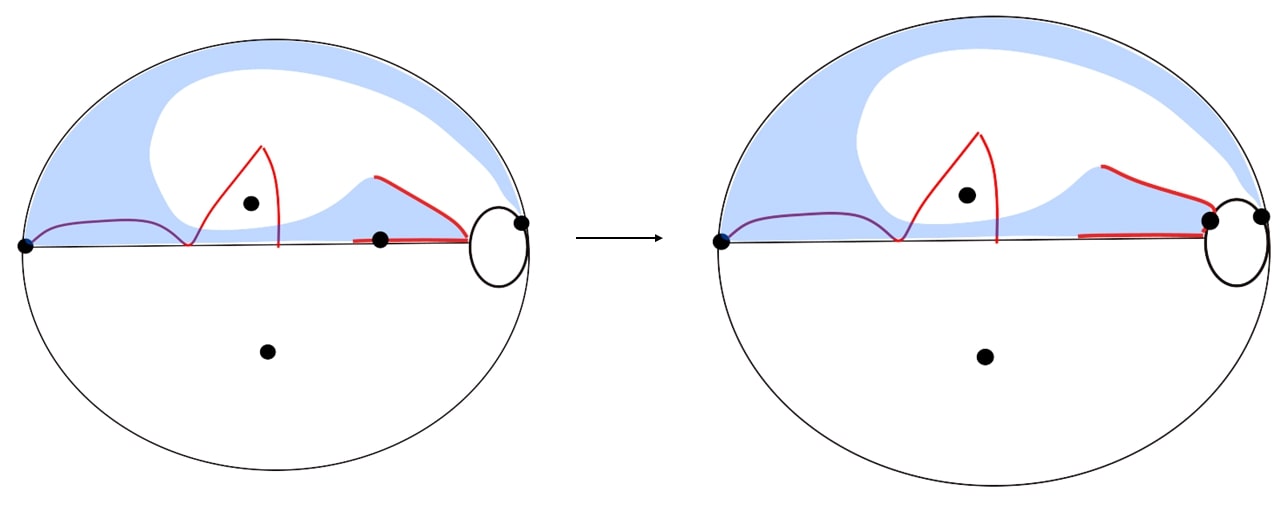}
\put(-90,200){$P_{Out}$}
\put(260,180){$x_1$}
\put(170,270){$P'_0$}
\put(130,100){$P'_1$}
\put(410,240){$P'_{In}$}
\put(500,160){$P_{Out}$}
\put(910,260){$x_1$}
\put(1000,200){$P'_{In}$}
\put(800,100){$P'_{1}$}
\put(730,270){$P'_0$}
\end{overpic}

\caption[Collapsing periodic orbits for $f_0$ on $l$ and $\delta$.]{\textit{The isotopy collapsing points on $f_0(\delta_1\cup\delta_2)$ like $x_1$ (or more generally, points $x_1\in l\cup f_0(\delta_2)$) to $P'_{In}$.}}
\label{DEFFF}
\end{figure}

    Now, let us remark that we can choose the isotopy from $f_0$ to $G$ s.t. any point on $\partial\delta'=\delta_1\cup\delta_2$ is eventually collapsed to an arc on the directions sphere $P_{In}$ along the way (see the illustration in Fig.\ref{DEFF}) - and since $P'_{In}$ is the only periodic orbit on the direction sphere $P_{In}$, it follows every periodic orbit on $\partial\delta'$ for $f_0$ can be destroyed or collapsed by the isotopy. Consequentially, as the periodic orbit $\{x_1,...x_n\}$ is both unremovable and uncollapsible by Th.\ref{stability}, it must lie away from the curves $\delta_1,\delta_2$ and their pre-images in $\overline{H_p}$ -  or, in other words, it lies away from the pre-images of the sector $\overline{\delta'}$ in $\overline{H_p}$.\\
    
    All in all, we conclude that save for one periodic orbit for $G$ which is deformed to $P'_{In}$, every other periodic orbit for $G$ is deformed to a periodic orbit for $f_0$, which lies away from $P'_{In}$ and the pre-images of $\overline{\delta}$ in $H_p$. Consequentially, since $f_p$ and $f_0$ are conjugate away from $\delta'$, its pre-images, and the directions sphere $P_{In}$ we conclude $\{x_1,...x_k\}$ is also a periodic orbit for $f_p$ which lies away from the curve $\delta$ and its pre-images under $f_p$ - as such, by Prop.\ref{lem33}, $f_p$ is continuous at the said periodic orbit. Finally, since this is true for every periodic orbit for the Smale Horseshoe map $G$ (save perhaps for one), we conclude $f_p$ generates infinitely many periodic orbits in $H_p$, of all minimal periods - and moreover, $f_p$ is continuous at every such periodic orbit.
\end{proof}

\begin{figure}[h]
\centering
\begin{overpic}[width=0.8\textwidth]{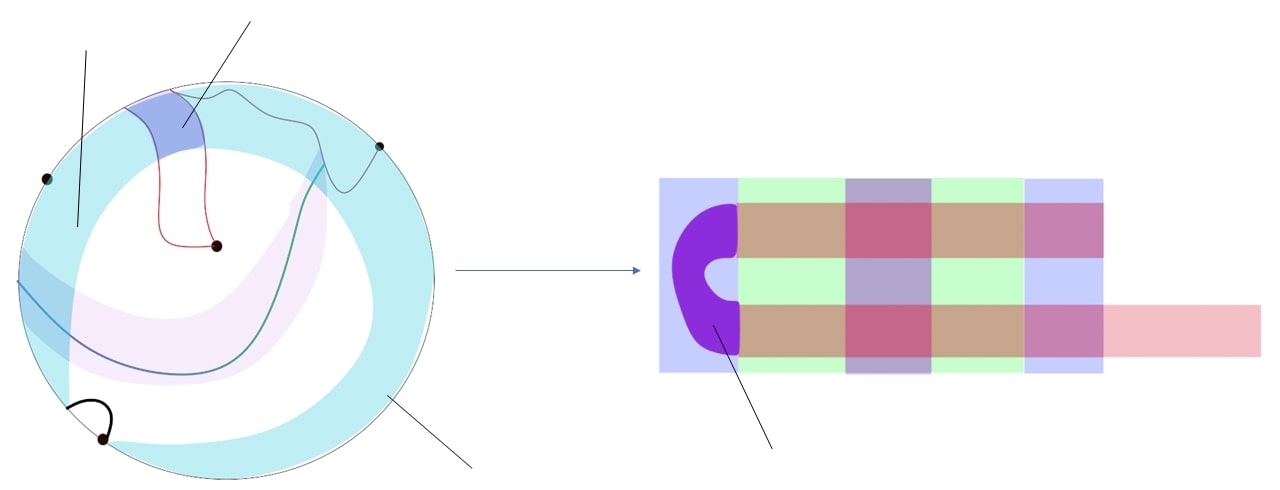}
\put(10,30){$P'_{In}$}
\put(170,200){$P'_{0}$}
\put(200,350){$f_0(\Delta)$}
\put(350,0){$f_0(D_1)$}
\put(170,250){$\delta$}
\put(100,100){$\Delta$}
\put(70,350){$f_0(D_2)$}
\put(300,270){$P_{Out}$}
\put(600,20){$G(L_3)$}
\put(550,50){$C$}
\put(550,250){$D$}
\put(670,270){$L_3$}
\put(800,250){$B$}
\put(800,50){$A$}
\put(730,270){$L_1$}
\put(600,270){$L_2$}
\end{overpic}
\caption[Deforming $\Delta$ to $L_3$.]{\textit{Similarly to Case $A$, the isotopy of $F|_{\overline{H_p}}=f_0$ on $H_p$ to $G:ABCD\to S$ deforms the purple region $\Delta$ to the sub-rectangle $L_3$.}}
\label{DEFF}
\end{figure}
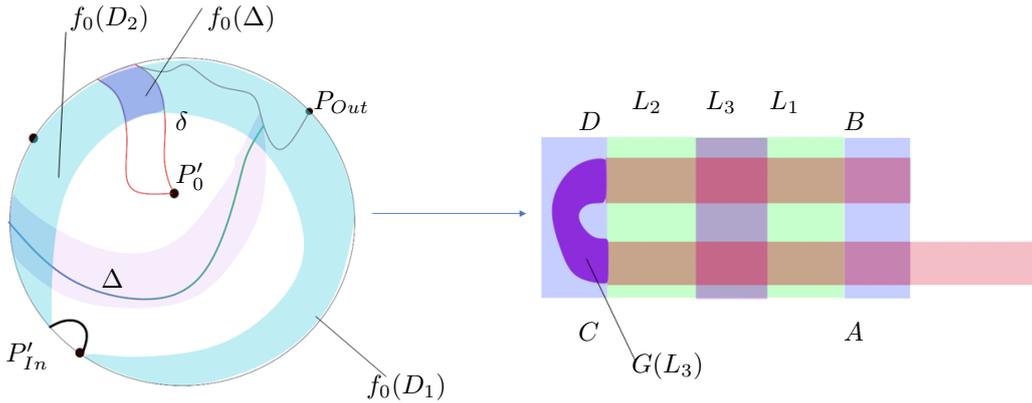

Having proven Lemma \ref{continu2}, we can almost conclude the proof in Case $B$. To do so, we must first prove we can define symbolic dynamics for $f_p$ - that, however, is immediate. By considering the curve $\delta$ and by recalling that when we deform $f_p$ to the Smale Horseshoe map $G$ we open $\delta$ to a sector and then continuously deform it to the $CD$ side (see the illustration in Fig.\ref{DEFF}), using a similar argument to the one used to prove Lemma \ref{sym1}, we conclude:

\begin{lemma}
    \label{sym2} Let $F_p$ be a trefoil parameter falling into Case $B$. Then, there exists a curve $\Delta\subseteq f^{-1}_p(\delta)$ s.t. $\overline{H_p}\setminus\Delta$ includes two components - $D_1$ and $D_2$, s.t. $P_{In}\in\partial D_1\setminus\Delta$ and $\delta\subseteq\partial D_2\setminus\Delta$ (see the illustration in Fig.\ref{DEF4}). Consequentially, we can define symbolic dynamics on the invariant set of $f_p$ in $\overline{H_p}\setminus f^{-1}_p(l)$, denoted by $I'$ - i.e., there exists a continuous $\pi':J\to\{1,2\}^\mathbf{N}$ s.t. $\pi'\circ f_p=\sigma\circ \pi'$.\\
\end{lemma}
\begin{proof}
By considering Fig.\ref{firstpic}, it immediately follows $\Delta$ exists and divides $\overline{H_p}$ to $D_1$ and $D_2$, as in Fig.\ref{DEF4}. Now, it only remains to define symbolic dynamics on $I'$. Let us note that the set $I'$, by definition, is a subset of $\overline{D_\alpha}\setminus(\cup_{n\geq1}f^{-n}_p(l))$, and it is also invariant under $f_p$ - as such, since $\Delta\subseteq f^{-1}_p(\delta)\subseteq f^{-2}_p(l)$, any component $C$ of $I'$ either lies in $D_1$ or $D_2$, and the same is true for $f_p(C)$. Furthermore, by Prop.\ref{lem33}, $f_p$ and all its iterates are continuous at initial conditions in $I'$ - as such, given $x\in I'$, there exists a function $\pi':I\to\{1,2\}^\mathbf{N}$ s.t. $\pi'(x)=\{i_0,i_1,i_2,...\}$ where $i_k=1$ when $f^k_p(x)\in D_1$ and $2$ otherwise. It immediately follows that as $f^k_p$ is continuous at $x$ for every $k$, so is $\pi'$ - and by definition, we also have $\pi'\circ f_p=\sigma\circ\pi'$, where $\sigma:\{1,2\}^\mathbf{N}\to\{1,2\}^\mathbf{N}$ is the one-sided shift. The proof of Lemma \ref{sym2} is now complete.
\end{proof}
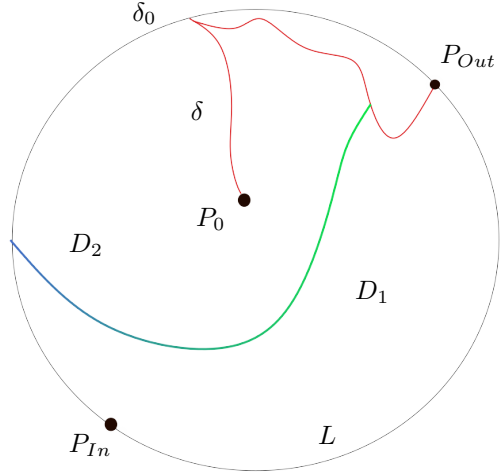
\begin{figure}[h]
\centering
\begin{overpic}[width=0.4\textwidth]{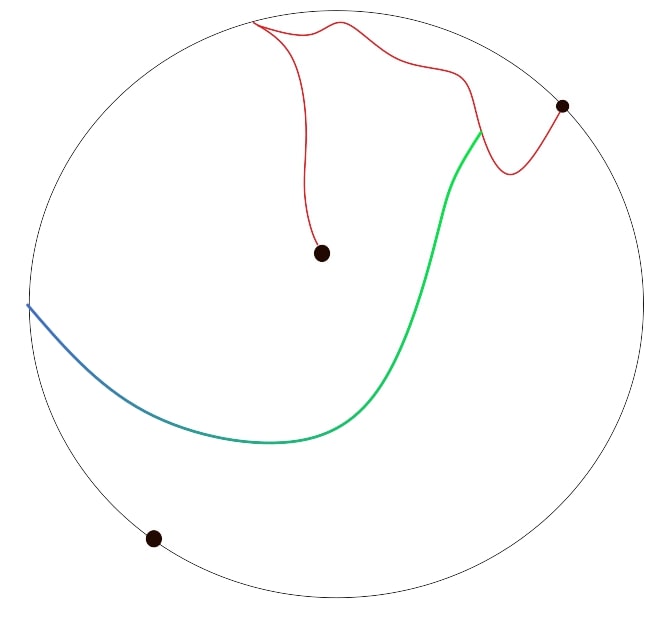}
\put(150,70){$P_{In}$}
\put(620,85){$L$}
\put(380,700){$\delta$}
\put(270,890){$\delta_0$}
\put(390,500){$P_0$}
\put(150,450){$D_2$}
\put(850,810){$P_{Out}$}
\put(690,360){$D_1$}
\end{overpic}
\caption[Symbolic dynamics in $H_p$ for Case $B$.]{\textit{The curve $\Delta$, partitioning $\overline{H_p}$ to $D_1$ and $D_2$.}}
\label{DEF4}
\end{figure}

Now, having proven Cor.\ref{sym2}, we now use it to conclude the following fact:

\begin{corollary}
    \label{conin} The periodic points given by Lemma \ref{continu2} all lie in the set $I'$,given by Cor.\ref{sym2} - and consequentially, we can define symbolic dynamics on them. Moreover, $\pi'(I)$ includes every periodic symbol in $\{1,2\}^\mathbf{N}$ that is not the constant $\{1,1,1...\}$.
\end{corollary}
\begin{proof}
Let $\{x_1,...,x_k\}$ be a periodic orbit given by Lemma \ref{continu2} - by the proof of Lemma \ref{continu2}, we already know $f_p$ is continuous on $\{x_1,...,x_k\}$, and that it lies away from $\delta$. Hence, we have $\{x_1,...,x_k\}\subseteq \overline{H_p}\setminus(\cup_{n\geq0}f^{-n}_p(\delta))$ and by $\Delta\subseteq f^{-1}_p(\delta)$ we conclude $\{x_1,...,x_k\}$ also lies in $I'$. Similarly to the argument in Lemma \ref{continu1} for Case $A$, again the arc $\Delta$ is continuously deformed into the sub-rectangle $L_3$ of $ABCD$, defined by $G(L_3)\cap ABCD=\emptyset$ (see the illustration in Fig.\ref{DEFF}). As such, similarly to the proof of Lemma \ref{continu1}, it follows that when we deform $\{x_1,...,x_k\}$ to $\{y_1,...y_k\}$, a periodic orbit for the Smale Horseshoe map $G$, the symbol corresponding to $\{x_1,...,x_k\}$ does not change. As such, since every periodic orbit for $G$ (save the one corresponding to the constant $\{1,1,1...\}$ - i.e., $P'_{In}$) can be continuously deformed to a periodic orbit for $f_p$, using a similar argument to the one used to prove Lemma \ref{continu1} Cor.\ref{conin} now follows.
\end{proof}
We now conclude the proof of Th.\ref{th31} for Case $B$. To do so, again, for every periodic $s\in\{1,2\}^\mathbf{N}$ of minimal period $k$ which is not the constant $\{1,1,1...\}$, set $D_s$ as the component of $\pi'^{-1}(s)$ in $I'$ s.t.:

\begin{itemize}
    \item $D_s$ contains a periodic orbit of minimal period $k$, $\{x_1,...,x_k\}$ - by Cor.\ref{conin}, it exists.
    \item $\{x_1,...,x_k\}$ can be deformed isotopically to a periodic orbit of $G$, the Smale Horseshoe map.
\end{itemize}

Now, similarly to Case $A$, when $s=\{1,1,1...\}$ set $D_s=\{P_{In},P_{Out}\}\cup L$ (where $L=\partial B_\alpha\cap U_p$ - see Def.\ref{def32}) - and let $Q$ denote the collection of such $D_s$, for periodic $s$. Setting $\pi=\pi'|_Q$ and summarizing our results, we conclude:

\begin{itemize}
    \item The first-return map $f_p$ is continuous on $Q$.
    \item Also, there exists a continuous $\pi:Q\to\{1,2\}^\mathbf{N}$ s.t. $\pi\circ f_p=\sigma\circ\pi$ (where $\sigma:\{1,2\}^\mathbf{N}\to\{1,2\}^\mathbf{N}$ denotes the one-sided shift).
    \item $\pi(Q)$ includes every periodic symbol in $\{1,2\}^\mathbf{N}$.
    \item If $s$ is periodic of minimal period $k$, $\pi^{-1}(s)$ includes a periodic orbit of minimal period $k$ for $f_p$.
\end{itemize}

All in all, we have completed the proof of Th.\ref{th31} for Case $B$.

\subsubsection{Stage $III$ - proving Th.\ref{th31} for Case $C$:}
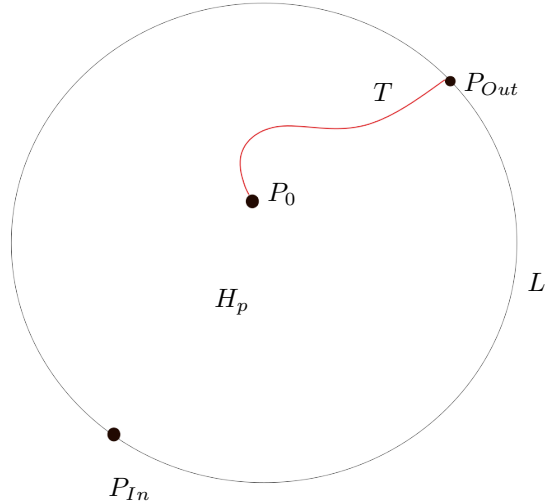
\begin{figure}[h]
\centering
\begin{overpic}[width=0.40\textwidth]{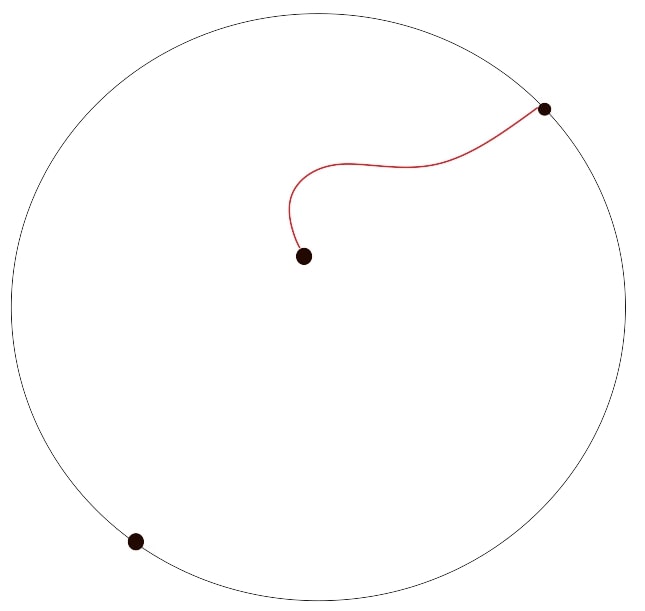}
\put(200,-10){$P_{In}$}
\put(500,550){$P_{0}$}
\put(400,350){$H_p$}
\put(870,760){$P_{Out}$}
\put(700,740){$T$}
\put(990,380){$L$}
\end{overpic}
\caption[Case $C$.]{\textit{Case $C$ - $\delta_0= P_{Out}$, and $\delta$ is the red curve connecting $P_0$ and $P_{Out}$.}}
\label{optc}
\end{figure}
Having proven Th.\
\ref{th31} holds for Cases $A$ and $B$, we now prove the same for Case $C$. In this case we have $\delta_0=P_{Out}$, as illustrated in Fig.\ref{optc}. Consequentially, in this scenario, by $f_p(P_{Out})=P_{Out}$ we have $\delta=f^{-1}_p(l)$, and as such $H_p$ is homeomorphic to a slit disc - see the illustration in Fig.\ref{optc}. However, this implies $f_p(l)=f^2_p(\delta)$ is a curve in ${D_\alpha}$, beginning at $P_{In}$ and terminating at $P_{Out}$ - see the illustration in Fig.\ref{first3}.\\

Since $f_p(\delta)=l$ for every $s_0\in \delta$ the limit $\lim_{s\to s_0}f_p(s)$ is either $f_p(s_0)$ or $f^2_p(s_0)$ - depending on the direction from which $s$ tends to $s_0$. This remains true even when $s_0=P_{Out}$ - it is therefore easy to see that when $s\in H_p$ tends to $P_{Out}$ from above $\delta$, $f_p(s)\to P_{Out}$, and when it does so from below $\delta$ (i.e., from the sector $T$ given by Cor.\ref{disconcurve} - see Fig.\ref{optc}), we'd have $f_p(s)\to P_0$.\\

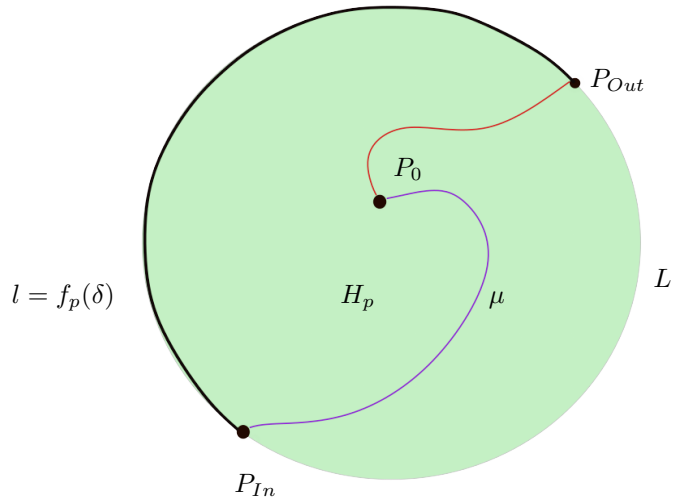
\begin{figure}[h]
\centering
\begin{overpic}[width=0.40\textwidth]{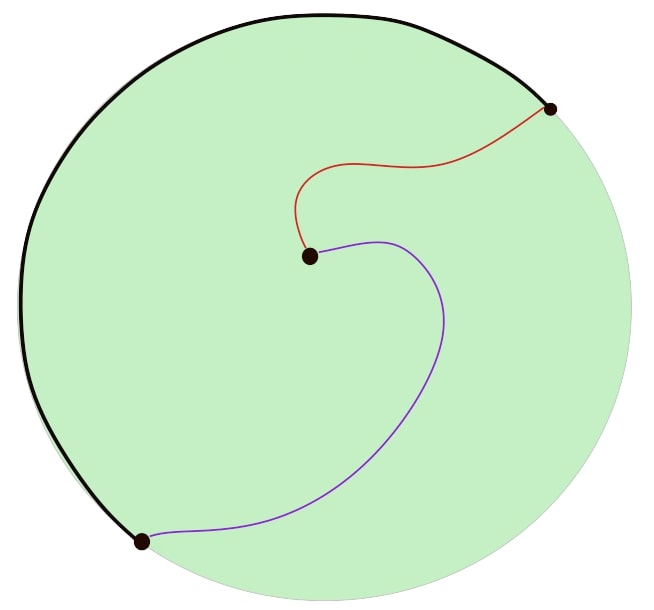}
\put(200,-10){$P_{In}$}
\put(500,590){$P_{0}$}
\put(400,350){$H_p$}
\put(-220,350){$l=f_p(\delta)$}
\put(680,350){$\mu$}
\put(870,760){$P_{Out}$}
\put(990,380){$L$}
\end{overpic}
\caption[The first-return map in Case $C$.]{\textit{The first-return map in Case $C$ - the curve $\mu$ denotes $f^2_p(\delta)=f_p(l)$. $\mu$ begins at $P_{In}$ and terminates at $P_0$.}}
\label{first3}
\end{figure}
By this entire discussion we see that Case $C$ is very similar to Case $B$ - only that in Case $B$ the point $P_0$ was away from $f_p(\overline{H_p})$. In this scenario this is not the case - as such, making $P_0$ a periodic point w.r.t. to any extension of $f_p$ to a homeomorphism of a disc would have to be done more carefully.\\

To do so, we deform the cross-section $H_p$ as follows - first, note that by Fig.\ref{first3} and Lemma \ref{firstret} (and by Remark \ref{extend}), every $s\in\delta$ which is not $P_0$ has a unique $n$-th pre-image in $\overline{D_\alpha}$. Now, begin deforming $H_p$ by opening $P_{Out}$ to a direction sphere  (see the illustration in Fig.\ref{DEF5}). Then, open up the curve $\delta$ into an open sector, as in Fig.\ref{DEF5} - and do the same simultaneously for every component $\cup_{n\geq0}f^{-n}_p(\delta)$ in $D_\alpha$. Much like Case $B$, after this deformation $\delta$ becomes a sector $\delta'$, whose boundary includes two curves, terminating at $P_0$ - $\delta_1$ and $\delta_2$, as illustrated in Fig.\ref{DEF5}.\\

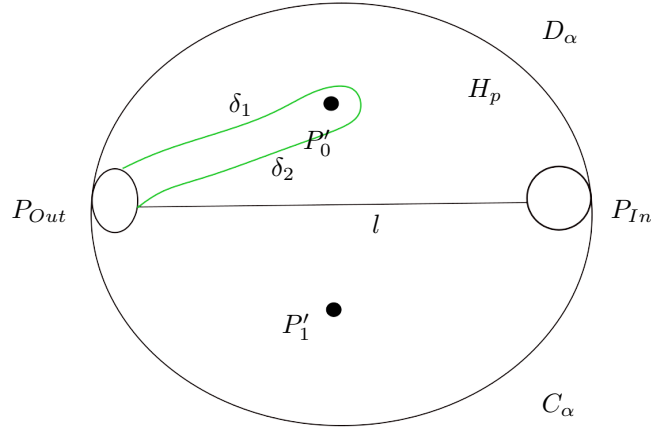
\begin{figure}[h]
\centering
\begin{overpic}[width=0.40\textwidth]{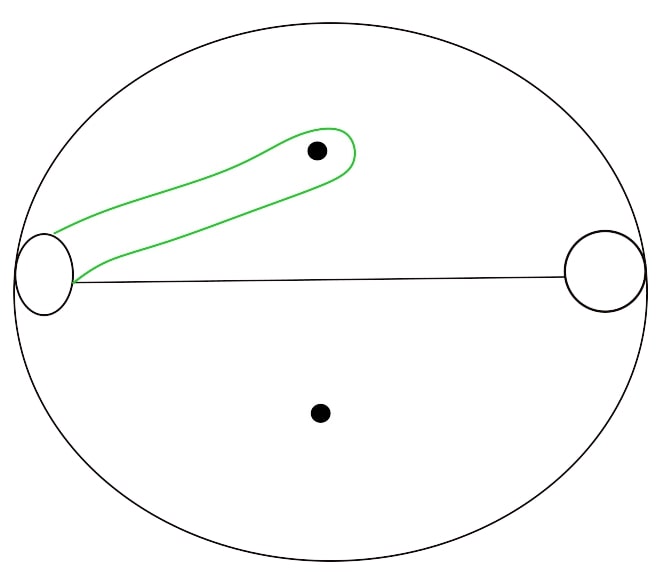}
\put(420,550){$P'_{0}$}
\put(280,620){$\delta_1$}
\put(360,500){$\delta_2$}
\put(730,650){$H_p$}
\put(550,390){$l$}
\put(-130,420){$P_{Out}$}
\put(380,200){$P'_1$}
\put(870,760){$D_\alpha$}
\put(870,50){$C_\alpha$}
\put(1000,420){$P_{In}$}
\end{overpic}
\caption[Deforming Case $C$.]{\textit{$H_p$ after the deformation - $\delta$ is blown to a sector, while $P_{In}$ and $P_{Out}$ are both opened to direction spheres. Similarly to cases $A$ and $B$, we glue $C_\alpha$ to $D_\alpha$ at $l$.}}
\label{DEF5}
\end{figure}

Now, open up $P_{In}$ to a directions sphere with a unique fixed-point on it, $P'_{In}$, and let $f_0$ denote some homeomorphism $f_0:\overline{H_p}\to\overline{D_\alpha}$ s.t. $f_0$ is conjugate to $f_p$ on $\overline{H_p}\setminus\cup_{n\geq0}f^{-n}_0(\delta')$, and behaves as in Fig.\ref{DEF6}. In particular, when $x\to\delta_i$, we have $f_0(x)\to f_0(\delta_i)$, $i=1,2$. Consequentially, $f_0(\delta)$ is a sector whose boundary intersects the directions sphere $P_{In}$ and includes $P_0$ on its vertex (see the illustration in Fig.\ref{DEF5}) - as such, because $f^2_p(\delta)$ is an arc beginning at $P_{In}$ and terminating at $P_0$, when we deform $f_p$ to $f_0$ we can arrange for the intersection $f_0(\delta')\cap\delta'$ to include an open disc $D$ as in Fig.\ref{DEF6}.\\

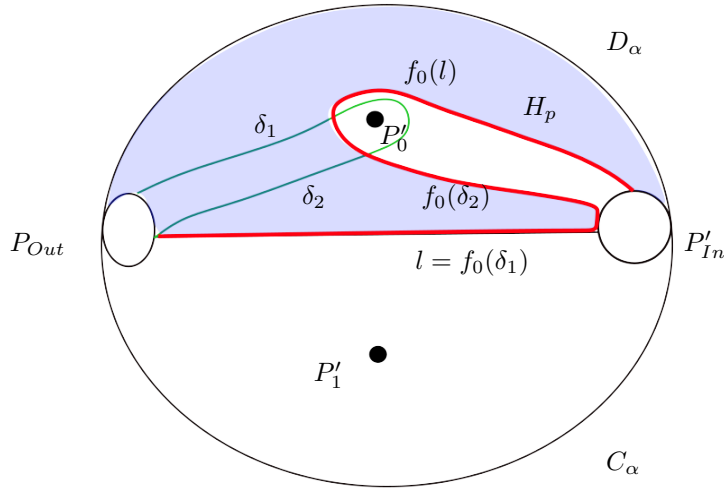
\begin{figure}[h]
\centering
\begin{overpic}[width=0.45\textwidth]{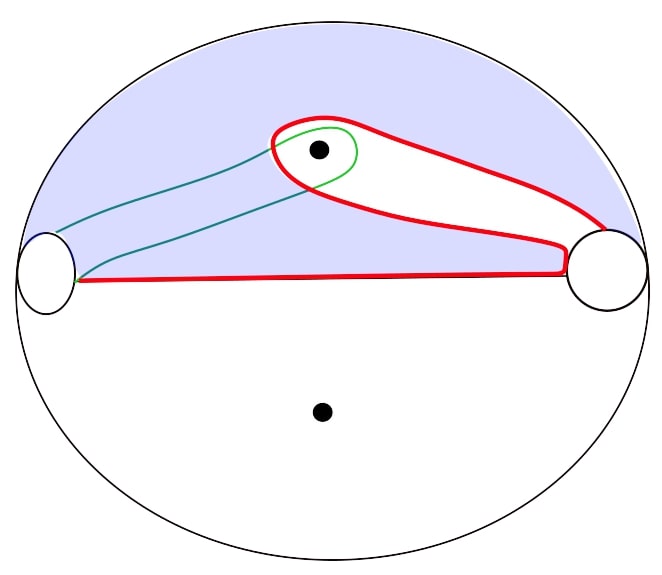}
\put(490,600){$P'_{0}$}
\put(280,620){$\delta_1$}
\put(360,500){$\delta_2$}
\put(730,650){$H_p$}
\put(530,710){$f_0(l)$}
\put(550,390){$l=f_0(\delta_1)$}
\put(560,500){$f_0(\delta_2)$}
\put(-130,420){$P_{Out}$}
\put(380,200){$P'_1$}
\put(870,760){$D_\alpha$}
\put(870,50){$C_\alpha$}
\put(1000,420){$P'_{In}$}
\end{overpic}
\caption[$f_0$ in Case $C$.]{\textit{$f_0(H_p)$ imposed on $H_p$ - the red curve denotes $f_0(l\cup\partial\delta$). $D$ is the region trapped between $f_0(\partial\delta)$ and $\partial \delta$.}}
\label{DEF6}
\end{figure}

Similarly to Cases $A$ and $B$, glue some half-disc $C_\alpha$ to $D_\alpha$ at $l$ to create a disc $V_\alpha$ as illustrated in Fig.\ref{DEF5} and Fig.\ref{DEF6}, and choose some $P'_0\in f_0(\delta)\cap\delta$ and another $P_1$, interior to $C_\alpha$. Finally, extend $f_0$ to a disc homeomorphism $F:V_\alpha\to V_\alpha$  s.t. $F$ coincides with $f_0$ on $\overline{H_p}$ and $F(P_{In})=P_{In}$, $F(P'_0)=P_1$ and $F(P_1)=P'_0$. Now, consider a graph $T$ inside $H_p$ connecting $P'_1,P_0$ and $P_{In}$ as in Fig.\ref{DEF7}, and then consider its image under $f_0$ (see Fig.\ref{DEF7}) - again, applying the same procedures used to prove Cases $A$ and $B$ we see the resulting graph map is essentially the same as the one obtained for Cases $A$ and $B$.\\

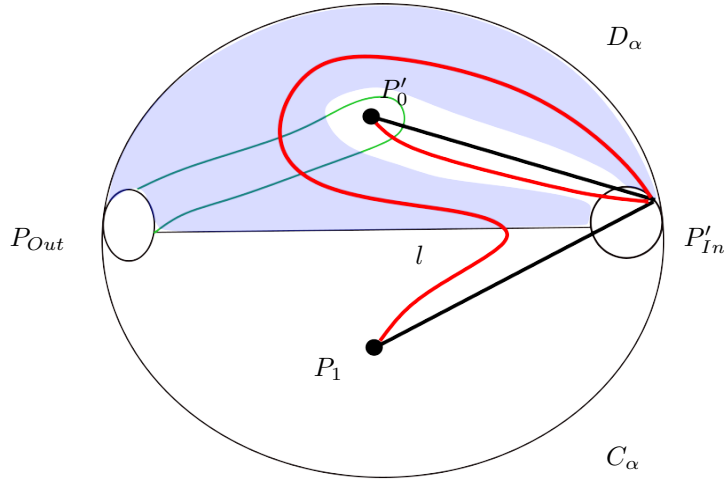
\begin{figure}[h]
\centering
\begin{overpic}[width=0.45\textwidth]{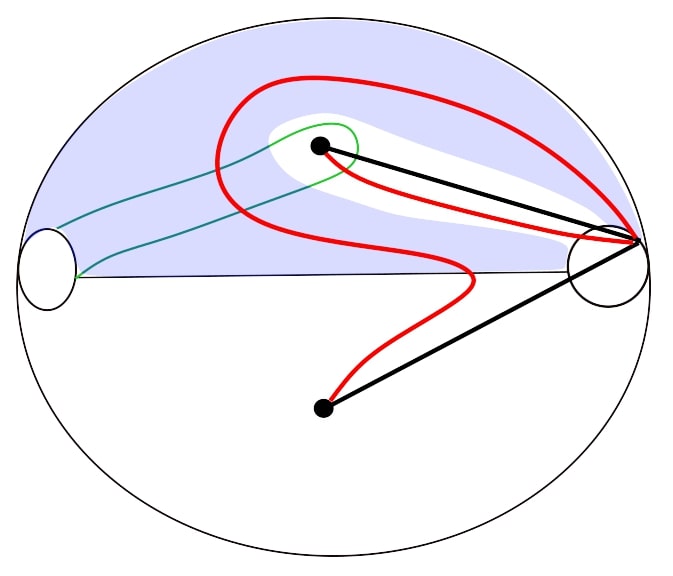}
\put(490,670){$P'_{0}$}
\put(550,390){$l$}
\put(-130,420){$P_{Out}$}
\put(380,200){$P_1$}
\put(870,760){$D_\alpha$}
\put(870,50){$C_\alpha$}
\put(1000,420){$P'_{In}$}
\end{overpic}
\caption[The graph map in Case $C$.]{\textit{The image of the black graph $T$ under $F$ is the red curve - the image of one edge connects $P'_{In}$ and $P'_0$, while the image of the second begins at $P'_{In}$, loops around $P_0$ and then enters $C_\alpha$ and connects with $P_1$.}}
\label{DEF7}
\end{figure}

As such, we have shown Case $C$ can be reduced to Case $B$ - now, recalling Cor.\ref{periodic2},  Lemma \ref{continu2}, Cor.\ref{sym2} and and their corollaries did not depend on the location of $\delta_0$ in the curve $\overline{l}$, we see their conclusions extend to Case $C$ as well. All in all, using similar arguments to those used in $B$, Th.\ref{th31} similarly follows for Case $C$. The proof of Th.\ref{th31} is now complete.
\end{proof}

At this point, we remark the proof of Th.\ref{th31} can probably be generalized to other heteroclinic parameters for the Rössler system, which generate knots more complex than a trefoil (see Def.\ref{def31}) - however, much like Th.\ref{th31}, the proof of any such generalization would greatly depend on the topology of the heteroclinic knot involved.\\

\begin{figure}[h]
\centering
\begin{overpic}[width=0.4\textwidth]{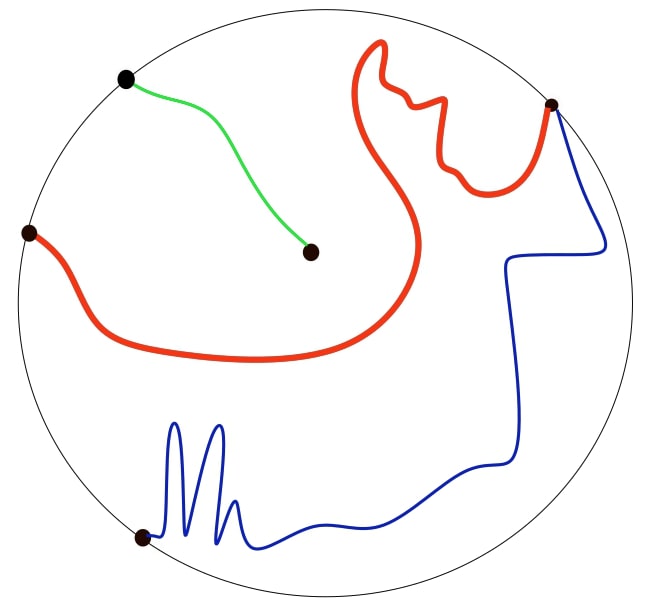}
\put(120,70){$P_{In}$}
\put(550,35){$L$}
\put(450,285){$U_{1}$}
\put(700,720){$\rho$}
\put(280,710){$\delta$}
\put(140,870){$\delta_0$}
\put(0,680){$l$}
\put(390,535){$P_0$}
\put(850,810){$P_{Out}$}
\put(270,500){$U_{2}$}
\end{overpic}
\caption[Splitting $U_p$ to $U_1$ and $U_2$.]{\textit{The cross-section $U_p$, sketched as a disc instead of a half-plane. The curve $\rho$ bisects $U_p$ to $U_1,U_2$, with $L$ denoting the curve $\partial B_\alpha\cap U_p$ (see Def.\ref{def32}).}}
\label{DALPHA23}
\end{figure}
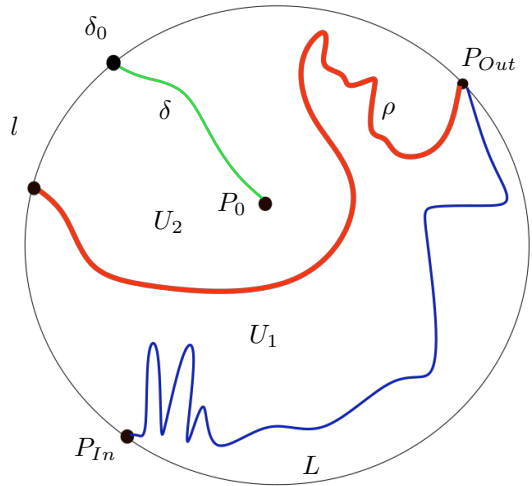

Let us also remark Th.\ref{th31} can be reformulated in a somewhat simpler way (the proof, of course, remains the same). To begin, recall the cross-section $U_p$ (see the discussion before Lemma \ref{obs}), the curve $\Delta$ and the regions $D_1,D_2$ from the proof of Th.\ref{th31} (as defined per Case $A,B$ or $C$). Let $\rho$ denote some curve in $\overline{D_\alpha}$ s.t. the following is satisfied (see the illustration in Fig.\ref{DALPHA23}):
\begin{itemize}
    \item $\Delta\subseteq\rho$.
    \item $\rho$ is homeomorphic to a closed interval - in particular, the endpoints of $\rho$ are in $\overline{l}$.
    \item $\overline{U_p}\setminus\rho$ is composed of two components, $U_1,U_2$ - indexed by $D_i\subseteq U_i,i=1,2$.
\end{itemize}

Now, set $I$ as the maximal invariant set of $f_p$ in $\overline{U_p}\setminus \rho$ - that is, $I$ is the collection of initial conditions in the cross-section $U_p$ whose trajectories both never hit $\rho$ and never escape to $\infty$. It is easy to see the set $Q$ constructed in the proof of Th.\ref{th31} is a subset of $I$. We now restate Th.\ref{th31} as follows:

\begin{corollary}\label{corsym}
 Let $F_p$ be a trefoil parameter for the Rössler system and let  $f_p:\overline{U_p}\to \overline{U_p}$ denote the first-return map. Let $I$ be as above and denote by $\sigma:\{1,2\}^\mathbf{N}\to\{1,2\}^\mathbf{N}$ the one-sided shift - then, there exists an $f_p$-invariant $Q\subseteq I$, s.t. the following holds:
 \begin{itemize}
     \item $f_p$ is continuous on $Q$.
     \item There exists a continuous $\pi:Q\to\{1,2\}^\mathbf{N}$ s.t. $\pi\circ f_p=\sigma\circ\pi$.
     \item $\pi(Q)$ includes every periodic $s\in\{1,2\}^\mathbf{N}$.
     \item If $s\in\{1,2\}^\mathbf{N}$ is periodic of minimal period $k>0$, $\pi^{-1}(s)$ includes at least one periodic point for $f_p$ of minimal period $k$.
\end{itemize}

\end{corollary}

A second remark is that the periodic orbits in $Q$ given by Th.\ref{th31} are isotopy stable and uncollapsible, in a sense which is somewhat different than that of Th.\ref{stability}. To state it, let us recall the curve $\Delta$ from the proof of Th.\ref{th31} and the curve $\delta$ from Prop.\ref{lem33}. Now, let $g:\overline{H_p}\setminus\delta\to\overline{D_\alpha}$ be a homeomorphism isotopic to $f_p:\overline{H_p}\setminus\delta\to \overline{D_\alpha}$, s.t. the following is satisfied (see the illustration in Fig.\ref{home}):

\begin{itemize}
\item $g(P_{In})=P_{In}$, $g(P_{Out})=P_{Out}$ and $g(L)=L$.
    \item  Whenever $P_0\in\partial H_p$, we have  $\lim_{s\to P_0}g(s)=P_{In}$. 
    \item $P_0\not\in\overline g(\overline{H_p}\setminus\delta)$.
\end{itemize}

Applying a similar logic used to prove Th.\ref{th31}, we have the following result:

\begin{corollary}
\label{deformation11}    Let $F_p$ be a trefoil parameter, and let $g:\overline{H_p}\setminus\delta\to\overline{D_\alpha}$ be as above. Then, if $s\in\{1,2\}^\mathbf{N}$ is periodic of minimal period $k$, we have the following:

    \begin{itemize}
        \item There exists some periodic $x\in D_s$ of minimal period $k$ which is continuously deformed to $y$, a periodic orbit of minimal period $k$ for $g$. In particular, $y$ is in the invariant set of $g$ in $\overline{H_p}\setminus\delta$.
        \item Let $s,\omega\in\{1,2\}^\mathbf{N}$ be periodic orbits of minimal periods $k_1$ and $k_2$, and let $x_1\in D_s$ and $x_2\in D_\omega$ be periodic for $f_p$ of minimal periods $k_1$ and $k_2$. If the isotopy deforms $x_i$ to $y_i$, $i=1,2$, then $y_1=y_2$ if and only if $s=\omega$ - that is, the periodic dynamics are also uncollapsible. 
    \end{itemize}
\end{corollary}
\begin{proof}
Recall the disc $S=V_\alpha\setminus\{P_0,P_1\}$ introduced in the proof of Th.\ref{th31}, and let us recall that the set $Q$ given by Th.\ref{th31} lies in the invariant set of $f_p$ in $\overline{H_p}\setminus\Delta$. Since $g$ is isotopic to $f_p$ as described above, similarly to how we extend $f_p$ to $F:S\to S$ we can extend $g$ to $G':S\to S$, and moreover, we can do so s.t. $F$ and $G'$ are isotopic on $S$ (see the illustration in Fig.\ref{home}). Since by Th.\ref{stability} the periodic orbits in $Q$ are unremovable and uncollapsible it follows that as we isotope $F$ to $G'$, the periodic orbits for $f_p=F|_{\overline{H_p}\setminus\Delta}$ are continuously deformed to the periodic orbits of $g=G'|_{\overline{H_p}\setminus\Delta}$ without changing their minimal periods or collapse into one another. Cor.\ref{deformation11} now follows.
\end{proof}

\begin{remark}
   The function $g$ from Cor.\ref{deformation11} need not necessarily be a first-return map for a flow. That is, it need only be a homeomorphism of $\overline{H_p}\setminus\delta$.
\end{remark}
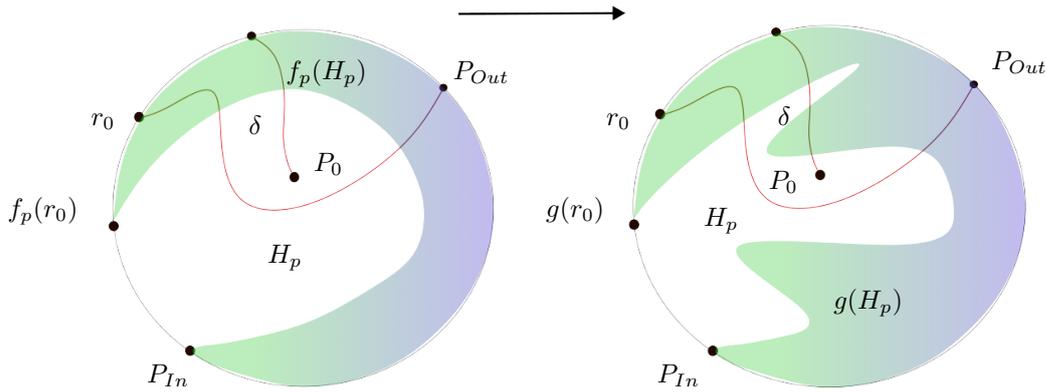
\begin{figure}[h]
\centering
\begin{overpic}[width=0.7\textwidth]{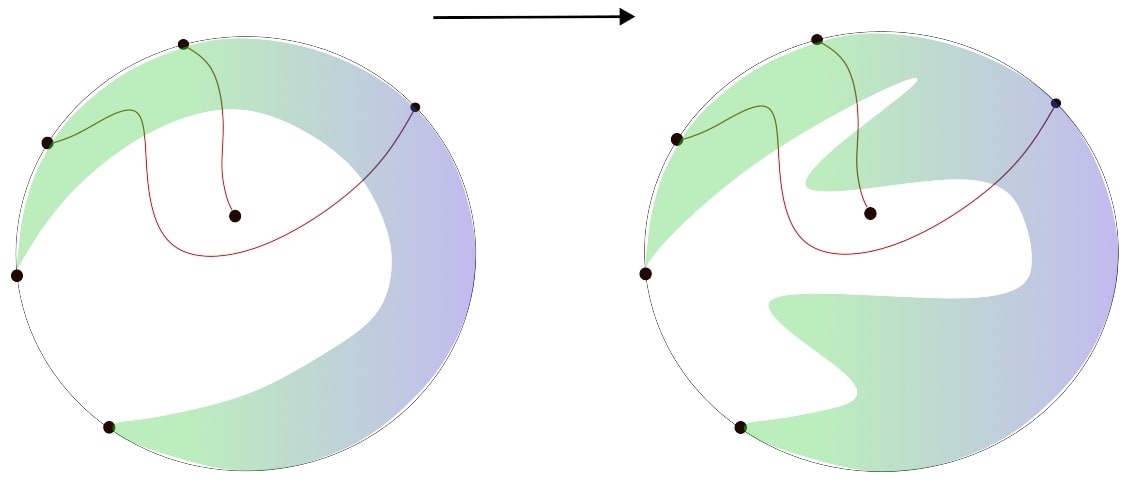}
\put(50,20){$P_{In}$}
\put(230,250){$P_{0}$}
\put(160,290){$\delta$}
\put(-10,300){$r_0$}
\put(545,300){$r_0$}
\put(-100,200){$f_p(r_0)$}
\put(500,200){$g(r_0)$}
\put(380,350){$P_{Out}$}
\put(720,230){$P_0$}
\put(730,300){$\delta$}
\put(600,20){$P_{In}$}
\put(180,150){$H_p$}
\put(420,150){$L$}
\put(-0,150){$l$}
\put(650,190){$H_p$}
\put(1000,150){$L$}
\put(560,150){$l$}
\put(960,360){$P_{Out}$}
\put(790,100){$g(H_p)$}
\put(200,350){$f_p(H_p)$}
\end{overpic}
\caption[An isotopy which preserves $Q$.]{\textit{Since $g$ (on the right) and $f_p$ (on the left) are isotopic and $g$ satisfies the assumptions of Cor.\ref{deformation11}, the periodic orbits in $Q$ all persist when $f_p$ is isotoped to $g$. As can be seen, in this scenario $F_p$ is a trefoil parameter falling into Case $A$ (see Fig.\ref{optA}).}}
\label{home}
\end{figure}
\section{ The implications:}
Th.\ref{th31} teaches us that at trefoil parameters, the dynamics of trefoil parameters are complex essentially like those of a Smale Horseshoe, suspended around a heteroclinic trefoil knot. Motivated by the theory of homoclinic bifurcations (and in particular, Shilnikov's Theorem - see \cite{LeS}), in this section we use Th.\ref{th31} to prove several results about the dynamical complexity of the Rössler system and its bifurcations around trefoil parameters.\\ 

To begin, recall $P$ always denote the parameter space introduced in page \pageref{eq:9}, and that given a parameter $v\in P$ we denote by $F_v$ the vector field corresponding to $v=(a,b,c)$ (see Eq.\ref{Field}). This section is organized as follows - first we study the question of hyperbolicity of the dynamics at trefoil parameter: namely, in Prop.\ref{nohyp} we prove the flow at trefoil parameters $F_p$ are non-hyperbolic. Following that, we prove that in the space of $C^\infty$ vector fields (on $\mathbf{R}^3$), trefoil parameters are inseparable from period-doubling and saddle-node bifurcation sets (for a more precise formulation, see Prop.\ref{spiraltheorem}). Finally, we conclude this Section with the following result: given a trefoil parameter $F_p$, a positive $n>0$ and $v\in P$, provided $v$ is sufficiently close to $p$ the vector field $F_v$ generates at least $n$ distinct periodic trajectories (see Th.\ref{conti}).\\

To begin, let us first recall the notion of Dominated Splitting, a weaker form of hyperbolicity:

\begin{definition}
    \label{hyperbolicvector}     Let $M$ be a Riemannian manifold and let $\phi_t:M\rightarrow M, t\in\mathbf{R}$ be a smooth flow. A compact, $\phi_t$-invariant set $\Lambda$ is said to satisfy a \textbf{dominated splitting} \textbf{condition} if the following conditions are satisfied:
\begin{itemize}
    \item The tangent bundle satisfies $T\Lambda=S\oplus U$, s.t. $S=\cup_{x\in\Lambda}S(x)\times\{x\}$, $U=\cup_{x\in\Lambda}U(x)\times\{x\}$, where $T_x M= S(x)\oplus U(x)$ and $S(x)$ denotes the stable directions while $U(x)$ the unstable directions.
    \item There exists some $ c>0,0<\lambda<1$, s.t. for all $ t>0$ and every $x\in\Lambda$, $(||D\phi_t|_{U(x)}||)(||D\phi_{-t}|_{S(\phi_t(x))}||)<c\lambda^t$.
    \item   For all $x\in \Lambda$, $S(x),U(x)$ vary smoothly with $t\in\mathbf{R}$ to $S(\phi_t(x)),U(\phi_t(x))$.
\end{itemize}
\end{definition}
Now, let us recall that in Cor.\ref{nosep} we defined $I$ as the maximal invariant set of the first-return map $f_p:\overline{D_\alpha}\setminus\{P_0\}\to\overline{D_\alpha}\setminus\{P_0\}$ in $\overline{D_\alpha}\setminus l$ - as proven in Th.\ref{th31}, the set $I$ includes infinitely many periodic orbits for $f_p$. We now prove:

\begin{proposition}\label{nohyp}
Let $F_p$ be a trefoil parameter. Then, the Rössler system at trefoil parameters does not satisfy any Dominated Splitting condition on $I$. 
\end{proposition}
\begin{proof}
We first note the set $I_n= \overline{D_\alpha}\setminus \cup_{i=0}^n \overline{f^{-i}_p(l)}$ is open and dense in $D_\alpha$, hence by the Baire Category Theorem so is $I=\cap_n I_n$. This implies $I$ is dense around the fixed points, hence its suspension w.r.t. the flow $\Lambda$ must also be dense around the heteroclinic trajectories. With these ideas in mind, we prove Prop.\ref{nohyp} by contradiction. To do so assume there exists an invariant set $\Lambda$ as above which includes the fixed points for $F_p$  - by the closure of $\Lambda$ and the continuity of the flow it follows $H\subseteq\Lambda$ Moreover, as we assume by contradiction we can decompose $T\Lambda=S\oplus U$ as in Def.\ref{hyperbolicvector}  it follows that at each fixed point $s\in\Lambda$ for $F$ we can write $T_sS^3=S(s)\oplus U(s)$ - where $U(s)$ corresponds to the unstable directions and $S(s)$ to the stable directions.\\

By assumption we know $F_p$ has precisely two fixed points in $S^3$, $P_{In}$ and $P_{Out}$ - both saddle foci of opposing indices, connected by a heteroclinic trajectory which lies in the intersection of their one-dimensional manifolds. In more detail, we know the said heteroclinic trajectory forms the one-dimensional unstable manifold of $P_{Out}$, while it is also the stable manifold of, say, $P_{In}$. Therefore, by the discussion above above we know $U(P_{Out})$ is one-dimensional while $S(P_{Out})$ is two dimensional - and similarly, $U(P_{In})$ is two dimensional while $S(P_{In})$ is one-dimensional. Moreover, by the continuity of the flow we know the same is true for all initial conditions $x\in\Lambda$ sufficiently close to the fixed points - that is, of $x\in\Lambda$ is sufficiently close to $P_{Out}$ the space $U(x)$ would be one dimensional while $S(x)$ would be two dimensional, and when $x$ is sufficiently close to $P_{In}$ we have the opposite.\\

Now, consider any given $s\in\Lambda$ sufficiently close to the bounded heteroclinic trajectory and denote the flow by $\phi_t,t\in\mathbf{R}$. Since the trajectory of $x$ tends to $P_{In}$ by $s\in\Lambda$ we conclude that for all sufficiently large $t>0$ the space $U(\phi_t(x))$ is two-dimensional while $S(\phi_t(x))$ is one dimensional - and since by the dominated splitting assumption these invariant subspaces vary smoothly with $t$ it follows the same is true for $s$, i.e. $U(s)$ is also two-dimensional while $S(s)$ is one dimensional. On the other hand, since the backwards trajectory of $s$ tends to $P_{Out}$ it follows that for all sufficiently small $t<0$ the space $U(\phi_t(s))$ would be one dimensional while $S(\phi_t(s))$ would be two dimensional - and similarly it would follow $U(s)$ is one dimensional while $S(s)$ is two dimensional. This is a contradiction, which implies there can be no such $\Lambda$ Prop.\ref{nohyp} follows.
\end{proof}

Having proven the non-hyperbolic nature of the flow, we proceed to study the bifurcations around trefoil parameters. To begin, let the vector field $F_p$ be a trefoil parameter and recall we denote the saddle-indices for the saddle-foci $P_{In},P_{Out}$ by $\nu_{In},\nu_{Out}$ - respectively (see the definition at page \pageref{eq:9}). Additionally, recall that we assume the saddle-indices corresponding to $v$ always satisfy either $\nu_{In}<1$ or $\nu_{Out}<1$ (see the discussion in page \pageref{eq:9}). Plugging in this assumption, we derive the following result:

\begin{proposition}\label{spiraltheorem}
Let $F_p$ be a trefoil parameter. Then, the vector field $F_p$ is an accumulation point for infinitely many period-doubling and saddle-node bifurcation sets in the space of $C^\infty$ vector fields on $\mathbf{R}^3$. 
\end{proposition}
\begin{proof}
For completeness, let us first recall Sections 3.2 and 3.1 of \cite{GKP}. In that paper, the setting is that of a smooth, two-parameter family of vector fields $\{X_{(a_1,a_2)}\}_{(a_1,a_2)\in O}$, with $(a_1,a_2)$ varying smoothly in some open set $O\subseteq\mathbf{R}^2$. Assume there exists a curve $\gamma\subseteq O$ corresponding to the existence of homoclinic trajectories. Also assume the saddle indices along $\gamma$ are all strictly lesser then $1$. As proven in \cite{GKP}, under these assumptions there exists a family of spirals, $\{\delta_o\}_{o\in\gamma}$ satisfying:
\begin{itemize}
    \item  $\delta_o$ is a curve in $O$, spiralling towards $o$.
    \item  $\delta_o$ is a subset of $O$ corresponding the existence of a periodic trajectory $T$ for the vector field $X_s$, $s\in\delta_0$.
    \item Provided the saddle index along $\gamma$ is bounded, there exists some $c>0$ s.t. $\forall o\in\gamma$, $diam(\delta_o)>c$ (with that diameter taken w.r.t. the Euclidean metric in $O$)
    \item As $s\in\delta_o$ goes to $o$ along $\delta_0$, $T$ undergoes a cascade of period-doubling and saddle node bifurcations. Consequentially, the period-doubling and saddle-node bifurcation sets in $O$ are dense around every point in $o\in\gamma$.
\end{itemize}

Now, back to the Rössler system. To begin, recall we denote by $F_p$ the corresponding vector field to the Rössler system at trefoil parameters, and assume first that $\nu_{In}<1$ - additionally, recall that since $F_p$ is a trefoil parameter, by Def.\ref{def32} we know the two dimensional manifolds $W^s_{Out},W^u_{In}$ coincide. Therefore, we can perturb $F_p$ at some small neighborhood of $P_{Out}$ and combine the heteroclinic trajectory $\Theta$ and some flow line in $W^u_{In}=W^s_{Out}$ to create a homoclinic trajectory to $P_{In}$ - and moreover, we can choose this perturbation to be arbitrarily close to $F_p$ in the $C^\infty$ metric. In other words, we have just proven there exists an arbitrarily small $C^\infty$ approximation of $F_p$, $V$ s.t. $V$ generates a homoclinic trajectory to the fixed point $P_{In}$. Consequentially, there exists a smooth curve $\gamma_{In}$ in the space of $C^\infty$ vector fields on $\mathbf{R}^3$, s.t. for every $\omega\in\gamma_{In}$, the vector field corresponding to $\omega$ generates a homoclinic curve to $P_{In}$. Moreover, by the description of these perturbations above, we can choose the curve $\gamma_{In}$ s.t. for every $\omega\in\gamma_{In}$, the vector fields $\omega, F_p$ coincide around the fixed point $P_{In}$ - consequentially, the saddle index $\nu_{In}$ is constant along $\gamma_{In}$.\\

As $\nu_{In}<1$ the discussion above implies the period-doubling and saddle node bifurcations sets are dense around every vector field on $\gamma_{In}$ (in the space of $C^\infty$ vector fields on $\mathbf{R}^3$) - and since $F_p\in\overline{\gamma_{In}}$, it follows these bifurcation sets are also dense around $F_p$. The argument for the case $\nu_{Out}<1$ is symmetric, and Prop.\ref{spiraltheorem} now follows.
\end{proof}
\begin{remark}
   Assume both $\nu_{In},\nu_{Out}<1$. In that case, the proof of Prop.\ref{spiraltheorem} essentially proves the bifurcations around trefoil parameters in the space of $C^\infty$ vector fields on $\mathbf{R}^3$ are doubly as complicated when compared to homoclinic bifurcations. This should be contrasted with the results of \cite{By}, where it was proven a two-dimensional parameter space cannot be used to completely describe the bifurcations around heteroclinic vector fields like $F_p$.
\end{remark}

\begin{figure}[h]
\centering
\begin{overpic}[width=0.65\textwidth]{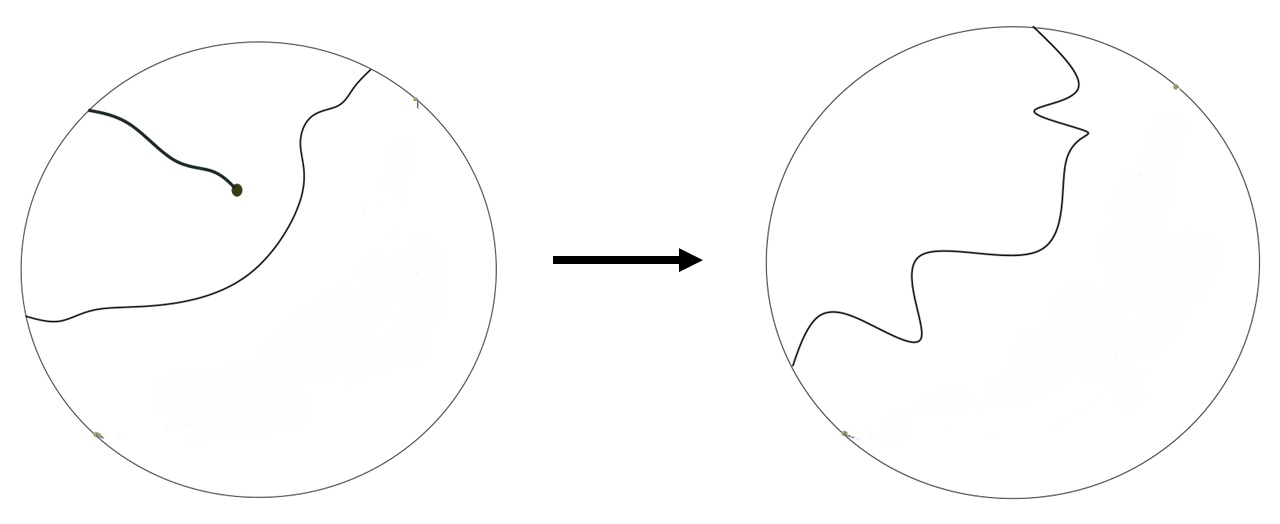}
\put(30,50){$P_{In}$}
\put(610,50){$P_{In}$}
\put(190,250){$P_0$}
\put(150,290){$\delta$}
\put(170,200){$\rho$}
\put(670,200){$\rho_v$}
\put(80,230){$U_2$}
\put(720,290){$U_{2,v}$}
\put(210,120){$U_1$}
\put(790,120){$U_{1,v}$}
\put(320,340){$P_{Out}$}
\put(930,340){$P_{Out}$}
\put(40,330){$\delta_0$}
\end{overpic}
\caption[Fig32]{the deformation of $U_p$ (on the left) to $U_v$ (on the right). The curve $\rho$ is deformed to $\rho_v$. For simplicity, $U_p$ and $U_v$ are sketched as discs rather than half-planes.}
\label{part}
\end{figure}
We now study how the dynamical complexity given by Th.\ref{th31} wears off on nearby vector fields in the parameter space $P$ (or more generally, on $C^1$ close vector fields). That is, we now prove that given a trefoil parameter $p$ in $P$ and given any $n>0$, provided a parameter $v$ in the parameter space $P$ generates dynamics sufficiently $C^1$ close to those of trefoil parameter, the Rössler system corresponding to $v$ generates at least $n$ periodic trajectories. This has the following heuristic meaning - the dynamical complexity given by Th.\ref{th31} "wear off" on nearby vector fields.\\

Despite the simple formulation, the proof would be very technical - as it would rely on constructing explicit (local) homotopies of the first-return map $f_p$, and proving these homotopies satisfy very specific topological conditions. However, the idea behind the proof is relatively simple, and based off the following intuition - if $F_p$ is a trefoil parameter, by Th.\ref{th31} the dynamics of the flow are essentially those of a suspended Smale Horseshoe map. As such, since suspended Smale Horseshoes are hyperbolic on their periodic trajectories (see, for example, \cite{S}), we would expect every periodic trajectory given by Th.\ref{th31} to persist under sufficiently small perturbations of the vector field.\\

In practice, due to Prop.\ref{nohyp} we cannot assume any hyperbolicity condition on any of the periodic trajectories given by Th.\ref{th31}. Therefore, to prove the persistence of periodic dynamics we will take a different route, and apply the notion of the Fixed-Point Index to the first-return map of the flow (see Ch.VII.5 in \cite{Dold}). The reason we do so is because the Fixed-Point Index allows us to "ignore" the non-hyperbolic nature of the flow, and treat the periodic dynamics as a homotopy-invariant of the first-return map. In more detail, whenever the Fixed-Point Index is non-zero, it serves as a homotopy-invariant indicator for the existence of periodic points - as such, it allows us to study the persistence of periodic dynamics in the absence of hyperbolicity conditions.\\  

To begin, let us first recall several facts and notations. First, given a trefoil parameter $F_p$, $p=(a,b,c)$, recall $f_p:\overline{U_p}\to\overline{U_p}$ always denotes the first-return map generated by the Rössler system on the half-plane $U_p\subseteq\{\dot{y}=0\}$ (wherever defined - see Lemma \ref{obs}) - in particular, recall $U_p=\{(x,-\frac{x}{a},z)|-z+\frac{x}{a}<0\}$. Additionally, recall the curve $\rho$ given by Cor.\ref{corsym} which partitions $\overline{U_p}\setminus\rho$ into two components, $U_1$ and $U_2$ - and further recall we denote the one-sided shift by $\sigma:\{1,2\}^\mathbf{N}\to\{1,2\}^\mathbf{N}$. Now, let $v\in P$, $v=(a',b',c')$ be some parameter - by the parameterization of $U_p$, as we smoothly deform the vector field $F_p$ to $F_v$, the half-plane $U_p$ is smoothly deformed to the half-plane $U_v$. Consequentially, the curve $\rho\subseteq\overline{U_p}$ is deformed to some curve $\rho_v\subseteq\overline{U_v}$ - which implies $U_1$ and $U_2$ are continuously deformed to $U_{1,v}$ and $U_{2,v}$, the components of $\overline{U_v}\setminus\rho_v$ (see the illustration in Fig.\ref{part}).\\

To continue, denote by $f_v:\overline{U_v}\to\overline{U_v}$ the first-return map corresponding to the vector field $F_v$ (wherever defined in $\overline{U_v}$ - see Lemma \ref{obs}), and set $I_v$ as the maximal invariant set of $f_v$ in $\overline{U_v}\setminus\rho_v$ - that is, $I_v$ is the maximal collection of initial conditions $x\in\overline{U_v}\setminus\rho_v$ whose forward trajectory never hits $\rho_v$ or diverges to $\infty$. Finally, since by definition $I_v\subseteq (U_{1,v}\cup U_{2,v})$ it follows there exists a symbolic coding $\pi_v:I_v\to\{1,2\}^\mathbf{N}$ s.t. $\pi_v\circ f_v=\sigma\circ\pi_v$. With these ideas and notations in mind, we now prove the following:

\begin{theorem}
    \label{conti}
Let $s\in\{1,2\}^\mathbf{N}$ be periodic of minimal period $k$ that is not the constant $\{1,1,1...\}$. Then, for all parameters $v\in P$ which generate a Rössler system sufficiently $C^1$-close to the dynamics of a trefoil parameters we have:

    \begin{itemize}
        \item $s\in\pi_v(I_v)$. 
        \item $\pi^{-1}_v(s)$ contains at least one periodic point $x_s$ for $f_v$, of minimal period $k$. 
        \item  The functions $f_v,f^2_v,...,f^k_v$ are all continuous at $x_s$.
        \item $\pi_v$ is continuous at $x_s,f_v(x_s),...,f^{k-1}_v(x_s)$.
    \end{itemize}

    Consequentially, given any $n$, whenever $v$ corresponds to a Rössler system sufficiently $C^1$-close to the dynamics at trefoil parameter, the Rössler system corresponding to $v$ generates at least $n$ distinct periodic trajectories. In other words, the more a given Rössler system is $C^1$-close to its idealized model, the more complex its dynamics will be.
\end{theorem}
\begin{proof}
From now on to the end of the proof, for any parameter $v$ in the parameter space $P$, denote by $F_v$ the corresponding vector field (see Eq.\ref{Field}) - moreover, from now on, $p$ would always denote a trefoil parameter (we often refer to $v$ and $p$ both as parameters and as vector fields). The proof of Th.\ref{conti} would be technical and based on direct topological analysis of the first-return map. Therefore, before giving an outline of the proof, we first recall several facts from Section $3$ (and in particular, the proof of Th.\ref{th31}). To begin, recall the topological disc $D_\alpha\subseteq U_p$, bounded by the curve $L\cup l\cup\{P_{In},P_{Out}\}$ - as proven in Prop.\ref{arccor}, $D_\alpha$ is a topological disc on $U_p$, and as shown in Lemma \ref{noempty}, the heteroclinic trefoil knot intersects $\overline{D_\alpha}$ in precisely one interior point, $P_0$, thus making $D_\alpha\setminus\{P_0\}$ homeomorphic to a punctured disc (see the illustration in Fig.\ref{trefint}). Moreover, as proven in Lemma \ref{firstret}, the first-return map $f_p:\overline{D_\alpha}\setminus\{P_0\}\to\overline{D_\alpha}\setminus\{P_0\}$ is well defined (by Remark \ref{extend}, we also know $f_p$ is surjective).\\

Now, recall there exists another cross-section $H_p\subseteq D_\alpha$, also a topological disc (as defined in Lemma \ref{lemhp}) - and that we proved Th.\ref{th31} by studying the map $f_p: \overline{H_p}\setminus\overline\delta\to D_\alpha$ (recall $f_p$ is continuous on $H_p$ - see Lemma \ref{lemhp}). In more detail, we proved the existence of an invariant set $Q$ for $f_p$ in $\overline{H_p}\setminus\overline\delta$ (where $\delta$ is as in Prop.\ref{lem33}) s.t. the following is satisfied: 

\begin{itemize}
    \item $Q\subseteq \overline{U_p}\setminus\rho$, and $f_p$ is continuous on $Q$ - in particular, $Q$ is a subset of the invariant set $I$ given by Cor.\ref{corsym}.
    \item $Q$ includes infinitely many periodic orbits for $f_p$.
    \item There exists a continuous map $\pi:Q\to\{1,2\}^\mathbf{N}$ s.t. $\pi\circ f_p=\sigma\circ\pi$ (where $\sigma:\{1,2\}^\mathbf{N}\to\{1,2\}^\mathbf{N}$ denotes the one-sided shift).
    \item $\pi(Q)$ includes every periodic sequence $s\in\{1,2\}^\mathbf{N}$ 
    \item  For every periodic $s\in\{1,2\}^\mathbf{N}$ of minimal period $k$, $D_s=\pi^{-1}(s)$ is connected - and includes a periodic point $x_s$ of minimal period $k$ for $f_v$. Moreover, whenever $s$ is not the constant $\{1,1,1...\}$, $x_s$ is not a fixed-point for the flow.
\end{itemize}
 
With these ideas in mind, we are now ready to give a general outline of the proof of Th.\ref{conti}. The general idea is as follows - we begin by constructing an isotopy of $f_p:\overline{H_p}\setminus\delta\to\overline{D_\alpha}$ to a Smale Horseshoe map, which we use to prove the Fixed-Point Index (as defined below) on periodic orbits in $Q$ is non-zero. Following that, we prove that for any parameter $v\in P$ s.t. the corresponding Rössler system is sufficiently $C^1-$ close to $p$, the perturbation of $p$ to $v$ induces a homotopy between the corresponding first-return maps - then, using the invariance of the Fixed-Point Index under homotopies Th.\ref{conti} would follow.\\

To begin, fix some periodic $s\in\{1,2\}^\mathbf{N}$ which is not the constant $\{1,1,1...\}$ - by the construction of $Q$ in the proof of Th.\ref{th31} in Stages $I$, $II$ and $III$, this implies $D_s$ lies away from the fixed-point $P_{In}$. Now, define $Per(s)=\{x\in D_s|f^k_p(x_s)=x_s\}$ - it is easy to see $Per(s)$ is compact in $\overline{H_p}$. Moreover, since for every $0\leq j< k$ we have $f^j_p(D_s)\cap D_s=\emptyset$, it is also immediate that if $x\in Per(s)$, then $x$ is periodic of minimal period $k$ for $f_p$. We first prove the following technical (yet useful) fact:

\begin{lemma}
    \label{isolat}
    Let $F_p$ be a trefoil parameter, and let $s\in\{1,2\}^\mathbf{N}$ be periodic of minimal period $k$ s.t. $s$ is not the constant $\{1,1,1...\}$. Then, there exists a connected open set $V_s\subseteq {H_p}$, $Per(s)\subseteq V_s$ s.t. the following is satisfied:
    \begin{itemize}
        \item   $f_p^k$ has no fixed points in $V_s\setminus Per(s)$.
        \item $Per(s)$ is compact in $V_s$.
        \item For every $0\leq i<j\leq k-1$, $f^{j}_p(V_s)\cap f^i_p(V_s)=\emptyset$.
        \item $f_p,...,f^k_p$ is continuous on $V_s$.
    \end{itemize}
    
\end{lemma}
\begin{proof}
We first make some general observations. Let us first recall we have $\partial D_\alpha=L\cup l\cup\{P_{In},P_{Out}\}$ (see Prop.\ref{arccor}), and that since the curve $L$ is the intersection between the two-dimensional invariant manifold $W^s_{Out}$ and the half-plane $U_p$ (see the illustration in Fig.\ref{trefint1}) - consequentially, there are no periodic orbits for $f_p$ in $L$. Additionally, using a similar argument to the one used to prove Prop.\ref{lem33}, it follows every component of $f^{-k}_p(l),k\geq0$ is a curve in $\overline{D_\alpha}$ with two endpoints on $\cup_{0\leq j\leq k-1}f^{-j}_p(l)$ - see the illustration in Fig.\ref{trefint1}.\\

\begin{figure}[h]
    \centering
\begin{overpic}[width=0.6\textwidth]{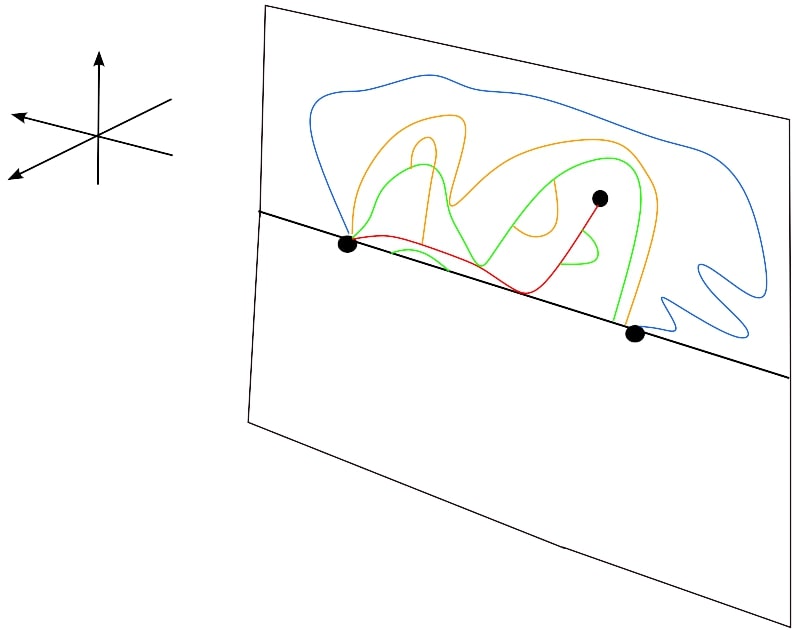}
\put(720,350){$P_{In}$}
\put(640,380){$\delta_0$}
\put(690,530){$P_0$}
\put(1050,200){$\{\dot{y}<0\}$}
\put(100,150){$\{\dot{y}>0\}$}
\put(350,450){$P_{Out}$}
\put(570,400){$l$}
\put(770,640){$L$}
\put(350,730){$U_p$}
\put(410,600){$D_\alpha$}
\put(700,200){$L_p$}
\put(-10,645){$x$}
\put(-20,550){$y$}
\put(110,730){$z$}
\end{overpic}
\caption[The disc $D_\alpha$.]{\textit{The cross-section $U_p$. $L$ is the blue curve, $l$ is the straight line connecting $P_{In}$ and $P_{Out}$, while $f^{-1}_p(l)$, $f^{-2}_p(l)$ and $f^{-3}_p(l)$ are the red, green and orange arcs (respectively). Every component in $f^{-j}_p(l)$ is a curve with endpoints on $\cup_{i=0}^{j-1} f^{-i}_p(l)$. In this scenario, $H_p$ is the sub-region of $D_\alpha$ trapped between $l,f^{-1}_p(l)$ and $L$. $\delta$ is the red arc connecting $P_0$ and $\delta_0$ - see Prop.\ref{lem33}.}}
\label{trefint1}
\end{figure}

We now claim that $Per(s)\cap(\cup_{n\geq0}f^{-n}_p(l))=\emptyset$. To see why, assume by contradiction this is not the case - which implies there exists some $x\in Per(s)$ and $j\geq0$ s.t. $f^j_p(x)\in l$. Using a similar argument to the proof of Lemma \ref{continu2}, this implies there exists an isotopy of $f_p:\overline{H_p}\to\overline{D_\alpha}$ which collapses $x$ to $P_{In}$, similarly to as illustrated in Fig.\ref{DEFFF}. Since all the periodic points in $Q$ given by Th.\ref{th31} are uncollapsible (per Th.\ref{stability} and Cor.\ref{deformation11}), by $Per(s)\subseteq Q$ we derive a contradiction. Consequentially, it follows $Per(s)\cap(\cup_{n\geq0}f^{-n}_p(l))=\emptyset$.\\

\begin{figure}[h]
\centering
\begin{overpic}[width=0.7\textwidth]{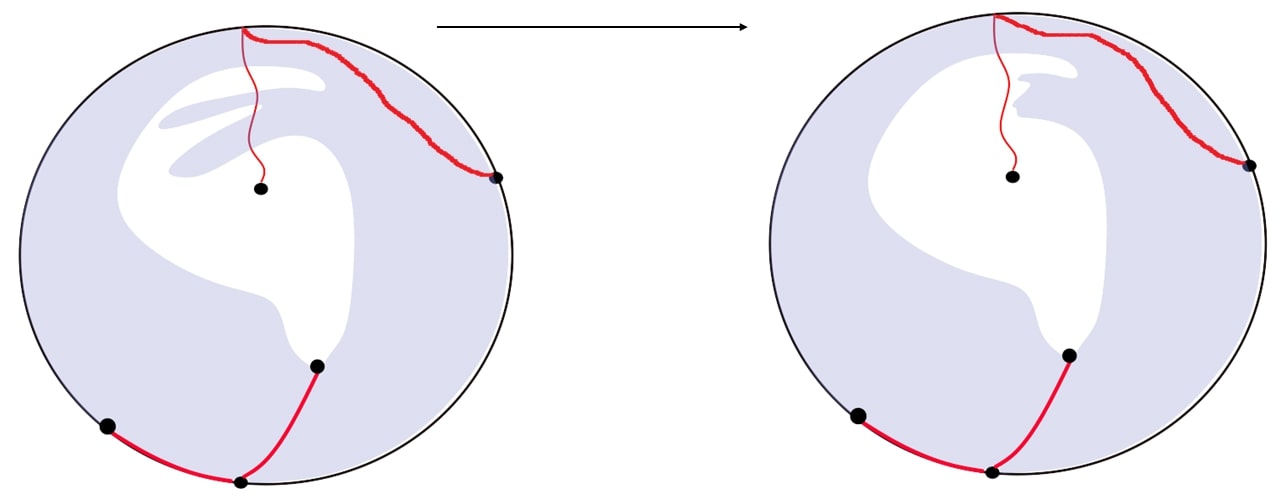}
\put(110,-20){$P_{In}$}
\put(200,130){$f^2_p(\delta_0)$}
\put(230,250){$f_p^{-1}(l)$}
\put(170,200){$P_0$}
\put(170,390){$\delta_0$}
\put(100,70){$f_p(\delta_0)$}
\put(410,240){$P_{Out}$}
\put(700,0){$P_{In}$}
\put(660,90){$g(\delta_0)$}
\put(990,250){$P_{Out}$}
\put(800,140){$g^2(\delta_0)$}
\put(830,250){$g^{-1}(l)$}
\put(740,250){$P_0$}
\put(730,400){$\delta_0$}
\end{overpic}

\caption[Collapsing periodic orbits for $f_0$ on $l$ and $\delta$.]{\textit{$l$ is the left-arc on the circle connecting $P_{In}$ and $P_{Out}$, the shaded region is $f_p(\overline{H_p}\setminus\delta)$ (where $\delta\subseteq f^{-1}_p(l)$ is the red curve connecting $P_0$ and $\delta_0$), while the red arcs are $f_p(\delta)$ and $f^2_p(\delta)$. As $\delta\subseteq f^{-1}_p(l)$, it is easy to see the isotopy of $f_p:\overline{H_p}\setminus\delta\to\overline{D_\alpha}$ to $g:\overline{H_p}\setminus\delta\to\overline{D_\alpha}$ above removes periodic orbits not in $Q$, while preserving $Q$ due to Cor.\ref{deformation11}. }}
\label{remove1}
\end{figure}

To continue, consider $\{x_n\}_n$, a sequence of periodic points in $\overline{D_\alpha}\setminus Per(s)$ which tends to some $x\in Per(s)$ - we now claim the minimal period of $f_p$ on elements in $\{x_n\}_n$ is unbounded. To do so, assume by contradiction the assertion is incorrect, i.e., that the period of $x_n$ is bounded - since $x\in Per(s)$, by the continuity and smoothness of the flow this would imply the minimal period of $x_n$ has to be $k$ (for any sufficiently large $n$). Since $\{1,2\}^\mathbf{N}$ only has a finite number of periodic symbols of minimal period $k$, and because for every periodic $\omega,s\in\{1,2\}^\mathbf{N}$ both $Per(s)$ and $Per(\omega)$ are compact and disjoint, it follows that for any sufficiently large $n$, $x_n\not\in Q$ - therefore, without any loss of generality, we may assume $\{x_n\}_n\cap Q=\emptyset$. Since $Per(s)\cap(\cup_{n\geq0}f^{-n}_p(l))=\emptyset$, by $\{x_n\}_n\cap Q=\emptyset$ it follows $x$ is separated from $\{x_n\}_n$ by the components of $\cup_{n\geq0}f^{-n}_p(l)$, hence the set $\cup_{n\geq0}f^{-n}_p(l)$ accumulates on $x$. \\

\begin{figure}[h]
\centering
\begin{overpic}[width=0.7\textwidth]{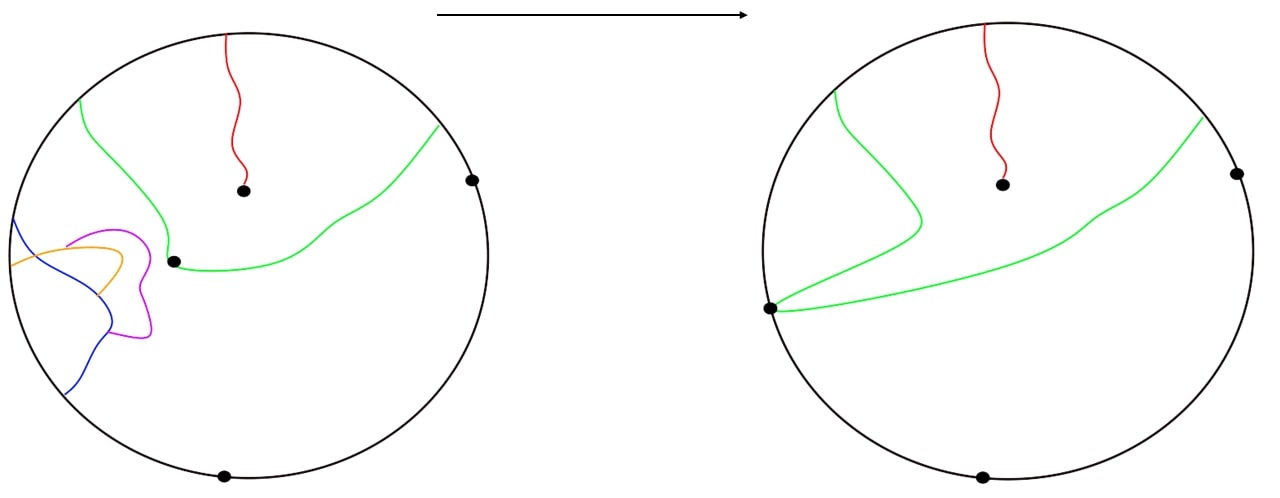}
\put(110,-20){$P_{In}$}
\put(200,130){$D_s$}
\put(200,240){$P_0$}
\put(170,390){$\delta_0$}
\put(150,200){$x$}
\put(30,130){$x_3$}
\put(60,180){$x_2$}
\put(80,150){$x_1$}
\put(390,240){$P_{Out}$}
\put(700,0){$P_{In}$}
\put(570,140){$x$}
\put(990,250){$P_{Out}$}
\put(800,140){$D_s$}
\put(740,250){$P_0$}
\put(730,400){$\delta_0$}
\end{overpic}

\caption[Collapsing periodic orbits for $f_0$ on $l$ and $\delta$.]{\textit{$l$ is the left-arc on the circle connecting $P_{In}$ and $P_{Out}$, the shaded region is $f_p(\overline{H_p}\setminus\delta)$ (where $\delta$ is the red curve connecting $P_0$ and $\delta_0$), while the red arcs are $f_p(\delta)$ and $f^2_p(\delta)$. As $\{x_n\}_n$ accumulates on $x$, the isotopy described in the proof of Lemma \ref{isolat} and depicted in Fig.\ref{remove1} pushes the sequence $\{x_n\}_n$ towards $l$, thus pushing $g^j(x)$ to $l$ as well (for some $j$ - in this illustration, we chose $j=1$). }}
\label{remove}
\end{figure}
    
Now, deform $f_p:\overline{H_p}\setminus\overline\delta\to\overline{D_\alpha}$ isotopically to some $g:\overline{H_p}\setminus\overline\delta\to\overline{D_\alpha}$ by destroying all the periodic orbits generated by $\{x_n\}$ as illustrated in Fig.\ref{remove1} and \ref{remove} - moreover, we do so s.t. the assumptions in Cor.\ref{deformation11} are satisfied. In more detail, we do so by collapsing to singletons the arcs in $\cup_{n\geq0}f^{-n}_p(l)$ and domains in $\overline{D_\alpha}\setminus(\cup_{n\geq0}f^{-n}_p(l))$ which include elements of $\{x_n\}_n$ (we can do so by $\{x_n\}_n\cap Q=\emptyset$). Because $x$ is separated from $\{x_n\}_n$ by the set $\cup_{n\geq0}f^{-n}_p(l)$, as $x$ is inseparable from $\{x_n\}_n$ it follows  the isotopy described above pushes $x$ to some arc $\gamma$ on $f^{-j}_p(l)$ - which, consequentially, implies $g^j(x)\in l$ (see the illustration in Fig.\ref{remove}). However, this implies we can isotopically deform $g:\overline{H_p}\setminus\delta\to\overline{D_\alpha}$ to $g':\overline{H_p}\setminus\delta\to\overline{D_\alpha}$ by collapsing $x$ to $P_{In}$ as in Fig.\ref{remove1} - and again, using an argument similar to Lemma \ref{continu2}, we have a contradiction. Therefore, by this contradiction we conclude that if $\{x_n\}_n$ is a sequence of periodic orbits for $f_p$ which approximates $x\in Per(s)$, the minimal periods of $\{x_n\}_n$ have to be unbounded.\\

Consequentially, since $Per(s)\cap\partial H_p=\emptyset$, it follows there exists an open, connected neighborhood $V_s\subseteq H_p$ s.t. $Per(s)\subseteq V_s$ s.t. there are no fixed points for $f^k_p$ in $V_s\setminus Per(s)$ - hence $Per(s)$ is compact in $V_s$. Now, recall $f_p,...,f^k_p$ are all continuous on $D_s$, and that for every $0\leq i<j\leq k-1$ we have $f^{j}_p(D_s)\cap f^i_p(D_s)=\emptyset$ - moreover, recall that $D_s$ is in the invariant set of $f_p$ in $\overline{H_p}\setminus\delta$, which implies it also lies away from $H_p\cap(\cup_{n\geq1}f^{-n}_p(l))$. By $Per(s)\subseteq D_s$ and because $D_s$ is connected, it follows we can choose $V_s$ s.t. the flow lines emanating from $V_s$ and it to $f_p(V_s),...,f^k_p(V_s)$ never hit $l$, i.e., they are all transverse to the half-plane $U_p$ - consequentially, by Prop.\ref{lem33} it follows $f_p,...,f^k_p$ are all continuous on $V_s$. Moreover, by $f^{j}_p(D_s)\cap f^i_p(D_s)=\emptyset$ we may also choose $V_s$ s.t. for every $0\leq i<j\leq k-1$ we have $f^{j}_p(V_s)\cap f^i_p(V_s)=\emptyset$. The proof of Lemma \ref{isolat} is now complete.
\end{proof}

To continue, we now recall the notion of the Fixed-Point Index (see \cite{Dold}). To do so, let $O\subseteq R^2$ be an open set, and let $f:O\to\mathbf{R}^2$ be continuous s.t. the set $\{x\in O|f(x)=x\}$ is compact in $O$. The \textbf{Fixed-Point Index} of $f$ would be defined as the degree of $f(x)-x$ on $O$ - it is well-known that when $f$ has no fixed points in $O$ its Fixed-Point Index is $0$ (see Lemma VII.5.5 in \cite{Dold}). Now, let $f_t:O\to\mathbf{R}^2$, $t\in[0,1]$ be a homotopy of continuous maps, and for every $t\in[0,1]$ set $F_t=\{x\in O|f_t(x)=x\}$. As proven in Th.VII.5.8 in \cite{Dold} (termed "homotopy invariance"), provided $\cup_{t\in[0,1]}F_t\times\{t\}$ is compact in $O\times[0,1]$, the fixed point index is invariant under the homotopy. Or in other words, given a homotopy as described above, if the Fixed-Point Index for $f_0$ is non-zero, the same is true for $f_1$ - that is, if $f_0$ has a fixed point in $O$ of non-zero index, $f_1$ also has at least one fixed-point in $O$. \\

We now proceed to compute the Fixed Point Index of $f^k_p$ on $V_s$. We begin by isotopically deforming $f_p:\overline{H_p}\to\overline{D_\alpha}$ by an isotopy $g_t:\overline{H_p}\to\overline{D_\alpha}$, $t\in[0,1]$ s.t. $g_0=f_p$, and $g_1$ is a Smale Horseshoe map $H:ABCD\to\overline{D_\alpha}$ as appears below in Fig.\ref{remove2} (i.e., we deform $H_p$ to a topological rectangle $ABCD$). Moreover, we do so s.t. no new fixed points for $f^k_p$ are added in $V_s\setminus Per(s)$ throughout the isotopy.\\

\begin{figure}[h]
\centering
\begin{overpic}[width=0.7\textwidth]{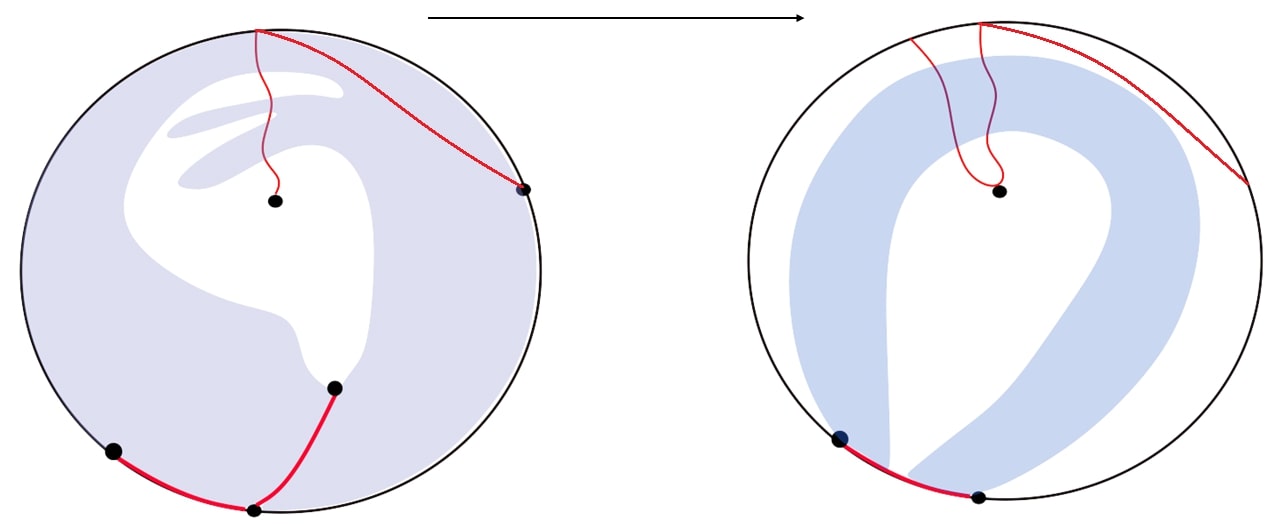}
\put(120,-20){$P_{In}$}
\put(200,130){$f^2_p(\delta_0)$}
\put(250,270){$f^{-1}_p(l)$}
\put(170,210){$P_0$}
\put(170,400){$\delta_0$}
\put(100,70){$f_p(\delta_0)$}
\put(420,250){$P_{Out}$}
\put(720,-20){$B$}
\put(610,40){$A$}
\put(640,90){$H(CD)$}
\put(750,60){$H(AB)$}
\put(740,240){$P_0$}
\put(740,405){$D$}
\put(680,390){$C$}
\end{overpic}

\caption[Collapsing periodic orbits for $f_0$ on $l$ and $\delta$.]{\textit{The isotopy deforming $f_p:\overline{H_p}\setminus\delta\to\overline{D_\alpha}$ to the Smale Horseshoe map $H:ABCD\to\overline{D_\alpha}$. $ABCD$ is obtained from $\overline{H_p}\setminus\delta$ by opening $P_{In}$ to the $AB$ side and $\delta$ to the $CD$ side.}}
\label{remove2}
\end{figure}

It is easy to see we can choose the isotopy s.t. the periodic orbits in $Per(s)$ are deformed continuously throughout the isotopy - i.e., for every $1\geq t>0$ the set $Per(s)=Per_0(s)$ is continuously deformed to $Per_t(s)=\{x\in V_s|g^k_t(x)=x\}$, a collection of periodic points of minimal period $k$ for $g_t$ - and moreover, by previous paragraph, for every $t\in[0,1]$, $Per_t(s)$ is compact in the open disc $V_s$. Therefore, we have the following:

    \begin{itemize}
        \item   $g_s^k$ has no fixed points in $V_s\setminus Per_t(s)$.
        \item $Per_t(s)$ is compact in $V_s$.
        \item For every $0\leq i<j\leq k-1$, $g^{j}_t(V_s)\cap g^i_t(V_s)=\emptyset$.
        \item $g_t,...,g^k_t$ is continuous on $V_s$.
    \end{itemize}

Consequentially, it follows there exist an isotopy of continuous maps $g^k_t:V_s\to\overline{D_\alpha}$, $t\in[0,1]$ s.t. for every $t$, $g^k_t$ has no fixed points in $V_s\setminus Per_t(s)$ - hence the set $\cup_{t\in[0,1]}Per_t(s)\times\{t\}$ is compact in $V_s\times[0,1]$. Applying the invariance of the Fixed-Point Index under homotopies to the isotopy $g^k_t:V_s\to\overline{D_\alpha}$, $t\in[0,1]$, we now conclude:

\begin{corollary}
\label{fixedpoint}    With the notations above, the fixed-point index of $f^k_p$ on $V_s$ is non-zero.
\end{corollary}
\begin{proof}
    By the invariance of the Fixed-Point Index under homotopies mentioned above, it immediately follows the Fixed-Point Index of $g^k_t=f^k_p$ on $V_s$ is the same as that of $g^k_1=H^k$ on $V_s$ (where $H$ is a Smale Horseshoe) - therefore, provided we prove the Fixed Point Index of $H^k$ on $V_s$ is negative, we are done. To do so, first note $H$ is a $\cap-$Horseshoe map (see Fig.\ref{remove2}), which implies it is orientation-preserving - hence its differential on its invariant set has one eigenvalue in $(1,\infty)$ and another at $(0,1)$. It therefore follows that given any periodic orbit for $H$ of minimal period $k$, the degree of $H^k-Id$ is $\pm 1$.\\
    
Since for every $k>0$ the Smale Horseshoe map $H$ has only a finite number of periodic orbits of period $k$ (minimal or not), it follows $Per_1(s)$ is finite. Consequentially, by the additivity of the Fixed-Point Index (see Lemma VII.5.6 in \cite{Dold}), it follows the Fixed-Point Index of $H^k$ in $V_s$ is given by $\sum_{x\in Per_1(s)}\pm1=I_1\ne0$. It follows $I_1$ is non zero, and Cor.\ref{fixedpoint} now follows.
\end{proof}

Having proven Cor.\ref{fixedpoint}, we now use it to conclude the proof of Th.\ref{conti}. To do so, recall we assume $p=(a,b,c)$ is a trefoil parameter, and consider $v$, a vector field sufficiently $C^1-$close to $p$. Whenever $v$ is sufficiently $C^1-$close to $p$, the half-plane $\overline{U_p}=\{(x,-\frac{x}{a},z)|-z+\frac{x}{a}\leq0\}$ is smoothly deformed the a set to $\overline{U_v}\subseteq\{\dot{y}=0\}$, where $\dot{y}$ is taken w.r.t. $v$ (when $v$ is some parameter in $P$, $U_v$ is the set given by Lemma \ref{obs}). Therefore, since $V_s$ is an open, connected set inside ${U_p}$ it follows that when the dynamics at trefoil parameters are smoothly deformed to those of $v$, the topological disc $V_s\subseteq U_p$ is continuously deformed to $V_v\subseteq U_v$, an open disc subset of the half-plane $U_v$ (see the illustration in Fig.\ref{remove3}). Moreover, as $V_s$ is transverse to the flow at $p$, $V_v$ is transverse  to the flow at $v$.\\

\begin{figure}[h]
\centering
\begin{overpic}[width=0.8\textwidth]{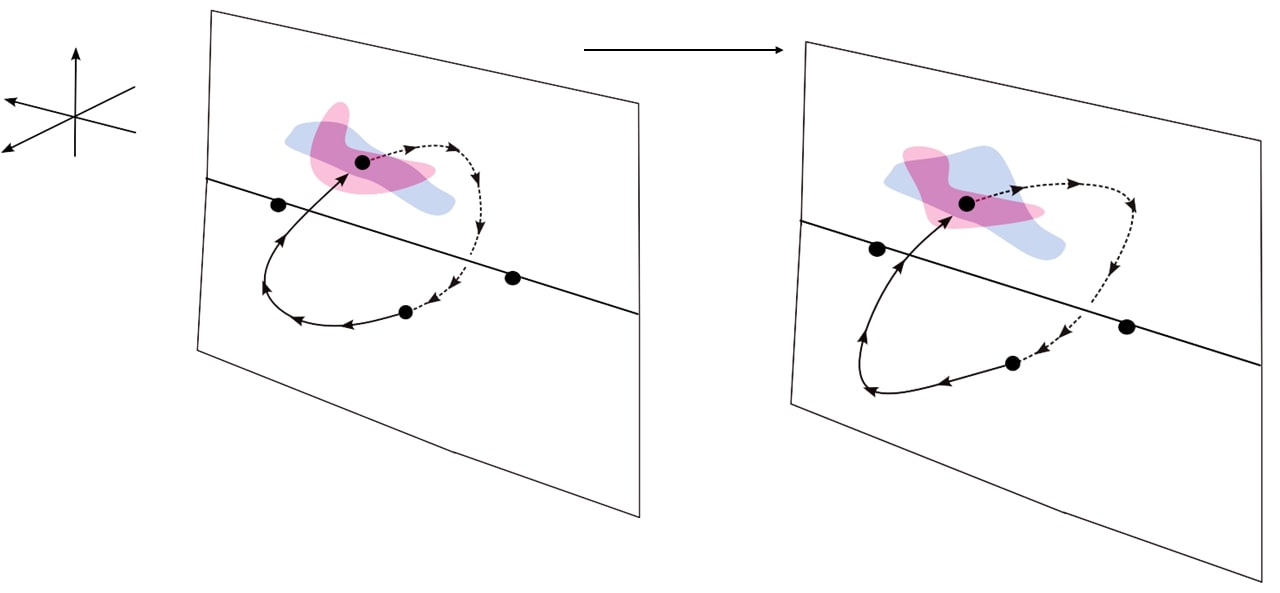}
\put(190,350){$V_s$}
\put(150,285){$P_{Out}$}
\put(270,250){$l$}
\put(290,160){$L_p$}
\put(190,400){$U_p$}
\put(250,400){$f_p(V_s)$}
\put(420,250){$P_{In}$}
\put(720,150){$L_v$}
\put(655,320){$V_v$}
\put(620,255){$P_{Out}$}
\put(850,185){$P_{In}$}
\put(740,220){$l$}
\put(740,365){$f_v(V_v)$}
\put(740,320){$x_s$}
\put(650,390){$U_v$}
\put(50,430){$z$}
\put(-15,380){$x$}
\put(-15,340){$y$}
\end{overpic}

\caption[Collapsing periodic orbits for $f_0$ on $l$ and $\delta$.]{\textit{The sets $V_s$ and $V_v$ in blue, with $f_p(V_s)$ and $f_v(V_v)$ in red (for simplicity, we assume $s\in\{1,2\}^\mathbf{N}$ is periodic of minimal period $1$). As the dynamics at trefoil parameters are smoothly deformed to those of $v$, we induce a homotopy between $f_p:V_s\to U_p$ and $f_v:V_v\to U_v$. }}
\label{remove3}
\end{figure}

To continue, recall we constructed $V_s$ in Lemma \ref{isolat} s.t. the flow lines connecting $V_s,f_p(V_s),...,f^{k}_p(V_s)$ are all transverse to $U_p$. Recalling the cross-section $H_p$ is bounded in $U_p$, it follows $V_s$ is pre-compact - and therefore, $F_v$ is a Rössler system sufficiently $C^1$- close to the dynamics of a trefoil parameter (see Eq.\ref{Field}), the first-hit map $f_v:V_v\to U_v$ w.r.t. the vector field $F_v$ is well defined, and so are its iterates $f^2_v,...,f^k_v$.  Consequentially, it follows that given any $x\in \overline{V_v}$ whenever the Rössler system corresponding to $v$ is sufficiently $C^1-$ close to $p$, the flow lines for $F_v$ connecting $V_v,f_v(V_v),...,f^k_v(V_v)$ are all transverse to $U_v$ (see the illustration in Fig.\
\ref{remove3}) - which implies $f_v,...,f^k_v$ are all continuous on $V_v$. Similarly, since by Lemma \ref{isolat} for every $0\leq i<j\leq k-1$ we have $f^j_p(V_s)\cap f^i_p(V_s)=\emptyset$, it follows that provided $v$ is sufficiently close to $p$, for every $0\leq i<j\leq k-1$ we also have $f^j_v(V_v)\cap f^i_v(V_v)=\emptyset$. Summarizing our results, we conclude:

\begin{corollary}
    \label{pert} Let $V_s$ be as in Lemma \ref{isolat}, and let $V_v$ be as above. Then, provided the dynamics of $v\in P$ are sufficiently $C^1$-close to those of trefoil parameters, we have the following: 
    \begin{itemize}
        \item     $f_v,...,f^k_v$ are all continuous on $V_v$.
        \item For every $0\leq i<j\leq k-1$, $f^i_v(V_v)\cap f^j_v(V_v)=\emptyset$.
        \item $f^k_p:V_s\to U_p$ and $f^k_v:V_v\to U_v$ are homotopic (see the illustration in Fig.\ref{remove3}).
    \end{itemize}

Consequentially, if $x\in V_v$ is periodic for $f_v$, its minimal period is at least $k$.
\end{corollary}

Our next goal is to compute the Fixed-Point Index of $f^k_v$ in $V_v$ - provided we show it is non-zero, by Cor.\ref{pert} it would immediately follow $f_v$ has a periodic point of minimal period $k$ in $V_v$. To this end, we first prove:

\begin{lemma}
\label{compactper}    Let $F_p$ be a trefoil parameter. Then, there exists a $C^1-$ neighborhood $O$, s.t. $F_p\in O$, and for every Rössler system $F_v\in O$ (where $v\in P$ is a parameter for the Rössler system - see Eq.\ref{Field}) the corresponding first-return map $f^k_v$ has no fixed points in $\partial V_v$. 
\end{lemma}
\begin{proof}
We prove Lemma \ref{compactper} by contradiction. To this end, assume the Lemma is incorrect - that is, assume that for any $O$, a $C^1-$ neighborhood of $p$, there exists a sequence $\{(x_n,v_n)\}_n\subseteq \cup_{v\in O}V_v\times \{v\}$ s.t. the following holds:
    \begin{itemize}
        \item $(x_n,v_n)\in\overline {\cup_{v\in O}V_{v_n}\times \{v\}}$ for every $n$ (where $V_{v_n}$ is defined analogously as above for $C^1-$perturbations of $p$). 
        \item The vector fields corresponding to $v_n$ tend to $p$ in the $C^1$-metric.
        \item $f^k_{v_n}(x_n)=x_n$ - where $f_{v_n}:V_n\to U_{v_n}$ is as defined above.
        \item $(x_n,v_n)\to\partial V_s\times\{p\}$.
        \end{itemize}
    Let $x\in\partial V_s$ be the limit (or partial limit) of $x_n$. Since by construction $\overline{V_s}\subseteq\overline{H_p}$, from $\overline{H_p}\subseteq\overline{D_\alpha}$ and by both Prop.\ref{arccor}  and Lemma \ref{firstret} it follows the trajectory of $x$ is bounded. Now, note the vector fields $\{{v_n}\}_n$ converge to $p$ in the $C^1$ metric on vector fields in $\mathbf{R}^3$ - therefore, because every $x_n$ lies on a periodic trajectory for ${v_n}$ which intersects the half-plane $U_{v_n}$ at most $k$ times, it follows $x$ must lie on $T$, a periodic trajectory for $F_p$, and satisfy $f^k_p(x)=x$. However, this contradicts Lemma \ref{isolat} where we proved that for all $x\in\partial V_s$, $f^k_p(x)\ne x$ - hence no such sequence exists and the assertion follows. 
\end{proof}
As a consequence of Lemma \ref{compactper}, Cor.\ref{fixedpoint}, and Cor.\ref{pert}, it follows that for every $v\in P$ whose corresponding Rössler system is sufficiently $C^1$ close to that of a trefoil parameter $p$, we have:
\begin{corollary}
    \label{fix3}
    Let $F_p$ be a trefoil parameter - then, given $v\in P$ s.t. the Rössler system corresponding to $v$ is sufficiently $C^1$-close to $F_p$ the Fixed Point Index of $f^k_v:V_v\to\overline{U_v}$ is non-zero - hence, $f_v$ has a periodic point $x_s$ in $V_v$ of minimal period $k$. Moreover, $f_v,..., f^k_v$ are continuous on $x_s$ (see the illustration in Fig.\ref{remove3}).
\end{corollary}

\begin{figure}[h]
\centering
\begin{overpic}[width=0.8\textwidth]{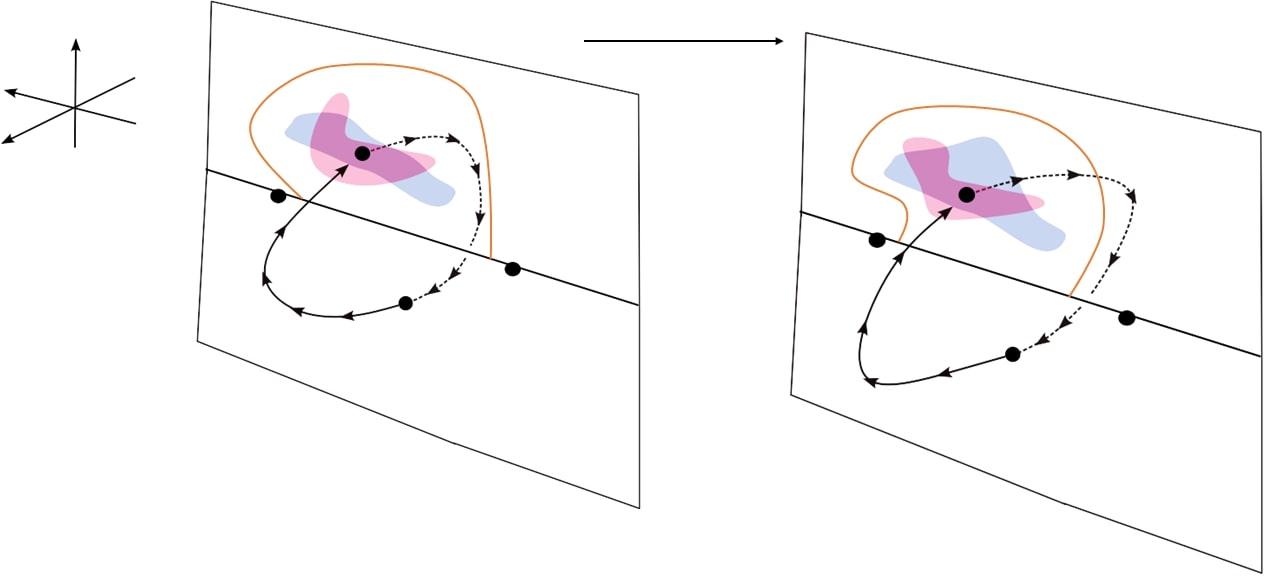}
\put(210,350){$V_s$}
\put(150,285){$P_{Out}$}
\put(270,250){$l$}
\put(290,160){$L_p$}
\put(170,420){$U_1$}
\put(300,285){$U_2$}
\put(260,370){$f_p(V_s)$}
\put(380,370){$\rho$}
\put(420,250){$P_{In}$}
\put(720,150){$L_v$}
\put(675,320){$V_v$}
\put(620,255){$P_{Out}$}
\put(850,185){$P_{In}$}
\put(740,220){$l$}
\put(740,335){$f_v(V_v)$}
\put(855,335){$\rho_v$}
\put(725,310){$x_s$}
\put(640,390){$U_{1,v}$}
\put(780,245){$U_{2,v}$}
\put(50,430){$z$}
\put(-15,380){$x$}
\put(-15,340){$y$}
\end{overpic}

\caption[Collapsing periodic orbits for $f_0$ on $l$ and $\delta$.]{\textit{The sets $V_s$ and $V_v$ in blue, with $f_p(V_s)$ and $f_v(V_v)$ in red (for simplicity, we assume $s\in\{1,2\}^\mathbf{N}$ is periodic of minimal period $1$). The curves $\rho$ and $\rho_v$ are sketched in orange - as $\rho$ is continuously deformed to $\rho_v$, $U_{i}$ is continuously deformed to $U_{i,v}$, $i=1,2$. Provided $v$ is sufficiently close to $p$, $\pi_v(x_s)=s$.}}
\label{symboll}
\end{figure}

We can now conclude the proof of Th.\ref{conti}. To do so, first recall we originally chose $s\in\{1,2\}^\mathbf{N}$ s.t. it is periodic of minimal period $k$, $s=\{i_0,i_1,i_2,...,i_k,i_0,...\}$ and not the constant $\{1,1,1...\}$. Additionally, recall the curve $\rho\subseteq U_p$ and the curve $\rho_v$ which divides the cross-section $U_v$ to $U_{1,v},U_{2,v}$ - and additionally, recall the invariant set $I_v$, and the map $\pi_v:I_v\to\{1,2\}^\mathbf{N}$  (see the discussion before Th.\ref{conti} and the illustration in Fig.\ref{part} and Fig.\ref{symboll}). Having proven Cor.\ref{fix3}, all that remains to conclude the proof of Th.\ref{conti} is to prove the point $x_s$ given by Cor.\ref{fix3} lies inside $I_v$, that $\pi_v(x_s)=s$, and that $\pi_v$ is continuous on the finite sequence $\{x_s,f_v(x_s),...,f^{k-1}_v(x_s)\}$. Note that since $f_v,...,f^k_v$ are continuous at $x_s$ by  Cor.\ref{fix3}, provided $x_s\in I_v$ then $\pi_v$ is automatically continuous at $x_s,f_v(x_s),...,f^{k-1}_v(x_s)$ - therefore, we need only prove that $x_s\in I_v$, and that $\pi_v(x_s)=s$. \\

However, that is immediate - recall that by Cor.\ref{corsym}, for every $0\leq j\leq k$ we have $f^j_p(D_s)\subseteq U_{i_j}$ - where $s=\{i_0,i_1,...,i_k,i_0,...\}$. It therefore follows by $Per(s)\subseteq D_s$ that we can choose $V_s$ to be a sufficiently small neighborhood of $Per(s)$ s.t. for every $0\leq j\leq k$, $f^j_p(V_s)\subseteq U_{i_j}$. Because $U_p$ is continuously deformed to $U_v$ (and consequentially, the curve $\rho$ is continuously deformed to $\rho_v$) it follows that provided the dynamics of the Rössler system corresponding to some $v\in P$ are sufficiently $C^1$- close to the trefoil parameter $F_p$, the set $V_v$ also satisfies $f^j_v(V_v)\subseteq U_{i_j,v}$ for all $0\leq j\leq k$ (see the illustration in Fig.\ref{symboll}). By $x_s\in V_v$, $f^k_v(x_s)=x_s$ and because the minimal period of $x_s$ w.r.t. $f_v$ is $k$, this implies $\pi(x_v)=s$ - and Th.\ref{conti} now follows.
\end{proof}
\begin{remark}
    Given $v\in P$, the set $\pi_v(I_v)$ is never empty. It is easy to see the curve $\rho$ from Cor.\ref{corsym} does not include the fixed point $P_{In}$. As such, the same is true for $\rho_v$ - or in other words, for every $v\in P$, $P_{In}\in U_{1,v}$. Since $P_{In}$ is a fixed-point, it follows $P_{In}\in I_v$ and consequentially that the constant $\{1,1,1,...\}$ is in $\pi_v(I_v)$.
\end{remark}
\begin{remark}
    In \cite{Zgli97}, the Fixed-Point Index was applied to prove the existence of infinitely many periodic trajectories in the Rössler system.
\end{remark}
\begin{remark}
    Given a trefoil parameter $F_p$, it is easy to see the proof of Th.\ref{conti} above can be generalized for sufficiently small $C^k$ perturbations of $F_p$, where $k\geq1$.
\end{remark}
\begin{figure}[h]
\centering
\begin{overpic}[width=0.35\textwidth]{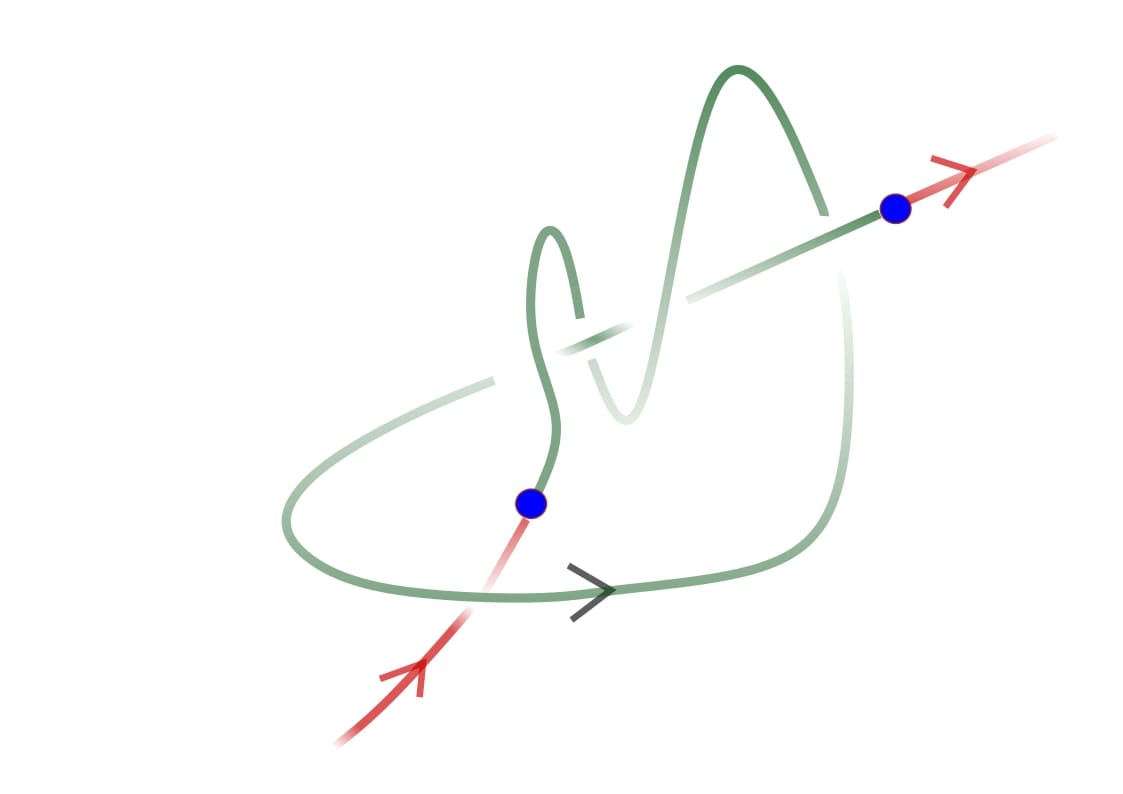}
\put(350,245){$P_{In}$}
\put(780,440){$P_{Out}$}
\end{overpic}
\caption[Fig39]{A heteroclinic knot more complex than a trefoil.}
\label{fig39}
\end{figure}
\section{Discussion:}

Before we conclude this paper, we would like to make several remark. To begin, first note both Th.\ref{th31} and Th.\ref{conti} are essentially facts about the periodic dynamics of the Rössler system - and that their proof is essentially two-dimensional. It is well-known in the theory of two-dimensional dynamics that one can study the topological dynamics of a surface homeomorphism by studying its class of essential periodic dynamics (see \cite{Bo} for a review of the topic). This leads us to ask the following question: can we describe the dynamics (and bifurcations) of the Rössler system in and around trefoil parameters by studying its periodic trajectories?\\

This question can and will be studied rigorously. In a series of future papers it will be proven that by applying ideas from \cite{BW}, it is possible to classify the knot-types of periodic trajectories in the set $Q$ - and prove their persistence under perturbations in the parameter space $P$. This would allow to reduce the flow in and around trefoil parameters to a one-dimensional model, which also acts as a geometrical model describing the flow.\\

Another natural question arising from both Th.\ref{conti} and Th.\ref{th31} is the following - can these results be extended to other heteroclinic knots in the Rössler system, or better yet, can we extend them to the actual Rössler system, as opposed to idealized models? At the moment, we do not have an answer for this question - and moreover, the our knowledge, no study (numerical or analytical) attempted to study which heteroclinic knots are generated by the Rössler system. Therefore we will be content with the following conjecture:

\begin{conj}
    Let $p\in P$ be a heteroclinic parameter for the Rössler system (see Def.\ref{def31}). Then, provided the heteroclinic knot it generates is complex at least like a trefoil knot, the flow includes infinitely many periodic trajectories.
\end{conj}

Finally, another interesting question which arises from both Th.\ref{th31} and Th.\ref{conti} is their possible connection with the theory of homoclinic bifurcations. The theory of homoclinic bifurcations had been extensively applied to numerically study and analyze the onset of chaos in the Rössler system (see \cite{MBKPS} and the references therein) - and yet, surprisingly, save for Prop.\ref{spiraltheorem} and Prop.\ref{cor216}, none of the results proven in this paper depend on that theory. Therefore, one is motivated to ask the following open questions with which we conclude this paper:

\begin{enumerate}
    \item Motivated by Cor.\ref{cor216} and Prop.\ref{spiraltheorem}, we first ask - what is the connection between heteroclinic trajectories in the Rössler system and homoclinic bifurcation phenomena in the Rössler system? 
    \item Let us now further note that in Th.\ref{th31} (and to a lesser extent, also in Th.\ref{conti}) the bounded heteroclinic trajectory serves a somewhat analogous role to that of the homoclinic trajectory in the proof of Shilnikov's Theorem. This motivates us to ask the following - does there exist a general theory of heteroclinic bifurcations, s.t. Th.\ref{th31} and Th.\ref{conti} are just special instances of facts within it? 
    \item Finally, if the answer to the above is positive, just how much of this hypothesized Heteroclinic Bifurcation Theory is analogous to the theory of homoclinic bifurcations?
\end{enumerate}

\printbibliography
\end{document}